\documentclass[twoside,11pt]{article}

\PassOptionsToPackage{numbers,sort&compress}{natbib}

\usepackage{blindtext}
\usepackage{placeins}
\usepackage{amsthm}
\usepackage[preprint]{formatting}
\usepackage{graphicx}
\usepackage{amsmath,amssymb,amsfonts,bm,dsfont}
\usepackage{tikz-cd}
\usepackage{tikz}
\usetikzlibrary{
  arrows.meta,
  positioning,
  calc
}

\usepackage[utf8]{inputenc}
\usepackage[T1]{fontenc}
\usepackage{url}
\usepackage{booktabs}
\usepackage{nicefrac}
\usepackage{microtype}
\usepackage{xspace}
\usepackage{array}
\usepackage{tabularx}
\usepackage{caption}
\usepackage{subcaption}
\usepackage{multirow}
\usepackage{mathrsfs}
\usepackage[dvipsnames,table]{xcolor}
\usepackage{mathtools}
\usepackage{float}
\usepackage{soul}
\usepackage{enumitem}
\usepackage{geometry}
\usepackage{lastpage}
\usepackage{hyperref}
\usepackage{cleveref}

\usepackage[normalem]{ulem}

\setcitestyle{numbers,square,sort&compress}
\hypersetup{
  colorlinks=true,
  linkcolor=blue,
  citecolor=blue,
  urlcolor=blue
}

\newcommand{\hlred}[2][red!7]{\colorbox{#1}{$\displaystyle #2$}}
\newcommand{\sssoap}{\textsc{SS--eSOAP}\xspace}

\newcommand{\wgray}[1]{#1}

\newcommand{\bvar}[1]{\textcolor{blue}{#1}}

\newcommand{\cures}{c_{u,\mathrm{res}}}

\newcommand{\Lean}{(\ensuremath{\checkmark}\ Formalized and verified in Lean)}

\definecolor{grayshade}{RGB}{240,247,255}

\definecolor{certiforange}{RGB}{180,68,2}

\definecolor{certifpurple}{RGB}{15,25,100}

\definecolor{certpipelinetextmuted}{RGB}{82,82,88}
\definecolor{certpipelinerulegray}{RGB}{140,140,148}
\definecolor{certpipelinerulemid}{RGB}{176,176,183}

\tikzset{
  certpipeline flow/.style={
    -{Stealth[length=1.15mm,width=0.80mm]},
    draw=certpipelinerulegray,
    line width=0.34pt
  },
  certpipeline transition/.style={
    midway,
    right=2.9pt,
    fill=white,
    inner xsep=0.75pt,
    inner ysep=0.35pt,
    text=certpipelinetextmuted,
    font=\fontsize{5.85}{6.2}\selectfont\itshape
  },
  certpipeline stage/.style={
    draw=certpipelinerulemid,
    fill=white,
    line width=0.31pt,
    minimum width=8.15cm,
    text width=7.60cm,
    minimum height=11.45mm,
    align=center,
    inner xsep=6.2pt,
    inner ysep=3.2pt,
    outer sep=0pt
  }
}

\newcommand{\certpipelinetitle}[2]{%
  {\fontsize{8.1}{8.55}\selectfont
    \bfseries\color{#1}\textsc{#2}\par}%
}

\newcommand{\certpipelinesubtitle}[1]{%
  {\vspace{0.45pt}%
    \fontsize{6.2}{6.55}\selectfont
    \color{certpipelinetextmuted}#1\par}%
}

\newcommand{\certpipelinephasetag}[2]{%
  {\fontsize{4.8}{5.15}\selectfont
    \bfseries\scshape\color{#1}#2}%
}

\newcommand{\certpipelinetagpanel}[6]{%
  \node[certpipeline stage,#6] (#1) {%
    \certpipelinetitle{#2}{#4}%
    \certpipelinesubtitle{#5}%
  };
  \draw[#2!78,line width=0.54pt]
    ([xshift=0.2pt]#1.north west) --
    ([xshift=-0.2pt]#1.north east);
  \node[
    anchor=west,
    fill=white,
    inner xsep=2.15pt,
    inner ysep=0.55pt,
    text=#2
  ] at ([xshift=6.4mm]#1.north west) {%
    \certpipelinephasetag{#2}{#3}%
  };
}

\definecolor{nsgreen}{RGB}{0,100,0}
\definecolor{lightproofgray}{gray}{0.55}
\newcommand{\nsgreen}[1]{\textcolor{nsgreen}{#1}}

\definecolor{darksalmon}{RGB}{240,80,70}
\newcommand{\dsalmon}[1]{\textcolor{darksalmon}{#1}}

\newcommand{\DLamL}{\ensuremath{\dsalmon{\Lambda_{\mathcal L}}}}
\newcommand{\DLamLlow}{\ensuremath{\dsalmon{\Lambda_{\mathcal L}^{\rm low}}}}
\newcommand{\DLamLhigh}{\ensuremath{\dsalmon{\Lambda_{\mathcal L}^{\rm high}}}}
\newcommand{\DLamDLhigh}{\ensuremath{\dsalmon{\Lambda_{D^\alpha\mathcal L}^{\rm high}}}}
\newcommand{\DLamDLlow}{\ensuremath{\dsalmon{\Lambda_{D^\alpha\mathcal L}^{\rm low}}}}
\newcommand{\DLamNthree}{\ensuremath{\dsalmon{\Lambda_{\mathcal N,3}}}}

\newcommand{\DLamDNthree}{\ensuremath{\dsalmon{\Lambda_{D^\alpha\mathcal N,3}}}}

\newcommand{\DLamThree}{\ensuremath{\dsalmon{\Lambda_{3}}}}

\newcommand{\DLamRPhi}{\ensuremath{\dsalmon{\Lambda_{\mathcal R,\Phi}}}}
\newcommand{\DLamDR}{\ensuremath{\dsalmon{\Lambda_{D^\alpha\mathcal R}}}}
\newcommand{\DLamR}{\ensuremath{\dsalmon{\Lambda_{\mathcal R}}}}
\newcommand{\DLamStab}{\ensuremath{\dsalmon{\Lambda_{\rm stab}}}}

\newcounter{certalg}
\newenvironment{certalgorithm}[1]{%
  \refstepcounter{certalg}%
  \par\medskip
  \noindent\rule{\linewidth}{0.4pt}\par\smallskip
  \noindent\textbf{Algorithm \thecertalg. #1}\par\smallskip
  \begin{enumerate}[label=\textbf{\arabic*.},leftmargin=2.2em,itemsep=0.35em,topsep=0.25em]
}{%
  \end{enumerate}
  \smallskip\noindent\rule{\linewidth}{0.4pt}
  \par\medskip
}

\makeatletter

\let\theproposition\relax
\makeatother

\theoremstyle{plain}
\newenvironment{grayproof}
  {\hfill \begin{proof} \hspace{0.5mm} \color{lightproofgray}}
  {\end{proof}}
\newtheorem*{remark*}{Remark}
\newtheorem{proposition}{Proposition}

\newcommand{\paperTitle}{Self-Similar Singularity of the Euler Equations on $\mathbb{R}^3$}
\newcommand{\paperAuthors}{%
  \name Adarsh Ganeshram* \email aganeshram@berkeley.edu
  \AND
  \name Valentin Duruisseaux* \email vduruiss@caltech.edu
  \AND
  \name Anima Anandkumar \email anima@caltech.edu
  \AND
  \\
  \addr Department of Mathematics, University of California Berkeley, Berkeley, CA, USA
  \AND
  \addr Department of Computing and Mathematical Sciences, California Institute of Technology, Pasadena, CA, USA%
}
\title{\paperTitle}
\author{\paperAuthors}
\date{}

\makeatletter
\newcommand{\supportingInformationAuthors}{%
  \noindent
  \begin{tabular*}{\textwidth}{@{}l@{\extracolsep{\fill}}r@{}}
    {\large\bfseries Adarsh Ganeshram*}
      & {\normalsize\scshape aganeshram@berkeley.edu} \\[0.48em]
    {\large\bfseries Valentin Duruisseaux*}
      & {\normalsize\scshape vduruiss@caltech.edu} \\[0.48em]
    {\large\bfseries Anima Anandkumar}
      & {\normalsize\scshape anima@caltech.edu}
  \end{tabular*}\par
  \vspace{1.95em}%
  \noindent{\normalsize\itshape
    Department of Mathematics, University of California Berkeley, Berkeley, USA\par}%
  \vspace{0.55em}%
  \noindent{\normalsize\itshape
    Department of Computing and Mathematical Sciences, California Institute of Technology, Pasadena, USA\par}%
}

\newcommand{\supportingInformationHeading}{%
  \begin{center}
    {\LARGE\bfseries Supporting Information for\par}%
    \vspace{0.55em}%
    {\Large\bfseries \paperTitle\par}%
  \end{center}%
  \vspace{1.65em}%
  \supportingInformationAuthors
  \vspace{0.9em}%
}

\newcommand{\printSupplementalContents}{%
  \section*{Contents}%
  \@starttoc{stoc}%
}
\newcommand{\startSupplementalTocCapture}{%
  \let\mainaddcontentsline\addcontentsline
  \renewcommand{\addcontentsline}[3]{%
    \def\supplementalTocSource{##1}%
    \def\mainTocSource{toc}%
    \ifx\supplementalTocSource\mainTocSource
      \mainaddcontentsline{stoc}{##2}{##3}%
    \else
      \mainaddcontentsline{##1}{##2}{##3}%
    \fi
  }%
}
\newcommand{\stopSupplementalTocCapture}{%
  \let\addcontentsline\mainaddcontentsline
}
\makeatother

\begin{document}

\twocolumn[
  \begin{@twocolumnfalse}
    \maketitle

    \begin{abstract}
      
\hfill 

We provide evidence of a finite-time singularity in the 3D Euler equations on the unbounded domain. Using a physics-informed neural network (PINN) with a self-similar ansatz, we find an approximate singular profile for the Euler system at the critical blowup rate of $0.5$ and certify it using a spline representation. The transport field associated with the obtained profile has local outgoing property throughout the domain that suggests linear damping, a key stabilizing mechanism for the candidate profile. We also establish a framework for proving nonlinear stability of the approximate self-similar profile, reducing the analysis to a large but finite collection of explicit estimates and computable constants.

\hfill  


\hfill \\ 

\vspace{1mm}

    \end{abstract}
  \end{@twocolumnfalse}
]

\section{Introduction}

A central problem in the analysis of partial differential equations (PDEs) is whether all solutions arising from smooth initial data remain regular at all times or whether they may lose regularity in finite time. Such a loss of regularity is known as \emph{singularity formation}, or \emph{blowup}. In particular, establishing whether solutions of the 3D incompressible Navier--Stokes equations can develop a finite-time blowup is among the most prominent open problems in fluid dynamics, and is one of the six unsolved Millennium problems~\citep{FeffermanNavierStokesCMI} and the $15^{{\textit th}}$ Smale problem~\citep{smale1998mathematical}. Closely related is the question of finite-time singularity formation in the Euler equations, which govern inviscid flows and differ from the Navier–Stokes equations by the absence of viscous dissipation.\\ 

Blowups have been established in related simpler problems or under modified assumptions. For instance, finite-time blowups have been established for the Euler equations in the presence of a boundary given smooth initial conditions~\citep{HouJiaje}, or under non-smooth initial conditions on the unbounded domain $\mathbb{R}^3$~\citep{nonSmoothProof,shkoller2026incompressible,JiajeIncompI,JiajeIncompII}.  
However, it is not clear how these additional assumptions can be removed. For instance, boundary-driven blowup mechanisms exploit the confinement imposed by solid walls to generate rapid vorticity growth~\citep{HouLuoBoundary,HouLuo1,HouJiaje,Chen2026Analysis}, which does not work when there is no boundary. In addition, there are non-existence results which rule out finite-time blowups in the Euler system on $\mathbb{R}^3$ under certain conditions \citep{chae2007nonexistence,chae2013blow}, which include strong decay and integrability assumptions on the vorticity profile. \\

\begin{figure*}[htbp] 
\centering
\includegraphics[width=\linewidth]{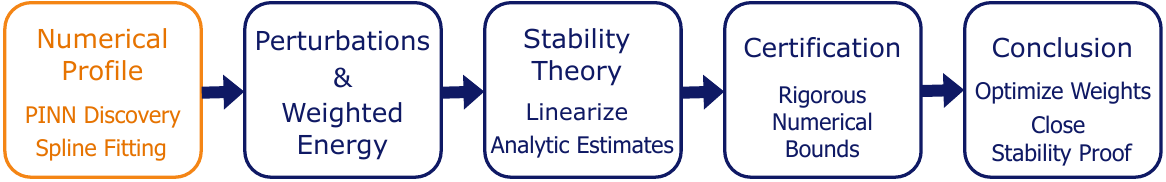}
\caption{\textbf{From numerical discovery toward certified stability.} An approximate self-similar profile is first discovered using PINNs and converted to a spline for analytic estimates and rigorous certification. Perturbations around this profile are introduced together with a parametrized weighted energy that captures the relevant singular behavior of the solution. The stability theory is then developed by linearizing around the spline profile and deriving analytic estimates to control the linear and nonlinear contributions. The analytic estimates depend on explicit numerical constants and quantitative properties of the approximate profile, which are rigorously enclosed using interval arithmetic together with the profile data, stability constants, and PDE residuals. Finally, the parametrized singular weights and other parameters are optimized so that the certified damping dominates the other terms and PDE residuals. This provides a route toward rigorous computer-assisted nonlinear stability, with the analytic and certified estimates also suitable for subsequent formalization in Lean. }\label{fig:flowchart} \vspace{2mm}
\end{figure*}

A common strategy for establishing finite-time blowup is to construct an approximate solution numerically, and then to use stability theory to prove the existence of a singular solution to the governing PDEs that remains sufficiently close to this approximation. A self-similar ansatz offers a tractable path to both the numerical search for such an approximate profile and its stability analysis. It allows us to continuously zoom in or out at the specified rate,  as the solution evolves, so that its overall ``shape'' is unchanged over time. More precisely,  a blowup solution is \emph{self-similar} if, under a suitable dynamic rescaling, it remains close to an approximate blowup profile all the way up to the blowup time. The dynamic rescaling equations thus serve as a starting point for stability analysis around a given approximate profile.\\

Since it is challenging to manually construct approximate self-similar profiles as blowup candidates, Physics-Informed Neural Networks (PINNs)~\citep{PINN_OG} have recently been adopted as a more flexible alternative. When the PINN computations are carried in high precision (FP64), they can discover approximate singular profiles \citep{deepMindPaper1,deepMindPaper2,wang2025high}, while (approximately) enforcing the governing equations together with symmetries, normalization, decay and other physical constraints. They have been successful in discovering profiles for simplified versions of the Euler equations, such as Burgers and Boussinesq equations.  However, PINN formulations have not yet succeeded in discovering potential singularities in the Euler and Navier--Stokes equations on the unbounded domain. It remains unclear how to choose suitable constraints and ansatz for such a PINN formulation, and how to overcome the resulting optimization challenges. \\

Recent theoretical results also impose additional constraints on the regime where such finite-time blowups may occur, further reducing flexibility for the selection of ansatz and PINN formulations.  \citet{constantin2026putative} recently establish that the critical self-similar scaling exponent of $\lambda=0.5$, known to be natural for Navier--Stokes equations, is also the most promising (or ``mathematically distinguished'') regime for self-similar Euler blowups. They essentially rule out  stable self-similar singularities for the Euler system when $\lambda<0.5$, and $\lambda>0.5$ for the commonly used axis-symmetric ansatz considered here.  Despite sharing this critical scaling exponent of $0.5$ for potential blowups, the Euler and Navier--Stokes equations cannot share a common singular profile, since the effect of viscosity cannot be ignored in this regime. Overall, establishing the existence of singularities in either system is highly non-trivial.\\

In addition to finding an approximate profile, establishing stability is crucial but complicated. This becomes simpler when a global outgoing property is satisfied~\citep{DGBlowup}, which intuitively pushes the self-similar flow away from the singular core. This expansive effect of the global outgoing flow is also closely connected to the favorable higher-order damping that underlies the stability analysis. However, for the Euler system with a self-similar ansatz, \citet{constantin2026putative} provides strong theoretical evidence that nontrivial fixed points of the meridional flow are likely to exist and have been ``well known to cause tremendous headaches" in other related problems. Given that the global outgoing property cannot hold for the Euler system with self-similar ansatz, they instead consider a weaker local property that locally pushes the self-similar flow away from these fixed points, preventing perturbations from concentrating there and providing a favorable structure for the stability analysis. Establishing stability with this weaker local outgoing property however requires a substantially more delicate analysis.\\

\begin{figure*}[htbp] 
\centering

\begin{subfigure}{\linewidth}
\centering
\includegraphics[width=0.495\linewidth]{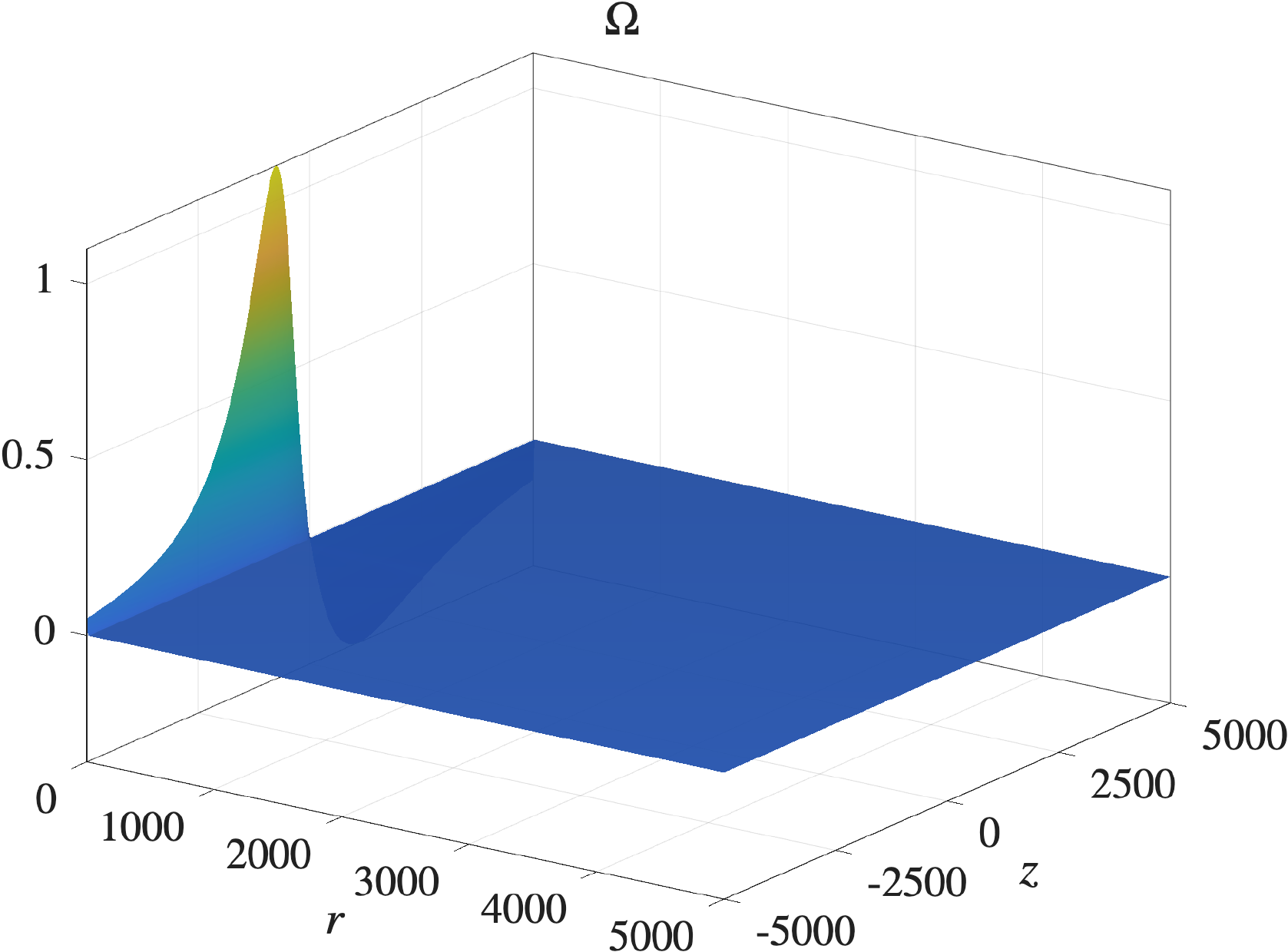}\hfill
\includegraphics[width=0.495\linewidth]{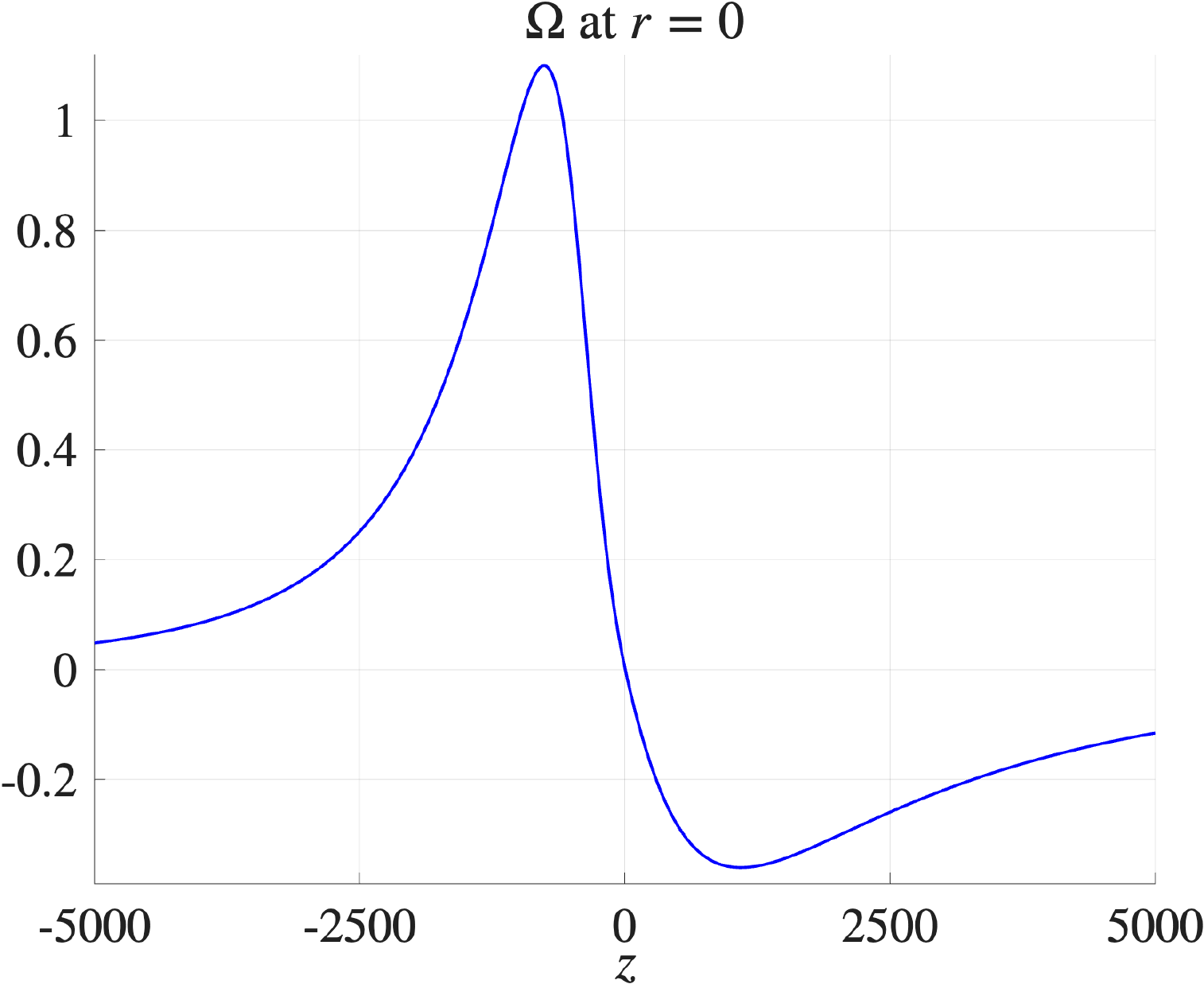}
\end{subfigure}\vspace{0mm}

\begin{subfigure}{\linewidth}
\centering
\includegraphics[width=0.495\linewidth]{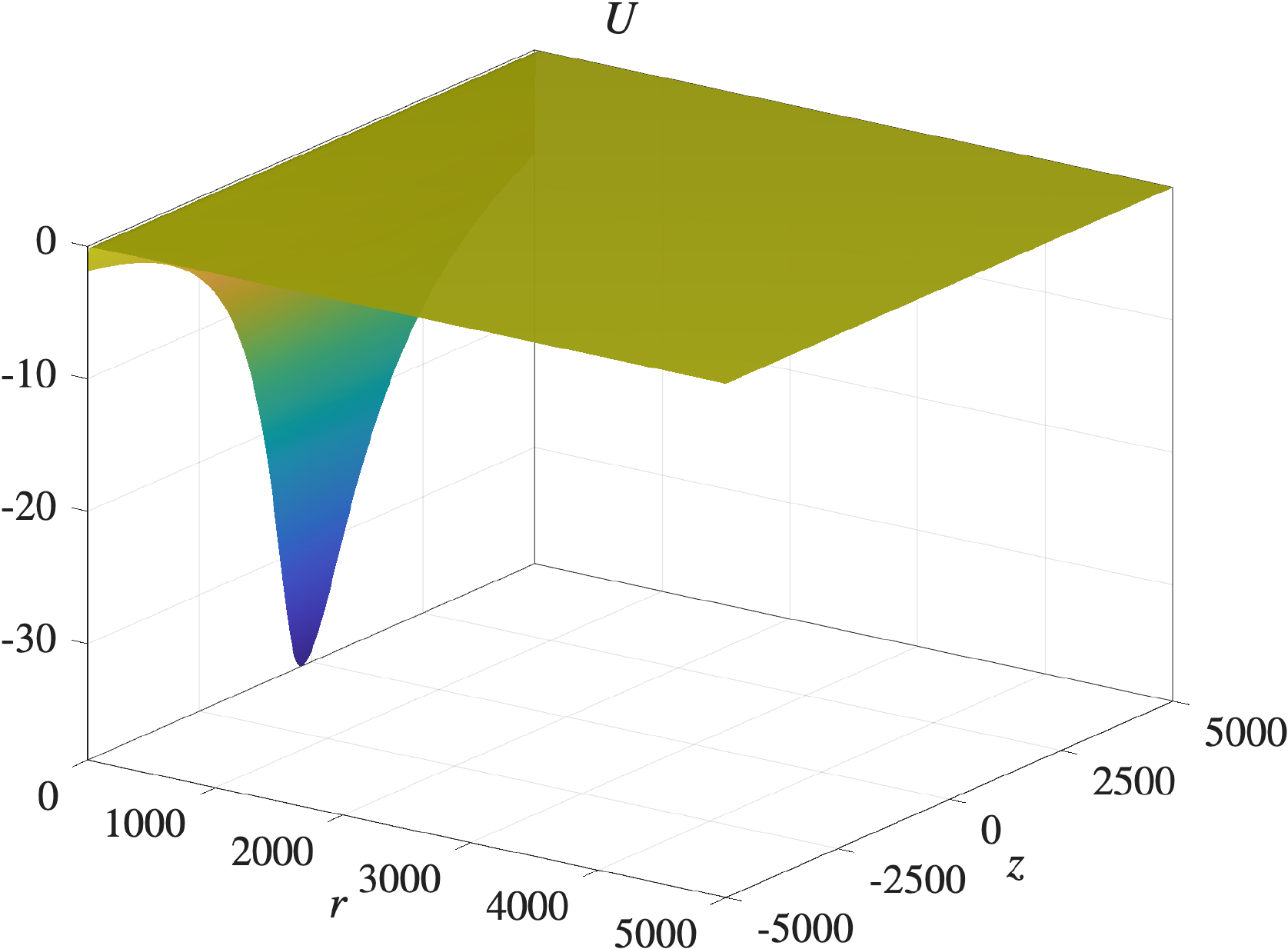}\hfill
\includegraphics[width=0.495\linewidth]{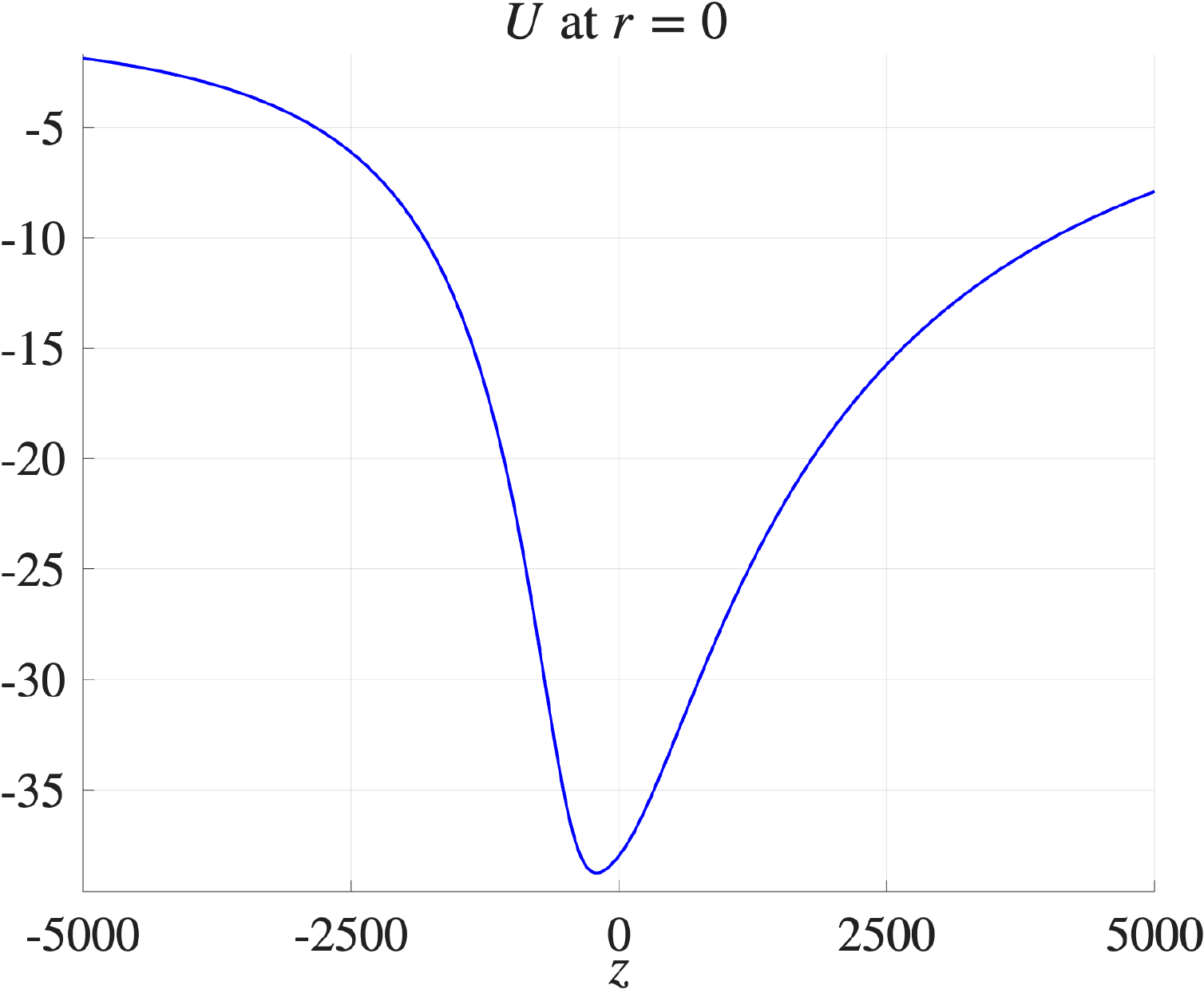}
\end{subfigure}\vspace{0mm}

\begin{subfigure}{\linewidth}
\centering
\includegraphics[width=0.495\linewidth]{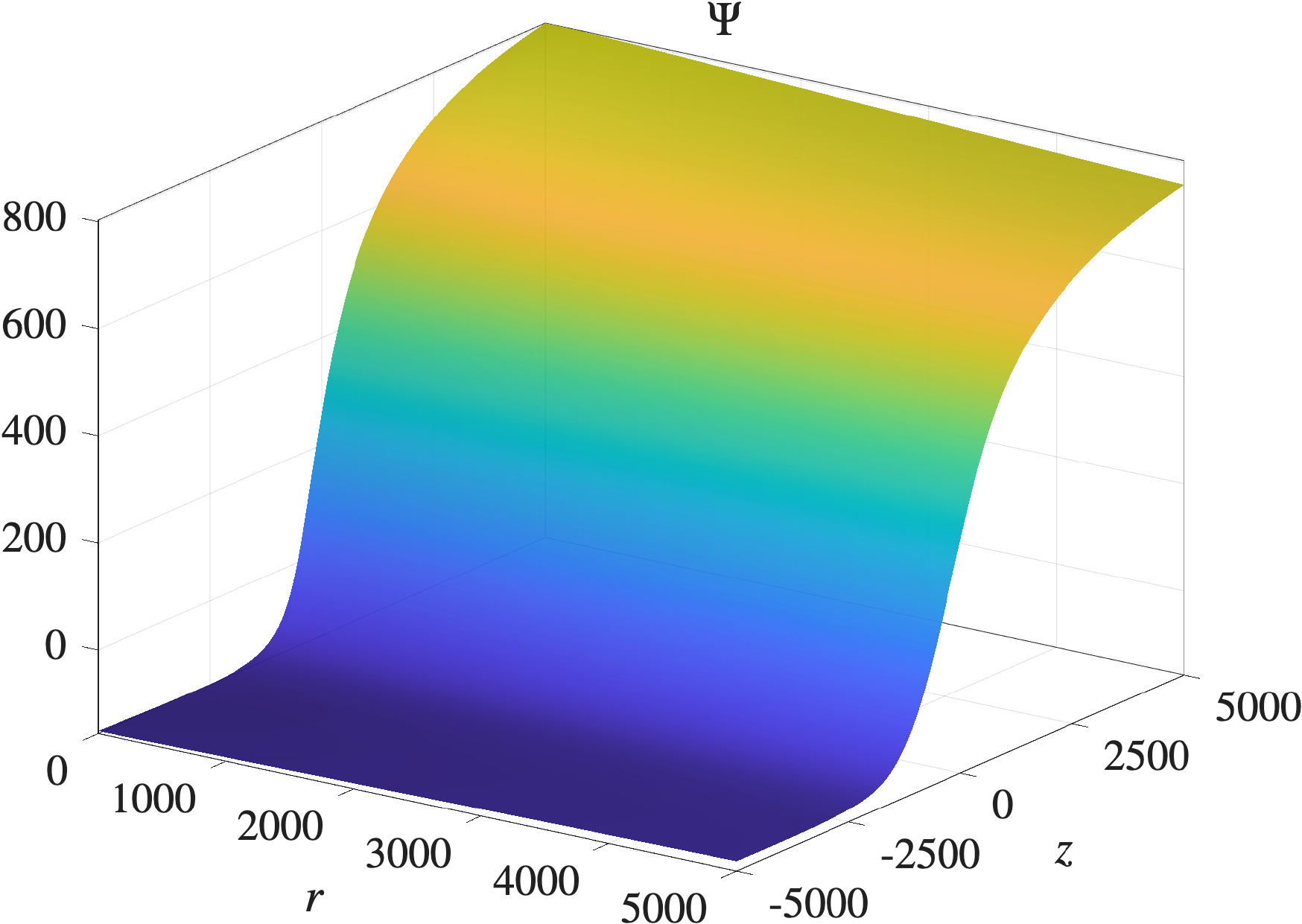}\hfill
\includegraphics[width=0.495\linewidth]{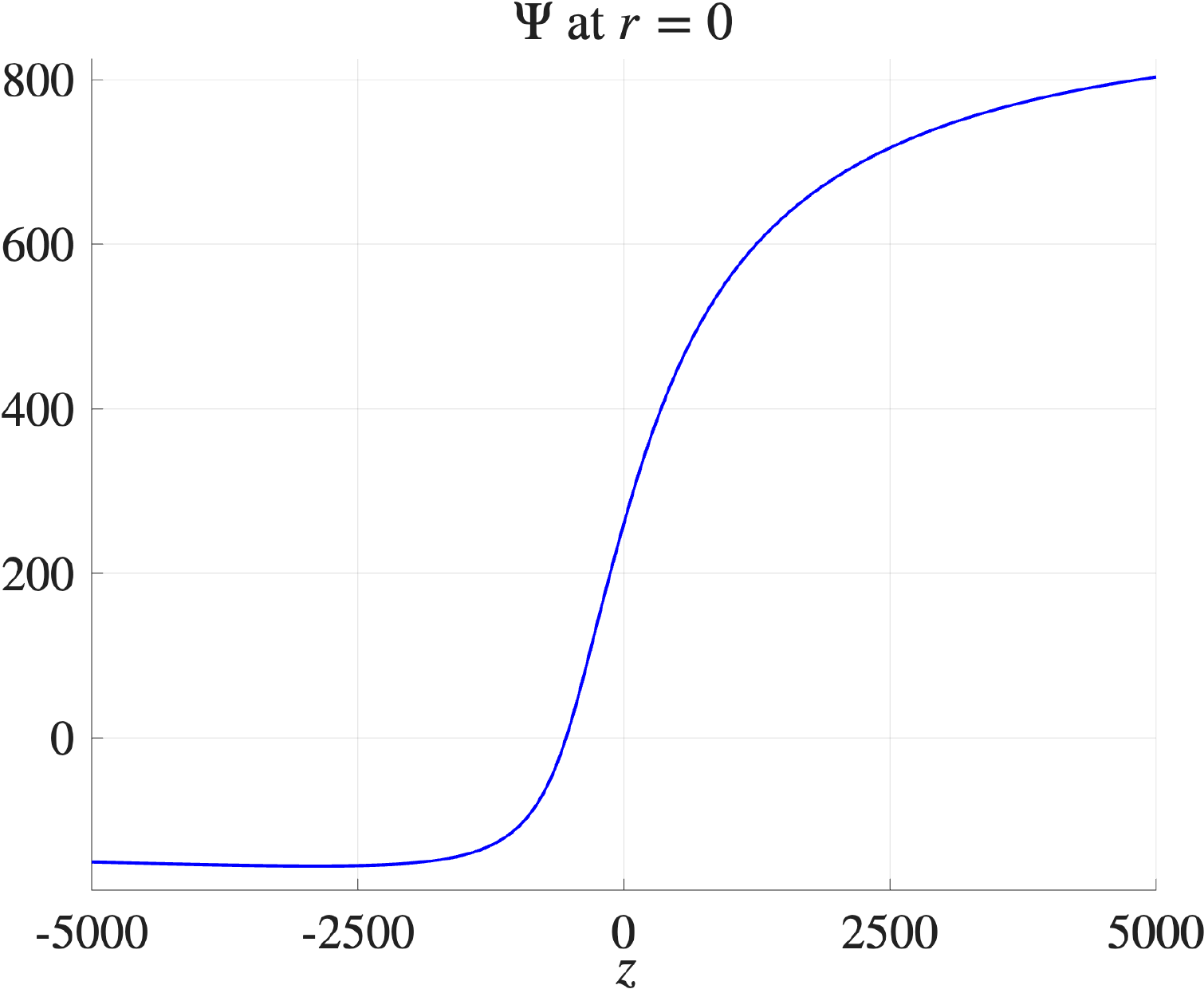}
\end{subfigure}

\vspace{3mm}

\caption{Visualization of the approximate profile $\Omega, U, \Psi$ as three-dimensional surface plots for the 3D axisymmetric Euler equations with their corresponding $r=0$ cross-sections.}
\label{fig: short_profiles} \vspace{1mm}
\end{figure*}

\hfill 

 \textbf{We provide the first evidence that the Euler equations admits a self-similar finite-time blowup on $\mathbb{R}^3$} at the critical scaling exponent $\lambda=0.5$, referred to as ``mathematically distinguished'' by \citet{constantin2026putative}. Our numerical results indicate substantial linear damping in the domain, which could be a key stabilizing mechanism. We also establish a proof framework for the nonlinear stability of the approximate blowup profiles in a separate stability paper~\citep{EulerBlowupStability}, reducing the argument to a finite collection of explicit, rigorously verifiable estimates. The spline profile will be made available together with code demonstrating that the Euler equations are satisfied to high precision. The Lean formal verification of the derivations presented in this manuscript is available at \href{https://github.com/lean-dojo/Euler}{github.com/lean-dojo/Euler}. Our methodology is shown in \Cref{fig:flowchart}.

\hfill \\

\noindent \emph{Numerical Profiles.} Using a  PINN formulation, we obtain a highly accurate approximate self-similar profile for the Euler system. A novel combination of complementary tools is essential to achieving this level of accuracy, ensuring that the formulation captures the correct physics while converging to a high-precision solution. Imposing appropriate soft and hard constraints is necessary to prevent convergence to trivial or degenerate solutions, a failure mode observed with other PINN formulations. 

\vspace{3mm}

Motivated by~\citet{pengfeiPaper}, we adopt a traveling-wave ansatz propagating along the $z$-axis in the axisymmetric setting. Although the traveling-wave and stationary formulations are equivalent in the Euler setting~\citep{constantin2026putative} since we can translate and scale the system appropriately, they induce substantially different optimization landscapes. In the traveling-wave formulation, the singularity propagation speed is treated as a free parameter, providing greater flexibility to the optimizer. Moreover, the traveling-wave ansatz does not impose the $z$-parity constraints that often accompany the stationary self-similar ansatz, and PINN optimization consistently fails to converge when these $z$-parity constraints are imposed. These features substantially alleviate the optimization difficulties encountered in previous PINN formulations and are essential to the success of our approach.\\

We employ optimization techniques to obtain a profile with high-accuracy. We use adaptive loss-weighting schedules to balance competing objectives, FP64 arithmetic to reduce round-off errors and improve gradient accuracy, PINN boosting~\citep{fang2024boosting,boostingYongji} to progressively correct residual errors, and adaptive collocation with residual-based resampling to concentrate training points in regions where the profile is most difficult to resolve. We also employ a $\sinh$ transformation to efficiently sample the unbounded domain. Finally, the choice of optimizer is also important, and we use efficient self-scaled curvature-aware optimizers.  While Adam~\citep{adamOptim} typically plateaus at a profile-residual MSE of order $10^{-4}$, the recently proposed SS-eSOAP~\citep{SS-SOAP} optimizer reaches the $10^{-6}$ regime with only a $15\%$ increase in per-step cost. We finally employ the more computationally expensive SS-Broyden~\citep{Al-Baali1998,SSBroyden2} for high-accuracy fine-tuning, reducing the profile-residual MSE to order $10^{-10}$.

\vspace{2.5mm}

The optimized PINN profile is converted into piecewise polynomial splines suitable for rigorous certification. The spline representation facilitates exact interval-wise differentiation and evaluation of profile-dependent quantities. We  use \texttt{Arb}~\citep{Arb}, an arbitrary-precision interval arithmetic library, to rigorously control numerical and rounding errors throughout the certification procedure. This allows in particular to bound the residual errors of the spline profile by $\mathcal{O}(10^{-5})$ in $L^2$ and $\mathcal{O}(10^{-3})$ in $L^\infty$.

\hfill 

\noindent \emph{Connections to theory.} We verify that our numerical profile does not violate non-existence results for singularities in the Euler system. In particular, our vorticity profile exhibits only polynomial spatial decay, whereas the non-existence theorems of~\citet{chae2007nonexistence,chae2013blow} require super-algebraic decay. We further investigate the case in which the scaling exponent $\lambda$ is treated as a free parameter and optimized jointly with the profile. Despite this added flexibility, the optimized scaling exponent converges to a value close to the critical exponent $\lambda=0.5$ predicted by \citet{constantin2026putative}. \\

Additionally,  we consider the family of convection-varied Euler equations proposed by \citet{pengfeiPaper}, where a parameter $\varepsilon$ allows convection to be progressively increased from the non-convective case at $\varepsilon=0$ to the full Euler system at $\varepsilon=1$. Using direct time evolution from a common numerically generated initial state, \citet{pengfeiPaper} observed a traveling self-similar singularity scenario at the tested values \(\varepsilon\in \{0,0.1,0.2,0.3\} \), whereas the evolution at \(\varepsilon=0.4\) did not exhibit the same behavior. Their dynamic-rescaling computations further indicated that the computed self-similar state became progressively less stable as \(\varepsilon\) increased. These results left open whether self-similar singular profiles persist for stronger convection but are not reached by the forward dynamics. Our PINN formulation finds high-accuracy approximate self-similar profiles throughout \(\varepsilon\in[0,1]\), providing numerical evidence that the disappearance of the earlier dynamically observed scenario should not be interpreted as the nonexistence of singular profiles. As $\varepsilon$ approaches $1$, the scaling exponent $\lambda$ decreases toward the critical value $\lambda=0.5$, providing further agreement with the theoretical prediction of \citet{constantin2026putative}. \\

Looking at the transport field of the  profiles, we identify two regions of low transport, which can be more challenging for stability analysis. The first  surrounds a fixed point of the meridional flow on the symmetry axis, while the second is centered at an off-axis meridional fixed point. The axial fixed point has a nearly degenerate locally outgoing property which agrees with the theory from~\citep{constantin2026putative} for $\lambda=0.5$ (with $c_*=0$ in their notation). In contrast, no such constraint is imposed on the local outgoing behavior of the off-axis fixed point in the theory of~\citep{constantin2026putative}, which is consistent with our numerical results where the off-axis fixed point is strictly locally outgoing. We  use these properties in our stability analysis. \\

\hfill 

\noindent \emph{Stability.}  Starting from the spline representation of the numerical profiles, we develop a stability framework, following the strategy from prior blowup stability works~\citep{DGBlowup,HouJiaje,Chen2026Analysis}. We linearize the appropriate dynamically rescaled equations around the spline  profiles and formulate a stability theorem in suitable weighted $L^2$ and Sobolev norms. We then derive complementary low-order and high-order energy estimates for both the linear and nonlinear terms, with explicit computable bounds. These estimates are incorporated into a single energy functional that provides the quantitative structure needed to determine whether the full stability argument can be closed. We refer to a separate paper~\citep{EulerBlowupStability} for the details of the analysis. The stability question is thereby reduced to determining whether a large finite collection of explicit estimates can be rigorously certified with sufficient margin using interval arithmetic or \texttt{Arb} computations.\\

We first focus on linear damping for stability analysis in this paper. Establishing damping at the linear level is a key prerequisite for the complete stability proof, since insufficient linear damping cannot generally be overcome by the remaining nonlinear estimates. Its certification thus isolates a potentially important stabilizing mechanism and provides quantitative evidence that a full nonlinear stability argument may be within reach. In contrast to prior works, the singular weight functions defining the stability energy norms, together with other tunable parameters entering the stability argument, are optimized numerically rather than selected through a lengthy and often highly nontrivial process of manual trial and error. Here, we represent the weight functions as a combination of locally supported Gaussian basis functions, a singular component near the origin, asymptotically growing neural networks in the far field, and a global neural-network correction, with all components parametrized and jointly optimized.\\

After optimizing the weights, we obtain low-order damping throughout most of the domain, with the exception of two neighborhoods around the meridional fixed points. For the fixed point on the axis, we can use singular weights and local analysis to obtain low-order damping. At the off-axis fixed point, the flow is strictly locally outgoing, enabling higher-order damping to be achieved through appropriate tuning of the weighting functions. In particular, the relevant higher-order derivatives have a favorable sign and contribute positively to the damping estimate for adequate Sobolev derivative order. Overall, we obtain low-order damping on most of the computational domain, and high-order damping in the problematic off-axis low transport region, and damping outside the computational domain analytically. The combined energy exploits these complementary mechanisms. The resulting damping margin provides a useful quantitative indication that the stability scheme may be feasible. \\

\noindent \emph{Remaining Work to Complete the Proof.} The remaining work within this framework is primarily computational for the full nonlinear stability in the Euler system: rigorously certifying the profile-dependent constants and margins required in the stability framework using spline representations, interval arithmetic, certified matrix bounds, and finite-dimensional optimization. Some estimates, presented in detail in a separate paper~\citep{EulerBlowupStability}, may require further refinement to obtain sufficient margin, but the modular structure of the argument allows such adjustments without changing the overall proof strategy. If these certifications yield the required positive margins and thereby establish nonlinear stability of the rescaled profile, this stability can then be transferred through the dynamic rescaling to finite-time singularity formation in the original variables. At the critical scaling $\lambda=1/2$, control of the modulation parameters near their self-similar values would yield stability of the candidate profile and, after reconstruction, an admissible solution that becomes singular in finite time.\\

\noindent \emph{Formalization in Lean.} In the \texttt{LeanPDE} paper, we formalize key components of the argument in Lean~\citep{Lean4}, using \texttt{Mathlib}~\citep{mathlib} to verify symbolic identities, derivations, and other algebraic steps that enter the proof. Rigorous numerical certificates are generated separately in \texttt{Julia} using \texttt{Arb} interval arithmetic~\citep{Julia,Arb}, then replayed and verified in Lean using \texttt{TorchLean} ~\citep{TorchLean} so that the verification pipeline assigns distinct roles to symbolic formalization and certified numerical computation. More broadly, this work is part of the \texttt{LeanPDE} effort to develop tools and infrastructure for machine-checked PDE analysis, including closer integration between proof assistants, computer algebra systems, and rigorous numerical methods.  \\

\noindent \emph{Broader Implications.} While much of the recent progress in AI-assisted mathematics has focused on the role of large language models (LLMs) in tasks such as theorem proving and conjecture testing, the central ingredient in this work is physics-informed optimization. The governing PDEs and physics knowledge are incorporated directly into the optimization problem, allowing candidate singularity structures to be found. In this paper, simple multilayer perceptron architectures are sufficient for PINNs to identify candidate singularity profiles. In other settings, however, more sophisticated architectures and frameworks, such as Neural Operators~\citep{FNO,duruisseaux2025FNOGuide,azizzadenesheli2024neural}, may be needed. While the approach considered here seeks an individual solution for a specific PDE instance, neural operators are designed to learn solution operators across entire families of parameterized PDEs. Physics-informed neural operators (PINOs) further incorporate the governing PDEs and other physical constraints directly into the operator-learning framework~\citep{PINO_OG,ganeshram2025fcpinohighprecisionphysicsinformed}. Such approaches could enable training across families of PDEs and support curriculum strategies that progress from simpler regimes toward the more challenging target problem. In the current context, for example, one could use a PINO to learn the solution operator across a family of Euler equations parameterized by convection strength, starting from weaker convection and progressively approaching the full-convection regime.  Such progressive learning may be helpful to discover blowups in other open problems. More broadly, such physics-informed and physics-centric AI are critical ingredients  across many areas of research involving physical systems~\citep{Anima-Daedulus}.\\

Closely related blowup results were developed concurrently with our results, including constructions for the incompressible porous-medium equation~\citep{alpoge2026ipm}, the inviscid Boussinesq system~\citep{alpoge2026boussinesq}, the 3D Euler equations with smooth forcing~\citep{alpoge2026euler}, the unforced 3D Euler equations~\citep{openai2026euler} and the forced 3D Navier--Stokes equations~\citep{openai2026navierstokes}. These results are almost entirely developed by LLMs, reusing and leveraging ideas from the existing mathematical literature, e.g. the forcing strategy from \citet{cordoba2023blowup}, using massive multi-agent searches requiring enormous computational resources. In our paper, we consider the unforced Euler equations and use a self-similar profile ansatz to guide PINN-based discovery of an approximate blowup profile, which we then seek to promote to an exact finite-time blowup solution by establishing its nonlinear stability under the dynamically rescaled Euler evolution. The central methodological components of our work, including the PINN framework, the choice of ansatz, and the optimization strategy, are developed entirely by the authors themselves. LLMs served in a limited supporting capacity for locating references, understanding relevant background material and definitions, checking derivations, and assisting with Lean formalization. In presenting the method, we make the reasoning behind the main methodological choices explicit, showing how they arise from the structure of the problem, why they are appropriate, and how they align with existing theoretical results. In this sense, our approach is more closely aligned with the view that the lasting value lies not only in obtaining a solution, but also in developing mathematical concepts, techniques, and intuition that can be assimilated and built upon by the broader community~\citep{fields2026misalignment}.

\clearpage

\section{The Euler Equations}

\subsection{The 3D Axisymmetric Equations}

Starting from the incompressible 3D Euler equations
\begin{equation}
\label{eq:short_Original_Euler_shortV}
\partial_t \tilde u + (\tilde u\cdot\nabla) \tilde u = -\nabla \tilde p,
\qquad
\nabla\cdot \tilde u = 0,
\end{equation}
we solve for the velocity field $\tilde u(x,t)\in\mathbb{R}^3$ and the scalar pressure $\tilde p(x,t)\in\mathbb{R}$, for $x\in\mathbb{R}^3$ and $t\ge 0$. Taking the curl yields the vorticity formulation
\begin{equation}
\label{eq:short_short_vorticity}
\partial_t \tilde \omega + \tilde u\cdot\nabla \tilde \omega = \tilde \omega\cdot\nabla \tilde u,
\qquad
\tilde \omega := \nabla\times \tilde u .
\end{equation}
\noindent For flows that are symmetric about a fixed spatial axis, such as the $\mathring z$-axis, we can rewrite~\eqref{eq:short_Original_Euler_shortV} in cylindrical coordinates $x=(\mathring{r}\cos\theta,\mathring{r}\sin\theta,\mathring{z})$ with $ \mathring{r} \ge 0$, $   \theta\in[0,2\pi)$, $\   \mathring{z}\in\mathbb{R}$. In these coordinates, 
\begin{align*}
\tilde u = u^r(\mathring{r},\mathring{z},t)\,e_r + u^\theta(\mathring{r},\mathring{z},t)\,e_\theta + u^z(\mathring{r},\mathring{z},t)\,e_z, \\
\tilde \omega = \omega^r(\mathring{r},\mathring{z},t)\,e_r + \omega^\theta(\mathring{r},\mathring{z},t)\,e_\theta + \omega^z(\mathring{r},\mathring{z},t)\,e_z,
\end{align*}
where $e_r, e_\theta, e_z$ are the cylindrical unit vectors. By definition of axisymmetry, we also have $\partial_\theta(\cdot)=0$. Following \citet{HouLuoBoundary,Chen2026Analysis}, we introduce the rescaled swirl and streamfunction variables
\begin{equation*}
u_{\bullet} := \frac{u^\theta}{\mathring{r}},\qquad
\omega_{\bullet} := \frac{\omega^\theta}{\mathring{r}},\qquad
\psi_{\bullet} := \frac{\psi^\theta}{\mathring{r}},
\end{equation*}
where $\psi^\theta$ denotes the $\theta$-component of the streamfunction $\psi$. This reformulation has the advantage of removing the $1/\mathring{r}$ singularity from the cylindrical coordinates. We also define the radial component $u^r$ and axial component $u^z$ of the velocity in cylindrical coordinates, which can be recovered from $\psi_{\bullet}$ via 
\begin{equation*}
u^{r} = -\,\mathring{r}\,\partial_{\mathring{z}} \psi_{\bullet},
\qquad
u^{z} = 2\psi_{\bullet} + \mathring{r}\,\partial_{\mathring{r}} \psi_{\bullet}.
\end{equation*}

 \hfill

The axisymmetric 3D Euler equations are
\begin{align*}
\partial_t u_{\bullet}
+ u^r\,\partial_{\mathring{r}} u_{\bullet}
+ u^z\,\partial_{\mathring{z}} u_{\bullet}
&=
2u_{\bullet}\,\partial_{\mathring{z}}\psi_{\bullet}, \\
\partial_t \omega_{\bullet}
+ u^r\,\partial_{\mathring{r}} \omega_{\bullet}
+ u^z\,\partial_{\mathring{z}} \omega_{\bullet}
&=
2u_{\bullet}\,\partial_{\mathring{z}}u_{\bullet}, \\
-\mathcal E_{\mathring r,\mathring z}\psi_{\bullet}
&=
\omega_{\bullet}
\end{align*}
where \(
\mathcal E_{\mathring r,\mathring z}
=\partial_{\mathring r}^2+\frac{3}{\mathring r}\partial_{\mathring r}
+\partial_{\mathring z}^2
\).
\noindent The 3D axisymmetric Euler equations can be viewed as the non-viscous counterpart of the 3D axisymmetric Navier--Stokes equations. \\

\subsection{Traveling-Wave Ansatz} \label{sec: short_Ansatz}

\noindent
We consider a \emph{traveling self-similar} ansatz. \emph{Self-similar} refers to the fact that the geometry of the solution, after a proper rescaling of space, keeps an approximately fixed shape. \emph{Traveling} means that the center of this singularity is not fixed. Instead, it moves along the symmetry ($\mathring{r} = 0$) axis while the solution converges to the self-similar geometry.\\

The translation is introduced only in the axial ($\mathring{z}$) variable. This is consistent with the axisymmetric setting, where $\mathring r=0$ is the fixed symmetry axis. A radial ($\mathring{r}$) translation would move the singularity away from the symmetry axis and would generally destroy the parity and regularity structure imposed at $\mathring r=0$. By contrast, an axial translation in $\mathring z$ preserves the cylindrical symmetry and provides a natural mechanism for the singular structure to drift while it concentrates.\\

\paragraph{The traveling-wave ansatz.}
\noindent
Given the original axisymmetric variables as $(u_\bullet,\omega_\bullet,\psi_\bullet)$, we introduce the \emph{rescaled profile variables} $(u,\omega,\psi)$ via
\begin{align}
u_\bullet(\mathring r,\mathring z,t)
&= (T-t)^{c_u}\,
 u\!\left(\frac{\mathring r}{(T-t)^\lambda},\frac{\mathring z-\mathring z_c(t)}{(T-t)^\lambda}\right),
\label{ansantz1_shortV}\\
\omega_\bullet(\mathring r,\mathring z,t)
&= (T-t)^{c_\omega}\,
 \omega\!\left(\frac{\mathring r}{(T-t)^\lambda},\frac{\mathring z-\mathring z_c(t)}{(T-t)^\lambda}\right),
\label{ansantz2_shortV}\\
\psi_\bullet(\mathring r,\mathring z,t)
&= (T-t)^{c_\psi}\,
 \psi\!\left(\frac{\mathring r}{(T-t)^\lambda},\frac{\mathring z-\mathring z_c(t)}{(T-t)^\lambda}\right).
\label{ansantz3_shortV}
\end{align}
\noindent In these equations, $T$ is the \emph{blowup time}, $\lambda$ is the \emph{spatial blowup rate}, $\mathring z_c(t)$ is the center of the traveling profile, and $c_u,c_\omega,c_\psi$ are \emph{amplitude exponents}. \\

\noindent The rescaled coordinates are $$    r=\mathring r/(T-t)^\lambda,$$ and $$   z=(\mathring z-\mathring z_c(t))/(T-t)^\lambda,$$ and
\begin{equation}
    \frac{d}{dt}\mathring z_c(t)=-C(T-t)^{\lambda-1}.
    \label{eq:short_center_speed_shortV}
\end{equation}
Thus $C$ is the \emph{traveling speed} of the profile in the rescaled axial coordinate. If $C=0$, the singularity is centered at a fixed axial point in the rescaled frame. If $C\neq0$, the singularity drifts in the original variables while remaining stationary in the moving rescaled frame.\\

We are considering a \emph{single-scale traveling ansatz}. The same length scale $(T-t)^\lambda$ is used to measure the radial size of the profile, the axial size of the profile, and the axial displacement relative to the traveling center. The purpose is to isolate the traveling self-similar regime in which one spatial scale, one traveling speed, and the corresponding amplitude exponents determine the leading-order profile. \\

The exponents are fixed by balancing the powers of $T-t$ in the axisymmetric system, to enforce that time differentiation, transport, stretching, and the stream function relation contribute at the same order:
\begin{equation}
    c_u=-1,
    \qquad
    c_\omega=-1-\lambda,
    \qquad
    c_\psi=-1+\lambda.
    \label{eq:short_scaling_exponents}
\end{equation}
\noindent

\noindent The traveling-wave ansatz provides more flexibility for the optimization than the more common \emph{stationary self-similar} ansatz~\citep{deepMindPaper1,deepMindPaper2,wang2025high,HouJiaje}, in which one exploits scaling invariance and seeks a time-independent backward self-similar profile
\begin{equation*}
    u_\bullet(x,t) = (T-t)^{c_u} \, u\!\left(\frac{x - x_0}{ (T-t)^{\lambda}}\right).
\end{equation*}
\noindent
In our paper, the singularity is allowed to \emph{travel} along the symmetry axis. This introduces an additional axial center $\mathring z_c(t)$, and the profile is written in coordinates centered at $\mathring z_c(t)$.

\hfill

\subsection{Steady-State Profile Equations}
\label{sec: short_Steady-State Profile Equations}

We next present the profile equations generated by the traveling self-similar ansatz. This reduction transforms the original singular evolution into a time-independent nonlinear system for the rescaled variables $(U,\Omega,\Psi)$, together with the drift speed $C$.  \\

Assume $(u_\bullet,\omega_\bullet,\psi_\bullet)$ is given by the ansatz \eqref{ansantz1_shortV}--\eqref{ansantz3_shortV}, that the center satisfies \eqref{eq:short_center_speed_shortV}, that the exponents are chosen as in \eqref{eq:short_scaling_exponents}. Then $(U,\Omega,\Psi)$ satisfies the profile equations
\begin{align}
& U + \bigl(\lambda r+U^r\bigr)\partial_rU + \bigl(C+\lambda z+U^z\bigr)\partial_zU \label{eq:short_ProfileEq1_eps1_shortV}
\\ & \qquad\qquad\qquad\qquad\quad\qquad\quad\qquad =2U\partial_z\Psi, \nonumber\\
& (1+\lambda)\Omega \!+\! \bigl(\lambda r+U^r\bigr)\partial_r\Omega
\!+\! \bigl(C+\lambda z+U^z\bigr)\partial_z\Omega \label{eq:short_ProfileEq2_eps1_shortV} \\
& \qquad\qquad\qquad\qquad\quad\qquad\quad\qquad =2U\partial_zU,  \nonumber
\\
& \qquad \qquad \qquad\qquad  -\, \mathcal E \,\Psi=\Omega, &
\label{eq:short_ProfileEq3_eps1_shortV}
\end{align}
with $U^r=-r\partial_z\Psi$, 
$U^z=2\Psi+r\partial_r\Psi,$ and the elliptic operator $
    \mathcal{E} \coloneqq \partial_r^2 + \frac{3}{r}\partial_r + \partial_z^2$.

Equations \eqref{eq:short_ProfileEq1_eps1_shortV}--\eqref{eq:short_ProfileEq3_eps1_shortV} are the steady-state equations for the traveling self-similar profile. Thus, the traveling-wave ansatz reduces the study of a possible finite-time singularity to an autonomous nonlinear profile problem for $(U,\Omega,\Psi)$ and the parameter $C$. \

\section{Discovery of a Finite-Time Blowup via Physics-Informed Optimization} \label{sec: short_Finding the blowup}

\subsection{Physics-Informed Optimization}

Traditional numerical methods for PDE are typically developed and implemented on a case-by-case basis, with finite element methods among the most widely used approaches, but their computational cost can limit their practical applicability and they struggle to identify new self-similar blowup profiles. \\ 

Physics-informed optimization provides a more flexible and computationally practical alternative. 
Physics-Informed Neural Networks (PINNs)~\citep{PINN_OG} are neural networks designed to approximate solutions to PDE. Their parameters are obtained by minimizing deviations from the governing PDE in an appropriate norm, taking a finite set of collocation points in the spatiotemporal domain as input and producing a parametrized approximation to the solution at each grid point. The PDE residual loss is often combined with additional penalty terms that encode boundary and initial conditions, smoothness priors, and other problem-specific constraints, resulting in a composite objective used for gradient-based optimization to train the neural network.  \\

Despite their success, PINNs can be difficult to train due to poorly conditioned optimization problems and increasingly complex loss landscapes as the physics loss is emphasized \citep{Krishnapriyan2021}. They may also converge to trivial or undesired solutions \citep{Leiteritz2021}, while imbalances among loss terms can lead to uneven gradient contributions during training \citep{Wang2021}. Consequently, obtaining reliable PINN solutions requires careful design and tuning of the network architecture, sampling strategy, loss weights, hyperparameters, and optimization procedure. These challenges are particularly important in our setting, where making PINN solutions amenable to computer-assisted proofs requires very high numerical accuracy when solving the Euler equations.\\

To address these difficulties, we combine several complementary mechanisms designed to improve optimization, enforce the desired solution structure, and increase accuracy. In particular, \\

\begin{itemize}[leftmargin=1em, itemsep=-2pt, topsep=1.5pt]
    \item \textbf{Soft/Hard constraints} to guide the optimization and encode known properties of the profile, together with smoothness losses.\\
    \item \textbf{Adaptive loss weighting schedules}, allowing the relative importance of the various objectives to evolve during training and helping mitigate imbalances between competing loss terms.\\
    \item \textbf{64-bit floating point precision} (FP64) to reduce round-off errors in autograd and achieve the highest fidelity available in PyTorch, improving gradient quality, convergence, and accuracy.\\
    \item \textbf{PINN boosting}~\citep{fang2024boosting,boostingYongji}, where new networks are added sequentially after the current model reaches an optimization plateau, with each new network trained to correct the residual errors left by the previous ones. This turns one difficult global optimization problem into a sequence of simpler correction problems and allows later stages to use specialized architectures and localized corrections to target the remaining error more effectively.\\
    \item \textbf{Self-scaled curvature-aware optimizers}. We first use $\sssoap$~\citep{SS-SOAP} as an efficient curvature-aware optimizer, followed by SS-Broyden~\citep{Al-Baali1998,SSBroyden2} for high-accuracy refinement. SS-Broyden approximates second-order curvature information without explicitly computing the Hessian, and its adaptive scaling can improve conditioning and convergence in challenging nonlinear optimization problems. In our case, this optimization strategy reduces the MSE from $\sim10^{-4}$ with Adam~\citep{adamOptim} to $\sim10^{-10}$.\\ 
    \item \textbf{A $\mathbf{sinh}$ transformation} to efficiently sample on the unbounded domain.\\ 
    \item \textbf{Collocation points resampling}, concentrated near the origin and regions where the solution varies rapidly, and residual-based adaptive resampling to target points with the largest PDE errors.
\end{itemize}

\subsection{Numerical Results}

\begin{figure*}[htbp]
\centering

\includegraphics[height=9.6cm]{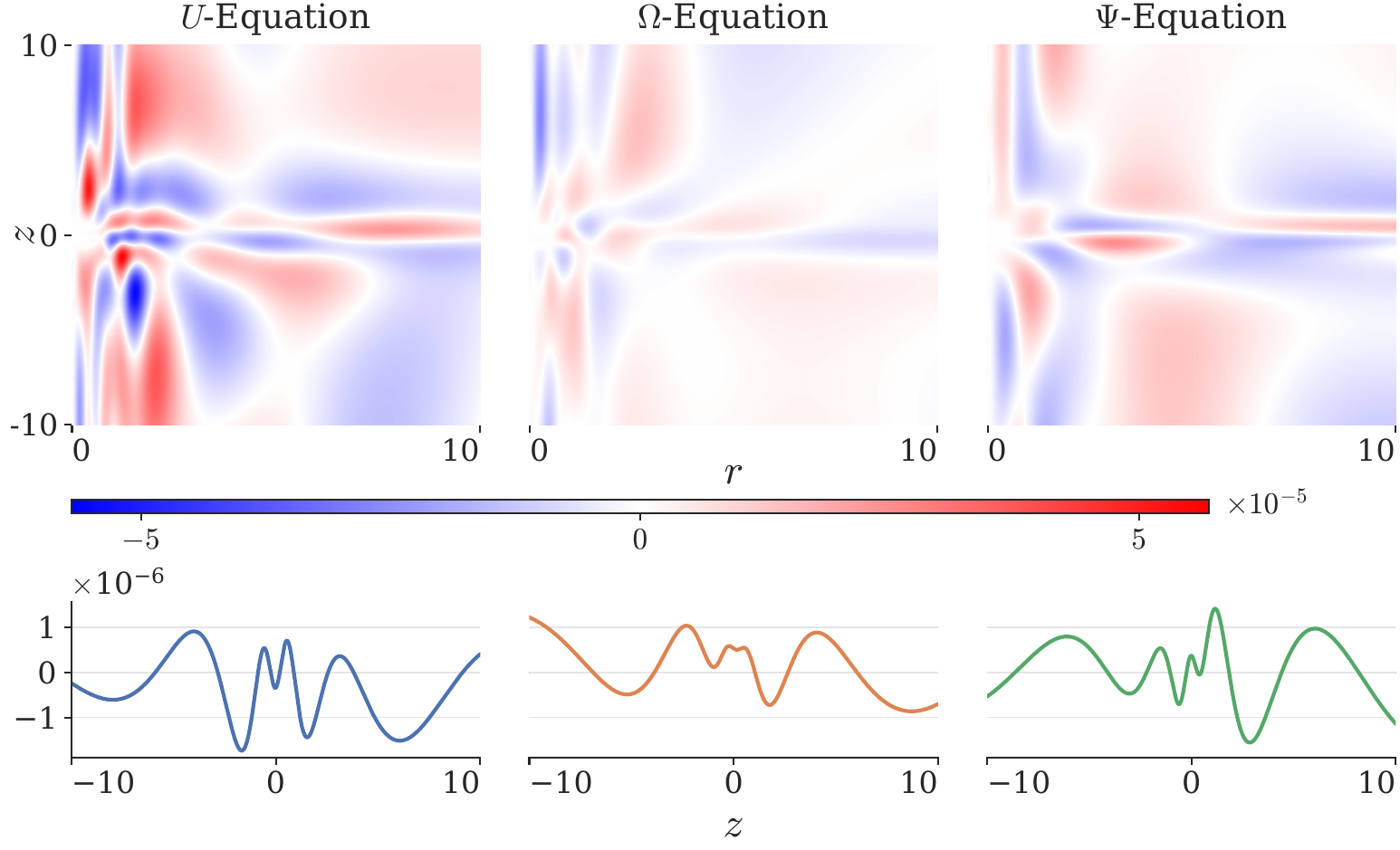}

\caption{\textbf{(Top)} Profile equation residuals on the box $(r,z)\in[0,10]\times[-10,10]$. \textbf{(Bottom)} Profile equation residuals on the symmetry-axis ($r=0$) with $z \in [-10,10]$.}
\label{fig:short_residuals} \vspace{7mm}

\end{figure*}

\noindent We optimize PINNs for the Euler system and visualize the resulting approximate profiles $U$,~$\Omega$,~$\Psi$ in Figure~\ref{fig: short_profiles}, which appear consistent with the profile geometry reported by~\citet{pengfeiPaper} in the low convection setting. The approximate profiles are obtained to high accuracy, as indicated by \Cref{tab:short_table_errors} and displayed in \Cref{fig:short_residuals}.

\begin{table*}[htbp]
\centering

\caption{Residual errors for the $U$, $\partial_r U$, $\partial_z U$, $\Omega$, and $\Psi$ profile equations for the Euler system, comparing PINN and spline representations. Note that the residuals are estimated empirically on an independent set of $100{,}000{,}000$ randomly sampled collocation points for the PINN representations, whereas they are evaluated analytically and exactly for the spline representations. in particular, cRMSE refers to the continuous root-mean squared error $ \sqrt{
\frac{1}{\left|\mathbb{D}\right|}
\int_{\mathbb{D}}
\mathsf{residual}\!\left(r_{\!\mathfrak{comp}}, z_{\mathfrak{comp}}\right)^2
\,\mathrm{d}r_{\mathfrak{comp}}\,\mathrm{d}z_{\mathfrak{comp}}
} $. Here, $\mathbb{D}$ denotes the full computational domain (corresponding to the half plane), $B_{0.1}(0^{\star})$, $B_1(0^{\star})$, and $B_{10}(0^{\star})$ denote balls of radii $0.1$, $1$, and $10$, respectively, centered at the origin $0^{\star}$ of the stability analysis (i.e. the on-axis meridional fixed point), $\mathsf{axis}$ denotes the $r=0$ axis, and $\mathsf{FF}$ denotes the far-field region $\left(r_{\!\mathfrak{comp}}, z_{\mathfrak{comp}}\right) \in [10, 30] \times [-30, -10] \cup [10, 30] \times [10, 30] $.}
\label{tab:short_table_errors}

\vspace{0.7mm}
\renewcommand{\arraystretch}{1.04}

\begin{tabularx}{0.9\linewidth}{
@{}
c@{\hspace{2mm}}
l
*{4}{>{\centering\arraybackslash}X}
@{}
}

&
&
\multicolumn{2}{c}{\textbf{PINNs}}
&
\multicolumn{2}{c}{\textbf{Splines}}
\\[-1.1mm]

\cmidrule(lr){3-4}
\cmidrule(lr){5-6}

&
\textbf{Domain}
& RMSE
& $\ell^\infty$
& cRMSE
& $L^\infty$
\\

\midrule
\addlinespace[0.7mm]


\multirow{6}{*}{
\rotatebox[origin=c]{90}{\textbf{\small $U$-equation}}
}
&
$B_{0.1}(0^{\star})$
& $2.8\!\cdot\!10^{-7}$
& $3.4\!\cdot\!10^{-7}$
& $2.8\!\cdot\!10^{-7}$
& $3.5\!\cdot\!10^{-7}$
\\

&
$B_1(0^\star)$
& $7.0\!\cdot\!10^{-6}$
& $3.3\!\cdot\!10^{-5}$
& $6.7\!\cdot\!10^{-6}$
& $3.4\!\cdot\!10^{-5}$
\\

&
$B_{10}(0^{\star})$
& $1.2\!\cdot\!10^{-5}$
& $5.7\!\cdot\!10^{-5}$
& $1.7\!\cdot\!10^{-5}$
& $5.8\!\cdot\!10^{-5}$
\\

&
$\mathsf{axis}$
& $9.0\!\cdot\!10^{-5}$
& $1.4\!\cdot\!10^{-3}$
& $9.1\!\cdot\!10^{-5}$
& $1.4\!\cdot\!10^{-3}$
\\

&
$\mathsf{FF}$
& $6.2\!\cdot\!10^{-6}$
& $2.2\!\cdot\!10^{-5}$
& $6.3\!\cdot\!10^{-6}$
& $2.3\!\cdot\!10^{-5}$
\\

&
\cellcolor{grayshade}$\mathbb{D}$
& \cellcolor{grayshade}$4.2\!\cdot\!10^{-5}$
& \cellcolor{grayshade}$3.6\!\cdot\!10^{-3}$
& \cellcolor{grayshade}$5.7\!\cdot\!10^{-5}$
& \cellcolor{grayshade}$3.7\!\cdot\!10^{-3}$
\\

\addlinespace[1.0mm]
\midrule
\addlinespace[1.0mm]


\multirow{6}{*}{
\rotatebox[origin=c]{90}{\textbf{\small $\partial_r U$-equation}}
}
&
\rule[-0.55ex]{0pt}{3.15ex}$B_{0.1}(0^{\star})$
& $1.5\!\cdot\!10^{-6}$
& $3.3\!\cdot\!10^{-6}$
& $1.6\!\cdot\!10^{-6}$
& $3.4\!\cdot\!10^{-6}$
\\

&
\rule[-0.55ex]{0pt}{3.15ex}$B_1(0^\star)$
& $5.2\!\cdot\!10^{-5}$
& $1.9\!\cdot\!10^{-4}$
& $5.0\!\cdot\!10^{-5}$
& $1.9\!\cdot\!10^{-4}$
\\

&
\rule[-0.55ex]{0pt}{3.15ex}$B_{10}(0^{\star})$
& $4.2\!\cdot\!10^{-5}$
& $3.1\!\cdot\!10^{-4}$
& $6.9\!\cdot\!10^{-5}$
& $3.2\!\cdot\!10^{-4}$
\\

&
\rule[-0.55ex]{0pt}{3.15ex}$\mathsf{axis}$
& $0$
& $0$
& $0$
& $0$
\\

&
\rule[-0.55ex]{0pt}{3.15ex}$\mathsf{FF}$
& $1.3\!\cdot\!10^{-10}$
& $2.4\!\cdot\!10^{-9}$
& $1.3\!\cdot\!10^{-10}$
& $2.5\!\cdot\!10^{-9}$
\\

&
\cellcolor{grayshade}\rule[-0.55ex]{0pt}{3.15ex}$\mathbb{D}$
& \cellcolor{grayshade}$1.2\!\cdot\!10^{-4}$
& \cellcolor{grayshade}$8.3\!\cdot\!10^{-3}$
& \cellcolor{grayshade}$1.5\!\cdot\!10^{-4}$
& \cellcolor{grayshade}$8.5\!\cdot\!10^{-3}$
\\

\addlinespace[1.0mm]
\midrule
\addlinespace[1.0mm]


\multirow{6}{*}{
\rotatebox[origin=c]{90}{\textbf{\small $\partial_z U$-equation}}
}
&
\rule[-0.55ex]{0pt}{3.15ex}$B_{0.1}(0^{\star})$
& $1.1\!\cdot\!10^{-6}$
& $1.9\!\cdot\!10^{-6}$
& $1.1\!\cdot\!10^{-6}$
& $2.3\!\cdot\!10^{-6}$
\\

&
\rule[-0.55ex]{0pt}{3.15ex}$B_1(0^\star)$
& $2.3\!\cdot\!10^{-5}$
& $1.0\!\cdot\!10^{-4}$
& $2.3\!\cdot\!10^{-5}$
& $1.1\!\cdot\!10^{-4}$
\\

&
\rule[-0.55ex]{0pt}{3.15ex}$B_{10}(0^{\star})$
& $1.3\!\cdot\!10^{-5}$
& $1.3\!\cdot\!10^{-4}$
& $2.8\!\cdot\!10^{-5}$
& $1.3\!\cdot\!10^{-4}$
\\

&
\rule[-0.55ex]{0pt}{3.15ex}$\mathsf{axis}$
& $1.1\!\cdot\!10^{-6}$
& $1.7\!\cdot\!10^{-5}$
& $1.2\!\cdot\!10^{-6}$
& $1.8\!\cdot\!10^{-5}$
\\

&
\rule[-0.55ex]{0pt}{3.15ex}$\mathsf{FF}$
& $1.5\!\cdot\!10^{-11}$
& $1.8\!\cdot\!10^{-10}$
& $1.6\!\cdot\!10^{-11}$
& $2.7\!\cdot\!10^{-10}$
\\

&
\cellcolor{grayshade}\rule[-0.55ex]{0pt}{3.15ex}$\mathbb{D}$
& \cellcolor{grayshade}$2.9\!\cdot\!10^{-6}$
& \cellcolor{grayshade}$1.3\!\cdot\!10^{-4}$
& \cellcolor{grayshade}$4.7\!\cdot\!10^{-6}$
& \cellcolor{grayshade}$1.3\!\cdot\!10^{-4}$
\\

\addlinespace[1.0mm]
\midrule
\addlinespace[1.0mm]


\multirow{6}{*}{
\rotatebox[origin=c]{90}{\textbf{\small $\Omega$-equation}}
}
&
$B_{0.1}(0^{\star})$
& $5.0\!\cdot\!10^{-7}$
& $5.8\!\cdot\!10^{-7}$
& $5.1\!\cdot\!10^{-7}$
& $5.9\!\cdot\!10^{-7}$
\\

&
$B_1(0^\star)$
& $3.4\!\cdot\!10^{-6}$
& $1.1\!\cdot\!10^{-5}$
& $3.3\!\cdot\!10^{-6}$
& $1.1\!\cdot\!10^{-5}$
\\

&
$B_{10}(0^{\star})$
& $4.4\!\cdot\!10^{-6}$
& $2.7\!\cdot\!10^{-5}$
& $5.2\!\cdot\!10^{-6}$
& $2.8\!\cdot\!10^{-5}$
\\

&
$\mathsf{axis}$
& $4.8\!\cdot\!10^{-5}$
& $3.7\!\cdot\!10^{-4}$
& $4.9\!\cdot\!10^{-5}$
& $3.8\!\cdot\!10^{-4}$
\\

&
$\mathsf{FF}$
& $3.1\!\cdot\!10^{-6}$
& $1.4\!\cdot\!10^{-5}$
& $3.2\!\cdot\!10^{-6}$
& $1.4\!\cdot\!10^{-5}$
\\

&
\cellcolor{grayshade}$\mathbb{D}$
& \cellcolor{grayshade}$1.7\!\cdot\!10^{-5}$
& \cellcolor{grayshade}$9.5\!\cdot\!10^{-4}$
& \cellcolor{grayshade}$1.9\!\cdot\!10^{-5}$
& \cellcolor{grayshade}$9.6\!\cdot\!10^{-4}$
\\

\addlinespace[1.0mm]
\midrule
\addlinespace[1.0mm]


\multirow{6}{*}{
\rotatebox[origin=c]{90}{\textbf{\small $\Psi$-equation}}
}
&
$B_{0.1}(0^{\star})$
& $3.6\!\cdot\!10^{-7}$
& $1.2\!\cdot\!10^{-6}$
& $3.7\!\cdot\!10^{-7}$
& $4.2\!\cdot\!10^{-7}$
\\

&
$B_1(0^\star)$
& $2.8\!\cdot\!10^{-6}$
& $8.9\!\cdot\!10^{-6}$
& $2.7\!\cdot\!10^{-6}$
& $9.0\!\cdot\!10^{-6}$
\\

&
$B_{10}(0^{\star})$
& $6.9\!\cdot\!10^{-6}$
& $2.6\!\cdot\!10^{-5}$
& $8.0\!\cdot\!10^{-6}$
& $2.7\!\cdot\!10^{-5}$
\\

&
$\mathsf{axis}$
& $8.5\!\cdot\!10^{-5}$
& $1.4\!\cdot\!10^{-3}$
& $8.6\!\cdot\!10^{-5}$
& $1.5\!\cdot\!10^{-3}$
\\

&
$\mathsf{FF}$
& $2.1\!\cdot\!10^{-6}$
& $8.3\!\cdot\!10^{-6}$
& $2.1\!\cdot\!10^{-6}$
& $8.3\!\cdot\!10^{-6}$
\\

&
\cellcolor{grayshade}$\mathbb{D}$
& \cellcolor{grayshade}$3.8\!\cdot\!10^{-5}$
& \cellcolor{grayshade}$1.4\!\cdot\!10^{-3}$
& \cellcolor{grayshade}$5.3\!\cdot\!10^{-5}$
& \cellcolor{grayshade}$1.5\!\cdot\!10^{-3}$
\\

\addlinespace[0.7mm]
\bottomrule

\end{tabularx}

\vspace{1mm}
\end{table*}

The PINN serves as a numerical discovery tool that allows us to locate an approximate self-similar profile with high accuracy in a practical and well-conditioned way, especially when direct optimization over analytic ansatz would be difficult or unstable. For certification, we replace the computed profile by a fixed piecewise polynomial spline, which provides an explicit analytic representation on which derivatives and $L^\infty$ bounds can be evaluated exactly. This avoids the more challenging task of certifying network evaluation, automatic differentiation, and $L^\infty$ bounds. We choose splines because their low-degree local structure supports simple and exact differentiation and residual estimates while remaining sufficiently expressive to accurately capture the computed profile.     Their locality also allows the fit to the PINN to be refined selectively in more important regions, e.g. near the profile center and the axis. The spline profile retains the small residuals of the PINN, with only a minor loss of accuracy.

\section{Stability}

A singularity is \emph{stable} if the associated blowup mechanism persists under small perturbations of the initial data, in the sense that the same mechanism occurs for an open neighborhood of nearby initial conditions. Conversely, an \emph{unstable} singularity requires infinitely precise initial conditions: small perturbations deflect the evolution away from the singularity mechanism.  

The dynamic rescaling formulation presented in \Cref{sec: short_Dynamic Rescaling} is the natural starting point for such stability analysis. In the original physical variables, the solution may grow and concentrate as $t\to T^-$, so there is no fixed object around which to linearize. After rescaling, the candidate blowup profile becomes a steady state of the rescaled equations, while the blowup rate, translation rate, and amplitude normalization are recorded by modulation parameters. Linear stability can therefore be studied by perturbing this steady rescaled profile and analyzing the resulting linearized rescaled dynamics.  

We follow the framework of Chen et al.~\citep{DGBlowup,HouJiaje,Chen2026Analysis} and related works on stability of blowups. We first compute the blowup profile numerically, then develop a framework for establishing linear stability of the dynamically rescaled equation in suitable weighted norms and for proving nonlinear stability using Sobolev embeddings and additional weighted estimates. The main distinction is that the weight functions and other tunable parameters entering the stability proof are optimized numerically to seek sufficient margin for the argument to close, rather than being chosen and refined through a potentially long and highly nontrivial process of trial and error.

\subsection{Dynamic Rescaling Equations}
\label{sec: short_Dynamic Rescaling}

In the dynamic formulation, the fixed self-similar parameters in the profile ansatz are replaced by time-dependent \emph{modulation parameters}. Rather than prescribing a blowup time, a spatial scale, and a traveling center in advance, it defines a long-time evolution in normalized variables, with the singularity kept at order-one size and near a fixed location in rescaled coordinates. This avoids the need to resolve a structure whose amplitude diverges and whose length scale collapses in the physical variables, while allowing the computation to identify am approximate self-similar profile and its associated scaling parameters. 

More precisely, let $(\mathring r,\mathring z,\mathring t)$ and $(r,z,t)$  denote the physical and dynamically rescaled variables. We introduce a common spatial scale $\mathrm{s}_r(t)>0$, amplitude scales $\mathrm{s}_u(t), \, \mathrm{s}_\omega(t)>0$, and an axial center $\mathring z_c(t)$ by \begin{align*} \mathring r=\mathrm{s}_r(t) \, r, \quad & \qquad \mathring z-\mathring z_c(t)=\mathrm{s}_r(t) \, z, \\ u_\bullet(\mathring r,\mathring z,\mathring t) & =\mathrm{s}_u(t) \, u(r,z,t), \\ \omega_\bullet(\mathring r,\mathring z,\mathring t) & =\mathrm{s}_\omega(t) \, \omega(r,z,t). \end{align*} 
We then define the rescaling rates and axial drift by \begin{align*}  &\partial_t\mathrm{s}_r(t)=-\lambda \,\mathrm{s}_r(t),  \quad \partial_t\mathrm{s}_u(t)=-c_u(t) \, \mathrm{s}_u(t), \\ & \qquad \qquad \quad 
\partial_t\mathrm{s}_\omega(t)=-c_\omega(t) \, \mathrm{s}_\omega(t) \end{align*} 
and $\partial_t\mathring z_c(t)=-C(t) \, \mathrm{s}_r(t)$. Overall, this is encoded by the coordinate drift $\left(\partial_t   r, \, \partial_t z\right)_{\mathring r,\mathring z}= (\lambda r, \lambda z+C(t) )$ where $C(t)$ is the \emph{instantaneous traveling speed} of the axial center in the same rescaled frame. The scalar functions $c_u(t)$ and $c_\omega(t)$, governing the amplitude scaling of $u$ and $\omega$, maintain an order-one amplitude in the rescaled variables. \\ 

 In the original variables, perturbations can change the blowup time, shift the center, or alter the amplitude normalization, producing apparent growth or drift even when the profile is stable. The normalization conditions remove these neutral directions by recentering, rescaling, and renormalizing the solution, so stability can be assessed directly in the rescaled variables. If the dynamic rescaling captures a stable self-similar blowup, the rescaled solution converges to a steady profile and the modulation parameters approach constants. The time derivatives then vanish, and for the traveling-wave ansatz the limiting amplitude rates are $c_u=-1$ and $c_\omega=-(1+\lambda)=-3/2$. Substituting these values recovers the profile equations. Thus, self-similar profiles are steady states of the dynamically rescaled evolution, with finite-time blowup in physical variables corresponding to convergence toward a stable steady state.

\hfill

\subsection{Linearization}

\noindent We derive the perturbation equations around a steady state of the dynamically-rescaled system, and separate linear and nonlinear terms in the perturbations.

\noindent We denote perturbations via \textcolor{blue}{$\delta$}'s in \textcolor{blue}{blue}, and use the subscript $_0$ to denote evaluation at the origin. We perturb the approximate steady state $(    \bar u,\bar\omega,\bar\psi,\bar C,\bar c_u,\bar c_\omega)$ of the dynamic rescaling equations, by writing
\begin{equation*}
    u = \bar u + \textcolor{blue}{\delta u},
    \quad \ \ 
    \omega=\bar\omega+\textcolor{blue}{\delta\omega},
    \quad \ \ 
    \psi = \bar\psi + \textcolor{blue}{\delta\psi},
\end{equation*}
\begin{equation*}
    C=\bar C+\textcolor{blue}{\delta C},
    \quad \ \ 
    c_u = \bar c_u + \textcolor{blue}{\delta c_u},
    \quad \ \ 
    c_\omega=\bar c_\omega+\textcolor{blue}{\delta c_\omega}.
\end{equation*}
\noindent For notational convenience, we will use the symbol $\textcolor{blue}{\delta}$ as shorthand for the full list of perturbation variables, $$
    \textcolor{blue}{\delta}
    :=
    \{ \textcolor{blue}{\delta u},\textcolor{blue}{\delta\omega},\textcolor{blue}{\delta\psi},\textcolor{blue}{\delta C},
    \textcolor{blue}{\delta c_u},\textcolor{blue}{\delta c_\omega} \}.$$
The perturbation variables measure the difference between the evolving rescaled solution and the steady blowup profile. The perturbations must inherit the symmetries of the underlying steady state and belong to the same functional class. In particular, the parity, far-field decay, and regularity conditions imposed on the rescaled variables are inherited by $\textcolor{blue}{\delta u}$, $\textcolor{blue}{\delta\omega}$, and $\textcolor{blue}{\delta\psi}$, with the latter two coupled by the elliptic equation relating vorticity and streamfunction. \\

The reduced variables $u$ and $\omega$ do not transform in the same way under rescaling. Instead, $\omega$ and $\nabla u$ are the natural pair of quantities in a scale-consistent linear stability analysis. \\

\noindent We collect these perturbations into the scale-consistent perturbation vector
\begin{equation*}
\textcolor{blue}{\delta v}
:=
(\textcolor{blue}{\delta\omega},\nabla\textcolor{blue}{\delta u})
=
(\textcolor{blue}{\delta\omega},\partial_r\textcolor{blue}{\delta u},\partial_z\textcolor{blue}{\delta u}),
\end{equation*}
\begin{equation*}
\textcolor{blue}{\delta v_\omega}:=\textcolor{blue}{\delta\omega},
\quad
\textcolor{blue}{\delta v_r}:=\partial_r\textcolor{blue}{\delta u},
\quad
\textcolor{blue}{\delta v_z}:=\partial_z\textcolor{blue}{\delta u}.
\end{equation*}
Linearization consists of expanding the perturbation equation around the approximate profile and retaining only linear terms in the perturbation. The linearized equation then gives the leading-order dynamics of small perturbations near the profile. The discarded terms record the nonlinear interactions among the perturbations. \\   

 Starting from the dynamic rescaling equations, perturbing all variables, subtracting the steady-state equations, and separating first-order terms from higher-order terms gives
\begin{align*}
    \partial_t\textcolor{blue}{\delta u}
    & =
    \mathcal L_u(\textcolor{blue}{\delta}) + \mathcal N_u(\textcolor{blue}{\delta}) \ + \ \mathcal R_u,
    \\ 
    \partial_t\textcolor{blue}{\delta\omega}
    &=
    \mathcal L_\omega(\textcolor{blue}{\delta}) + \mathcal N_\omega(\textcolor{blue}{\delta}) \ + \ \mathcal R_\omega.
\end{align*}
Here $\mathcal L_u$ and $\mathcal L_\omega$ are linear in the perturbation variables, while $\mathcal N_u$ and $\mathcal N_\omega$ contain the nonlinear products of perturbations, and $\mathcal R_u$ and $\mathcal R_\omega$ are the PDE residuals accounting for the fact that the numerical profiles are not exact. 

Since the scale-consistent variables are $\omega$ and $\nabla u$, we differentiate the $\delta u$-equation:
\begin{align*}
\partial_t\partial_r\textcolor{blue}{\delta u}
&=\mathcal L_r(\textcolor{blue}{\delta})+
\mathcal N_{r}(\textcolor{blue}{\delta}) \ +\ \mathcal R_r, \\ 
\partial_t\partial_z\textcolor{blue}{\delta u}
& =\mathcal L_z(\textcolor{blue}{\delta})+
\mathcal N_{z}(\textcolor{blue}{\delta})  \ +\ \mathcal R_z,
\end{align*}
where we write $\mathcal L_r:=\partial_r\mathcal L_u$, $\mathcal N_r:=\partial_r\mathcal N_u$, $\mathcal L_z:=\partial_z\mathcal L_u$, $\mathcal N_z:=\partial_z\mathcal N_u$,
$
    \mathcal R_z:=\partial_z\mathcal R_u$, $\mathcal R_z:=\partial_z\mathcal R_u$. \\

\noindent For stability, we use a two-head fit for $U$ and $\Psi$, while $\Omega$ is obtained by applying the elliptic operator exactly to $\Psi$. Hence, the elliptic residual vanishes by construction, removing one source of error and simplifying substantially the stability analysis. \\  

The modulation parameters $
    \textcolor{blue}{\delta C},
    \textcolor{blue}{\delta c_u},\textcolor{blue}{\delta c_\omega}$ are fixed by normalization conditions, which express them in terms of the profile and its perturbations evaluated at the origin
\begin{equation*}
     C(t)+u^z_0(t)=0,
    \qquad
    u_0(t)=\bar u_0.
\end{equation*}
For the perturbations, these conditions give
\begin{equation*} 
    \textcolor{blue}{\delta C}+2\textcolor{blue}{\delta\psi}_0=0,
    \qquad
    \textcolor{blue}{\delta u}_0 = 0,
    \qquad
    (\partial_t\textcolor{blue}{\delta u})_0=0.
\end{equation*}

These normalization conditions remove the artificial freedom associated with axial translation, spatial rescaling, amplitude normalization. They are essential for formulating linear stability in the modulated variables, since otherwise the perturbation could drift along symmetry directions rather than measuring a genuine change in the profile.\\

\subsection{Weighted Norms and Energy}

The stability estimates are measured in weighted norms adapted to the blowup profile and damping structure of the linearized operator. We define positive weight functions $$\wgray{\Phi_\omega},\wgray{\Phi_r},\wgray{\Phi_z}
:
\mathbb{D} \to (0,\infty)$$ on the  half-plane $$\mathbb{D}=\{(r,z):r\ge0,\ z\in\mathbb R\}.$$  The purpose of the weighted energy is to produce damping from selected linear terms, such as the linear transport terms. We choose weights that are singular at the origin since they strengthen the linear damping mechanism at the origin where the transport terms vanish, while still allowing the weighted energy to remain finite for the admissible perturbation class.  \\ 

More explicitly, $\wgray{\Phi_\omega}$ weights the vorticity perturbation $\textcolor{blue}{\delta\omega}$, $\wgray{\Phi_r}$ weights $\partial_r\textcolor{blue}{\delta u}$, and $\wgray{\Phi_z}$ weights $\partial_z\textcolor{blue}{\delta u}$. We denote the collection of weights by $\wgray{\Phi}
:=
(\wgray{\Phi_\omega},\wgray{\Phi_r},\wgray{\Phi_z}).$  For a positive scalar weight $\wgray{\phi}$, define the associated inner product and norm by $$\langle f,g\rangle_{\wgray{\phi}} := \int_\mathbb{D} fg\,\wgray{\phi}\,dr\,dz$$ and $\|f\|_{\wgray{\phi}}^2 := \langle f,f\rangle_{\wgray{\phi}}$. Given vector quantities $q := (q_\omega,q_r,q_z)$, and   $\tilde q := (\tilde q_\omega,\tilde q_r,\tilde q_z)$, their $\wgray{\Phi}$-weighted inner product is defined as 
\begin{equation*}
\langle q,\tilde q\rangle_{\wgray{\Phi}} := \langle q_\omega,\tilde q_\omega\rangle_{\wgray{\Phi_\omega}} + \langle q_r,\tilde q_r\rangle_{\wgray{\Phi_r}} + \langle q_z,\tilde q_z\rangle_{\wgray{\Phi_z}}. \end{equation*}

\noindent The low-order energy is defined as the $\wgray{\Phi}$-norm of $\textcolor{blue}{\delta v}$:
\begin{equation*}
  \mathfrak{E}_0(t)^2 = \sum_{i\in \{\omega,r,z\}}\|\textcolor{blue}{\delta v_i}(t)\|_{\wgray{\Phi_i}}^2.
\end{equation*}

\noindent The three weights $\wgray{\Phi_\omega}$, $\wgray{\Phi_r}$, $\wgray{\Phi_z}$ are kept separate because the damping calculation acts on the three components of $\textcolor{blue}{\delta v}$ in different ways.  \\ 

While the low-order energy remains the main source of damping, it is not sufficient by itself to close the estimates since the linear and nonlinear terms contain terms which cannot be controlled using a weighted $L^2$ norm, such as terms involving $L^\infty$ norms or pointwise quantities evaluated at the origin.  To control these terms, we supplement $\mathfrak E_0$ with a high-order weighted energy $\mathfrak H_k$, chosen so that the top derivatives provide the required pointwise, interpolation, elliptic, and modulation bounds.
\\

The high-order part uses the adjusted weights
\begin{equation*} 
\wgray{\Phi_{i,k}}:=1+\rho^{\,k}\wgray{\Phi_i},
\quad \rho(r,z):=r^2+z^2,
\quad i\in\{\omega,r,z\}. 
\end{equation*}
\noindent We then define the \emph{high-order weighted energy}
\begin{equation}
\label{eq:short_ns-Hk}
\mathfrak H_k^2
:={}\sum_{i\in \{\omega,r,z \}}\sum_{|\alpha|=k}\|D^\alpha\bvar{\delta v_i}\|_{\wgray{\Phi_{i,k}}}^2. 
\end{equation}

\hfill 

The full energy used for the stability argument is finally obtained as
\begin{equation*} 
\mathfrak E_k^2:=\mathfrak E_0^2+\mu_k\mathfrak H_k^2,
\qquad 0<\mu_k\le 1,
\end{equation*}
where the parameter $\mu_k$ is tuned to allow both the low-order and high-order damping to absorb the low-order and high-order losses in the linear estimates. \\ 

The order of the top derivative $k$ should be chosen large enough so that the energy controls all point values, products, commutators, and streamfunction terms used in the analytic estimates. With our proof construction, the lowest admissible choice is $k=4$. Larger choices of $k$ are possible, but increase combinatorially the number of derivatives and analytic constants that must be certified.  \\

No intermediate derivatives are needed in the argument: every derivative level $0<j<k$ is controlled by the low-order piece $\mathfrak E_0^2$, the top-order piece $\mathfrak H_k^2$, and explicit interpolation constants $I_{i,j,k}(\eta)$:
\begin{equation*}
\|D^j f\|_{\wgray{\Phi_{i,j}}}
\le \eta\|D^k f\|_{\wgray{\Phi_{i,k}}}+I_{i,j,k}(\eta)\|f\|_{\wgray{\Phi_i}}.
\end{equation*}

\hfill

\subsection{Stability Estimates}

Linear stability asks whether perturbations of solutions of the linearized system remain bounded in a suitable norm. The desired estimate is
\begin{equation*}
\|\textcolor{blue}{\delta v}(t)\|_{\wgray{\Phi}} \leq C e^{-\gamma t} \|\textcolor{blue}{\delta v}(0)\|_{\wgray{\Phi}}, \qquad \gamma>0.
\end{equation*}

\noindent A standard rigorous way to prove this is to establish the differential energy estimate
\begin{equation*}
\frac{1}{2}\frac{d}{dt}\mathfrak{E}_0(t)^2
=
\left\langle
\mathcal L\textcolor{blue}{\delta v},
\textcolor{blue}{\delta v}
\right\rangle_{\wgray{\Phi}}
\leq
-\DLamL \mathfrak{E}_0(t)^2.
\end{equation*}

\noindent For a blowup profile, this is the relevant notion of stability: the original physical solution may still blow up, but the normalized spatial profile is stable in the rescaled variables, i.e. the perturbations converge to zero under the linearized rescaled flow. \\ 

However, the low-order energy is not sufficient to close the linear estimates since the linear operators contain terms which cannot be controlled using a weighted $L^2$ norm. It is thus supplemented with the high-order energy $\mathfrak H_k$, to get an estimate of the form
\begin{equation*}  
\frac{1}{2}\frac{d}{dt}\mathfrak{E}_0(t)^2
\leq -\DLamLlow\mathfrak E_0^2+\DLamLhigh\mathfrak H_k^2.
\end{equation*}

The nonlinear stability argument upgrades this statement to the full perturbation equation, to prove that, if the initial perturbation is sufficiently small in a suitable energy, then the perturbation remains small for as long as the estimates apply. In the case of an approximate numerical profile, the solution is allowed to remain in a small neighborhood whose size also depends on the certified residual error.

\hfill

\subsection{Stability Theorem}

\noindent The stability argument is reduced to an $\mathfrak E_k$ energy estimate. Differentiating $\mathfrak E_k^2$ in time and inserting the perturbation equations separates the argument into distinct components: 
\begin{itemize}[leftmargin=1em, itemsep=-2pt, topsep=1.5pt]
    \item \emph{the residual error}, which records the defect of the approximate steady profile
    \item \emph{low- and high-order linear parts} which contain the stabilizing mechanism: the low-order estimate gives the main source of damping, while the high-order estimate controls the differentiated equations, allowing only losses that can be absorbed by the low-order energy after choosing $\mu_k$.
    \item \emph{low- and high-order nonlinear parts}, bounded by cubic powers of the full energy $\mathfrak E_k$. 
\end{itemize}

\vspace{3mm}

\noindent These estimates culminate in the stability theorem below, with the understanding that every certification constant entering the analytic bounds must be rigorously validated for the proof to be complete. The detailed estimates and derivations are given in our separate paper~\citep{EulerBlowupStability}.

\begin{theorem}[Euler Stability]
Assume that the approximate profile has a sufficiently small residual and that the certified low- and high-order linear estimates can be combined, after a suitable choice of the high-order weight $\mu_k$, to produce strict damping of the full energy $\mathfrak E_k$. Assume further that the nonlinear terms contribute only cubic powers of $\mathfrak E_k$. \\

\noindent These estimates yield an inequality of the form
\begin{equation*}
\frac{d}{dt}\mathfrak E_k^2
\le 
-\DLamStab\mathfrak E_k^2
+\DLamThree\mathfrak E_k^3
+\DLamR \mathfrak E_k,
\end{equation*}
where $\DLamStab>0$ represents the net linear damping, $\DLamThree$ controls the nonlinear terms, and $\DLamR$ measures the residual.\\

\noindent If there exists $\delta_*>0$ such that
\begin{equation}
\DLamThree\delta_*
+\frac{\DLamR}{\delta_*}
<
\DLamStab, \label{eq:short_ns-smallness_shortV}
\end{equation}
then the damping dominates both the nonlinear effects and the residual forcing on the boundary $\mathfrak E_k=\delta_*$. \\

\noindent Consequently, the ball $\{\mathfrak E_k<\delta_* \}$
is forward invariant: every solution satisfying $
\mathfrak E_k(0)<\delta_*$
continues to satisfy $\mathfrak E_k(t)<\delta_*$ for as long as the solution exists and the certified estimates remain valid. Thus, sufficiently small perturbations remain uniformly controlled in the full energy $\mathfrak E_k$.
\end{theorem}

\noindent The single-radius certificate can be relaxed by prescribing separate target bounds for the low-order and top-order energies, i.e. $\mathfrak E_0(t)<\delta_0$ and $\mathfrak H_k(t)<\delta_k$.

\hfill 

\subsection{Rescaled Stability to Physical Blowup}

\noindent
The stability analysis is carried out in dynamically rescaled variables, so a reconstruction is needed to recover the physical behavior. If an exact admissible rescaled solution exists for all $t\ge0$ and remains in the stability regime, then the amplitude scale $\mathrm{s}_u(t)$ grows exponentially whenever $c_u(t)$ stays negative. Since the physical clock satisfies $d\mathring t /dt=1/\mathrm{s}_u(t)$, this growth makes the total remaining physical time finite. At the same time, the fixed normalization at the rescaled origin implies that the physical axial vorticity at the moving center grows proportionally to $\mathrm{s}_u(t)$ and thus becomes unbounded. The center itself converges to a finite physical location. Hence the rescaled stability regime reconstructs to a finite-time physical singularity.\\

This reconstruction step is conditional on the existence of the controlled rescaled trajectory for all $t\ge0$, but does not require convergence to the numerical profile. An unconditional blowup result would additionally require a local well-posedness and continuation argument ruling out breakdown at finite rescaled time.

\subsection{Computer-Assisted Proof Strategy}
\label{sec: short_Computer-Assisted Strategy for Proving Stability}

\noindent The stability theorems reduce stability to a finite-dimensional optimization problem.  The aim is to choose the weights so that the linear damping margin is larger than the nonlinear estimates and residual error.  For each candidate weight, we compute the certified constants, optimize the radius, and update the weights until the closing inequality becomes strict. \\

\noindent The optimization loop proceeds as follows:

\vspace{2mm}

\emph{1. Parametrize the weight functions} with parameters \(\wgray{\vartheta}\), while ensuring that the weights have the required positivity, regularity, and symmetries. Choose an initial admissible candidate \((\wgray{\vartheta},\mu_k)\).
 
\vspace{2mm}

\emph{2. Compute the linear constants} $
\DLamLlow$, $
\DLamLhigh$, $
\DLamDLlow$,
$\DLamDLhigh$ for the current candidate. If the damping margin available $\DLamStab\le0$, this candidate cannot close the theorem, so we move to a new candidate.

\vspace{2mm}

\emph{3. Compute the nonlinear and residual constants.}  If \(\DLamStab>0\), compute $\DLamThree$ and $\DLamR$. 

\vspace{2mm}

\emph{4. Optimize the radius and evaluate the gap} $\mathcal J(\wgray{\vartheta},\mu_k)$ which determines whether the candidate succeeds. When successful, choose an admissible radius \(\delta_*\) such that \eqref{eq:short_ns-smallness_shortV} holds.

\vspace{2mm}

\emph{5. Iterate} steps 2--4 by  update \((\wgray{\vartheta},\mu_k)\) to decrease \(\mathcal J(\wgray{\vartheta},\mu_k)\), until a candidate satisfying the numerical closing condition is obtained. 

\vspace{2mm}

\emph{6. Certify} rigorously and analytically that the candidate weights and radius satisfy the closing condition and all required estimates.

\subsection{Partial Results: Linear Damping}

A central question in the stability argument is whether the leading linear terms provide a positive damping margin before the remaining terms are added. This matters because missing damping cannot generally be repaired by taking smaller perturbations. We evaluate the signed linearized operator directly in the weighted energies and obtain strongly favorable results. Low-order damping controls most of the half-plane, high-order damping supplies the localized repair near the off-axis low transport regions, singular weights can control the on-axis low transport region, the exterior region is closed analytically, and the profile residual is small in the weighted norms. \\

The low-order transport terms and selected couplings form a symmetric matrix $\mathbb M_{\mathcal L}(r,z)$. Pointwise damping is equivalent to $\lambda_{\max}\!\left(\mathbb M_{\mathcal L}(r,z)\right)<0.$ The weights $\wgray{\Phi_i}$ enter through the transport contribution and are optimized to decrease the largest eigenvalue of $\mathbb M_{\mathcal L}$. The largest eigenvalue is negative on the entire computational region, except in two small neighborhoods where the transport terms become very small. This defect is highly localized, and the high-order linear equations provide the natural repair off-axis because derivatives of the transport field need not be small when the transport itself vanishes. The numerical evaluation reveals higher-order damping in precisely the off-axis neighborhood where the low-order margin deteriorates. In addition, singular weights can be used to control the on-axis low transport region. 

Outside the spline support (i.e. the finite computational box), the profile-dependent coefficients vanish exactly, the signed matrix reduces to an explicit diagonal tail form, and the weights can be continued analytically so that the tail inequalities close without additional numerical computation.

\subsection{Remaining quantitative certification}

The signed matrix provides damping, but the final linear estimate must also include the remaining linear terms. The relevant linear stability margin is
\begin{equation*}
\DLamStab
=
\min\!\left\{
\DLamLlow-\mu_k\DLamDLlow,
\ \DLamDLhigh-\frac{\DLamLhigh}{\mu_k}
\right\}>0.
\end{equation*}
A positive $\DLamStab$ gives a strict linear damping margin. The certified signed matrix bounds show that the central linear mechanism required for this margin is already present. The remaining task is to certify the terms treated outside the signed matrices and verify that their combined losses remain below the available damping. If additional margin is needed, the same pipeline can be improved by refining the profile in the regions that contribute most to the residuals, optimizing the weights and other parameters, or sharpening the analytic estimates. \\

The immediate remaining task within this framework is quantitative: rigorously certify and optimize the remaining terms and estimates of the proof. Some selected analytic estimates may still require sharpening to recover additional margin, but such changes can be made within the existing modular proof framework. If the remaining constants are rigorously certified with a positive closing margin, thereby establishing nonlinear stability of the candidate profile, then the nonlinear stability argument is fully closed. Combined with the finite-time reconstruction from the dynamically rescaled variables, this would produce an admissible solution converging to the approximate self-similar profile while developing a singularity in finite physical time. Thus, remaining work is rigorous quantitative certification and optimization within the proposed stability and singularity framework.

\section*{Conclusion}

We provided the first evidence that the Euler equations admit a self-similar blowup on $\mathbb{R}^3$ at the critical scaling exponent $\lambda$ of $0.5$, in agreement with prior theory. We obtained a highly accurate approximate profile using PINNs with a novel combination of complementary tools, including soft and hard constraints, adaptive loss weighting, FP64 arithmetic, boosting, adaptive resampling, and self-scaled curvature aware optimizers. The PINN profile was then converted into piecewise polynomial spline surrogates suitable for rigorous certification. \\

Complementary low- and high-order energy estimates for the linear and nonlinear terms are derived in a separate stability paper~\citep{EulerBlowupStability}, with every constant and estimate made explicit or reduced to a rigorously computable quantity. These estimates are incorporated into a single energy functional that organizes the quantitative conditions required to establish stability of the profiles. The stability question is thereby reduced to determining whether a large finite collection of explicit inequalities can be rigorously certified with sufficient margin using interval arithmetic or \texttt{Arb} computations. We initiated this computational program by focusing on linear damping, and overall obtained low-order damping on most of the domain and high-order damping in the problematic off-axis low-transport region. This damping provides an important quantitative indication that the stability scheme may be feasible, and hence that the profiles could be stable.

\newpage

The remaining work within the proposed stability framework is primarily computational: rigorously certifying the profile-dependent constants and margins in the stability framework using spline representations, interval arithmetic, certified matrix bounds, and finite-dimensional optimization. Some estimates may require refinement to achieve sufficient margin, but the modular structure of the argument permits such adjustments without altering the overall proof strategy. If the required constants and margins can be certified with sufficient positive margin, thereby establishing nonlinear stability of the rescaled profile, this stability can then be transferred through dynamic rescaling to finite-time singularity formation in the original variables. When $\lambda=1/2$, control of the modulation parameters near their self-similar values would yield stability of the candidate profile and, after reconstruction, an admissible solution developing a finite-time singularity.\\

In parallel, the \texttt{LeanPDE} paper develops a formalization in Lean~\citep{Lean4} of key components of the argument, with particular emphasis on the symbolic derivations and proof steps. Numerical certificates are generated separately in \texttt{Julia} using \texttt{Arb} interval arithmetic~\citep{Julia,Arb}, while Lean verifies the algebraic identities and symbolic relations underlying the certified estimates. The resulting pipeline assigns complementary roles to its different components: physics-informed optimization is used for discovery, splines and interval arithmetic for rigorous numerical certification, and Lean for formal symbolic verification and, ultimately, machine-checked analysis. 

\newpage 

\onecolumn 

\section*{Author Contributions}

\textbf{A.G.} performed all numerical experiments and implementations, generated the figures, and contributed substantially to the development of the theoretical results in close collaboration with A.A. and V.D. A.G. also proofread the manuscript and provided suggestions on the exposition.\\

\noindent \textbf{V.D.} established most of the derivations and theoretical results, closely supervised the design, development, and execution of the numerical experiments, and carried out most of the writing and overall preparation of the manuscript.\\

\noindent \textbf{A.A.} led the project, closely supervised numerical experiments and theoretical developments, and made substantial contributions to the framing and presentation of the paper, including writing the abstract, introduction, and conclusion with V.D.

\hfill

\section*{Acknowledgments}

A. Anandkumar is supported in part by the Bren endowed chair, ONR (MURI grant N00014-18-12624), Darpa ExpMath HR0011 and by the AI2050 senior fellow program at Schmidt Sciences. V. Duruisseaux is supported in part by the AI2050 program and SciDAC P-240000887. A. Ganeshram was supported by the Caltech SURF program during the summer.\\

\noindent  We are especially grateful to Robert J. George for his generous contributions to this work, including his help with the certification framework and his formal verification of the derivations in Lean as part of the \texttt{LeanPDE} paper. We are also grateful to Peicong Song for his generous time, insightful discussions, and valuable suggestions. We also benefited from discussions with Changhe Yang and Yixuan Wang during the early stages of this work. \\

\noindent We acknowledge the use of LLMs for general brainstorming, verification of derivations, and Lean autoformalization. Importantly, LLMs were not used in the central component of this work, namely, the discovery of the approximate profiles. The PINN formulation, including the choice of soft and hard constraints and efficient sampling and optimization strategies, were all developed solely by the authors.

\clearpage
\onecolumn
\newgeometry{margin=0.8in}

\setcounter{section}{0}
\setcounter{subsection}{0}
\setcounter{subsubsection}{0}
\setcounter{figure}{0}
\setcounter{table}{0}
\setcounter{theorem}{0}
\setcounter{proposition}{0}
\renewcommand{\thefigure}{S\arabic{figure}}
\renewcommand{\thetable}{S\arabic{table}}
\renewcommand{\thetheorem}{S\arabic{theorem}}
\renewcommand{\theproposition}{S\arabic{proposition}}
\providecommand*{\theHfigure}{}
\providecommand*{\theHtable}{}
\providecommand*{\theHtheorem}{}
\providecommand*{\theHproposition}{}
\providecommand*{\theHsection}{}
\providecommand*{\theHsubsection}{}
\providecommand*{\theHsubsubsection}{}
\renewcommand*{\theHfigure}{SI.\arabic{figure}}
\renewcommand*{\theHtable}{SI.\arabic{table}}
\renewcommand*{\theHtheorem}{SI.\arabic{theorem}}
\renewcommand*{\theHproposition}{SI.\arabic{proposition}}
\renewcommand*{\theHsection}{SI.\arabic{section}}
\renewcommand*{\theHsubsection}{SI.\arabic{section}.\arabic{subsection}}
\renewcommand*{\theHsubsubsection}{SI.\arabic{section}.\arabic{subsection}.\arabic{subsubsection}}

\setcounter{equation}{0}
\renewcommand{\theequation}{S\arabic{equation}}
\providecommand*{\theHequation}{}
\renewcommand*{\theHequation}{SI.\arabic{equation}}

\supportingInformationHeading

\hfill \\

\printSupplementalContents
\startSupplementalTocCapture

\clearpage

\begin{figure}[p]
\centering
\begingroup
\definecolor{eulerstaborange}{RGB}{230,120,20}
\definecolor{eulerstabpurple}{RGB}{15,25,100}
\definecolor{eulerstabsalmon}{RGB}{240,80,70}
\definecolor{eulerstabteal}{RGB}{45,125,135}
\definecolor{eulerstabtextmuted}{RGB}{82,82,88}
\definecolor{eulerstabrulegray}{RGB}{140,140,148}
\definecolor{eulerstabrulemid}{RGB}{176,176,183}

\newcommand{\profileorange}[1]{\textcolor{eulerstaborange}{#1}}
\newcommand{\certpurple}[1]{\textcolor{eulerstabpurple}{#1}}
\newcommand{\salmonlambda}[1]{\textcolor{eulerstabsalmon}{#1}}

\newcommand{\LamR}{\ensuremath{\salmonlambda{\Lambda_{\mathcal R}}}}
\newcommand{\LamStab}{\ensuremath{\salmonlambda{\Lambda_{\rm stab}}}}

\tikzset{
  certflow/.style={
    -{Stealth[length=1.85mm,width=1.38mm]},
    draw=eulerstabrulegray,
    line width=0.50pt
  },
  methodflow/.style={
    -{Stealth[length=2.15mm,width=1.45mm]},
    draw=eulerstabrulegray,
    line width=0.60pt,
    shorten >=0.5pt,
    shorten <=0.5pt
  },
  certtransition/.style={
    fill=white,
    inner xsep=1.65pt,
    inner ysep=0.75pt,
    text=eulerstabtextmuted,
    font=\fontsize{10.0}{10.9}\selectfont\itshape,
    align=left
  },
  certstage/.style={
    draw=eulerstabrulemid,
    fill=white,
    line width=0.40pt,
    minimum width=16.55cm,
    text width=15.55cm,
    align=left,
    inner xsep=8.2pt,
    inner ysep=7.3pt,
    outer sep=0pt
  },
  methodboxleft/.style={
    draw=eulerstaborange!82,
    fill=white,
    rounded corners=5pt,
    line width=0.5pt,
    minimum width=2.46cm,
    text width=1.92cm,
    minimum height=0.78cm,
    align=center,
    inner xsep=2.8pt,
    inner ysep=2.0pt,
    outer sep=0pt,
    text=eulerstaborange,
    font=\fontsize{8.6}{9.1}\selectfont\bfseries
  },
  methodboxright/.style={
    draw=eulerstaborange!82,
    fill=white,
    rounded corners=5pt,
    line width=0.5pt,
    minimum width=2.46cm,
    text width=2.5cm,
    minimum height=0.78cm,
    align=center,
    inner xsep=2.8pt,
    inner ysep=2.0pt,
    outer sep=0pt,
    text=eulerstaborange,
    font=\fontsize{8.6}{9.1}\selectfont\bfseries
  }
}

\newcommand{\certtitle}[2]{%
  {\fontsize{13.9}{14.8}\selectfont\bfseries\color{#1}#2\par}%
}
\newcommand{\certbody}[1]{%
  {\vspace{5.0pt}\fontsize{10.45}{11.65}\selectfont\color{black!90}#1\par}%
}
\newcommand{\certphasetag}[2]{%
  {\fontsize{8.1}{8.7}\selectfont\bfseries\scshape\color{#1}#2}%
}

\newcommand{\certfullpanel}[6]{%
  \node[certstage,#6] (#1) {%
    \certtitle{#2}{#4}%
    \certbody{#5}%
  };
  \draw[#2!78,line width=0.72pt]
    ([xshift=0.2pt]#1.north west) --
    ([xshift=-0.2pt]#1.north east);
  \node[
    anchor=west,
    fill=white,
    inner xsep=2.5pt,
    inner ysep=1.0pt,
    text=#2
  ] at ([xshift=7.6mm]#1.north west) {\certphasetag{#2}{#3}};
}

\makebox[\linewidth][c]{%
\begin{tikzpicture}[node distance=5.6mm,font=\rmfamily]

\node[align=center,text=eulerstabpurple] (figtitle) at (0,0) {%
  {\fontsize{15.8}{16.6}\selectfont\bfseries
  From Numerical Profile Discovery Toward Certified Nonlinear Stability}%
};

\certfullpanel
  {profile}
  {eulerstaborange}
  {Numerics}
  {Approximate Numerical Profile}
  {\textbf{Profile discovery.}
   Starting from a parametrized traveling wave self-similar profile ansatz, obtain a blowup profile  $(u, \, \omega, \, \psi)$ using physics-informed neural networks (PINNs). 

   \par\vspace{5.2pt}\noindent
   \textbf{Analytic representation.}
   Fit the PINN profile with splines to obtain a representation suitable for rigorous certification, enabling exact estimates, and analytic expressions without introducing further numerical approximation errors.}
  {below=3.8mm of figtitle,
   minimum height=3.10cm,
   minimum width=9.80cm,
   text width=9.79cm}

\node[methodboxleft,anchor=east] (m1)
    at ([xshift=-7mm,yshift=1.25cm]profile.west)
    {Soft\\Constraints};

\node[methodboxleft,anchor=east] (m2)
    at ([xshift=-7mm,yshift=0cm]profile.west)
    {Hard\\Constraints};

\node[methodboxleft,anchor=east] (m3)
    at ([xshift=-7mm,yshift=-1.25cm]profile.west)
    {Boosting};

\node[methodboxright,anchor=west] (m4)
    at ([xshift=4.8mm,yshift=1.25cm]profile.east)
    {Curature-Aware\\Optimizer};

\node[methodboxright,anchor=west] (m5)
    at ([xshift=4.8mm,yshift=0cm]profile.east)
    {High Precision};

\node[methodboxright,anchor=west] (m6)
    at ([xshift=4.8mm,yshift=-1.25cm]profile.east)
    {Specialized\\Sampling};

\draw[methodflow] (m1.east) -- ([yshift=10.4mm]profile.west);
\draw[methodflow] (m2.east) -- ([yshift=0mm]profile.west);
\draw[methodflow] (m3.east) -- ([yshift=-10.4mm]profile.west);
\draw[methodflow] (m4.west) -- ([yshift=10.4mm]profile.east);
\draw[methodflow] (m5.west) -- ([yshift=0mm]profile.east);
\draw[methodflow] (m6.west) -- ([yshift=-10.4mm]profile.east);

\certfullpanel
  {energy}
  {blue}
  {Perturbations}
  {Perturbation Equations and Stability Energy}
  {\textbf{Linearize} around the spline profile and introduce \bvar{perturbations}
   $\bvar{\delta=(\delta\omega,\delta u)}$, with
   $\bvar{\delta v_\omega}:=\bvar{\delta\omega}$,
   $\bvar{\delta v_r}:=\partial_r\bvar{\delta u}$, and
   $\bvar{\delta v_z}:=\partial_z\bvar{\delta u}$. Schematically,
   \[
   \partial_t\bvar{\delta\omega}
   =\mathcal L_\omega(\bvar{\delta})+\mathcal N_\omega(\bvar{\delta})+\mathcal R_\omega,
   \qquad
   \partial_t\partial_i\bvar{\delta u}
   =\partial_i\mathcal L_u(\bvar{\delta})+\mathcal N_i(\bvar{\delta})+\mathcal R_i,
   \quad i\in\{r,z\}.
   \]

   \par\vspace{4pt}\noindent
   \textbf{Weights and energy.}
   Define the full energy
   $\mathfrak E_k^2:=\mathfrak E_0^2+\mu_k\mathfrak H_k^2$, where
   \[
   \mathfrak E_0^2:=\|\bvar{\delta\omega}\|_{\wgray{\Phi_\omega}}^2
   +\|\partial_r\bvar{\delta u}\|_{\wgray{\Phi_r}}^2
   +\|\partial_z\bvar{\delta u}\|_{\wgray{\Phi_z}}^2,
   \qquad
   \mathfrak H_k^2:=\sum_{i\in\{\omega,r,z\}}\ \sum_{|\alpha|=k}
   \|D^\alpha\bvar{\delta v_i}\|_{\wgray{\Phi_{i,k}}}^2, \vspace{-3mm}
   \]
   with $\wgray{\Phi_{i,k}}:=1+\rho^{\,k}\wgray{\Phi_i}$ and
   $\rho:=r^2+z^2$.}
  {below=6.6mm of profile,
   minimum height=5.55cm,
   minimum width=16.55cm,
   text width=15.55cm,
   inner ysep=1.6pt}

\certfullpanel
  {theory}
  {eulerstabpurple}
  {Stability Theory}
  {Analytic Estimates and Conditional Stability Proof}
  {\textbf{Analytic results and estimates.} Leverage Sobolev embeddings, weighted transport estimates, Leibniz estimates, elliptic bounds, and modulation estimates to derive the linear and nonlinear energy estimates required to close the stability argument, conditional on obtaining the necessary numerical constants in these estimates and suitable properties of the numerical profile.

   \par\vspace{4.8pt}\noindent
   The derivations, analytic estimates, and the bootstrap argument, can later be formalized in Lean.}
  {below=6.6mm of energy,
   minimum height=3.38cm,
   minimum width=16.55cm,
   text width=15.55cm,
   inner ysep=0.0pt}

\certfullpanel
  {certbox}
  {eulerstabpurple}
  {Certification}
  {Certified Constants and Estimates}
  {\textbf{Rigorous certification.}
   Starting from the spline profile, interval arithmetic rigorously encloses the profile data, the operator and stability constants, and the PDE residual bound, so that every numerical ingredient is converted into an explicit, machine-verifiable bound.

   \par\vspace{4.8pt}\noindent
   \textbf{Certified proof.} The certified quantities can then be incorporated into the conditional stability argument, with the aim of closing the proof. Subsequently, these estimates can be encoded within a Lean formalization framework to obtain a fully computer-verified proof.}
  {below=6.6mm of theory,
   minimum height=2.95cm,
   minimum width=16.55cm,
   text width=15.55cm}

\certfullpanel
  {stability}
  {eulerstabpurple}
  {Conclusion}
  {Computer-Assisted Nonlinear Stability}
  {\textbf{Optimize}
   the admissible singular weight functions \wgray{$\Phi$}, $\mu_k$, and stability radius   $\delta_*$ so that the certified damping dominates the nonlinear and residual contributions.

   \par\vspace{4.8pt}\noindent
   \textbf{Stability conclusion.}
   If $\mathfrak E_k(0)<\delta_*$, then
   $\mathfrak E_k(t)<\delta_*$ for every later admissible time.
   Therefore $\{\mathfrak E_k<\delta_*\}$ is forward invariant,
   yielding the desired nonlinear stability.}
  {below=6.6mm of certbox,
   minimum height=3.20cm,
   minimum width=16.55cm,
   text width=15.55cm}

\draw[eulerstabpurple!62,line width=0.32pt]
  ([xshift=1.9mm,yshift=-1.6mm]stability.north west)
  rectangle
  ([xshift=-1.9mm,yshift=1.6mm]stability.south east);

\draw[certflow]
  ([xshift=-1.75cm]profile.south) --
  ([xshift=-1.75cm]energy.north)
  node[midway,right=1.2pt,certtransition] {linearize around the numerical profile};

\draw[certflow]
  ([xshift=-1.75cm]energy.south) --
  ([xshift=-1.75cm]theory.north)
  node[midway,right=1.2pt,certtransition] {derive analytic estimates and form the stability proof};

\draw[certflow]
  ([xshift=-1.75cm]theory.south) --
  ([xshift=-1.75cm]certbox.north)
  node[midway,right=1.2pt,certtransition] {certify the numerical inputs and explicit estimates};

\draw[certflow]
  ([xshift=-1.75cm]certbox.south) --
  ([xshift=-1.75cm]stability.north)
  node[midway,right=1.2pt,certtransition] {optimize the singular weights and parameters to conclude};

\end{tikzpicture}
}%
\endgroup
\end{figure}

\clearpage

\section{The Euler Equations} \label{sec: Euler}

\subsection{The 3D Axisymmetric Equations} \label{sec: The 3D Axisymmetric Equations}

Starting from the incompressible 3D Euler equations written by the great \citet{EulerOriginalPaper},
\begin{equation}
\label{eq:Original_Euler}
\partial_t \tilde u + (\tilde u\cdot\nabla) \tilde u = -\nabla \tilde p,
\qquad
\nabla\cdot \tilde u = 0,
\end{equation}
we solve for the velocity field $\tilde u(x,t)\in\mathbb{R}^3$ and the scalar pressure $\tilde p(x,t)\in\mathbb{R}$, for $x\in\mathbb{R}^3$ and $t\ge 0$. Taking the curl of the above equation yields the vorticity formulation
\begin{equation}
\label{eq:vorticity}
\partial_t \tilde \omega + \tilde u\cdot\nabla \tilde \omega = \tilde \omega\cdot\nabla \tilde u,
\qquad
\tilde \omega := \nabla\times \tilde u .
\end{equation}
In $\mathbb{R}^3$, one can recover $\tilde u$ from $\tilde \omega$ via the Biot--Savart law. Equivalently, one can introduce a vector streamfunction $\tilde \psi(x,t)\in\mathbb{R}^3$ satisfying
\begin{equation}
\label{eq:vector_potential}
-\Delta \tilde \psi = \tilde \omega,
\qquad
\tilde u = \nabla\times \tilde \psi,
\end{equation}
together with a gauge condition (e.g.\ $\nabla\cdot \tilde \psi=0$) to fix $\tilde \psi$. \\ 

\noindent For flows that are symmetric about a fixed spatial axis, such as the $\mathring z$-axis, it is convenient to rewrite equation~\eqref{eq:Original_Euler} in cylindrical coordinates
\begin{equation}
x=(\mathring{r}\cos\theta,\mathring{r}\sin\theta,\mathring{z}),\qquad \mathring{r} \ge 0,\ \  \theta\in[0,2\pi),\ \  \mathring{z}\in\mathbb{R}.
\end{equation}
In these coordinates, the velocity and vorticity decompose as
\begin{align}
\tilde u = u^r(\mathring{r},\mathring{z},t)\,e_r + u^\theta(\mathring{r},\mathring{z},t)\,e_\theta + u^z(\mathring{r},\mathring{z},t)\,e_z, \\
\tilde \omega = \omega^r(\mathring{r},\mathring{z},t)\,e_r + \omega^\theta(\mathring{r},\mathring{z},t)\,e_\theta + \omega^z(\mathring{r},\mathring{z},t)\,e_z,
\end{align}
where $e_r, e_\theta, e_z$ are the unit vectors in cylindrical coordinates. In addition, by definition of axisymmetry, we have that $\partial_\theta(\cdot)=0$. \\ 

\noindent Following \citet{HouLuoBoundary,Chen2026Analysis}, it is convenient to introduce the rescaled swirl and streamfunction variables
\begin{equation}
\label{eq:scaled_vars}
u_{\bullet} := \frac{u^\theta}{\mathring{r}},\qquad
\omega_{\bullet} := \frac{\omega^\theta}{\mathring{r}},\qquad
\psi_{\bullet} := \frac{\psi^\theta}{\mathring{r}},
\end{equation}
where $\psi^\theta$ denotes the $\theta$-component of the vector potential $\psi$ in~\eqref{eq:vector_potential}. This reformulation has the advantage of removing the $1/\mathring{r}$ singularity from the cylindrical coordinates. \\

\noindent We also define the radial component $u^r$ and axial component $u^z$ of the velocity in cylindrical coordinates. These meridian velocity components $(u^r,u^z)$ can be recovered from $\psi_{\bullet}$ via the (axisymmetric) streamfunction relations
\begin{equation}
\label{eq:stream}
u^{r} = -\,\mathring{r}\,\partial_{\mathring{z}} \psi_{\bullet},
\qquad
u^{z} = 2\psi_{\bullet} + \mathring{r}\,\partial_{\mathring{r}} \psi_{\bullet}.
\end{equation}

\hfill 

\paragraph{3D Axisymmetric Euler Equations.} The axisymmetric 3D Euler equations can be written as
\begin{align}
\label{eq:Euler1}
\partial_t u_{\bullet} + u^r\,\partial_{\mathring{r}} u_{\bullet} + u^z\,\partial_{\mathring{z}} u_{\bullet}
&= 2u_{\bullet}\,\partial_{\mathring{z}} \psi_{\bullet}, \\
\label{eq:Euler2}
\partial_t \omega_{\bullet} + u^r\,\partial_{\mathring{r}} \omega_{\bullet} + u^z\,\partial_{\mathring{z}} \omega_{\bullet}
&= 2u_{\bullet}\,\partial_{\mathring{z}}u_{\bullet}, \\
\label{eq:Euler3}
-\Bigl[\partial_{\mathring{r}}^2 + \frac{3}{\mathring{r}}\partial_{\mathring{r}} + \partial_{\mathring{z}}^2\Bigr]\psi_{\bullet}
&= \omega_{\bullet}.
\end{align}

\noindent Equations \eqref{eq:stream}, \eqref{eq:Euler1}--\eqref{eq:Euler3} form the reduced axisymmetric-with-swirl system in $(\mathring{r},\mathring{z})$ for $(u_{\bullet},\omega_{\bullet},\psi_{\bullet})$.

\clearpage

\noindent \textbf{Axis interpretation.}
Every occurrence of \(\partial_{\mathring r}^2+3\mathring r^{-1}\partial_{\mathring r} +\partial_{\mathring z}^2\) at \(\mathring r=0\) denotes its axis-regular continuation, not literal division by zero.  For a \(C^2\) profile even in
\(\mathring r\),
\begin{equation}
\left(\partial_{\mathring r}^2+\frac{3}{\mathring r}
\partial_{\mathring r}+\partial_{\mathring z}^2\right)f(0,\mathring z)
=4\,\partial_{\mathring r\mathring r}f(0,\mathring z)
+\partial_{\mathring z\mathring z}f(0,\mathring z).
\end{equation}

\hfill  

\subsection{Traveling-Wave Ansatz} \label{sec: Ansatz}

\noindent
We consider a \emph{traveling self-similar} ansatz. 
\begin{itemize}
    \item \emph{Self-similar} refers to the fact that the geometry of the solution, after a proper rescaling of space, keeps an approximately fixed shape. 
    \item \emph{Traveling} means that the center of this singularity is not fixed. Instead, it moves along the symmetry ($\mathring{r} = 0$) axis while the solution converges to the self-similar geometry.
\end{itemize} 

\noindent Since the ansatz is written in terms of powers of $\tau = T-t$, it is also a \emph{backward self-similar ansatz}: the possible singular profile is described by rescaling backward from the blowup time $T$.\\

\noindent
The translation is introduced only in the axial ($\mathring{z}$) variable. This is consistent with the axisymmetric setting, where $\mathring r=0$ is the fixed symmetry axis and $\mathring r\geq 0$ measures distance to that axis. A radial ($\mathring{r}$) translation would move the singularity away from the symmetry axis and would generally destroy the parity and regularity structure imposed at $\mathring r=0$. By contrast, an axial translation in $\mathring z$ preserves the cylindrical symmetry and provides a natural mechanism for a coherent structure to drift while it concentrates.\\

\paragraph{The traveling-wave ansatz.}
\noindent
Given the original axisymmetric variables as $(u_\bullet,\omega_\bullet,\psi_\bullet)$, we introduce the \emph{rescaled profile variables} $(u,\omega,\psi)$ via
\begin{align}
u_\bullet(\mathring r,\mathring z,t)
&= (T-t)^{c_u}\,
 u\!\left(\frac{\mathring r}{(T-t)^\lambda},\frac{\mathring z-\mathring z_c(t)}{(T-t)^\lambda}\right),
\label{ansantz1}\\
\omega_\bullet(\mathring r,\mathring z,t)
&= (T-t)^{c_\omega}\,
 \omega\!\left(\frac{\mathring r}{(T-t)^\lambda},\frac{\mathring z-\mathring z_c(t)}{(T-t)^\lambda}\right),
\label{ansantz2}\\
\psi_\bullet(\mathring r,\mathring z,t)
&= (T-t)^{c_\psi}\,
 \psi\!\left(\frac{\mathring r}{(T-t)^\lambda},\frac{\mathring z-\mathring z_c(t)}{(T-t)^\lambda}\right).
\label{ansantz3}
\end{align}
\noindent
Here $T$ is the \emph{blowup time}, $\lambda$ is the \emph{spatial blowup rate}, $\mathring z_c(t)$ is the center of the traveling profile, and $c_u,c_\omega,c_\psi$ are \emph{amplitude exponents}. \\

\noindent The rescaled coordinates are
\begin{equation}
    r=\frac{\mathring r}{(T-t)^\lambda},
    \qquad
    z=\frac{\mathring z-\mathring z_c(t)}{(T-t)^\lambda},
\label{eq:rescaled_coordinates_ansatz}
\end{equation}
\noindent
and the rescaled meridian velocity components $u^r$ and $u^z$, can be recovered from $\psi$ via
\begin{equation}
    u^r=-r\partial_z\psi,
    \qquad
    u^z=2\psi+r\partial_r\psi.
    \label{eq:rescaled_stream_relations}
\end{equation}

\noindent The traveling speed is encoded by the constant $C$ through
\begin{equation}
    \frac{d}{dt}\mathring z_c(t)=-C(T-t)^{\lambda-1}.
    \label{eq:center_speed}
\end{equation}
\noindent
With this convention, the moving coordinate $z$ satisfies
\begin{equation}
    \left(\partial_t z\right)_{\mathring r,\mathring z}
    =\frac{\lambda z+C}{T-t}.
    \label{eq:z_coordinate_speed_ansatz}
\end{equation}
\noindent
Thus $C$ is the \emph{traveling speed} of the profile in the rescaled axial coordinate. If $C=0$, the singularity is centered at a fixed axial point in the rescaled frame. If $C\neq0$, the singularity drifts in the original variables while remaining stationary in the moving rescaled frame.\\

\noindent
The traveling-wave ansatz considered here provides more flexibility during the optimization than the more common \emph{stationary self-similar} dynamics studied in \citep{deepMindPaper1,deepMindPaper2,wang2025high,HouJiaje}. In the stationary case, one exploits scaling invariance and seeks a time-independent backward self-similar profile of the form
\begin{equation}
    u_\bullet(x,t) = (T-t)^{c_u} \, u\!\left(\frac{x - x_0}{ (T-t)^{\lambda}}\right).
    \label{normalselfsimilar}
\end{equation}
\noindent
Here $x$ denotes the relevant spatial variable, $(T-t)^\lambda$ is the concentrating length scale, and $(T-t)^{c_u}$ gives the amplitude growth. \\

\noindent The profile is \emph{stationary} because its center is fixed in physical space: after rescaling, the same spatial point remains the center of the core. In the present axisymmetric setting, the singularity is allowed to \emph{travel} along the symmetry axis. This introduces an additional axial center $\mathring z_c(t)$, and the profile is written in coordinates centered at $\mathring z_c(t)$ rather than at a fixed point.\\ 

\hfill 

\paragraph{Scaling exponents.}
\noindent
The exponents are fixed by balancing the powers of $T-t$ in the axisymmetric system. This balance enforces that time differentiation, transport, stretching, and the elliptic streamfunction relation all contribute at the same leading order. It gives
\begin{equation}
    c_u=-1,
    \qquad
    c_\omega=-1-\lambda,
    \qquad
    c_\psi=-1+\lambda.
    \label{eq:scaling_exponents}
\end{equation}
\noindent
These exponents have a simple interpretation. The variable $u_\bullet$ grows like $(T-t)^{-1}$. The vorticity variable $\omega_\bullet$ contains one additional spatial derivative at scale $(T-t)^\lambda$, and therefore grows like $(T-t)^{-1-\lambda}$. The streamfunction is two spatial derivatives smoother than $\omega_\bullet$, which leads to the exponent $-1+\lambda$.\\

\hfill

\paragraph{Energy scaling constraint.}
\noindent The traveling-wave ansatz also has consequences for energy quantities in the original variables. For the convection-varied Euler equations, the relevant energy identity is recalled from Theorem~4.2 of~\citep{PengfeiThesis}. The scaling calculation below applies whenever the corresponding energy is finite and remains bounded as $t\to T^-$. \\

\noindent Written in the original variables $(\mathring r,\mathring z)$, the relevant energy functional is
\begin{equation}
E
=\int_{-\infty}^{\infty}\int_{0}^{\infty}
\left(\lvert u_\bullet\rvert^2
+ \lvert \nabla_{\mathring r,\mathring z} \psi_\bullet\rvert^2\right)
\, \mathring r^3\,d\mathring r\,d\mathring z .
\label{eq:energy_functional_epsilon}
\end{equation}
\noindent
This identity links the self-similar scaling exponent $\lambda$ to a quantity expressed in the original physical variables. The following proposition shows how $E$ transforms under the traveling-wave ansatz and provides a simple but important consistency condition: a finite-energy self-similar profile cannot have an arbitrary concentration rate $\lambda$. 

\newpage

\begin{proposition}[Energy scaling under the self-similar ansatz]
\label{thm:energy-scaling}
Let $\tau=T-t$, and assume that
\begin{equation}
u_\bullet(\mathring r,\mathring z,t)
=
\tau^{-1}
u\!\left(
\frac{\mathring r}{\tau^\lambda},
\frac{\mathring z-\mathring z_c(t)}{\tau^\lambda}
\right),
\qquad
\psi_\bullet(\mathring r,\mathring z,t)
=
\tau^{\lambda-1}
\psi\!\left(
\frac{\mathring r}{\tau^\lambda},
\frac{\mathring z-\mathring z_c(t)}{\tau^\lambda}
\right).
\label{eq:energy_scaling_ansatz}
\end{equation}
Equivalently, introduce the rescaled variables
\begin{equation}
r = \frac{\mathring r}{\tau^\lambda},
\qquad
z = \frac{\mathring z - \mathring z_c(t)}{\tau^\lambda}.
\label{eq:energy_scaling_coordinates}
\end{equation}
If the rescaled profile energy is finite and nonzero, then boundedness
of $E(t)$ as $t\to T^-$ requires
\begin{equation}
\lambda\ge \frac{2}{5}.
\label{eq:energy_scaling_lambda_bound}
\end{equation}
\end{proposition}

\begin{grayproof}
    \ \ The proof is provided in Appendix~\ref{Energy Derivations_appx}.
\end{grayproof}

\hfill 

\noindent
\noindent
The proposition shows that the concentration rate $\lambda$ is constrained by the scaling of the energy in the original variables. If the rescaled profile energy is finite and nonzero, boundedness of $E(t)$ as $t\to T^-$ requires $\lambda\geq 2/5$. When $\lambda>2/5$, the contribution of an exact self-similar core to $E(t)$ tends to zero as $t\to T^-$. When $\lambda=2/5$, the energy contribution is invariant under the self-similar scaling. \\

\hfill 

\begin{remark*}
\citet{constantin2026putative} investigate the admissible scaling exponents for a hypothetical self-similar blowup of the Euler equations. They establish the lower bound $\lambda \geq 2/5$ and further show that, under a local outflow condition essential for the corresponding linear stability analysis, one must in fact have $\lambda \geq 1/2$. Prior to this work, it was hoped that an Euler singularity with $2/5 \leq \lambda<1/2$ might persist in the Navier--Stokes setting, since viscosity would be asymptotically lower order in this regime. In that case, a self-similar Euler blowup implies a Navier--Stokes singularity. The result of \citet{constantin2026putative}, however, rules out this scenario for profiles satisfying the local outflow condition. The exponent $\lambda = 1/2$ is especially important in both Euler and Navier--Stokes flows. For Euler, its importance is tied to the invariance of circulation, and for Navier--Stokes, this is the unique exponent at which self-similar blowup can occur with viscosity present as the blowup time is approached. 
\end{remark*}

\hfill \\

\subsection{Steady-State Profile Equations}
\label{sec: Steady-State Profile Equations}

We next present the profile equations generated by the traveling self-similar ansatz. This reduction transforms the original singular evolution into a time-independent nonlinear system for the rescaled variables $(U,\Omega,\Psi)$, together with the scaling exponent $\lambda$ and the drift speed $C$. \\

\noindent For notation simplicity, we define the elliptic operator
\begin{equation}
    \mathcal{E} \coloneqq \partial_r^2 + \frac{3}{r}\partial_r + \partial_z^2.
\end{equation}

\hfill 

\begin{proposition}[Traveling-wave profile equations]
\label{thm:traveling_wave_profile_eps1}
Assume that $(u_\bullet,\omega_\bullet,\psi_\bullet)$ is given by the ansatz \eqref{ansantz1}--\eqref{ansantz3}, that the center satisfies \eqref{eq:center_speed}, and that the exponents are chosen as in \eqref{eq:scaling_exponents}.  Set
\begin{equation} \lambda=\frac12. \label{eq:profile_eps1_fixed_lambda}
\end{equation}
Then $(U,\Omega,\Psi)$ satisfies
\begin{align}
U + \bigl(\lambda r+U^r\bigr)\partial_rU + \bigl(C+\lambda z+U^z\bigr)\partial_zU
&=2U\partial_z\Psi,
\label{eq:ProfileEq1_eps1}\\
(1+\lambda)\Omega + \bigl(\lambda r+U^r\bigr)\partial_r\Omega
+ \bigl(C+\lambda z+U^z\bigr)\partial_z\Omega
&=2U\partial_zU,
\label{eq:ProfileEq2_eps1}\\
-\, \mathcal E \,\Psi&=\Omega,
\label{eq:ProfileEq3_eps1}
\end{align}
where
\begin{equation}
U^r=-r\partial_z\Psi,
\qquad
U^z=2\Psi+r\partial_r\Psi.
\label{eq:profile_velocity_recovery_eps1}
\end{equation}
Equivalently,
\begin{align}
U+r\bigl(\lambda-\partial_z\Psi\bigr)\partial_rU
+\bigl(C+\lambda z+2\Psi+r\partial_r\Psi\bigr)\partial_zU
&=2U\partial_z\Psi,
\label{eq:ProfileEq1_expanded_eps1}\\
(1+\lambda)\Omega+r\bigl(\lambda-\partial_z\Psi\bigr)\partial_r\Omega
+\bigl(C+\lambda z+2\Psi+r\partial_r\Psi\bigr)\partial_z\Omega
&=2U\partial_zU,
\label{eq:ProfileEq2_expanded_eps1}\\
-\, \mathcal E \,\Psi&=\Omega.
\label{eq:ProfileEq3_expanded_eps1}
\end{align}
\end{proposition}

\begin{grayproof}
The derivation is provided in Appendix~\ref{appx: Traveling-wave profile equations derivation}.
\end{grayproof}

\hfill \\ 

\noindent In these equations,
\begin{itemize}
    \item The terms $\lambda r\partial_r$ and $\lambda z\partial_z$ come from observing the solution in a shrinking coordinate system.
    \item The term $C\partial_z$ comes from following the center of the profile as it travels along the symmetry axis. 
    \item The terms involving $ U^r$ and $ U^z$ are the convection terms.
    \item  The terms $2U\partial_z\Psi$ and $2U\partial_z U$ on the right-hand side are the stretching and coupling terms inherited from the axisymmetric formulation.
\end{itemize}

\hfill  

\noindent Equations \eqref{eq:ProfileEq1_eps1}--\eqref{eq:ProfileEq3_eps1} are the steady-state equations for the traveling self-similar profile. Thus, the traveling-wave ansatz reduces the study of a possible finite-time singularity to an autonomous nonlinear profile problem for $(U,\Omega,\Psi)$ and the parameters $\lambda$ and $C$. The parameter $\lambda$ determines the rate at which the spatial scale collapses, while $C$ determines the speed of the axial drift in the rescaled frame. 

\clearpage

\section{Discovering a Finite-Time Blowup via Physics-Informed Optimization} \label{sec: Finding the blowup}

Classical numerical methods for solving PDE typically rely on problem-specific discretizations and solver configurations, which can make their implementation and adaptation across different equations, geometries, and parameter regimes challenging. Among these approaches, finite element methods are among the most widely used due to their versatility and well-established mathematical foundation. However, high-fidelity finite element simulations can be computationally demanding, particularly for large-scale systems, problems characterized by rapidly varying dynamics, or applications requiring repeated PDE solutions. These limitations can become particularly acute for nonlinear PDE, where general-purpose finite element solvers may require weeks of CPU time and can struggle to discover previously unknown blowup profiles, especially when little prior information about their structure is available. \\

\noindent Physics-informed optimization offers a flexible and computationally tractable framework for solving PDE, with the potential to transfer across problem settings and facilitate the discovery of previously unknown solution profiles. Here, we focus on a particular class of such methods, physics-informed neural networks (PINNs), which recast the solution of a given PDE as an optimization problem whose objective is constructed directly from the governing equations and associated constraints. From this perspective, PINNs can be viewed as an alternative class of numerical methods for PDE, in which a parameterized representation of the solution is obtained through gradient-based optimization rather than through traditional algorithms applied to a prescribed discretization of the domain.\\

\subsection{Methods}
\label{sec: methodology}

\subsubsection{Physics-Informed Optimization}

Physics-Informed Neural Networks (PINNs)~\citep{PINN_OG}
are neural networks designed to approximate solutions to PDE. Their parameters are obtained by minimizing deviations from the governing PDE in an appropriate norm, taking a finite set of collocation points in the spatiotemporal domain as input and producing a parametrized approximation~$u_\theta(x,t)$ to the solution at each grid point $(x,t)$. Physics-Informed Neural Operators (PINOs) \citep{PINO_OG,ganeshram2025fcpinohighprecisionphysicsinformed} extend the same principle to learning a \emph{solution operator} mapping to an entire family of solutions~\citep{azizzadenesheli2024neural,Berner2025Principles}, for instance with Fourier Neural Operators~\citep{FNO,duruisseaux2025FNOGuide}. \\

\noindent The PDE residual loss is often combined with additional penalty terms that encode boundary and initial conditions, smoothness priors, and other problem-specific constraints. This is necessary because minimizing the residual may not, in general, guarantee that the optimized solution satisfies the prescribed constraints or exhibits the desired regularity, especially when training data are sparse or the optimization landscape is ill-conditioned. A typical resulting composite objective is given by
\begin{equation}
\mathcal{L}
= \mathcal{L}_{\text{PDE}}
+ \alpha\,\mathcal{L}_{\text{BC}}
+ \beta\,\mathcal{L}_{\text{Smoothness}}
+ \ldots.
\end{equation}
This composite loss is then used as the objective for gradient-based optimization to train the neural network $u_\theta $. Here, $\alpha,\beta,\ldots \in \mathbb{R}^{+}$ denote positive hyperparameters that weight the respective loss terms.

\hfill

\noindent While PINNs have proven extremely successful in numerous practical applications, the corresponding optimization task can be challenging and prone to numerical issues. The training loss can have poor conditioning as it involves differential operators that can be ill-conditioned. In particular, \citet{Krishnapriyan2021} showed that the loss landscape becomes increasingly complex and harder to optimize as the physics loss coefficient increases. The model could also converge to a trivial or non-desired solution that satisfies the physics laws on the set of points where the physics loss is computed \citep{Leiteritz2021}. There could also be conflicts between the multiple loss terms when their gradients point in opposite directions. Even without such conflicts, the losses can vary significantly in magnitude, leading to unbalanced gradients during training \citep{Wang2021}, and thus to the contribution and reduction of certain losses being relatively negligible during training. Consequently, obtaining reliable PINN solutions can require substantial care throughout the training process, including the choice of network architecture, collocation sampling strategy, loss weights, training hyperparameters, optimizer, and optimization strategy. In practice, these components must be tuned and monitored carefully to ensure that the optimized solution reflects the intended physical behavior rather than artifacts of the optimization procedure. \\ 

\hfill

\subsubsection{Achieving High-Precision with PINNs}

\noindent To make PINN solutions amenable to computer-assisted proofs, we must achieve very high accuracy when solving the Euler and Navier--Stokes equations. In this paper, we leverage several mechanisms that substantially improve the precision attainable with PINNs and enable the construction of sufficiently accurate approximate solutions.

\hfill 

\paragraph{Higher Floating Point Precision. }  64-bit floating point precision (FP64) with PINNs helps preserve accuracy when computing automatic-differentiation derivatives and PDE residuals, where FP32 round-off can accumulate and destabilize training. With cleaner gradients and more reliable physics losses, FP64 often improves convergence and yields higher-precision solutions that better satisfy the governing equations and physics constraints~\citep{xufp64}. We carry out PINN training in FP64 to ensure our solution losses are evaluated with the highest numerical fidelity that is possible in PyTorch.

\hfill 

\paragraph{Boosting. } Our methodology is inspired by the existing literature on PINN boosting~\citep{fang2024boosting,boostingYongji}. We write our solution approximation as a stage-wise additive model,
\begin{equation}
\label{eq:boosted_ansatz}
u(x,t) \approx f_{\theta_0}(x,t) + \alpha_1 f_{\theta_1}(x,t) + \alpha_2 f_{\theta_2}(x,t) + ... \ ,
\end{equation}
where each newly added PINN $f_{\theta_k}$ is introduced only after the current surrogate has reached an optimization plateau. The stage-$k$ PINN $f_{\theta_k}$ is trained \emph{conditioned} on the previously optimized sum: the parameters of the earlier PINNs $f_{\theta_0},\ldots,f_{\theta_{k-1}}$ are frozen, and only $f_{\theta_k}$ is optimized so that it serves as a correction targeting the residual structure left. In effect, boosting replaces a single, difficult global fitting task with a sequence of easier residual-fitting subproblems. The successive correction networks can continue to reduce the PDE residual beyond the saturation point of a single PINN, yielding higher-accuracy solutions. Because each stage solves a distinct residual-correction problem, one can vary the PINN architecture and training strategy across stages, tailoring $f_{\theta_k}$ and its optimization to the specific residual features it is intended to capture. This flexibility also allows later-stage residual networks to be localized to small subdomains or patches where errors are concentrated, or to be combined naturally with domain-decomposition strategies that assign different correction networks to different regions of the computational domain.

\hfill

\paragraph{Self-Scaled Curvature-Aware Optimizers.} Adam~\citep{adamOptim} is widely used for neural-network training, but its reliance on first-order gradient information can limit convergence in PINNs, whose loss landscapes are often stiff and poorly conditioned. Curvature-aware optimization can alleviate these difficulties and continue reducing the objective after first-order methods begin to stagnate, often leading to substantially higher solution accuracy, though this comes at an increased computational cost. Quasi-Newton methods achieve this by constructing an approximation to the inverse Hessian from successive parameter and gradient updates. However, these approximations may themselves become poorly conditioned, degrading the quality of the resulting search directions.\\

\noindent Self-scaling addresses this issue by adaptively rescaling the inverse-Hessian approximation to better capture the local geometry of the objective landscape. In PINNs, such self-scaled quasi-Newton updates have been shown to improve conditioning and accelerate convergence~\citep{SSBroyden2}.\\

\noindent More recently, structured curvature-aware optimizers such as Shampoo~\citep{Shampoo} and SOAP~\citep{SOAP} have provided a more scalable means of exploiting curvature information. Rather than maintaining a dense inverse-Hessian approximation, they construct structured preconditioners that capture important correlations in the gradient while remaining computationally practical for larger neural networks. In particular, SOAP performs adaptive optimization in the eigenbasis of the Shampoo preconditioner, yielding updates that are better aligned with the local geometry of the objective. Standard Shampoo and SOAP, however, do not have the self-scaling mechanism.\\

\noindent Recently, $\sssoap$~\citep{SS-SOAP} was introduced to augment SOAP with self-scaling and an adaptive eigenbasis update. This retains much of the computational efficiency of first-order optimization while providing a more effective curvature-aware update. In our experiments, $\sssoap$ incurs only a $15\%$ increase in computational cost relative to Adam, while substantially improving the PDE MSE residual from $\mathcal{O}(10^{-4})$ with Adam to $\mathcal{O}(10^{-6})$. For further refinement, we employ the self-scaled Broyden (SS-Broyden) method~\citep{Al-Baali1998,SSBroyden2}. SS-Broyden maintains a more expressive approximation of the inverse Hessian and updates both its scaling and Broyden parameters using local curvature information. This richer curvature model comes at a substantially higher computational and memory cost than $\sssoap$, making SS-Broyden less attractive as a general-purpose optimizer throughout training. However, once the solution has entered a sufficiently accurate regime, the improved curvature approximation can continue reducing the objective beyond the level at which the cheaper optimizer begins to saturate. We use SS-Broyden as a final high-accuracy refinement stage, reducing the PDE MSE residual in our experiments from $\mathcal{O}(10^{-6})$ with $\sssoap$ to $\mathcal{O}(10^{-10})$.

\hfill \\

\subsubsection{Dealing with Unbounded Domains}  \label{sec: sinhTransform}

\noindent Sampling on an unbounded domain is a fundamental challenge for numerical methods. Similar to \citep{deepMindPaper1,wang2025high}, we use the sinh transformation (and the chain rule for derivatives), to sample over the half-plane
\begin{equation}
    \mathbb{D}=\{(r,z): r\geq 0,\ z\in\mathbb{R}\}.
\end{equation}

\hfill 

\noindent More precisely, we sample the computational variables on the finite box
\begin{equation}
    \big(r_{\!\mathfrak{comp}}, z_{\mathfrak{comp}}\big)
    \in [0,30]\times[-30,30],
\end{equation}
and map them to the physical variables by
\begin{equation}
    (r,z)
    =
    \big(\sinh r_{\!\mathfrak{comp}}, \, \sinh z_{\mathfrak{comp}}\big).
\end{equation}

\noindent Under the $\sinh$ mapping, the transformation is approximately linear near the origin, since $\sinh(x)\sim x$ as $x\to 0$. Uniform sampling in the computational variables therefore produces nearly uniform, finely resolved sampling in the physical domain near the origin, where the blowup profile is expected to vary most rapidly. Farther from the origin, the mapping progressively stretches the sampling points, enabling the same computational grid to cover a large portion of the physical domain without unnecessarily high resolution in the far field. This distribution is well suited to the anticipated solution structure, which exhibits sharp variation near the blowup region and a much slower decay toward zero in the far field. \\

\noindent  Under the $\sinh$ mapping, since $
    \sinh(30)\approx 5.34\times 10^{12}$, the computational box $[0,30]\times[-30,30]$ maps roughly to $
    [0,\,5.34\times 10^{12}]
    \times
    [-5.34\times 10^{12},\,5.34\times 10^{12}]$, which is effectively unbounded for practical purposes. The restriction to nonnegative values of $r$ is consistent with the assumed even parity in $r$, so the solution can be mirrored across the axis $r=0$ (or enforced as a hard constraint in the PINN architecture). The odd symmetry of the $\sinh$ map is particularly convenient here, as it preserves the underlying parity structure when expressed in the computational variables. \\

\noindent The choice of $\sinh$ also avoids some of the drawbacks associated with alternative transformations. For example, a compactifying map such as $r=\mathrm{arctanh}\,\rho$, with $0\leq \rho<1$, maps the infinite physical domain to the finite interval $\rho\in[0,1)$, but becomes singular as $\rho\to 1$, since $dr/d\rho=1/(1-\rho^2)$. Polynomial maps avoid this singular behavior, but their algebraic growth is typically too slow to sample the physical far field efficiently. The $\sinh$ map therefore provides a natural compromise: it reaches large physical scales rapidly while remaining smooth and preserving parity through its odd symmetry. \\

\noindent We represent the PINN in the computational variables and obtain derivatives with respect to the physical variables via the chain rule. When the computational and physical coordinates are related by the $\sinh$ map introduced above, we use the same notation $u_\Theta$ for the induced representation in physical space:
\begin{equation}
    u_\Theta(r,z)
    :=
    u_\Theta\big(r_{\!\mathfrak{comp}},z_{\mathfrak{comp}}\big).
\end{equation}
The first derivatives in physical coordinates are then
\begin{equation}
    \frac{\partial u_\Theta}{\partial r}
    =
    \frac{\partial u_\Theta}{\partial r_{\!\mathfrak{comp}}}
    \mathrm{sech}\big(r_{\!\mathfrak{comp}}\big),
    \qquad
    \frac{\partial u_\Theta}{\partial z}
    =
    \frac{\partial u_\Theta}{\partial z_{\mathfrak{comp}}}
    \mathrm{sech}\big(z_{\mathfrak{comp}}\big).
\end{equation}
The same principle applies to higher-order derivatives: the transformed expressions involve higher-order chain-rule coefficients, built from smooth bounded functions such as sech and tanh on the computational box, and are evaluated using automatic differentiation~\citep{Autograd}, as is standard in PINN implementations~\citep{Baydin2017}. We give the full transformed PDE residual operator in Appendix~\ref{app:sinhDerivation}. 

\hfill 

\subsubsection{PINN Objectives and Constraints}

\noindent Let \(u_\Theta\) be a neural-network approximation to the PDE solution with parameters~\(\Theta\), being the collection of all the network parameters to be optimized and the traveling wave speed $C$. To optimize the PINN representation, we minimize a weighted combination of loss terms
\begin{equation}
\mathcal{L}(\Theta)
\ = \  \mathcal{L}_{\text{PDE}}(\Theta)
\ + \  \alpha\,\mathcal{L}_{\mathsf{FF}}(\Theta)
\ + \  \beta\,\mathcal{L}_{\text{Smoothness}}(\Theta)
\ + \  \gamma \,\mathcal{L}_{\text{Non-triviality}}(\Theta),
\end{equation}
where 
\begin{itemize}
\item $ \mathcal{L}_{\text{PDE}}(\Theta) =    \int_{\mathbb{D}} \left\lVert \mathcal{R}\big(u_\Theta)(x) \right\rVert^{2}\, dx$ denotes the PDE residual loss, where $\mathcal{R}$ is the PDE residual operator. 
\item $\mathcal{L}_{\mathsf{FF}}(\Theta)$ denotes the far field loss conditions that $u,\,\omega, \psi \to 0$ as $r,z \to \infty.$ We emphasize that these are decay conditions at infinity rather than classical boundary conditions (i.e Dirichlet, Neumann, Robin) at a finite boundary.
\item $
\mathcal{L}_{\text{Smoothness}}(\Theta) = \int_{\mathbb{D}} \left\lVert \nabla\mathcal{R}\big(u_\Theta)(x) \right\rVert^{2}\, dx$  penalizes the gradient of the PDE residual. This promotes smoothness of the optimized representation. Heuristically, controlling $\nabla \mathcal{R}$ discourages highly oscillatory residuals and biases the optimized solution toward higher regularity on the region where the penalty is enforced. The gradient is computed via automatic differentiation.
\item $\mathcal{L}_{\text{Non-triviality}}(\Theta)$ is an additional loss term designed to discourage the optimization from converging to trivial solutions.

\end{itemize}

\hfill

\noindent In contrast to the preceding constraints, which are encouraged softly through the loss function, we enforce the following symmetry exactly as a hard constraint. Because our self-similar ansatz travels along the $z$-axis, no symmetry condition is imposed in $z$. Instead, we impose even parity in $r$:
\begin{equation}
u(r,z) = u(-r,z), \quad \omega(r,z) = \omega(-r,z), \quad \psi(r,z) = \psi(-r,z). \label{eqn}
\end{equation}
This parity condition is physically motivated by \citet{HouLuoBoundary}, where, in the presence of a boundary, an analogous symmetry serves as a trapping mechanism that promotes vorticity growth toward singularity formation.

\hfill 

\subsubsection{Sampling Collocation Points}

\noindent To improve generalization and avoid overfitting to a fixed set of collocation points, we periodically resample training points and bias sampling toward the origin, where the solution varies most rapidly. Since the solution is even in $r$, we restrict sampling to $r\geq 0$. We also use residual-based adaptive sampling~\citep{Wu2023}, progressively adding points in regions where the PDE residual is largest. Separate samples are used for the PDE, smoothness, and far-field losses, with additional points placed along the symmetry axis $r=0$, where the formally singular radial term is evaluated using its regular limiting form.\\

\noindent For SS-Broyden, we use the same overall strategy but resample less frequently, since full-batch second-order methods are more sensitive to changes in the training set through their use of past iterates in Hessian approximations. We also find that sufficiently large batches are important for good generalization across the computational domain.

\hfill \\

\subsection{Numerical PINN Profiles}

\begin{figure}[t]
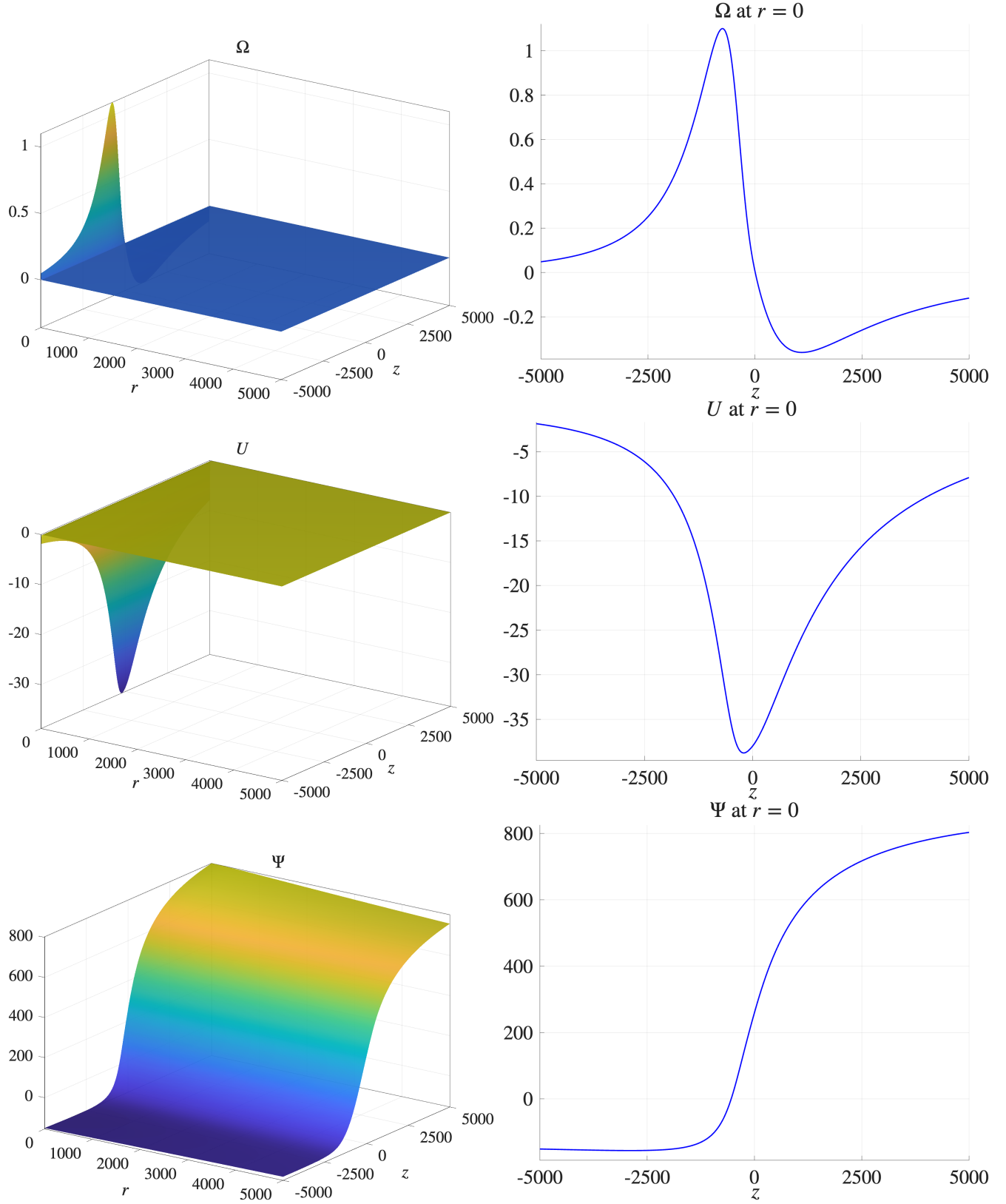
 
\centering

\begin{subfigure}{\linewidth}
\centering
\includegraphics[width=0.495\linewidth]{omega3D.pdf}\hfill
\includegraphics[width=0.495\linewidth]{omegaSlice.pdf}
\end{subfigure}

\begin{subfigure}{\linewidth}
\centering
\includegraphics[width=0.495\linewidth]{u3D.pdf}\hfill
\includegraphics[width=0.495\linewidth]{uSlice.pdf}
\end{subfigure}

\begin{subfigure}{\linewidth}
\centering
\includegraphics[width=0.495\linewidth]{psi3D.pdf}\hfill
\includegraphics[width=0.495\linewidth]{psiSlice.pdf}
\end{subfigure}

\caption{Visualization of the approximate profile $\Omega, U, \Psi$ as three-dimensional surface plots for the 3D axisymmetric Euler equations with their corresponding $r=0$ cross-sections.} 
\label{fig: profiles}
\end{figure}

\begin{figure*}[htbp] 
\centering

\includegraphics[height=10.4cm]{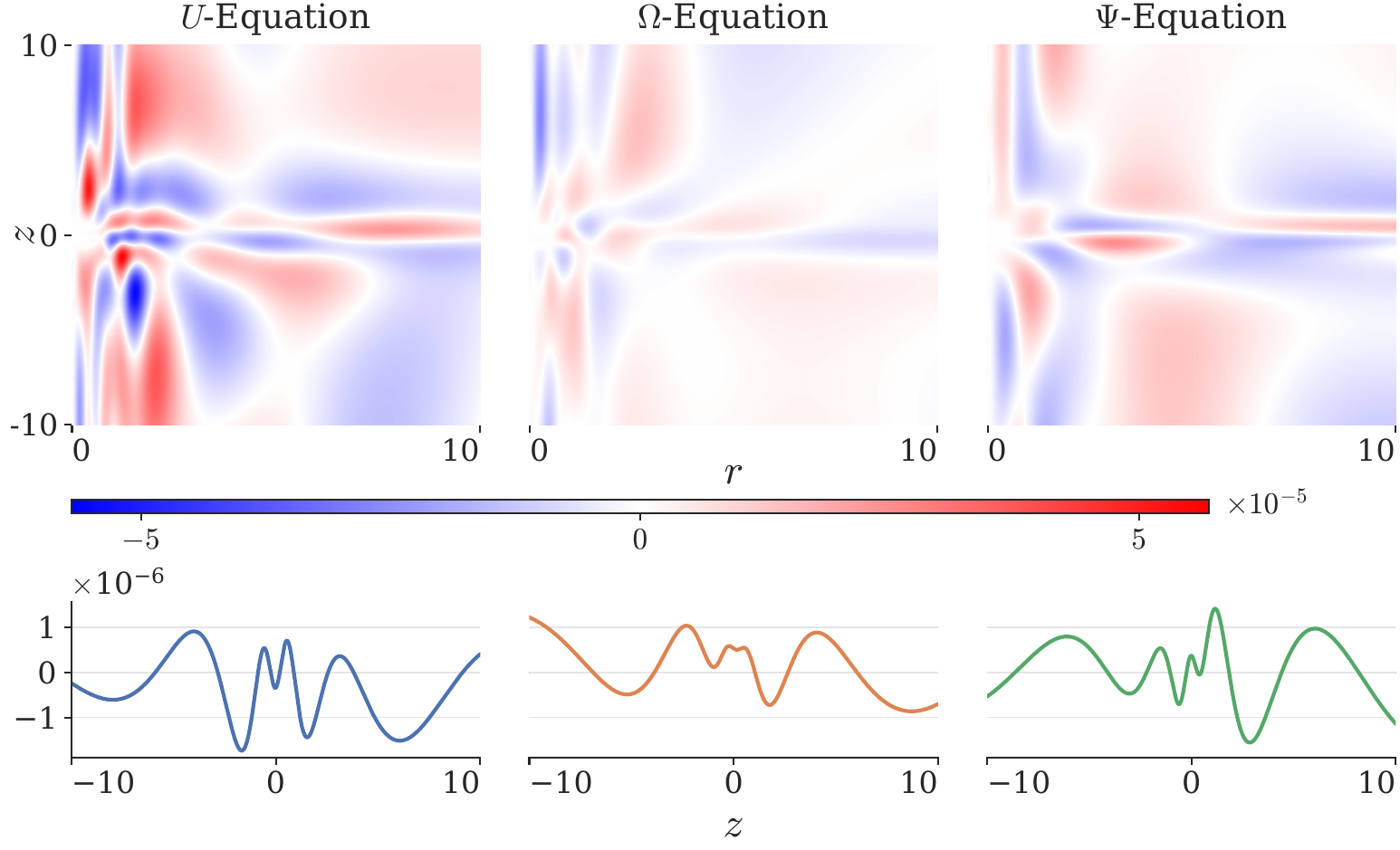}

\caption{\textbf{(Top)} Profile equation residuals on the box $(r,z)\in[0,10]\times[-10,10]$. \textbf{(Bottom)} Profile equation residuals on the symmetry-axis ($r=0$) with $z \in [-10,10]$.}
\label{fig:residuals} \vspace{8mm}

\end{figure*}

\noindent Leveraging the methods described in \Cref{sec: methodology}, we parametrize and train a PINN solution for the Euler self-similar profile system. We visualize the optimized approximate profiles $u$, $\omega$, and $\psi$ as three-dimensional surfaces and display the corresponding $r=0$ cross-sections in Figure~\ref{fig: profiles}. The optimized solution reproduces the qualitative profile structure reported by \citet{pengfeiPaper}, showing that the self-similar profile remains robust even when convection is fully amplified. These profiles were obtained to high accuracy, as can be seen from the profile equations residuals reported in \Cref{tab:table_errors} and displayed in \Cref{fig:residuals}.

\hfill

\begin{remark*}
    We also investigate in Appendix~\ref{appx: convection-varied} the family of self-similar profiles for the convection-varied Euler system from \citet{pengfeiPaper} where a parameter $\varepsilon$ allows convection to be progressively increased from the non-convective case at $\varepsilon=0$ to the full Euler system at $\varepsilon=1$. Using direct time evolution from a common numerically generated initial state, \citet{pengfeiPaper} observed a traveling self-similar singularity scenario at the tested values \(\varepsilon\in \{0,0.1,0.2,0.3\} \), whereas the evolution at \(\varepsilon=0.4\) did not exhibit the same behavior. Their dynamic-rescaling computations further indicated that the computed self-similar state became progressively less stable as \(\varepsilon\) increased. These results left open whether self-similar singular profiles persist for stronger convection but are not reached by the forward dynamics. Our PINN formulation finds high-accuracy approximate self-similar profiles throughout \(\varepsilon\in[0,1]\), providing numerical evidence that the disappearance of the earlier dynamically observed scenario should not be interpreted as the nonexistence of singular profiles. We observe from \Cref{fig:lambdaC_vs_eps} that the value of $\lambda$ decreases steadily towards $\lambda=1/2$ as $\varepsilon$ increases towards $\varepsilon=1$, as predicted by \citet{constantin2026putative}.
\end{remark*}

\clearpage

\begin{table*}[htbp]
\centering

\caption{Residual errors for the $U$, $\partial_r U$, $\partial_z U$, $\Omega$, and $\Psi$ profile equations for the Euler system, comparing PINN and spline representations. Note that the residuals are estimated empirically on an independent set of $100{,}000{,}000$ randomly sampled collocation points for the PINN representations, whereas they are evaluated analytically and exactly for the spline representations. in particular, cRMSE refers to the continuous root-mean squared error $ \sqrt{
\frac{1}{\left|\mathbb{D}\right|}
\int_{\mathbb{D}}
\mathsf{residual}\!\left(r_{\!\mathfrak{comp}}, z_{\mathfrak{comp}}\right)^2
\,\mathrm{d}r_{\mathfrak{comp}}\,\mathrm{d}z_{\mathfrak{comp}}
} $. Here, $\mathbb{D}$ denotes the full computational domain (corresponding to the half plane), $B_{0.1}(0^{\star})$, $B_1(0^{\star})$, and $B_{10}(0^{\star})$ denote balls of radii $0.1$, $1$, and $10$, respectively, centered at the origin $0^{\star}$ of the stability analysis (i.e. the on-axis meridional fixed point), $\mathsf{axis}$ denotes the $r=0$ axis, and $\mathsf{FF}$ denotes the far-field region $\left(r_{\!\mathfrak{comp}}, z_{\mathfrak{comp}}\right) \in [10, 30] \times [-30, -10] \cup [10, 30] \times [10, 30] $.}
\label{tab:table_errors}

\vspace{0.7mm}
\renewcommand{\arraystretch}{1.04}

\begin{tabularx}{0.9\linewidth}{
@{}
c@{\hspace{2mm}}
l
*{4}{>{\centering\arraybackslash}X}
@{}
}

&
&
\multicolumn{2}{c}{\textbf{PINNs}}
&
\multicolumn{2}{c}{\textbf{Splines}}
\\[-1.1mm]

\cmidrule(lr){3-4}
\cmidrule(lr){5-6}

&
\textbf{Domain}
& RMSE
& $\ell^\infty$
& cRMSE
& $L^\infty$
\\

\midrule
\addlinespace[0.7mm]


\multirow{6}{*}{
\rotatebox[origin=c]{90}{\textbf{\small $U$-equation}}
}
&
$B_{0.1}(0^{\star})$
& $2.8\!\cdot\!10^{-7}$
& $3.4\!\cdot\!10^{-7}$
& $2.8\!\cdot\!10^{-7}$
& $3.5\!\cdot\!10^{-7}$
\\

&
$B_1(0^\star)$
& $7.0\!\cdot\!10^{-6}$
& $3.3\!\cdot\!10^{-5}$
& $6.7\!\cdot\!10^{-6}$
& $3.4\!\cdot\!10^{-5}$
\\

&
$B_{10}(0^{\star})$
& $1.2\!\cdot\!10^{-5}$
& $5.7\!\cdot\!10^{-5}$
& $1.7\!\cdot\!10^{-5}$
& $5.8\!\cdot\!10^{-5}$
\\

&
$\mathsf{axis}$
& $9.0\!\cdot\!10^{-5}$
& $1.4\!\cdot\!10^{-3}$
& $9.1\!\cdot\!10^{-5}$
& $1.4\!\cdot\!10^{-3}$
\\

&
$\mathsf{FF}$
& $6.2\!\cdot\!10^{-6}$
& $2.2\!\cdot\!10^{-5}$
& $6.3\!\cdot\!10^{-6}$
& $2.3\!\cdot\!10^{-5}$
\\

&
\cellcolor{grayshade}$\mathbb{D}$
& \cellcolor{grayshade}$4.2\!\cdot\!10^{-5}$
& \cellcolor{grayshade}$3.6\!\cdot\!10^{-3}$
& \cellcolor{grayshade}$5.7\!\cdot\!10^{-5}$
& \cellcolor{grayshade}$3.7\!\cdot\!10^{-3}$
\\

\addlinespace[1.0mm]
\midrule
\addlinespace[1.0mm]


\multirow{6}{*}{
\rotatebox[origin=c]{90}{\textbf{\small $\partial_r U$-equation}}
}
&
\rule[-0.55ex]{0pt}{3.15ex}$B_{0.1}(0^{\star})$
& $1.5\!\cdot\!10^{-6}$
& $3.3\!\cdot\!10^{-6}$
& $1.6\!\cdot\!10^{-6}$
& $3.4\!\cdot\!10^{-6}$
\\

&
\rule[-0.55ex]{0pt}{3.15ex}$B_1(0^\star)$
& $5.2\!\cdot\!10^{-5}$
& $1.9\!\cdot\!10^{-4}$
& $5.0\!\cdot\!10^{-5}$
& $1.9\!\cdot\!10^{-4}$
\\

&
\rule[-0.55ex]{0pt}{3.15ex}$B_{10}(0^{\star})$
& $4.2\!\cdot\!10^{-5}$
& $3.1\!\cdot\!10^{-4}$
& $6.9\!\cdot\!10^{-5}$
& $3.2\!\cdot\!10^{-4}$
\\

&
\rule[-0.55ex]{0pt}{3.15ex}$\mathsf{axis}$
& $0$
& $0$
& $0$
& $0$
\\

&
\rule[-0.55ex]{0pt}{3.15ex}$\mathsf{FF}$
& $1.3\!\cdot\!10^{-10}$
& $2.4\!\cdot\!10^{-9}$
& $1.3\!\cdot\!10^{-10}$
& $2.5\!\cdot\!10^{-9}$
\\

&
\cellcolor{grayshade}\rule[-0.55ex]{0pt}{3.15ex}$\mathbb{D}$
& \cellcolor{grayshade}$1.2\!\cdot\!10^{-4}$
& \cellcolor{grayshade}$8.3\!\cdot\!10^{-3}$
& \cellcolor{grayshade}$1.5\!\cdot\!10^{-4}$
& \cellcolor{grayshade}$8.5\!\cdot\!10^{-3}$
\\

\addlinespace[1.0mm]
\midrule
\addlinespace[1.0mm]


\multirow{6}{*}{
\rotatebox[origin=c]{90}{\textbf{\small $\partial_z U$-equation}}
}
&
\rule[-0.55ex]{0pt}{3.15ex}$B_{0.1}(0^{\star})$
& $1.1\!\cdot\!10^{-6}$
& $1.9\!\cdot\!10^{-6}$
& $1.1\!\cdot\!10^{-6}$
& $2.3\!\cdot\!10^{-6}$
\\

&
\rule[-0.55ex]{0pt}{3.15ex}$B_1(0^\star)$
& $2.3\!\cdot\!10^{-5}$
& $1.0\!\cdot\!10^{-4}$
& $2.3\!\cdot\!10^{-5}$
& $1.1\!\cdot\!10^{-4}$
\\

&
\rule[-0.55ex]{0pt}{3.15ex}$B_{10}(0^{\star})$
& $1.3\!\cdot\!10^{-5}$
& $1.3\!\cdot\!10^{-4}$
& $2.8\!\cdot\!10^{-5}$
& $1.3\!\cdot\!10^{-4}$
\\

&
\rule[-0.55ex]{0pt}{3.15ex}$\mathsf{axis}$
& $1.1\!\cdot\!10^{-6}$
& $1.7\!\cdot\!10^{-5}$
& $1.2\!\cdot\!10^{-6}$
& $1.8\!\cdot\!10^{-5}$
\\

&
\rule[-0.55ex]{0pt}{3.15ex}$\mathsf{FF}$
& $1.5\!\cdot\!10^{-11}$
& $1.8\!\cdot\!10^{-10}$
& $1.6\!\cdot\!10^{-11}$
& $2.7\!\cdot\!10^{-10}$
\\

&
\cellcolor{grayshade}\rule[-0.55ex]{0pt}{3.15ex}$\mathbb{D}$
& \cellcolor{grayshade}$2.9\!\cdot\!10^{-6}$
& \cellcolor{grayshade}$1.3\!\cdot\!10^{-4}$
& \cellcolor{grayshade}$4.7\!\cdot\!10^{-6}$
& \cellcolor{grayshade}$1.3\!\cdot\!10^{-4}$
\\

\addlinespace[1.0mm]
\midrule
\addlinespace[1.0mm]


\multirow{6}{*}{
\rotatebox[origin=c]{90}{\textbf{\small $\Omega$-equation}}
}
&
$B_{0.1}(0^{\star})$
& $5.0\!\cdot\!10^{-7}$
& $5.8\!\cdot\!10^{-7}$
& $5.1\!\cdot\!10^{-7}$
& $5.9\!\cdot\!10^{-7}$
\\

&
$B_1(0^\star)$
& $3.4\!\cdot\!10^{-6}$
& $1.1\!\cdot\!10^{-5}$
& $3.3\!\cdot\!10^{-6}$
& $1.1\!\cdot\!10^{-5}$
\\

&
$B_{10}(0^{\star})$
& $4.4\!\cdot\!10^{-6}$
& $2.7\!\cdot\!10^{-5}$
& $5.2\!\cdot\!10^{-6}$
& $2.8\!\cdot\!10^{-5}$
\\

&
$\mathsf{axis}$
& $4.8\!\cdot\!10^{-5}$
& $3.7\!\cdot\!10^{-4}$
& $4.9\!\cdot\!10^{-5}$
& $3.8\!\cdot\!10^{-4}$
\\

&
$\mathsf{FF}$
& $3.1\!\cdot\!10^{-6}$
& $1.4\!\cdot\!10^{-5}$
& $3.2\!\cdot\!10^{-6}$
& $1.4\!\cdot\!10^{-5}$
\\

&
\cellcolor{grayshade}$\mathbb{D}$
& \cellcolor{grayshade}$1.7\!\cdot\!10^{-5}$
& \cellcolor{grayshade}$9.5\!\cdot\!10^{-4}$
& \cellcolor{grayshade}$1.9\!\cdot\!10^{-5}$
& \cellcolor{grayshade}$9.6\!\cdot\!10^{-4}$
\\

\addlinespace[1.0mm]
\midrule
\addlinespace[1.0mm]


\multirow{6}{*}{
\rotatebox[origin=c]{90}{\textbf{\small $\Psi$-equation}}
}
&
$B_{0.1}(0^{\star})$
& $3.6\!\cdot\!10^{-7}$
& $1.2\!\cdot\!10^{-6}$
& $3.7\!\cdot\!10^{-7}$
& $4.2\!\cdot\!10^{-7}$
\\

&
$B_1(0^\star)$
& $2.8\!\cdot\!10^{-6}$
& $8.9\!\cdot\!10^{-6}$
& $2.7\!\cdot\!10^{-6}$
& $9.0\!\cdot\!10^{-6}$
\\

&
$B_{10}(0^{\star})$
& $6.9\!\cdot\!10^{-6}$
& $2.6\!\cdot\!10^{-5}$
& $8.0\!\cdot\!10^{-6}$
& $2.7\!\cdot\!10^{-5}$
\\

&
$\mathsf{axis}$
& $8.5\!\cdot\!10^{-5}$
& $1.4\!\cdot\!10^{-3}$
& $8.6\!\cdot\!10^{-5}$
& $1.5\!\cdot\!10^{-3}$
\\

&
$\mathsf{FF}$
& $2.1\!\cdot\!10^{-6}$
& $8.3\!\cdot\!10^{-6}$
& $2.1\!\cdot\!10^{-6}$
& $8.3\!\cdot\!10^{-6}$
\\

&
\cellcolor{grayshade}$\mathbb{D}$
& \cellcolor{grayshade}$3.8\!\cdot\!10^{-5}$
& \cellcolor{grayshade}$1.4\!\cdot\!10^{-3}$
& \cellcolor{grayshade}$5.3\!\cdot\!10^{-5}$
& \cellcolor{grayshade}$1.5\!\cdot\!10^{-3}$
\\

\addlinespace[0.7mm]
\bottomrule

\end{tabularx}

\vspace{1mm}
\end{table*}

\clearpage

\subsection{Analytic Spline Representation of the Steady State} \label{sec: analytic representation splines}

The PINN computation serves as a numerical discovery tool, providing a practical and well-conditioned means of obtaining an accurate approximate self-similar profile when direct optimization over analytic ansatzes or alternative basis expansions would be considerably more difficult. We denote steady-state quantities by an overhead bar, and let $$ (\bar u_{\rm PINN}, \, \bar\psi_{\rm PINN}, \, \bar\omega_{\rm PINN})$$ denote the profile produced by the PINN. \\

\noindent In practice, however, the PINN solution is not the steady state we work over for our stability analysis and numerical certification. This is due to the fact that the basis optimized by the PINN can be highly nonlinear and oscillatory, making it difficult to interpret. We instead freeze an explicit analytic profile $(\bar u, \, \bar \psi, \, \bar \omega)$, together with certified intervals for the fitted parameters. All residuals, energy estimates, and nonlinear constants are then evaluated using this fitted profile rather than the original neural network. This different representation allows for certified residuals and norms, which is essential in our stability argument as it involves the repeatedly differentiating, calculating weighted perturbations, and finding $L^\infty$ bounds of the numerical profile. These operations can be certified for explicit analytic representations (via interval arithmetic), whereas certifying them directly for the PINN would require a separate verification of the network evaluation, automatic differentiation procedure, and $L^\infty$ bounds.  \\

\noindent We must now choose an analytic representation that is sufficiently expressive to reproduce the computed profile with high accuracy, while remaining sufficiently explicit for rigorous certification. For this purpose, we use \emph{piecewise polynomial splines}. This choice is common in numerical analysis and is particularly well suited to verification because, on each subinterval, the approximation is given by an explicit low-degree polynomial, so that its derivatives, extrema, and residuals can be evaluated and bounded using standard rigorous tools, including interval arithmetic. Moreover, the local nature of the representation allows the mesh to be refined only where the profile exhibits rapid variation without unnecessarily increasing the complexity of the approximation elsewhere. In contrast, global polynomial or spectral representations have to find a tradeoff between global and local behavior which may require a relatively high degree of polynomials, leading to larger expressions and less convenient rigorous bounds. Piecewise splines avoid this, achieving high local accuracy with moderate polynomial degree, preserve smoothness across subintervals, and yield a representation whose errors and derivatives can be certified in an more interpretable and computationally efficient manner. 

\hfill 

\noindent Overall, the certification pipeline has four separate layers, as described in the diagram below, where only the first arrow is a numerical approximation step: 

\begin{figure}[h]
  \centering

  \resizebox{0.73\linewidth}{!}{%
  \begin{tikzpicture}[node distance=5.0mm,font=\rmfamily]

    \certpipelinetagpanel
      {certprofile}
      {certiforange}
      {Numerics}
      {Approximate PINN Profile}
      {Numerical Discovery}
      {}

    \certpipelinetagpanel
      {certspline}
      {certifpurple}
      {Analytic}
      {Spline Basis Representation}
      {Analytic Representation}
      {below=of certprofile}

    \certpipelinetagpanel
      {certinterval}
      {certifpurple}
      {Certification}
      {Interval certificates}
      {Rigorous \texttt{Arb} Enclosures}
      {below=of certspline}

    \certpipelinetagpanel
      {certlean}
      {certifpurple}
      {Proof}
      {Computer-Assisted Stability Proof}
      {}
      {below=of certinterval}

    \draw[certifpurple!62,line width=0.22pt]
      ([xshift=1.6mm,yshift=-1.4mm]certlean.north west)
      rectangle
      ([xshift=-1.6mm,yshift=1.4mm]certlean.south east);

    \draw[certpipeline flow]
      (certprofile.south) --
      node[certpipeline transition] {fit splines}
      (certspline.north);

    \draw[certpipeline flow]
      (certspline.south) --
      node[certpipeline transition] {analytic evaluation}
      (certinterval.north);

    \draw[certpipeline flow]
      (certinterval.south) --
      node[certpipeline transition] {incorporate in formal argument}
      (certlean.north);

  \end{tikzpicture}%
  }

\end{figure}

\hfill \\

\noindent We define $\mathcal R_u^{\mathsf{raw}}$ and $\mathcal R_\omega^{\mathsf{raw}}$ to be the PDE residuals accounting for the fact that the numerical profiles do not satisfy exactly the Euler equations:
\begin{align}
\mathcal R_u^{\mathsf{raw}}
&:=-(\lambda r+\bar u^r)\partial_r\bar u-(\lambda z+\bar C+\bar u^z)\partial_z\bar u+2\bar u\partial_z\bar\psi+\bar c_u\bar u,\\
\mathcal R_\omega^{\mathsf{raw}}
&:=-(\lambda r+\bar u^r)\partial_r\bar\omega-(\lambda z+\bar C+\bar u^z)\partial_z\bar\omega+2\bar u\partial_z\bar u+(\bar c_u-\lambda)\bar\omega .
\end{align}

\noindent The spline fit preserves the small profile-equation residuals observed in the numerical solution, so the passage to an analytic representation incurs only a minor loss of accuracy. For the stability analysis, we further construct the surrogate through a two-head fit: \((\bar u,\bar\psi)\) are fitted separately, while \(\bar\omega\) is defined by applying the elliptic operator exactly to \(\bar\psi\). Consequently, the residual in the elliptic equation is exactly 0, eliminating one source of certification error and significantly simplifying the subsequent stability analysis.

\clearpage

\section{Stability}

A singularity is \emph{stable} if the associated blowup mechanism persists under small perturbations of the initial data, in the sense that the same mechanism occurs for an open neighborhood of nearby initial conditions. Conversely, an \emph{unstable} singularity requires infinitely precise initial conditions: small perturbations deflect the evolution away from the singularity mechanism.

\hfill 

\noindent The dynamic rescaling formulation presented in \Cref{sec: Dynamic Rescaling} is the natural starting point for such a stability analysis. In the original physical variables, the solution may grow and concentrate as $t\to T^-$, so there is no fixed object around which to linearize. After rescaling, the approximate blowup profile becomes a steady state of the rescaled equations, while the blowup rate, translation rate, and amplitude normalization are recorded by modulation parameters. Linear stability can therefore be studied by perturbing this steady rescaled profile and analyzing the resulting linearized rescaled dynamics. \\

\noindent We follow the same overall structure as \citep{DGBlowup,HouJiaje,Chen2026Analysis} and related works on stability of blowups. We first compute the blowup profile numerically, then develop a framework for establishing linear stability of the dynamically rescaled equation in suitable weighted norms and for proving nonlinear stability using Sobolev embeddings and weighted estimates. The main distinction is that the weight functions and other tunable parameters entering the stability proof are optimized numerically to seek sufficient margin for the argument to close, rather than being chosen and refined through expert intuition and a potentially long and technically challenging process of trial and error. The overall proof architecture is summarized in the diagram on page~\pageref{diagram: Stability Proof}, while the optimization strategy used to close the proof is described in detail in \Cref{sec: Computer-Assisted Strategy for Proving Stability}.

\hfill

\begin{remark*}
An alternative sufficient criterion is to prove invertibility of the full linearized profile operator~\citep{elgindi2025dynamics}. In this case, the self-similar profile equation is viewed perturbatively around an approximate profile, and the goal is to invert the linearized operator after removing the neutral scaling direction. This is stronger than the criterion used here. In contrast, the numerical eigenvalue analysis in \citep[Section~4.4]{PengfeiThesis} gives a more limited, finite-dimensional test of the same stability mechanism: it consists of discretizing the dynamically rescaled linearized equation around the computed profile and inspecting the spectrum of the resulting Jacobian. This provides useful evidence for stability, but is not by itself a rigorous substitute for the required linear estimate. It only concerns a finite-dimensional discretization and therefore controls, at best, resolved low-frequency perturbations. As noted in \citep[Section~1.3.1]{chen2022singularity}, the proof must control the infinite-dimensional operator, including high-frequency modes, nonlocal velocity terms, possible neutral directions, and transient growth caused by non-normality. Thus, negative real parts of the computed eigenvalues do not yield the uniform estimate required for nonlinear stability.
\end{remark*}

\subsection{Dynamic Rescaling Equations}
\label{sec: Dynamic Rescaling}

\subsubsection{Introduction}

\noindent
The steady profile equations above describe an exact traveling self-similar solution. In computations and stability arguments, however, one usually does not know the correct profile, or traveling speed in advance. The \emph{dynamic rescaling formulation} addresses this by letting the rescaled solution evolve while continuously adjusting its location, scale, and amplitude. The key idea is to replace the finite-time singular behavior in physical variables by a long-time evolution in normalized variables. The singularity is kept at order-one size and near a fixed location in the rescaled coordinates. In this frame, blowup no longer appears primarily as unbounded growth or shrinking length scale. Instead, a stable blowup scenario should appear as convergence toward a steady state. That steady-state is precisely the traveling-wave profile described in \Cref{thm:traveling_wave_profile_eps1}. \\

\paragraph{Modulated coordinates and amplitudes.}
\noindent
In the dynamic formulation, the fixed self-similar parameters in the ansatz are replaced by time-dependent \emph{modulation parameters}. Rather than prescribing a blowup time, a spatial scale, and a traveling center in advance, the rescaled frame is adjusted as the solution evolves. More precisely, let $(\mathring r,\mathring z,\mathring t)$ denote the physical variables and let $(r,z,t)$ denote the dynamically rescaled variables.  We introduce a common spatial scale $\mathrm{s}_r(t)>0$, amplitude scales $\mathrm{s}_u(t), \, \mathrm{s}_\omega(t)>0$, and an axial center $\mathring z_c(t)$ by
\begin{equation}
\mathring r=\mathrm{s}_r(t) \, r, \qquad\quad  \mathring z-\mathring z_c(t)=\mathrm{s}_r(t) \, z,
\label{eq:dynamic_spatial_rescaling_setup}
\end{equation}
and
\begin{equation}
u_\bullet(\mathring r,\mathring z,\mathring t)
=\mathrm{s}_u(t) \, u(r,z,t), \qquad \quad 
\omega_\bullet(\mathring r,\mathring z,\mathring t) =\mathrm{s}_\omega(t) \, \omega(r,z,t).
\label{eq:dynamic_u_omega_rescaling_setup}
\end{equation}
Then,
\begin{equation}
\psi_\bullet(\mathring r,\mathring z,\mathring t)
=\mathrm{s}_\omega(t) \, \mathrm{s}_r(t)^2 \, \psi(r,z,t).
\label{eq:dynamic_psi_rescaling_setup}
\end{equation}
The normalization of the transport, stretching, and time-derivative terms is fixed by
\begin{equation}
\mathrm{s}_u(t)=\mathrm{s}_\omega(t) \, \mathrm{s}_r(t),
\qquad
\frac{dt}{d\mathring t}=\mathrm{s}_\omega(t) \, \mathrm{s}_r(t).
\label{eq:dynamic_amplitude_time_rescaling_setup}
\end{equation}
We then define the rescaling rates by
\begin{equation}
\partial_t\mathrm{s}_r(t)=-\lambda \,\mathrm{s}_r(t),
\qquad
\partial_t\mathrm{s}_u(t)=-c_u(t) \, \mathrm{s}_u(t),
\qquad
\partial_t\mathrm{s}_\omega(t)=-c_\omega(t) \, \mathrm{s}_\omega(t).
\label{eq:dynamic_rate_definitions_setup}
\end{equation}
and define the axial drift via
\begin{equation}
\partial_t\mathring z_c(t)=-C(t) \, \mathrm{s}_r(t).
\label{eq:dynamic_axial_drift_setup}
\end{equation}

\noindent Overall, this adjustment is encoded by the coordinate drift
\begin{equation}
    \left(\partial_t r\right)_{\mathring r,\mathring z}=\lambda r,
    \qquad
    \left(\partial_t z\right)_{\mathring r,\mathring z}=\lambda z+C(t).
    \label{eq:dynamic_coordinate_drift}
\end{equation}
where
\begin{itemize}
    \item $\lambda$ is the \emph{spatial rescaling rate}. It measures, in the current rescaled time variable, how quickly the computational frame expands
    \item the function $C(t)$ is the \emph{instantaneous traveling speed} of the axial center in the same rescaled frame.
    \item  the term $\lambda z$ is the axial part of the dilation, while $C(t)$ is the additional axial translation. 
\end{itemize}

\noindent The amplitude of the solution must also be normalized. This is achieved via two scalar functions $c_u(t)$ and $c_\omega(t)$, which are the logarithmic amplitude rates for $u$ and $\omega$. The $c_u(t)u$ and $c_\omega(t)\omega$ terms appearing in the rescaled evolution equations do not represent new forcing, they only record the rate at which the normalization of the rescaled variables is changed in order to keep the profile at order-one amplitude. \\

\paragraph{Use and connection with stability.}
\noindent
The profile formulation and the dynamic rescaling formulation play complementary roles. The profile equations describe the limiting singularity once the blowup rate, traveling speed, and amplitude rates have settled to constants, while the dynamic equations describe the evolution toward that limiting object. Numerically, this is essential: instead of resolving a structure whose amplitude diverges and whose length scale collapses in physical variables, one follows a normalized solution whose main profile remains order one in a fixed computational frame. The modulation parameters provide direct diagnostics for the blowup scenario: if $C(t)$, $c_u(t)$, and $c_\omega(t)$ converge, their limits give the axial traveling speed, and the amplitude growth rates. If the rescaled profiles converge at the same time, the computation gives more than evidence of singular growth. It gives a candidate self-similar profile, the values of the scaling parameters, and the profile equations that the limiting object should satisfy.\\

\noindent
The dynamic rescaling formulation is the natural setting for stability. In the original physical variables, a perturbation of a blowup solution can change the blowup time, shift the center location, or alter the amplitude normalization. Such changes may appear as growth or drift even when the underlying profile is stable. The normalization conditions remove these neutral directions by continuously recentering, rescaling, and renormalizing the solution. After this gauge freedom has been fixed, stability becomes a question about the evolution of genuine perturbations in the dynamically rescaled variables. A stable traveling self-similar singularity is represented by a stable steady state of the dynamic rescaling equations: solutions starting close to the profile should remain controlled, the modulation parameters should approach limiting constants, and the rescaled solution should converge to the steady profile. In this formulation, finite-time blowup in the original variables is linked to convergence toward a stable fixed point of the evolution.\\

\subsubsection{Dynamic Rescaling Equations for the Euler Equations}

\noindent 
As before, define the elliptic operator
\begin{equation}
    \mathcal{E} \coloneqq \partial_r^2 + \frac{3}{r}\partial_r + \partial_z^2.
\end{equation}

\hfill  

\noindent The proposition below provides the axisymmetric Euler equations in the dynamically rescaled variables.  \\

\begin{proposition}[Dynamic rescaling equations for the Euler]
\label{thm:dynamic_rescaling_equations eps1}
In the dynamically rescaled variables, the axisymmetric system takes the form
\begin{align}
\partial_t u
+\bigl(\lambda r+ u^r\bigr)\partial_r u
+\bigl(\lambda z+C(t)+ u^z\bigr)\partial_z u
&=2u\partial_z\psi+c_u(t)u,
\label{eq:dynamic_rescaled_u eps1}\\
\partial_t\omega
+\bigl(\lambda r+ u^r\bigr)\partial_r\omega
+\bigl(\lambda z+C(t)+ u^z\bigr)\partial_z\omega
&=2u\partial_z u+c_\omega(t)\omega,
\label{eq:dynamic_rescaled_omega_eps1}\\
- \mathcal{E} \psi
&=\omega,
\label{eq:dynamic_rescaled_psi eps1}
\end{align}
with
\begin{equation}
    u^r=-r\partial_z\psi,
    \qquad
    u^z=2\psi+r\partial_r\psi,
    \label{eq:dynamic_rescaled_stream_intro eps1}
\end{equation}
and $\lambda = 1/2$. 
\end{proposition}
\begin{grayproof}
The derivation is provided in Appendix~\ref{appx: dynamic rescaling equations derivation}.
\end{grayproof}

\hfill

\paragraph{Normalization conditions.}

\noindent
When $\lambda=1/2$ is fixed, the independent modulation parameters are
$C(t)$ and $c_u(t)$. The vorticity amplitude rate is determined by
\begin{equation}
    c_\omega(t)=c_u(t)-\frac{1}{2}.
\end{equation}
Before the stability analysis, write the radial and axial transport coefficients of the numerically obtained profile as
\begin{equation}
G_r(r,z)=\lambda r+u^r_{\mathsf{num}}(r,z),
\qquad
G_z(r,z)=\lambda z+C_{\mathsf{num}}+u^z_{\mathsf{num}}(r,z).
\label{eq:numerical-transport-intro-add}
\end{equation}
We choose an axial shift $z_{\mathsf{shift}}$ so that the point $(0,z_{\mathsf{shift}})$ of the numerical profile becomes the analysis origin. For this point to be a fixed point of the meridional flow, we require
\begin{equation}
G_r(0,z_{\mathsf{shift}})=0,
\qquad
G_z(0,z_{\mathsf{shift}})=0.
\label{eq:numerical-transport-stagnation-intro-add}
\end{equation}
The radial condition is automatic because $u^r_{\mathsf{num}}(0,z)=0$ on the axis. Using $u^z_{\mathsf{num}}(0,z)=2\psi_{\mathsf{num}}(0,z)$, the axial condition $G_z(0,z_{\mathsf{shift}})=0$ becomes the scalar equation
\begin{equation}
C_{\mathsf{num}}+\lambda z_{\mathsf{shift}}+2\psi_{\mathsf{num}}(0,z_{\mathsf{shift}})=0.
\label{eq:shift-root-intro-add}
\end{equation}
We impose this condition and choose $z_{\mathsf{shift}}$ as a numerically certified root. Having fixed $z_{\mathsf{shift}}$, we recenter each numerical profile field by
\begin{equation}
\bar f(r,z):=f_{\mathsf{num}}(r,z+z_{\mathsf{shift}}).
\label{eq:profile-recentering-intro-add}
\end{equation}
Denoting the constant term in the recentered axial transport by $\bar C$, we continue to write the reference transport coefficients as
\begin{equation}
G_r(r,z)=\lambda r+\bar u^r(r,z),
\qquad
G_z(r,z)=\lambda z+\bar C+\bar u^z(r,z).
\label{eq:recentered-transport-intro-add}
\end{equation}
With this notation, the imposed root condition gives
\begin{equation}
G_r(0,0)=0,
\qquad
G_z(0,0)=\bar C+2\bar\psi_0=0.
\label{eq:recentered-stagnation-intro-add}
\end{equation}

\noindent After applying this shift, set $\bar u_0:=\bar u(0,0)$. We preserve the meridional fixed point and the amplitude at the origin by imposing
\begin{equation}
C(t)+u^z(0,0,t)=0,
\qquad
u(0,0,t)=\bar u_0.
\label{eq:normalization_conditions_intro_eps1}
\end{equation}
The first condition keeps the reference meridional fixed point at the origin. The second keeps the amplitude fixed there and therefore gives $\partial_tu(0,0,t)=0$.

\hfill \\

\paragraph{Limiting behaviour.}
\noindent
If the dynamic rescaling captures a stable traveling self-similar blowup, the rescaled solution should converge as rescaled time advances to a time-independent profile,
\begin{equation}
    u(r,z,t)\to u_\infty(r,z),
    \qquad
    \omega(r,z,t)\to \omega_\infty(r,z),
    \qquad
    \psi(r,z,t)\to \psi_\infty(r,z),\label{eq:profile_convergence_intro}
\end{equation}
\noindent
and the modulation parameters should converge to constants,
\begin{equation}
    C(t)\to C_\infty,
    \qquad
    c_u(t)\to c_{u,\infty},
    \qquad
    c_\omega(t)\to c_{\omega,\infty}.
\label{eq:modulation_convergence_intro}
\end{equation}
\noindent
In that limit, the time derivatives in the dynamic rescaling equations disappear. The steady dynamic equations are
\begin{align}
\Bigl(\lambda \, r+ u^r_\infty\Bigr)\partial_r u_\infty
+\Bigl(\lambda \, z+ C_\infty + u^z_\infty \Bigr)\partial_z u_\infty
&=2 u_\infty \, \partial_z \psi_\infty+ c_{u,\infty} \, u_\infty ,
\label{eq:steady_dynamic_u}\\
\Bigl(\lambda \, r+ u^r_\infty\Bigr)\partial_r\omega_\infty
+\Bigl(\lambda \, z+C_\infty+ u^z_\infty\Bigr)\partial_z\omega_\infty
&=2 u_\infty \, \partial_z u_\infty+ c_{\omega,\infty} \, \omega_\infty,
\label{eq:steady_dynamic_w}
\end{align}
\noindent
with $u^r_\infty =-r\partial_z \psi_\infty$ and $u^z_\infty=2 \psi_\infty+r\partial_r\psi_\infty$. \\ 

\noindent
For the traveling-wave ansatz in \eqref{ansantz1}--\eqref{ansantz3}, the limiting amplitude rates are
\begin{equation}
    c_{u,\infty}=-1,
    \qquad
    c_{\omega,\infty} = -(1+\lambda)=-3/2.
    \label{eq:limiting_amplitude_rates}
\end{equation}

\noindent By writing
\begin{equation}
\bigl(u_\infty,  \, \omega_\infty \,,\psi_\infty  ,\, C_\infty  \bigr) = \bigl(U, \, \Omega, \, \Psi, \, C \bigr),
\end{equation} 
and substituting \eqref{eq:limiting_amplitude_rates} into \eqref{eq:steady_dynamic_u}--\eqref{eq:steady_dynamic_w}, we see that the limiting dynamic equations agree with the corresponding profile equations. In this sense, the traveling-wave profiles are fixed points of the appropriate dynamically rescaled evolution. \\

\subsection{Linearization}

\noindent In the present section, we derive the perturbation equations around a steady state of the dynamically rescaled system. We then separate the terms that are linear in the perturbation from the higher-order remainders. The linearized equations identify the operator whose stability determines whether the candidate blowup profile is attracting to first order in the rescaled variables.

\subsubsection{Perturbations}

\noindent We denote perturbations using \textcolor{blue}{$\delta$}'s in \textcolor{blue}{blue}. We use the subscript $_0$ to denote evaluation at the origin $(r,z)=(0,0)$. We perturb a given steady state
\begin{equation}
    \bar u,\bar\omega,\bar\psi,\bar C,\bar c_u,\bar c_\omega
\end{equation}
of the dynamic rescaling equations, by writing
\begin{equation}
    u = \bar u + \textcolor{blue}{\delta u},
    \qquad
    \omega=\bar\omega+\textcolor{blue}{\delta\omega},
    \qquad
    \psi = \bar\psi + \textcolor{blue}{\delta\psi},
\label{eq:pert_ansatz_1}
\end{equation}
\begin{equation}
    C=\bar C+\textcolor{blue}{\delta C},
    \qquad
    c_u = \bar c_u + \textcolor{blue}{\delta c_u},
    \qquad
    c_\omega=\bar c_\omega+\textcolor{blue}{\delta c_\omega}. \label{eq:pert_ansatz_2}
\end{equation}

\hfill 

\noindent The perturbation variables measure the difference between the evolving rescaled solution and the steady blowup profile. The goal of linear stability analysis is to determine whether these perturbations remain bounded, or preferably decay, under the dynamics obtained by linearizing around the steady state. \\

\noindent The perturbations must inherit the symmetries of the underlying steady state and therefore belong to the same functional class. In particular, the parity and regularity conditions imposed on the rescaled variables are inherited by $\textcolor{blue}{\delta u}$, $\textcolor{blue}{\delta\omega}$, and $\textcolor{blue}{\delta\psi}$, with the latter two coupled by the elliptic equation relating vorticity and streamfunction. The perturbations must also decay at the far-field. The admissible perturbation class includes the weighted regularity, integrability, trace, and boundary-flux assumptions used in the energy estimates. In particular, the cutoff boundary-flux terms vanish as the cutoffs are removed. This ensures that the perturbation equations remain in the same functional class as the steady profile. \\

\noindent For notational convenience, we will use $\textcolor{blue}{\delta}$ as shorthand for the full list of perturbation variables,
\begin{equation}
\label{eq:delta_shorthand_intro}
    \textcolor{blue}{\delta}
    :=
    \{ \textcolor{blue}{\delta u},\textcolor{blue}{\delta\omega},\textcolor{blue}{\delta\psi},\textcolor{blue}{\delta C},
    \textcolor{blue}{\delta c_u},\textcolor{blue}{\delta c_\omega} \}.
\end{equation}

\noindent We also write the meridional velocity as
\begin{equation} \label{eq: meridional}
    u^r=\bar u^r+\textcolor{blue}{\delta u^r},
    \qquad
    u^z=\bar u^z+\textcolor{blue}{\delta u^z},
\end{equation}
where
\begin{equation} \label{eq:velocity_perturbation_intro}
\textcolor{blue}{\delta u^r} = -r\,\partial_z\textcolor{blue}{\delta\psi},
\qquad
\textcolor{blue}{\delta u^z} = 2\,\textcolor{blue}{\delta\psi} + r\,\partial_r\textcolor{blue}{\delta\psi}, 
\end{equation}
and also have
\begin{equation}
\bar c_\omega=\bar c_u-\lambda,
\qquad
\textcolor{blue}{\delta c_\omega}=\textcolor{blue}{\delta c_u}.
\end{equation}

\hfill \\

\noindent \textbf{Modulations and Normalization Conditions.} \, The modulation parameters are fixed by normalization conditions at the origin. We impose
\begin{equation}
    C(t)+u^z_0(t)=0,
    \qquad
    u_0(t)=\bar u_0.
\end{equation}

\hfill 

\noindent For the perturbations, these conditions give
\begin{equation} \label{eq: perturbation normalization}
    \textcolor{blue}{\delta C}+2\textcolor{blue}{\delta\psi}_0=0,
    \qquad
    \textcolor{blue}{\delta u}_0 = 0,
    \qquad
    (\partial_t\textcolor{blue}{\delta u})_0=0.
\end{equation}

\hfill

\noindent These normalization conditions fix the axial translation and amplitude normalization associated with $C(t)$ and $c_u(t)$. The spatial rate $\lambda=1/2$ has already been fixed. They are essential for formulating linear stability in the modulated variables, since otherwise the perturbation could drift along symmetry directions rather than measuring a genuine change in the profile.

\hfill \\

\subsubsection{Choice of Variables for the Linear Stability Analysis}

The reduced variables $u$ and $\omega$ do not transform in the same way under rescaling (see Appendix~\ref{appx: linear stability variables}). Indeed, the scaling of the full three-dimensional fields induces the following scaling of the reduced variables:
\begin{equation}
    u_\eta(r,z,t)=\eta^{\alpha+1}u(\eta r,\eta z,\eta^{\alpha+1}t),
    \qquad
    \omega_\eta(r,z,t)=\eta^{\alpha+2}\omega(\eta r,\eta z,\eta^{\alpha+1}t).
    \label{eq:euler-main-scaling}
\end{equation}

\noindent The relative scaling of $u$, $\omega$, and $\nabla u$ is unchanged. Thus $\omega$ carries one extra factor of $\eta$ relative to $u$, reflecting the one-derivative relation between the underlying full velocity and vorticity fields. By contrast,
\begin{equation}
    \omega_\eta(r,z,t)=\eta^{\alpha+2}\omega(\eta r,\eta z,\eta^{\alpha+1}t),
    \qquad
    \nabla u_\eta(r,z,t)=\eta^{\alpha+2}(\nabla u)(\eta r,\eta z,\eta^{\alpha+1}t).
\end{equation}

\hfill

\noindent Therefore $\omega$ and $\nabla u$ are the natural pair of quantities in a scale-consistent linear stability analysis. Quantities with different scaling correspond to different differential orders, so comparing them directly mixes objects of different analytic strength. For this reason, instead of using the equations for $(\partial_t\textcolor{blue}{\delta \omega},\partial_t\textcolor{blue}{\delta u})$, the stability analysis is formulated using the equations for $    \partial_t\textcolor{blue}{\delta\omega},
 \,
    \partial_t\partial_z\textcolor{blue}{\delta u},
 \,
    \partial_t\partial_r\textcolor{blue}{\delta u}$. We collect the associated perturbations into the scale-consistent perturbation vector
\begin{equation}
\label{eq:ns-deltav-class}
\textcolor{blue}{\delta v}
:=
(\textcolor{blue}{\delta\omega},\nabla\textcolor{blue}{\delta u})
=
(\textcolor{blue}{\delta\omega},\partial_r\textcolor{blue}{\delta u},\partial_z\textcolor{blue}{\delta u}),
\end{equation}
and for compact notation, we also set
\begin{equation}
\textcolor{blue}{\delta v_\omega}:=\textcolor{blue}{\delta\omega},
\qquad
\textcolor{blue}{\delta v_r}:=\partial_r\textcolor{blue}{\delta u},
\qquad
\textcolor{blue}{\delta v_z}:=\partial_z\textcolor{blue}{\delta u}.
\end{equation}

\hfill  

\subsubsection{Perturbation Equations}

Linearization consists of expanding the perturbation equation around the candidate blowup profile and retaining only the terms that are first order in the perturbation. The resulting linearized equation gives the leading-order dynamics of small disturbances near the profile. The discarded terms are collected in nonlinear remainders, which record the higher-order interactions among the perturbations. The full derivations of the linearized and modulation equations are given in Appendix~\ref{appx: linearization derivation}. \\

\noindent We start from the dynamic rescaling equations~\eqref{eq:dynamic_rescaled_u eps1}--\eqref{eq:dynamic_rescaled_psi eps1}:
\begin{align}
\partial_t u
+\bigl(\lambda r+ u^r\bigr)\partial_r u
+\bigl(\lambda z+C(t)+ u^z\bigr)\partial_z u
&=2u\partial_z\psi+c_u(t)u ,
\label{eq:dynamic_rescaled_u stab}\\
\partial_t\omega
+\bigl(\lambda r+ u^r\bigr)\partial_r\omega
+\bigl(\lambda z+C(t)+ u^z\bigr)\partial_z\omega
&=2u\partial_z u+c_\omega(t)\omega,
\label{eq:dynamic_rescaled_omega stab}\\
- \mathcal{E} \psi
&=\omega.
\label{eq:dynamic_rescaled_psi stab}
\end{align}
Here,
\begin{equation}
    u^r=-r\partial_z\psi,
    \qquad
    u^z=2\psi+r\partial_r\psi,
    \label{eq:dynamic_rescaled_stream_intro stab}
\end{equation}
and $\mathcal{E}$ is the elliptic operator
\begin{equation}
    \mathcal{E} \coloneqq \partial_r^2 + \frac{3}{r}\partial_r + \partial_z^2.
\end{equation}

\hfill \\

\noindent Substituting \eqref{eq:pert_ansatz_1}--\eqref{eq:pert_ansatz_2}, subtracting the steady-state equations, re-centering the modulation coefficients as specified below, and separating linear terms from nonlinear terms in the perturbations gives
\begin{equation}
\label{eq:linearized_system_intro}
    \partial_t\textcolor{blue}{\delta u}
    =
    \mathcal L_u(\textcolor{blue}{\delta}) + \mathcal N_u(\textcolor{blue}{\delta}) \ + \ \mathcal R_u,
    \qquad
    \partial_t\textcolor{blue}{\delta\omega}
    =
    \mathcal L_\omega(\textcolor{blue}{\delta}) + \mathcal N_\omega(\textcolor{blue}{\delta}) \ + \ \mathcal R_\omega.
\end{equation}

\hfill 

\noindent The raw profile residuals $\mathcal R_u^{\mathsf{raw}}$ and $\mathcal R_\omega^{\mathsf{raw}}$ account for the fact that the numerical profiles are not exact:
\begin{align}
\mathcal R_u^{\mathsf{raw}}
&=-(\lambda r+\bar u^r)\partial_r\bar u
-(\lambda z+\bar C+\bar u^z)\partial_z\bar u
+2\bar u\partial_z\bar\psi+\bar c_u\bar u,\\
\mathcal R_\omega^{\mathsf{raw}}
&=-(\lambda r+\bar u^r)\partial_r\bar\omega
-(\lambda z+\bar C+\bar u^z)\partial_z\bar\omega
+2\bar u\partial_z\bar u+(\bar c_u-\lambda)\bar\omega.
\end{align}

\hfill

\noindent We also define the residual amplitude offset
\begin{equation}
c_{u,\mathrm{res}}
:=-\frac{\mathcal R_u^{\mathsf{raw}}(0)}{\bar u_0},
\label{eq:residual_phase_offsets_intro}
\end{equation}
 under the nondegeneracy condition $\bar u_0\ne0$. \\
 
 \noindent We absorb this fixed offset into the reference amplitude rate. From this point onward, $\bvar{\delta C}$ and $\bvar{\delta c_u}$ denote the variable modulation coefficients, so that
\begin{equation}
C=\bar C+\bvar{\delta C},
\qquad
c_u=\bar c_u+c_{u,\mathrm{res}}+\bvar{\delta c_u},
\qquad
c_\omega=\bar c_\omega+c_{u,\mathrm{res}}+\bvar{\delta c_u}.
\label{eq:recentered_modulation_variables_intro}
\end{equation}

\hfill 

\noindent The residuals used in the perturbation equations are
\begin{align}
\mathcal R_u
&:=\mathcal R_u^{\mathsf{raw}}+c_{u,\mathrm{res}}\,\bar u,\\
\mathcal R_\omega
&:=\mathcal R_\omega^{\mathsf{raw}}+c_{u,\mathrm{res}}\,\bar\omega.
\label{eq:recentered_profile_residuals_intro}
\end{align}

\hfill  \\ 

\noindent The linearized $u$-operator and $\omega$-operator are
\vspace{2mm}

\begin{align}
\mathcal L_u(\textcolor{blue}{\delta})
:={}&
\bigl(2\partial_z\bar\psi+\bar c_u+c_{u,\mathrm{res}}\bigr)\,\textcolor{blue}{\delta u}
-\bigl(\lambda r+\bar u^r\bigr)\,\partial_r \textcolor{blue}{\delta u}
-\bigl(\lambda z + \bar C + \bar u^z\bigr)\,\partial_z \textcolor{blue}{\delta u}
+\bigl(2\bar u+ r\,\partial_r\bar u\bigr)\partial_z\textcolor{blue}{\delta\psi} \nonumber\\
&\quad
-\bigl( r\,\partial_z\bar u\bigr)\partial_r\textcolor{blue}{\delta\psi}
-\bigl(2\partial_z\bar u\bigr)\textcolor{blue}{\delta\psi}
-(\partial_z\bar u)\,\bvar{\delta C}
+\bar u\,\bvar{\delta c_u},
\label{eq:Lu_Euler_definition}
\end{align}
\begin{align}
\mathcal L_\omega(\textcolor{blue}{\delta})
:={}&
(\bar c_u+c_{u,\mathrm{res}}-\lambda)\,\textcolor{blue}{\delta\omega}
-\bigl(\lambda r + \bar u^r \bigr)\partial_r\textcolor{blue}{\delta\omega}
-\left(\lambda z+\bar C+ \bar u^z \right)\partial_z\textcolor{blue}{\delta\omega}
+2(\partial_z\bar u)\,\textcolor{blue}{\delta u}
+2\bar u\,\partial_z\textcolor{blue}{\delta u}
\nonumber\\
&\quad
-2(\partial_z\bar\omega)\,\textcolor{blue}{\delta\psi}
+ r\,(\partial_r\bar\omega)\,\partial_z\textcolor{blue}{\delta\psi}
- r\,(\partial_z\bar\omega)\,\partial_r\textcolor{blue}{\delta\psi}
-(\partial_z\bar\omega)\,\bvar{\delta C}
+\bar\omega\,\bvar{\delta c_u}.
\label{eq:Lw_Euler_definition}
\end{align}

\clearpage

\noindent The nonlinear remainders are
\vspace{1mm}

\begin{align}
\mathcal L_\omega(\textcolor{blue}{\delta})
={}&
 r\,\partial_z\textcolor{blue}{\delta\psi}\,\partial_r\textcolor{blue}{\delta u}
-\left(\bvar{\delta C}+\left(2\textcolor{blue}{\delta\psi}+r\partial_r\textcolor{blue}{\delta\psi}\right)\right)\partial_z\textcolor{blue}{\delta u}
+2\,\textcolor{blue}{\delta u}\,\partial_z\textcolor{blue}{\delta\psi}
+\bvar{\delta c_u}\,\textcolor{blue}{\delta u},
\label{eq:Nu_Euler_definition} \\ 
\mathcal N_\omega(\textcolor{blue}{\delta})
={}&
- \textcolor{blue}{\delta u^r} \partial_r\textcolor{blue}{\delta\omega}
-\left(\bvar{\delta C}+\textcolor{blue}{\delta u^z}\right)\partial_z\textcolor{blue}{\delta\omega}
+2\,\textcolor{blue}{\delta u}\,\partial_z\textcolor{blue}{\delta u}
+\bvar{\delta c_u}\,\textcolor{blue}{\delta\omega},
\label{eq:Nw_Euler_definition} 
\end{align}

\hfill

\noindent The elliptic relation~\eqref{eq:dynamic_rescaled_psi}   is already linear, so
\begin{equation}
\label{eq:linearized_elliptic_intro}
 -\mathcal E  \textcolor{blue}{\delta\psi}  \coloneqq   -\left(\partial_r^2 + \frac{3}{r}\partial_r + \partial_z^2\right)\textcolor{blue}{\delta\psi}
    = \textcolor{blue}{\delta\omega}.
\end{equation}

\hfill 

\noindent Consequently, $\textcolor{blue}{\delta\psi}$ is not an independent dynamical unknown: it is recovered from $\textcolor{blue}{\delta\omega}$ through this elliptic equation, and then $\textcolor{blue}{\delta u^r}$ and $\textcolor{blue}{\delta u^z}$ are recovered from \eqref{eq:velocity_perturbation_intro}. 

\hfill \\

\subsubsection{Differentiated Equations for the Scale-Consistent Variables}

Since the scale-consistent variables are $\omega$ and $\nabla u$, we differentiate the $\delta u$-equation:
\begin{align}
\label{eq: differentiated linearized}
\partial_t\partial_r\textcolor{blue}{\delta u}
&=\mathcal L_r(\textcolor{blue}{\delta})+
\mathcal N_{r}(\textcolor{blue}{\delta}) \ +\ \mathcal R_r,
\qquad \quad 
\partial_t\partial_z\textcolor{blue}{\delta u}
=\mathcal L_z(\textcolor{blue}{\delta})+
\mathcal N_{z}(\textcolor{blue}{\delta})  \ +\ \mathcal R_z,
\end{align}
where we denote
\begin{equation}
\mathcal L_r(\textcolor{blue}{\delta}):=\partial_r\mathcal L_u(\textcolor{blue}{\delta}),
\qquad
\mathcal N_r(\textcolor{blue}{\delta}):=\partial_r\mathcal N_u(\textcolor{blue}{\delta}), \qquad \mathcal L_z(\textcolor{blue}{\delta}):=\partial_z\mathcal L_u(\textcolor{blue}{\delta}),
\qquad
\mathcal N_z(\textcolor{blue}{\delta}):=\partial_z\mathcal N_u(\textcolor{blue}{\delta}),
\end{equation}
and 
\begin{equation}
    \mathcal R_z:=\partial_z\mathcal R_u,
\qquad\quad 
\,\mathcal R_r:=\,\partial_r\mathcal R_u.
\end{equation}

\hfill  

\noindent Differentiating the $\mathcal{L}_u$ expressions gives
\vspace{1mm}

\begin{align}
\mathcal L_r(\textcolor{blue}{\delta})
={}&
-(\lambda r+\bar u^r)\,\partial_{rr}\textcolor{blue}{\delta u}
-(\lambda z+\bar C+\bar u^z)\,\partial_{zr}\textcolor{blue}{\delta u}
+\bigl(2\partial_z\bar\psi+\bar c_u+c_{u,\mathrm{res}}-\lambda-\partial_r\bar u^r\bigr)
\partial_r\textcolor{blue}{\delta u}
\nonumber\\
&
-(\partial_r\bar u^z)\,\partial_z\textcolor{blue}{\delta u}  +2\partial_{rz}\bar\psi\,\textcolor{blue}{\delta u} +(2\bar u+ r\partial_r\bar u)\,\partial_{zr}\textcolor{blue}{\delta\psi}
+\bigl(3\partial_r\bar u+ r\partial_{rr}\bar u\bigr)
\partial_z\textcolor{blue}{\delta\psi}
- r(\partial_z\bar u)\,\partial_{rr}\textcolor{blue}{\delta\psi}
\nonumber\\
& -\bigl(3\partial_z\bar u+ r\partial_{rz}\bar u\bigr)
\partial_r\textcolor{blue}{\delta\psi}
-2\partial_{rz}\bar u\,\textcolor{blue}{\delta\psi}
-(\partial_{rz}\bar u)\,\bvar{\delta C}
+(\partial_r\bar u)\,\bvar{\delta c_u},
\label{eq:Lr_Euler}
\end{align}

\begin{align}
\mathcal L_z(\textcolor{blue}{\delta})
={}&
-(\lambda r+\bar u^r)\,\partial_{rz}\textcolor{blue}{\delta u}
-(\lambda z+\bar C+\bar u^z)\,\partial_{zz}\textcolor{blue}{\delta u}
+\bigl(2\partial_z\bar\psi+\bar c_u+c_{u,\mathrm{res}}-\lambda-\partial_z\bar u^z\bigr)
\partial_z\textcolor{blue}{\delta u}
\nonumber\\
& -(\partial_z\bar u^r)\,\partial_r\textcolor{blue}{\delta u} +2\partial_{zz}\bar\psi\,\textcolor{blue}{\delta u}
+(2\bar u+ r\partial_r\bar u)\,\partial_{zz}\textcolor{blue}{\delta\psi}
+\bigl(r\partial_{rz}\bar u\bigr)
\partial_z\textcolor{blue}{\delta\psi}
- r(\partial_z\bar u)\,\partial_{rz}\textcolor{blue}{\delta\psi}
\nonumber\\
& - r(\partial_{zz}\bar u)\,\partial_r\textcolor{blue}{\delta\psi}
-2\partial_{zz}\bar u\,\textcolor{blue}{\delta\psi}
-(\partial_{zz}\bar u)\,\bvar{\delta C}
+(\partial_z\bar u)\,\bvar{\delta c_u}.
\label{eq:Lz_Euler_explicit}
\end{align}

\hfill

\noindent The differentiated nonlinear remainders are

\begin{align}
\mathcal N_r(\textcolor{blue}{\delta})
={}&
\left(\partial_z\textcolor{blue}{\delta\psi}+r\partial_{rz}\textcolor{blue}{\delta\psi}\right)
\partial_r\textcolor{blue}{\delta u}
+ r\partial_z\textcolor{blue}{\delta\psi}\,\partial_{rr}\textcolor{blue}{\delta u}
-\left(3\partial_r\textcolor{blue}{\delta\psi}+r\partial_{rr}\textcolor{blue}{\delta\psi}\right)
\partial_z\textcolor{blue}{\delta u}
\label{eq:Nr_Euler}\\
& -\left(\bvar{\delta C}+\left(2\textcolor{blue}{\delta\psi}+r\partial_r\textcolor{blue}{\delta\psi}\right)\right)
\partial_{zr}\textcolor{blue}{\delta u}
+2(\partial_r\textcolor{blue}{\delta u})(\partial_z\textcolor{blue}{\delta\psi})
+2\textcolor{blue}{\delta u}\,\partial_{zr}\textcolor{blue}{\delta\psi}
+\bvar{\delta c_u}\,\partial_r\textcolor{blue}{\delta u},
\nonumber \end{align}
\begin{align} 
\mathcal N_z(\textcolor{blue}{\delta})
={}&
 r(\partial_{zz}\textcolor{blue}{\delta\psi})\,\partial_r\textcolor{blue}{\delta u}
+ r(\partial_z\textcolor{blue}{\delta\psi})\,\partial_{rz}\textcolor{blue}{\delta u}
-\left(2\partial_z\textcolor{blue}{\delta\psi}
+r\partial_{rz}\textcolor{blue}{\delta\psi}\right)\partial_z\textcolor{blue}{\delta u}
\label{eq:Nz_Euler} \\
& -\left(\bvar{\delta C}+\left(2\textcolor{blue}{\delta\psi}+r\partial_r\textcolor{blue}{\delta\psi}\right)\right)
\partial_{zz}\textcolor{blue}{\delta u}
+2(\partial_z\textcolor{blue}{\delta u})(\partial_z\textcolor{blue}{\delta\psi})
+2\textcolor{blue}{\delta u}\,\partial_{zz}\textcolor{blue}{\delta\psi}
+\bvar{\delta c_u}\,\partial_z\textcolor{blue}{\delta u}, \nonumber
\end{align}

\clearpage

\subsubsection{Modulation Equations}

The modulation equations are obtained in Appendix~\ref{appx: linearization derivation} by imposing the normalization conditions \eqref{eq: perturbation normalization}. They express the perturbations of the rescaling parameters in terms of the profile and its perturbations evaluated at the origin. The translation formula is the perturbation of the stagnation condition $C+u_0^z=0$. For the amplitude formula, we differentiate the normalization $u_0=\bar u_0$ in time, evaluate the $u$-equation at the meridional fixed point, and solve the resulting identity for $\bvar{\delta c_u}$.  \\

 \noindent With the fixed residual amplitude offset already absorbed into the reference coefficient, the modulation coefficients are
\begin{align}
\bvar{\delta C}
&=-2\bvar{\delta\psi}_0, \qquad\qquad  \bvar{\delta c_u} =-2(\partial_z\bvar{\delta\psi})_0.
\label{eq:ns_deltaC_deltacu}
\end{align}
Both $\bvar{\delta C}$ and $\bvar{\delta c_u}$ are linear in the perturbation. Products such as $\bvar{\delta C}\,\partial_z\bvar{\delta u}$ and $\bvar{\delta c_u}\,\bvar{\delta u}$ remain ordinary quadratic nonlinear terms. These equations determine the two modulation variables from the perturbations. 

\hfill

\subsubsection{Closed Linearized Variables}

 In terms of $\textcolor{blue}{\delta v}$, the closed linearized system may be written schematically as
\begin{equation}
\label{eq:closed-linearized-deltav}
\partial_t\textcolor{blue}{\delta v}
=
\mathcal L_{\nabla}\textcolor{blue}{\delta v},
\end{equation}
or equivalently in components:
\begin{equation}
\label{eq:closed-linearized-components}
\partial_t
\begin{pmatrix}
\textcolor{blue}{\delta\omega}\\
\partial_r\textcolor{blue}{\delta u}\\
\partial_z\textcolor{blue}{\delta u}
\end{pmatrix}
=
\mathcal L_{\nabla}
\begin{pmatrix}
\textcolor{blue}{\delta\omega}\\
\partial_r\textcolor{blue}{\delta u}\\
\partial_z\textcolor{blue}{\delta u}
\end{pmatrix}.
\end{equation}
The operator $\mathcal L_{\nabla}$ includes transport by the steady rescaled flow, coupling between $\textcolor{blue}{\delta u}$ and $\textcolor{blue}{\delta\omega}$, elliptic recovery of $\textcolor{blue}{\delta\psi}$, velocity reconstruction, and the linearized modulation terms. 

\hfill \\

\subsection{Weighted Norms and Energy}

\subsubsection{Weighted Norms}

\noindent The stability estimates are measured in weighted norms adapted to the candidate blowup profile and to the damping structure of the linearized operator. For this purpose, define
\begin{equation}
\wgray{\Phi_\omega},\ \wgray{\Phi_r},\ \wgray{\Phi_z}:\mathbb D\to(0,\infty)
\label{eq:low-order-weight-functions-main}
\end{equation}
to be the positive weight functions assigned respectively to the three components of $\textcolor{blue}{\delta v}$, on the half-plane
\begin{equation}
\mathbb D=\{(r,z):r\ge0,\ z\in\mathbb R\}.
\end{equation}

\hfill 

\noindent More explicitly, $\wgray{\Phi_\omega}$ weights the vorticity perturbation $\textcolor{blue}{\delta\omega}$, $\wgray{\Phi_r}$ weights $\partial_r\textcolor{blue}{\delta u}$, and $\wgray{\Phi_z}$ weights $\partial_z\textcolor{blue}{\delta u}$. We denote the collection of weights by
\begin{equation}
\label{eq:Phi-vector}
\wgray{\Phi}
:=
(\wgray{\Phi_\omega},\wgray{\Phi_r},\wgray{\Phi_z}).
\end{equation}

\hfill

\noindent For a positive scalar weight $\wgray{\phi}$, define the associated inner product and norm by
\begin{equation}
\label{eq:ns-weighted-inner}
\|f\|_{\wgray{\phi}}^2 := \int_\mathbb{D} |f|^2 \wgray{\phi}\,dr\,dz,
\qquad
\langle f,g\rangle_{\wgray{\phi}}:= \int_\mathbb{D} fg\,\wgray{\phi}\,dr\,dz.
\end{equation}

\hfill 

\noindent Given vector quantities $q := (q_\omega,q_r,q_z)$, and   $\tilde q := (\tilde q_\omega,\tilde q_r,\tilde q_z)$, we define the corresponding $\wgray{\Phi}$-weighted norm and $\wgray{\Phi}$-weighted inner product by 
\begin{equation}
\label{eq:Phi-vector-norm}
\|q\|_{\wgray{\Phi}}^2 := \|q_\omega\|_{\wgray{\Phi_\omega}}^2 + \|q_r\|_{\wgray{\Phi_r}}^2 + \|q_z\|_{\wgray{\Phi_z}}^2. \end{equation}
\begin{equation}
\label{eq:Phi-vector-inner}
\langle q,\tilde q\rangle_{\wgray{\Phi}} := \langle q_\omega,\tilde q_\omega\rangle_{\wgray{\Phi_\omega}} + \langle q_r,\tilde q_r\rangle_{\wgray{\Phi_r}} + \langle q_z,\tilde q_z\rangle_{\wgray{\Phi_z}}.
\end{equation}

\hfill

\subsubsection{Low-Order Weighted Energy.} 

\noindent The low-order energy uses the three singular weights $\wgray{\Phi_\omega}$, $\wgray{\Phi_r}$, $\wgray{\Phi_z}$, and is defined as the $\wgray{\Phi}$-norm of $\textcolor{blue}{\delta v}$:
\begin{equation}
\label{eq:ns-E0}
\mathfrak{E}_0(t) := \|\textcolor{blue}{\delta v}(t)\|_{\wgray{\Phi}}, \qquad \text{i.e.} \quad  \mathfrak{E}_0(t)^2 = \|\textcolor{blue}{\delta\omega}(t)\|_{\wgray{\Phi_\omega}}^2 + \|\partial_r\textcolor{blue}{\delta u}(t)\|_{\wgray{\Phi_r}}^2 + \|\partial_z\textcolor{blue}{\delta u}(t)\|_{\wgray{\Phi_z}}^2.
\end{equation}

\noindent The three weights $\wgray{\Phi_\omega}$, $\wgray{\Phi_r}$, $\wgray{\Phi_z}$ are kept separate because the damping calculation acts on the three components of $\textcolor{blue}{\delta v}$ in different ways. The purpose of the weighted energy is to make the coercive part of the linear part produce a positive damping constant $\DLamL>0$. A discussion of the constraints that the weight functions $\wgray{\Phi_\omega}$, $\wgray{\Phi_r}$, and $\wgray{\Phi_z}$ must satisfy in our setting is provided in \Cref{sec: Constraints on phi}. \\

\noindent While the low-order energy remains the main source of damping, it is not sufficient by itself to close the estimates since the linear and nonlinear terms contain terms which cannot be controlled using a weighted $L^2$ norm, such as terms that must be controlled in $L^\infty$, as well as modulation factors and pointwise quantities evaluated at the origin, such as $(\partial_z\bvar{\delta u})_0$ and $(\partial_{z}\bvar{\delta\omega})_0$.

\hfill

\subsubsection{High-Order Weighted Energy.} 

To control the linear and nonlinear terms which cannot be controlled in the low-order weighted energy $\mathfrak E_0$, we supplement $\mathfrak E_0$ with a high-order weighted energy $\mathfrak H_k$, chosen so that the top derivatives provide the required pointwise, interpolation, elliptic, and modulation bounds. \\

\noindent The high-order part uses the \emph{higher-order weight functions}
\begin{equation} 
  \label{eq:ns-compensated-weights-intro}
 \wgray{\Phi_{i,k}}:=1+\rho^{\,k}\wgray{\Phi_i},
\qquad \rho(r,z):=r^2+z^2,
\qquad i\in\{\omega,r,z\}. 
\end{equation}

\noindent The factor $\rho^{\,k}=|(r,z)|^{2k}$ compensates for the loss of $k$ powers under differentiation near the origin. If $f$ has local order $|(r,z)|^m$, then $D^k f$ has local order $|(r,z)|^{m-k}$; consequently, the factors $|D^k f|^2\rho^{,k}\wgray{\Phi_i}$ and $|f|^2\wgray{\Phi_i}$ have the same local power of $|(r,z)|$. The added $1$ prevents the high-order weight from degenerating at the origin and gives the direct bound
\begin{equation}
\wgray{\Phi_{i,k}}=1+\rho^k\wgray{\Phi_i}\ge 1,
\qquad
\|D^k f\|_{L^2(\mathbb D)}\le \|D^k f\|_{\wgray{\Phi_{i,k}}}.
\label{eq:high-weight-unit-floor-main}
\end{equation}

\noindent We then define the \emph{high-order weighted energy}
\begin{equation}
\label{eq:ns-Hk}
\mathfrak H_k^2
:={}\sum_{|\alpha|=k}\|D^\alpha\bvar{\delta\omega}\|_{\wgray{\Phi_{\omega,k}}}^2
+\sum_{|\alpha|=k}\|D^\alpha\partial_r\bvar{\delta u}\|_{\wgray{\Phi_{r,k}}}^2 +\sum_{|\alpha|=k}\|D^\alpha\partial_z\bvar{\delta u}\|_{\wgray{\Phi_{z,k}}}^2. 
\end{equation}

\hfill  

\subsubsection{Full Energy.} The full energy used for the stability argument is obtained as
\begin{equation} 
\label{eq:ns-Ek}
\mathfrak E_k^2:=\mathfrak E_0^2+\mu_k\mathfrak H_k^2,
\qquad 0<\mu_k\le 1.
\end{equation}
The parameter $\mu_k$ is chosen only small enough so the low-order damping absorbs the low-order loss coming from the differentiated linear estimate, but also large enough so that the differentiated damping absorbs the top-order loss left by the low-order linear estimate. \\

\paragraph{Choosing $k$.} \noindent The order of the top derivative $k$ should be chosen large enough so that the energy controls all point values, products, commutators, and elliptic or streamfunction terms used in the estimates. The modulation formulas require the origin values
\begin{equation}
\begin{aligned}
&\bvar{\delta\psi}_0,
\qquad
(\partial_z\bvar{\delta\psi})_0,
\nsgreen{\qquad
(\mathcal E\bvar{\delta u})_0.}
\end{aligned}
\label{eq:modulation-origin-values-k-choice}
\end{equation}
\noindent In two dimensions, a point value of a derivative of order \(m\) is controlled by Sobolev norms through order \(m+2\).  The largest origin derivative in \eqref{eq:modulation-origin-values-k-choice} has order one, so the modulation formulas require \(k\ge3\). Independently, the reverse allocation in the top-order transport commutators places a factor of derivative order at most two in \(L^\infty\), and the corresponding \(H^2\)-to-\(L^\infty\) estimate again uses derivatives through order four. \textbf{The lowest admissible choice is thus} $\boldsymbol{k=4}$. Larger choices of $k$ are possible, but they increase combinatorially the number of profile derivatives and analytic constants that must be certified. \\

\paragraph{Remark on intermediate derivative orders.} \noindent Intermediate derivative orders $0<j<k$ do not need to be included explicitly in the definition of the full energy. Indeed, for each such $j$, we derive in a separate paper~\citep{EulerBlowupStability} the certified interpolation bound \begin{equation} \sum_{|\alpha|=j}\|D^\alpha f\|_{\wgray{\Phi_{i,j}}} \le \eta_{i,j}\sum_{|\alpha|=k}\|D^\alpha f\|_{\wgray{\Phi_{i,k}}} +I_{i,j,k}(\eta_{i,j})\|f\|_{\wgray{\Phi_i}}. \label{eq:intermediate-order-roadmap-main} \end{equation} \noindent Thus, every intermediate-order term is controlled by a combination of the top-order and low-order energy levels. In each product estimate, these two contributions are kept separate until Young's inequality is applied: the top-order term is absorbed by the available $\mathfrak H_k^2$ damping, while the low-order term is controlled by $\mathfrak E_0^2$. The finite family of parameters $\eta_{i,j}$ is chosen so that the total top-order interpolation contribution remains strictly below the available differentiated damping.

\hfill \\

\subsubsection{Constraints on the Weight Functions} \label{sec: Constraints on phi}

\noindent Each weight function \(\wgray{\Phi_i}\) for \(i\in\{\omega,r,z\}\), is required to satisfy the following conditions.
\begin{enumerate}[leftmargin=2em]
\vspace{3mm}

\item \textbf{Positivity and coercivity.}
\noindent There is a constant \(C_{\rm floor}>0\) such that
\begin{equation}
\wgray{\Phi_i}(r,z)\ge C_{\rm floor}
\qquad\text{for almost every }(r,z)\in\mathbb D.
\label{eq:weight_coercivity}
\end{equation}
\noindent This prevents the weighted norm from becoming weak in a region where the perturbation must still be controlled.

\vspace{5mm}

\item \textbf{Near-origin admissibility.}
\noindent Each weight function $\wgray{\Phi_i}$ may have its own singularity rate at the origin. We require only local integrability, without imposing a common exponent or a comparison between different components. For certification, one convenient sufficient condition is
\begin{equation}
\Phi_i(r,z)\le C_i(r^2+z^2)^{-p_i},
\qquad
0<r^2+z^2\le \rho_i,
\qquad
0\le p_i<1,
\label{eq:origin-weight-power-model-revision}
\end{equation}
\noindent where $C_i$, $\rho_i$, and $p_i$ are chosen separately for each component. The condition $p_i<1$ guarantees local integrability in two dimensions. This power-law model is only a sufficient criterion, and a direct local-integrability certificate may be used instead.

\vspace{5mm}

\item \textbf{Regularity.}
\noindent Away from the origin, $C^1$ regularity is enough for the basic weighted integrations by parts. When a later estimate differentiates $\wgray{\Phi_i}$ or $\Phi_{i,k}$ further, the required derivatives are assumed to exist on that region and to satisfy the corresponding certified bounds. In the computer-assisted implementation, we use smooth parametrizations with $C^\infty$ activation functions, which meet these requirements and are convenient for gradient-based optimization.

\vspace{5mm}

\item \textbf{Far-field growth needed for velocity control.}
\noindent There are constants \(\rho_{\mathsf{growth}},\mathsf c_{\mathsf{growth}}>0\) such that
\begin{equation}
\min\{\Phi_r,\Phi_z\}\ge \mathsf c_{\mathsf{growth}} \, \rho
\qquad\text{whenever }\rho\ge\rho_{\mathsf{growth}},
\qquad
\rho=r^2+z^2.
\label{eq:delta-u-tail-growth-main}
\end{equation}
\noindent This lower growth condition is used to recover the unweighted \(L^2\) control of \(\bvar{\delta u}\) needed in the global pointwise estimate. 

\vspace{5mm}

\item (Optional) \textbf{Normalization.}
\noindent When the weight functions are optimized via gradient-based optimization, it can be useful to fix a normalization for each weight function to remove the irrelevant scaling freedom.

\vspace{5mm}

\item (Optional) \textbf{Evenness in \(r\).}
\noindent Since the pointwise estimates use reflection across the axis, it is convenient to choose the weights with the same even symmetry in $r$:
\begin{equation}
\wgray{\Phi_i}(r,z)=\wgray{\Phi_i}(-r,z).
\end{equation}
\noindent This symmetry can be built directly into the weight parametrization. 
\end{enumerate}

\hfill \\

\subsection{Linear Stability Estimates} \label{sec: meaning linear stability}

\subsubsection{Low-Order Energy Linear Stability and the Need for Higher-Order Energy} 

Linear stability asks whether perturbations of solutions of the closed linearized system remain bounded or decay in a suitable norm. In terms of $\textcolor{blue}{\delta v}$, the desired estimate is
\begin{equation}
\label{eq:semigroup_decay_intro_deltav}
\|\textcolor{blue}{\delta v}(t)\|_{\wgray{\Phi}} \leq C e^{-\gamma t} \|\textcolor{blue}{\delta v}(0)\|_{\wgray{\Phi}}, \qquad \gamma>0.
\end{equation}

\hfill 

\noindent A standard rigorous way to prove \eqref{eq:semigroup_decay_intro_deltav} is to establish the differential energy estimate
\begin{equation}
\label{eq:energy_decay_intro_E0}
\frac{1}{2}\frac{d}{dt}\mathfrak{E}_0(t)^2
=
\left\langle
\mathcal L_{\nabla}\textcolor{blue}{\delta v},
\textcolor{blue}{\delta v}
\right\rangle_{\wgray{\Phi}}
\leq
-\DLamL \mathfrak{E}_0(t)^2.
\end{equation}

\noindent This estimate gives a quantitative estimate for every solution of the closed linearized perturbation equation in the specified weighted space. Integrating \eqref{eq:energy_decay_intro_E0} gives exponential decay of $\mathfrak{E}_0(t)$, and hence of the scale-consistent perturbation $\textcolor{blue}{\delta v}$, rather than merely a numerical or visual indication of convergence. \\ 

\noindent For a blowup profile, this is the relevant notion of stability: the original physical solution may still blow up, but the normalized spatial profile is stable in the rescaled variables. Linear stability means that $\textcolor{blue}{\delta u}(t) \rightarrow0$ and $\textcolor{blue}{\delta\omega}(t) \rightarrow 0$ under the linearized rescaled flow. Note that the blowup time, blowup location, and scaling parameters may change under perturbation, but these changes are absorbed by the modulation parameters. \\

\noindent In spectral language, the idea is that the relevant spectrum of the closed linearized operator $\mathcal L_{\nabla}$ should lie in the stable half-plane. In a numerical stability calculation, one discretizes $\mathcal L_{\nabla}$ and studies the eigenvalues of the resulting matrix. Eigenvalues with negative real part correspond to decaying linear modes, while eigenvalues with positive real part correspond to growing directions. Neutral modes may arise from translation or amplitude normalization. The modulation conditions are designed to fix these directions, so that the stability question concerns perturbations transverse to the artificial symmetry directions. Thus, the role of the linearization is to produce the operator $\mathcal L_{\nabla}$ governing the first-order evolution of the scale-consistent perturbation $\textcolor{blue}{\delta v}$. The role of the linear stability analysis is to show that this operator has no growing directions in the weighted energy $\mathfrak{E}_0$, once the elliptic constraint, the velocity reconstruction, and the modulation conditions are imposed. If such an estimate is proved, then the proposed rescaled profile is linearly stable as a blowup profile. \\

\noindent However, the low-order energy is not sufficient by itself to close the linear estimates since the linear operators contain terms which cannot be controlled using a weighted $L^2$ norm. It is thus supplemented with the high-order weighted energy $\mathfrak H_k$, to get an estimate of the form
\begin{equation}  
\frac{1}{2}\frac{d}{dt}\mathfrak{E}_0(t)^2
=
\left\langle
\mathcal L_{\nabla}\textcolor{blue}{\delta v},
\textcolor{blue}{\delta v}
\right\rangle_{\wgray{\Phi}}
\leq -\DLamLlow\mathfrak E_0^2+\DLamLhigh\mathfrak H_k^2.
\end{equation}

\hfill \\

\subsubsection{Low-Order Linear Estimates}

\noindent The low-order linear estimate to certify in the Euler setting is
\begin{equation}  \ \ 
\label{eq:ns-CL}
\langle \mathcal L_\omega,\bvar{\delta\omega}\rangle_{\wgray{\Phi_\omega}}
+\langle \partial_r\mathcal L_u,\partial_r\bvar{\delta u}\rangle_{\wgray{\Phi_r}}
+\langle \partial_z\mathcal L_u,\partial_z\bvar{\delta u}\rangle_{\wgray{\Phi_z}}
\le -\DLamLlow \, \mathfrak E_0^2+\DLamLhigh \, \mathfrak H_k^2. \ \ 
\end{equation}
Here the constant \(\DLamLlow>0\) is the net low-order linear damping and \(\DLamLhigh\ge0\) is the top-order loss from direct estimates that require \(\mathfrak H_k\).

\hfill \\

\subsubsection{High-Order Linear Estimates}

\noindent We need to certify a top-order damping constant $\DLamDLhigh>0$ and a low-order constant $\DLamDLlow\ge0$ such that
\begin{equation} \ \ 
\label{eq:ns-high-lin}
\begin{aligned}
&\sum_{|\alpha|=k}\langle D^\alpha\mathcal L_\omega,D^\alpha\bvar{\delta\omega}\rangle_{\wgray{\Phi_{\omega,k}}}
+\sum_{|\alpha|=k}\langle D^\alpha\partial_r\mathcal L_u,D^\alpha\partial_r\bvar{\delta u}\rangle_{\wgray{\Phi_{r,k}}}
+\sum_{|\alpha|=k}\langle D^\alpha\partial_z\mathcal L_u,D^\alpha\partial_z\bvar{\delta u}\rangle_{\wgray{\Phi_{z,k}}}
 \\
&\qquad\le \DLamDLlow \, \mathfrak E_0^2 -\DLamDLhigh \, \mathfrak H_k^2.
\end{aligned} \ \ 
\end{equation}

\hfill \\

\subsection{Nonlinear Stability Estimates}

The linear stability analysis aims to show that the rescaled profile is stable under the first-order perturbation dynamics, measured in the weighted low-order norm associated with $\wgray{\Phi_\omega}$, $\wgray{\Phi_r}$, and $\wgray{\Phi_z}$. The nonlinear stability argument upgrades this statement to the full perturbation equation. Its goal is to prove that, if the initial perturbation is sufficiently small in a suitable energy, then the perturbation remains small for as long as the estimates apply. In the case of an approximate numerical profile, the solution is allowed to remain in a small neighborhood whose size also depends on the certified residual error. \\ 

\noindent The nonlinear stability argument is reduced to an energy estimate for the full energy $\mathfrak E_k$. Differentiating $\mathfrak E_k^2$ in time and inserting the perturbation equations separates the argument into five distinct components: the residual error, the low-order linear part, the high-order linear part, the low-order nonlinear part, and the high-order nonlinear part. The residual records the defect of the approximate steady profile. The linear terms contain the stabilizing mechanism: the low-order estimate gives the principal damping, while the high-order estimate controls the differentiated equations, allowing only losses that can be absorbed by the low-order energy after choosing $\mu_k$. The nonlinear terms are then estimated perturbatively, with the low-order and high-order contributions bounded by cubic and quartic powers of the full energy.

\hfill

\subsubsection{Low-Order Nonlinear Estimates}

We start from the expressions \eqref{eq:Nw_Euler_definition}, \eqref{eq:Nr_Euler}, \eqref{eq:Nz_Euler}, for the nonlinear parts $\mathcal{N}_i$ for $i\in \{\omega,r,z\}$. \\

\noindent We substitute the exact modulation formulas before separating linear and nonlinear terms:
\begin{equation}
\bvar{\delta C}=-2\bvar{\delta\psi}_0,
\qquad \qquad
\bvar{\delta c_u}
=-2(\partial_z\bvar{\delta\psi})_0.
\end{equation}

\hfill

\noindent The precise low-order nonlinear estimate to prove is
\begin{equation}  \ \ 
\label{eq:ns-low-nonlinear-final}
\left|\langle\mathcal N_\omega,\bvar{\delta\omega}\rangle_{\wgray{\Phi_\omega}}
+\langle\mathcal N_{r},\partial_r\bvar{\delta u}\rangle_{\wgray{\Phi_r}}
+\langle\mathcal N_{z},\partial_z\bvar{\delta u}\rangle_{\wgray{\Phi_z}}
\right|
\le \DLamNthree \, \mathfrak E_k^3. 
\end{equation}

\hfill \\ 

\subsubsection{High-Order Nonlinear Estimates}

\hfill 

\noindent The precise high-order nonlinear estimate to prove is

\begin{equation}
\label{eq:ns-high-nonlinear}
\mu_k\left|
\sum_{i\in\{\omega,r,z\}}\sum_{|\alpha|=k}
\left\langle D^\alpha\mathcal N_i,D^\alpha\bvar{\delta v_i}\right\rangle_{\wgray{\Phi_{i,k}}}
\right|
\le \DLamDNthree \, \mathfrak E_k^3.
\end{equation}

\hfill \\

\subsection{PDE Residuals}

If the $(\bar u,\bar\omega,\bar\psi,\bar\lambda,\bar C,\bar c_u)$ profile is exact, the raw residuals below vanish. If the profile is numerical, these are the defects obtained by substituting the profile into the steady equations:
\begin{align}
\label{eq:ns-Ru-raw}
\mathcal R_u^{\mathsf{raw}}
&:=-(\lambda r+\bar u^r)\partial_r\bar u
-(\lambda z+\bar C+\bar u^z)\partial_z\bar u
+2\bar u\partial_z\bar\psi+\bar c_u\bar u,\\
\label{eq:ns-Romega-raw}
\mathcal R_\omega^{\mathsf{raw}}
&:=-(\lambda r+\bar u^r)\partial_r\bar\omega
-(\lambda z+\bar C+\bar u^z)\partial_z\bar\omega
+2\bar u\partial_z\bar u+(\bar c_u-\lambda)\bar\omega.
\end{align}

\hfill

\noindent The residuals entering the perturbation equations are
\begin{align}
\label{eq:ns-Ru}
\mathcal R_u
&:=\mathcal R_u^{\mathsf{raw}}
+\cures\,\bar u,\\
\label{eq:ns-Romega}
\mathcal R_\omega
&:=\mathcal R_\omega^{\mathsf{raw}}
+\cures\,\bar\omega.
\end{align}

\hfill

\noindent The elliptic residual is not listed here because the spline representation of 
\(\bar\omega\) is obtained by applying the elliptic operator to \(\bar\psi\) exactly, as discussed in \Cref{sec: analytic representation splines}. As a result, the steady-state profile satisfies the elliptic equation exactly, and the elliptic residual is exactly 0.

\hfill \\

\noindent Define
\begin{equation}
\label{eq:ns-Rrz}
\mathcal R_r:=\partial_r\mathcal R_u,
\qquad
\mathcal R_z:=\partial_z\mathcal R_u.
\end{equation}

\hfill 

\noindent Then the low-order perturbation equations are
\begin{align}
\label{eq:ns-component-eqns}
\partial_t\bvar{\delta\omega}
&=\mathcal L_\omega(\bvar{\delta})+
\mathcal N_\omega(\bvar{\delta})+\mathcal R_\omega,\\
\partial_t\partial_r\bvar{\delta u}
&=\partial_r\mathcal L_u(\bvar{\delta})+
\mathcal N_{r}(\bvar{\delta})+\mathcal R_r,\\
\partial_t\partial_z\bvar{\delta u}
&=\partial_z\mathcal L_u(\bvar{\delta})+
\mathcal N_{z}(\bvar{\delta})+\mathcal R_z.
\end{align}

\hfill 

\noindent Define the residual constants by separating the low-order and top-order pieces:
\begin{align}
\label{eq:ns-eps0}
\DLamRPhi^2
&:=\|\mathcal R_\omega\|_{\wgray{\Phi_\omega}}^2
+\|\mathcal R_r\|_{\wgray{\Phi_r}}^2
+\|\mathcal R_z\|_{\wgray{\Phi_z}}^2,\\
\label{eq:ns-epsk}
\DLamDR^2
&:=\sum_{|\alpha|=k}
\Bigl(
\|D^\alpha\mathcal R_\omega\|_{\wgray{\Phi_{\omega,k}}}^2
+\|D^\alpha\mathcal R_r\|_{\wgray{\Phi_{r,k}}}^2
+\|D^\alpha\mathcal R_z\|_{\wgray{\Phi_{z,k}}}^2
\Bigr),\\
\label{eq:ns-eps-total}
\DLamR^2
&:=\DLamRPhi^2+\mu_k\DLamDR^2.
\end{align}

\hfill 

\noindent Applying the Cauchy--Schwarz inequality gives
\begin{equation} 
\label{eq:ns-res-pair}
\left| \,
\sum_{i\in\{\omega,r,z\}}
\langle\mathcal R_i,\bvar{\delta v_i}\rangle_{\wgray{\Phi_i}}
+\mu_k\sum_{i\in\{\omega,r,z\}}\sum_{|\alpha|=k}
\langle D^\alpha\mathcal R_i,D^\alpha\bvar{\delta v_i}\rangle_{\wgray{\Phi_{i,k}}} \,
\right|
\le \DLamR \, \mathfrak E_k. 
\end{equation}

\hfill \\

\subsection{Stability Theorems}

\hfill

\noindent We are now ready to state the stability results for the Euler equations. Importantly, every named certification constant attached to an analytic estimate is understood to be a rigorously validated finite numerical value. Conceptually, the analytic argument identifies the constants that are needed, while the numerical step certifies them once for the chosen profiles, weights, and parameters. The stability proofs are not complete until appropriate constants have been certified numerically. \\

\subsubsection{Single-Radius Stability Theorems}

\noindent We now state the full stability theorem.\\

\begin{theorem}[Single-radius stability]\label{thm:main-stability}
Assume that the following perturbation equations hold for admissible perturbations:
\begin{align}
\partial_t\bvar{\delta\omega}
&=\mathcal L_\omega(\bvar{\delta})+
\mathcal N_\omega(\bvar{\delta})+\mathcal R_\omega, \label{eq: Stability theorem perturbation 1}\\
\partial_t\partial_r\bvar{\delta u}
&=\partial_r\mathcal L_u(\bvar{\delta})+
\mathcal N_{r}(\bvar{\delta})+\mathcal R_r,\label{eq: Stability theorem perturbation 2}\\
\partial_t\partial_z\bvar{\delta u}
&=\partial_z\mathcal L_u(\bvar{\delta})+
\mathcal N_{z}(\bvar{\delta})+\mathcal R_z. \label{eq: Stability theorem perturbation 3}
\end{align}

\hfill 

\noindent Assume that the following estimates have been certified with finite constants:

\begin{equation}
\begin{aligned}
&\langle \mathcal L_\omega,\bvar{\delta\omega}\rangle_{\wgray{\Phi_\omega}}
+\langle \partial_r\mathcal L_u,\partial_r\bvar{\delta u}\rangle_{\wgray{\Phi_r}}
+\langle \partial_z\mathcal L_u,\partial_z\bvar{\delta u}\rangle_{\wgray{\Phi_z}}
\le -\DLamLlow \, \mathfrak E_0^2+\DLamLhigh \, \mathfrak H_k^2.
\end{aligned}
\label{eq:ns-thm-low-linear}
\end{equation}

\vspace{-1mm}

\begin{equation}
\begin{aligned}
&\sum_{i\in\{\omega,r,z\}}\sum_{|\alpha|=k}
\langle D^\alpha\mathcal L_i,D^\alpha\bvar{\delta v_i}\rangle_{\wgray{\Phi_{i,k}}}
\le \DLamDLlow \, \mathfrak E_0^2-\DLamDLhigh \, \mathfrak H_k^2.
\end{aligned}
\label{eq:ns-thm-high-linear}
\end{equation}

\vspace{-1mm}

\begin{equation}
\begin{aligned}
&\left|\langle\mathcal N_\omega,\bvar{\delta\omega}\rangle_{\wgray{\Phi_\omega}}
+\langle\mathcal N_{r},\partial_r\bvar{\delta u}\rangle_{\wgray{\Phi_r}}
+\langle\mathcal N_{z},\partial_z\bvar{\delta u}\rangle_{\wgray{\Phi_z}}
\right|
\le \DLamNthree \, \mathfrak E_k^3.
\end{aligned}
\label{eq:ns-thm-low-nonlinear}
\end{equation}

\vspace{-1mm}

\begin{equation}
\mu_k \, \Bigg|
\sum_{i\in\{\omega,r,z\}}\sum_{|\alpha|=k}
\langle D^\alpha\mathcal N_i,D^\alpha\bvar{\delta v_i}\rangle_{\wgray{\Phi_{i,k}}}
\Bigg|
\le \DLamDNthree \, \mathfrak E_k^3.
\label{eq:ns-thm-high-nonlinear}
\end{equation}

\vspace{-1mm}

\begin{equation}
\begin{aligned}
&\Bigg|
\sum_{i\in\{\omega,r,z\}}\langle \mathcal R_i,\bvar{\delta v_i}\rangle_{\wgray{\Phi_i}}
+\mu_k\sum_{i\in\{\omega,r,z\}}\sum_{|\alpha|=k}
\langle D^\alpha\mathcal R_i,D^\alpha\bvar{\delta v_i}\rangle_{\wgray{\Phi_{i,k}}}
\Bigg|
\le \DLamR \, \mathfrak E_k.
\end{aligned}
\label{eq:ns-thm-residual}
\end{equation}

\hfill \\

\noindent Choose $\mu_k$ so that $0<\mu_k\le 1$ and
\begin{equation}
\label{eq:ns-mu-choice}
\DLamStab(\mu_k):=
\min\left\{
\DLamLlow-\mu_k\DLamDLlow,
\ \DLamDLhigh-\frac{\DLamLhigh}{\mu_k}
\right\}>0.
\end{equation}

\hfill

\noindent After this choice of $\mu_k$ is fixed, write $\DLamStab$ for the positive number $\DLamStab(\mu_k)$.  Define the combined nonlinear constants by
\begin{equation}
\label{eq:ns-N2N3}
\DLamThree:=\DLamNthree+\DLamDNthree.
\end{equation}

\hfill 

\noindent Assume that $\delta_*>0$ satisfies

\begin{equation}
\label{eq:ns-smallness}
\DLamThree\delta_*+\frac{\DLamR}{\delta_*}<\DLamStab.
\end{equation}

\hfill 

\noindent Then every solution with $\mathfrak E_k(0)<\delta_*$ remains in the ball $\mathfrak E_k(t)<\delta_*$ for as long as the solution exists and the certified estimates apply.
\end{theorem}

\begin{grayproof}
The proof is presented in Appendix~\ref{appx: nonlinear stability proof}. 
\end{grayproof}

\clearpage 

\noindent The diagram below summarizes the strategy used for proving nonlinear stability.

\begin{figure}[h]
\centering
\begingroup

\definecolor{eulerstaborange}{RGB}{230,120,20}
\definecolor{eulerstabpurple}{RGB}{15,25,100}
\definecolor{eulerstabsalmon}{RGB}{240,80,70}
\definecolor{eulerstabtextmuted}{RGB}{82,82,88}
\definecolor{eulerstabrulegray}{RGB}{140,140,148}
\definecolor{eulerstabrulemid}{RGB}{176,176,183}

\newcommand{\eulerstabLamLlow}{%
  \ensuremath{\textcolor{eulerstabsalmon}
  {\Lambda_{\mathcal L}^{\rm low}}}%
}

\newcommand{\eulerstabLamLhigh}{%
  \ensuremath{\textcolor{eulerstabsalmon}
  {\Lambda_{\mathcal L}^{\rm high}}}%
}

\newcommand{\eulerstabLamDLhigh}{%
  \ensuremath{\textcolor{eulerstabsalmon}
  {\Lambda_{D^\alpha\mathcal L}^{\rm high}}}%
}

\newcommand{\eulerstabLamDLlow}{%
  \ensuremath{\textcolor{eulerstabsalmon}
  {\Lambda_{D^\alpha\mathcal L}^{\rm low}}}%
}

\newcommand{\eulerstabLamNthree}{%
  \ensuremath{\textcolor{eulerstabsalmon}
  {\Lambda_{\mathcal N,3}}}%
}

\newcommand{\eulerstabLamNfour}{%
  \ensuremath{\textcolor{eulerstabsalmon}
  {\Lambda_{\mathcal N,4}}}%
}

\newcommand{\eulerstabLamDNthree}{%
  \ensuremath{\textcolor{eulerstabsalmon}
  {\Lambda_{D^\alpha\mathcal N,3}}}%
}

\newcommand{\eulerstabLamDNfour}{%
  \ensuremath{\textcolor{eulerstabsalmon}
  {\Lambda_{D^\alpha\mathcal N,4}}}%
}

\newcommand{\eulerstabLamR}{%
  \ensuremath{\textcolor{eulerstabsalmon}
  {\Lambda_{\mathcal R}}}%
}

\newcommand{\eulerstabLamStab}{%
  \ensuremath{\textcolor{eulerstabsalmon}
  {\Lambda_{\rm stab}}}%
}

\tikzset{
  eulerstab flow/.style={
    -{Stealth[length=1.55mm,width=1.08mm]},
    draw=eulerstabrulegray,
    line width=0.46pt
  },
  eulerstab transition/.style={
    midway,
    right=3.5pt,
    fill=white,
    inner xsep=1.10pt,
    inner ysep=0.65pt,
    text=eulerstabtextmuted,
    font=\fontsize{9.6}{10.6}\selectfont\itshape
  },
  eulerstab stage/.style={
    draw=eulerstabrulemid,
    fill=white,
    line width=0.40pt,
    minimum width=15.75cm,
    text width=14.90cm,
    minimum height=17.0mm,
    align=center,
    inner xsep=8.0pt,
    inner ysep=5.8pt,
    outer sep=0pt
  }
}

\newcommand{\eulerstabtitle}[2]{%
  {\fontsize{12.4}{13.4}\selectfont
   \bfseries\color{#1}\textsc{#2}\par
   \vspace{6.6pt}}%
}

\newcommand{\eulerstabtitlecompact}[2]{%
  {\fontsize{12.4}{13.4}\selectfont
   \bfseries\color{#1}\textsc{#2}\par
   \vspace{5.0pt}}%
}

\newcommand{\eulerstabsubtitle}[1]{%
  {\fontsize{9.8}{11.0}\selectfont
   \color{eulerstabtextmuted}#1\par
   \vspace{-0.8pt}}%
}

\newcommand{\eulerstabphasetag}[2]{%
  {\fontsize{7.4}{8.2}\selectfont
   \bfseries\scshape\color{#1}#2}%
}

\newcommand{\eulerstabtagpanel}[6]{%
  \node[eulerstab stage,#6] (#1) {%
    \eulerstabtitle{#2}{#4}%
    \eulerstabsubtitle{#5}%
  };
  \draw[#2!78,line width=0.70pt]
    ([xshift=0.2pt]#1.north west) --
    ([xshift=-0.2pt]#1.north east);
  \node[
    anchor=west,
    fill=white,
    inner xsep=2.7pt,
    inner ysep=1.0pt,
    text=#2
  ] at ([xshift=7.5mm]#1.north west) {%
    \eulerstabphasetag{#2}{#3}%
  };
}

\newcommand{\eulerstabtagpanelcompact}[6]{%
  \node[eulerstab stage,#6] (#1) {%
    \eulerstabtitlecompact{#2}{#4}%
    \eulerstabsubtitle{#5}%
  };
  \draw[#2!78,line width=0.70pt]
    ([xshift=0.2pt]#1.north west) --
    ([xshift=-0.2pt]#1.north east);
  \node[
    anchor=west,
    fill=white,
    inner xsep=2.7pt,
    inner ysep=1.0pt,
    text=#2
  ] at ([xshift=7.5mm]#1.north west) {%
    \eulerstabphasetag{#2}{#3}%
  };
}

\hyphenpenalty=10000
\exhyphenpenalty=10000
\emergencystretch=1em

\resizebox{0.87\linewidth}{!}{%
\begin{tikzpicture}[node distance=5.1mm,font=\rmfamily]

  \eulerstabtagpanel
    {eulersetup}
    {blue}
    {Perturbations}
    {Perturbation Equations}
    {$\partial_t\bvar{\delta\omega}
      =\mathcal L_\omega(\bvar{\delta})
      +\mathcal N_\omega(\bvar{\delta})
      +\mathcal R_\omega$,\quad
     $\partial_t\partial_r\bvar{\delta u}
      =\partial_r\mathcal L_u(\bvar{\delta})
      +\mathcal N_r(\bvar{\delta})
      +\mathcal R_r$,\quad
     $\partial_t\partial_z\bvar{\delta u}
      =\partial_z\mathcal L_u(\bvar{\delta})
      +\mathcal N_z(\bvar{\delta})
      +\mathcal R_z$.\\[3pt]
     $-\mathcal E\bvar{\delta\psi} \!\!\!\!\!\! = \!\!\!\!\!\! \bvar{\delta\omega}$, \ \ 
     with modulation coefficients
     $\bvar{\delta C}$ and $\bvar{\delta c_u}$}
    {}

  \eulerstabtagpanel
    {weightedenergy}
    {eulerstabpurple}
    {Weights and Energy}
    {Stability Energy}
    {$\|f\|_{\phi}^{2} \!
      := \! \displaystyle\int_{\mathbb{D}}|f|^{2}\phi\,dr\,dz$,\quad  
     $\langle f,g\rangle_{\phi} \! 
      := \! \displaystyle\int_{\mathbb{D}}fg\,\phi\,dr\,dz$,\qquad 
     $\wgray{\Phi_{i,k}} \! 
      := \! 1+\rho^{\,k}\wgray{\Phi_i}$,\quad 
     $\rho \! := \! r^2+z^2$.\\[4pt]
     The singular base weights are
     $\wgray{\Phi_\omega},\wgray{\Phi_r},\wgray{\Phi_z}$,
     while $\wgray{\Phi_{i,k}}$ are the compensated top-order weights.\\[3pt]
     $\begin{aligned} \text{Low-Order Energy } \ \ 
      \mathfrak E_0^2
      &:=\|\bvar{\delta\omega}\|_{\wgray{\Phi_\omega}}^2
      +\|\partial_r\bvar{\delta u}\|_{\wgray{\Phi_r}}^2
      +\|\partial_z\bvar{\delta u}\|_{\wgray{\Phi_z}}^2,
      \\[1.8pt]
     \text{Top-Order Energy } \ \  \mathfrak H_k^2
      &:=\sum_{|\alpha|=k}
      \|D^\alpha\bvar{\delta\omega}\|_{\wgray{\Phi_{\omega,k}}}^2
      +\sum_{|\alpha|=k}
      \|D^\alpha\partial_r\bvar{\delta u}\|_{\wgray{\Phi_{r,k}}}^2
      +\sum_{|\alpha|=k}
      \|D^\alpha\partial_z\bvar{\delta u}\|_{\wgray{\Phi_{z,k}}}^2,
      \\[1.8pt]
      \text{Full Energy } \ \  \mathfrak E_k^2
      &:=\mathfrak E_0^2+\mu_k\mathfrak H_k^2,
      \qquad 0<\mu_k\le1
     \end{aligned}$}
    {below=of eulersetup,
     text width=15.35cm,
     inner xsep=4.5pt}

  \eulerstabtagpanel
    {eulertoolbox}
    {eulerstabpurple}
    {Theory}
    {Analytic Results and Estimates}
    {Sobolev extension and reflection theorems,
     Sobolev embeddings and interpolation,
     weighted transport, and Leibniz estimates control
     lower-order factors and intermediate derivatives.\\[1.8pt]
     Elliptic and streamfunction PDE bounds,
     modulation estimates, and certified operator and residual
     constants reduce all pairings to explicit bounds. }
    {below=of weightedenergy}

  \eulerstabtagpanel
    {eulerlinear}
    {eulerstabsalmon}
    {Linear}
    {Low-Order and High-Order Linear Estimates}
    {$\begin{aligned}
      \langle\mathcal L_\omega,\bvar{\delta\omega}\rangle_
        {\wgray{\Phi_\omega}}
      +\langle\partial_r\mathcal L_u,
        \partial_r\bvar{\delta u}\rangle_{\wgray{\Phi_r}}
      +\langle\partial_z\mathcal L_u,
        \partial_z\bvar{\delta u}\rangle_{\wgray{\Phi_z}}
      &\le
      -\eulerstabLamLlow\mathfrak E_0^2
      +\eulerstabLamLhigh\mathfrak H_k^2,
      \\[1.8pt]
      \sum_{i\in\{\omega,r,z\}}\sum_{|\alpha|=k}
      \langle D^\alpha\mathcal L_i,
        D^\alpha\bvar{\delta v_i}\rangle_{\wgray{\Phi_{i,k}}}
      &\le
      \eulerstabLamDLlow\mathfrak E_0^2
      -\eulerstabLamDLhigh\mathfrak H_k^2
    \end{aligned}$}
    {below=of eulertoolbox}

  \eulerstabtagpanel
    {eulernonlinear}
    {eulerstabsalmon}
    {Nonlinear}
    {Low-Order and High-Order Nonlinear Estimates}
    {$\begin{aligned}
      \bigl|
      \langle\mathcal N_\omega,\bvar{\delta\omega}\rangle_
        {\wgray{\Phi_\omega}}
      +\langle\mathcal N_r,\partial_r\bvar{\delta u}\rangle_
        {\wgray{\Phi_r}}
      +\langle\mathcal N_z,\partial_z\bvar{\delta u}\rangle_
        {\wgray{\Phi_z}}
      \bigr|
      &\le
      \eulerstabLamNthree\mathfrak E_k^3
      \\[1.8pt]
      \mu_k\Bigg|
      \sum_{i\in\{\omega,r,z\}}\sum_{|\alpha|=k}
      \langle D^\alpha\mathcal N_i,
        D^\alpha\bvar{\delta v_i}\rangle_{\wgray{\Phi_{i,k}}}
      \Bigg|
      &\le
      \eulerstabLamDNthree\mathfrak E_k^3
    \end{aligned}$}
    {below=2.0mm of eulerlinear}

  \eulerstabtagpanel
    {eulerresidual}
    {eulerstaborange}
    {PDE Residual}
    {Profile PDE Residual Estimates}
    {$\left|
      \sum_{i\in\{\omega,r,z\}}
      \langle\mathcal R_i,\bvar{\delta v_i}\rangle_{\wgray{\Phi_i}}
      +\mu_k
      \sum_{i\in\{\omega,r,z\}}\sum_{|\alpha|=k}
      \langle D^\alpha\mathcal R_i,
        D^\alpha\bvar{\delta v_i}\rangle_{\wgray{\Phi_{i,k}}}
      \right|
      \!\! \le \!\! \eulerstabLamR\mathfrak E_k$}
    {below=2.0mm of eulernonlinear}

  \eulerstabtagpanel
    {eulerclosure}
    {eulerstabpurple}
    {Closure}
    {Choose $\mu_k$ and the Stability Radius  $\delta_*$}
    {$\eulerstabLamStab(\mu_k) \!\!\!\!
      := \!\!\!\! \min\!\left\{
      \eulerstabLamLlow-\mu_k\eulerstabLamDLlow,\,
      \eulerstabLamDLhigh
      -\dfrac{\eulerstabLamLhigh}{\mu_k}
      \right\} \!\!\!\!\! > \!\!\!\!\! 0$,\\[1.8pt]
     choose \, $\delta_* \!\!\!\!\!\!\! > \!\!\!\!\!\!\!0$ \,  with \,
     $\bigl(\eulerstabLamNthree \!\!+\!\!\eulerstabLamDNthree\bigr) \,\delta_*
      \!\! + \!\! \dfrac{\eulerstabLamR }{\delta_*} \!\!\!\!\!
      < \!\!\!\!\! \eulerstabLamStab$}
    {below=of eulerresidual}

  \eulerstabtagpanelcompact
    {eulerproof}
    {eulerstabpurple}
    {Stability}
    {Stability}
    {If $\mathfrak E_k(0) \!\!\!\!<\!\!\!\! \delta_*$, then
     $\mathfrak E_k(t) \!\!\!\! < \!\!\!\! \delta_*$ for every later admissible time.\\[1.8pt]
     Thus the ball $\{\mathfrak E_k<\delta_*\}$ is forward invariant,
     which is the desired nonlinear stability statement}
    {below=of eulerclosure,
     inner ysep=6.8pt}

  \draw[eulerstabpurple!62,line width=0.28pt]
    ([xshift=1.9mm,yshift=-1.7mm]eulerproof.north west)
    rectangle
    ([xshift=-1.9mm,yshift=1.7mm]eulerproof.south east);

  \draw[eulerstab flow]
    (eulersetup.south) --
    node[eulerstab transition] {build the weighted energy}
    (weightedenergy.north);

  \draw[eulerstab flow]
    (weightedenergy.south) --
    node[eulerstab transition] {derive the analytic ingredients}
    (eulertoolbox.north);

  \draw[eulerstab flow]
    (eulertoolbox.south) --
    node[eulerstab transition] {obtain linear, nonlinear, and residual bounds}
    (eulerlinear.north);

  \draw[eulerstab flow]
    (eulerresidual.south) --
    node[eulerstab transition]
      {choose coupling and bootstrap radius}
    (eulerclosure.north);

  \draw[eulerstab flow]
    (eulerclosure.south) --
    node[eulerstab transition] {conclude the stability proof}
    (eulerproof.north);

\end{tikzpicture}%
}

\endgroup
\vspace{-18mm}
\label{diagram: Stability Proof}
\end{figure}

\clearpage

\subsubsection{Two-Radii Stability Theorems}

\noindent The single-radius certificate can be relaxed by prescribing separate target bounds for the low-order and top-order energies.  The same weighted energy argument still applies, provided the weighted ball is chosen small enough to lie inside the desired two-radii region,
\begin{equation}
\label{eq:two-radii-target-radii}
\mathfrak E_0(t)<\delta_0,
\qquad
\mathfrak H_k(t)<\delta_k,
\end{equation}
with $\delta_0>0$ and $\delta_k>0$.  These bounds are obtained from the same full energy $
\mathfrak E_k^2=
\mathfrak E_0^2+\mu_k\mathfrak H_k^2$. \\

\noindent For a chosen weight $\mu_k>0$, the invariant set is the weighted ball
\begin{equation}
\label{eq:two-radii-internal-ball}
\{ \mathfrak E_k(t)<\delta \}.
\end{equation}
This ball is contained in the target set \eqref{eq:two-radii-target-radii} provided that
\begin{equation}
\label{eq:two-radii-containment}
0<\delta\le \min\{\delta_0,\sqrt{\mu_k}\,\delta_k\},
\end{equation}
since this implies
\begin{equation}
\mathfrak E_0(t)\le \mathfrak E_k(t)<\delta\le\delta_0,
\qquad
\mathfrak H_k(t)\le \frac{\mathfrak E_k(t)}{\sqrt{\mu_k}}
<\frac{\delta}{\sqrt{\mu_k}}\le\delta_k.
\end{equation}

\hfill \\ 

\begin{theorem}[Two-radii nonlinear stability]
\label{thm:two-radii-stability}
Assume that the perturbation equations \eqref{eq: Stability theorem perturbation 1}--\eqref{eq: Stability theorem perturbation 3} hold for admissible perturbations. Fix $\delta_0>0$ and $\delta_k>0$.  Choose $\mu_k>0$ and $\delta>0$ satisfying
\begin{equation}
0<\delta\le \min\{\delta_0,\sqrt{\mu_k}\,\delta_k\}.
\end{equation}
Assume that the same certified estimates \eqref{eq:ns-thm-low-linear}--\eqref{eq:ns-thm-residual} of \Cref{thm:main-stability} hold, with the weighted energy
\begin{equation}
\mathfrak E_k^2=\mathfrak E_0^2+\mu_k\mathfrak H_k^2.
\end{equation}

\noindent Assume that $\mu_k$ satisfies the stability condition
\begin{equation}
\label{eq:two-radii-mu-choice}
\DLamStab(\mu_k):=
\min\left\{
\DLamLlow-\mu_k\DLamDLlow,
\ \DLamDLhigh-\frac{\DLamLhigh}{\mu_k}
\right\}>0.
\end{equation}
After this choice of $\mu_k$ is fixed, write $\DLamStab$ for the positive number $\DLamStab(\mu_k)$.  Define
\begin{equation}
\label{eq:two-radii-N2N3}
\DLamThree:=\DLamNthree+\DLamDNthree.
\end{equation}
Assume that
\begin{equation}
\label{eq:two-radii-smallness}
\DLamThree\delta+\frac{\DLamR}{\delta}<\DLamStab.
\end{equation}
Then every solution with
\begin{equation}
\label{eq:two-radii-initial-data}
\mathfrak E_0(0)^2+\mu_k\mathfrak H_k(0)^2<\delta^2
\end{equation}
satisfies
\begin{equation}
\label{eq:two-radii-conclusion}
\mathfrak E_0(t)<\delta_0,
\qquad
\mathfrak H_k(t)<\delta_k
\end{equation}
for as long as the solution exists and the certified estimates apply.
\end{theorem}

\begin{grayproof}
\noindent The weighted-energy part of the argument is identical to the single-radius argument in \Cref{thm:main-stability}. Thus, $\mathfrak E_k(t)<\delta$ for all times under consideration. It remains only to translate this invariant weighted ball into the two separate target bounds.  By the containment condition
\eqref{eq:two-radii-containment},
\begin{equation}
    \mathfrak E_0(t)\le \mathfrak E_k(t)<\delta\le\delta_0,
    \qquad
    \mathfrak H_k(t)\le\frac{\mathfrak E_k(t)}{\sqrt{\mu_k}}
    <\frac{\delta}{\sqrt{\mu_k}}\le\delta_k.
\end{equation}
\end{grayproof}

\clearpage

\subsection{From Rescaled Stability to Physical Blowup}
\label{sec:stability-to-physical-blowup}

\noindent The nonlinear stability theorems provide control of an exact solution in dynamically rescaled variables, where the concentrating core is kept at order-one scale and amplitude. By the exact dynamic-rescaling change of variables, every finite rescaled time corresponds to an equivalent physical solution. \\

\noindent To obtain finite-time blowup from a trajectory that exists for all rescaled times, two facts remain to be checked: the physical clock must have a finite limit, and undoing the normalization must force a physical quantity to become unbounded. This subsection provides exactly this reconstruction. \\

\noindent Both checks are controlled by the amplitude scale $\mathrm{s}_u(t)$. The normalization fixes $u(0,0,t)=\bar u_0\ne0$, where $\bar u_0$ is a constant, while
\begin{equation}
 \frac{d\mathring t}{dt}=\frac{1}{\mathrm{s}_u(t)}.
\label{eq:reconstruction-physical-clock}
\end{equation}
Thus growth of $\mathrm{s}_u$ makes the remaining physical time small and, after undoing the amplitude normalization, forces growth of the physical axial vorticity at the moving center. A uniform negative upper bound on $c_u(t)$ is enough to produce this growth. \\

\noindent As defined in \Cref{sec: Euler}, $\mathring t$ denotes physical time, $\mathring z_c(t)$ the physical axial center, $\tilde\omega$ the full physical vorticity, and $\omega^z$ its axial component.

\hfill 

\begin{proposition}[Finite-time physical blowup from a signed amplitude rate]
\label{thm:physical-reconstruction-signed-rate}
Consider the Euler equations in the dynamic rescaling formulation above, with $\lambda=1/2$. Assume that the exact admissible rescaled solution exists for every $t\ge0$, remains in the corresponding stability ball, and satisfies the modulation bounds throughout. Assume in addition that
\begin{equation}
 c_u(t)\le -\gamma_u<0,
 \qquad t\ge0,
\label{eq:reconstruction-cu-sign}
\end{equation}
for some $\gamma_u>0$. Then the physical time $\mathring t(t)$ converges to a finite limit $T$, with
\begin{equation}
0<T-\mathring t(t)
\le
\frac{e^{-\gamma_u t}}{\gamma_u\,\mathrm{s}_u(0)}.
\label{eq:reconstruction-time-tail-main}
\end{equation}
The traveling center $\mathring z_c(t)$ converges to a finite physical location. Moreover, at the moving center the physical axial vorticity satisfies
\begin{equation}
\left|\omega^z\bigl(0,\mathring z_c(t),\mathring t(t)\bigr)\right|
=2|\bar u_0|\,\mathrm{s}_u(t)
\longrightarrow\infty.
\label{eq:reconstruction-vorticity-main}
\end{equation}
Consequently, since $\omega^z$ is a component of the full physical vorticity $\tilde\omega$,
\begin{equation}
\|\tilde\omega(\cdot,\mathring t(t))\|_{L^\infty(\mathbb R^3)}\longrightarrow\infty,
\end{equation}
and the reconstructed physical solution cannot be continued smoothly through $T$.
\end{proposition}

\begin{grayproof}
The proof is given in Appendix~\ref{appx:physical-reconstruction}.
\end{grayproof}
\paragraph{Certification of the sign condition.}
\noindent Since $\bar c_u=-1$ and $c_u=\bar c_u+c_{u,\mathrm{res}}+\bvar{\delta c_u}$, the modulation bounds reduce \eqref{eq:reconstruction-cu-sign} to scalar checks. It is enough to verify
\begin{equation}
1-c_{u,\mathrm{res}}-M_{c_u}^{(1)}\delta_*>0,
\label{eq:reconstruction-euler-sign-test}
\end{equation}
where $M_{c_u}^{(1)}$ are constants that measure the sensitivity of the amplitude modulation parameter $c_u$ to the perturbation. Its precise definition is given in a separate stability paper~\citep{EulerBlowupStability}. \\

\begin{remark*}
\noindent This assumes that an exact admissible rescaled solution exists for every $t\ge0$ and remains in the stability ball. The stability estimates provide control of such a solution only for as long as it exists. It does not by itself prove all-time existence in rescaled time. To obtain an unconditional blowup theorem, the stability estimate must be combined with a local well-posedness and continuation result showing that the controlled rescaled solution cannot break down at any finite rescaled time. The residual of the numerical profile does not cause any issue because its certified residual is already included in the perturbation equations.
\end{remark*}

\hfill 

\subsection{Computer-Assisted Strategy for Proving Stability}
\label{sec: Computer-Assisted Strategy for Proving Stability}

\noindent Theorems~\ref{thm:main-stability} and~\ref{thm:two-radii-stability} reduce the nonlinear stability problem to a finite-dimensional optimization problem.  The aim is to choose the weights so that the linear damping margin is larger than the optimized nonlinear and residual error.  For each candidate choice of weights, we compute the certified constants, optimize the radius, and then update the weights until the closing inequality becomes strict. Once the weights are fixed, the remaining optimization is only over the scalar radius \(\delta\).  Decreasing \(\delta\) reduces the cubic contributions, while increasing \(\delta\) reduces the residual contribution \(\DLamR/\delta\).  Thus the best radius balances these two effects, and the outer search changes the weights to improve this balance.

\hfill

\noindent The optimization loop proceeds as follows.

\hfill 

\noindent\textbf{Step 1. Choose the class of admissible weight functions.} \\ 

\noindent Start with a finite-dimensional parametrization of a class of admissible weight functions, chosen so that the weights have the required positivity, regularity, and symmetries. \\

\noindent For each parameter value~\(\wgray{\vartheta}\), write the low-order weights as
\begin{equation}
\wgray{\Phi_{\vartheta}}
=
\bigl(\wgray{\Phi_{\omega,\vartheta}},
\wgray{\Phi_{r,\vartheta}},
\wgray{\Phi_{z,\vartheta}}\bigr),
\end{equation}
and define the order-\(k\) weights by
\begin{equation}
\wgray{\Phi_{i,k,\vartheta}}
=
1+\rho^k\wgray{\Phi_{i,\vartheta}},
\qquad
\rho(r,z)=r^2+z^2,
\qquad
 i\in\{\omega,r,z\}.
\end{equation}
The optimization variables are the weight parameter \(\wgray{\vartheta}\) and the relative top-order weight \(\mu_k\).  In the single-radius case of \Cref{thm:main-stability}, one imposes \(0<\mu_k\le1\).  In the two-radii case of \Cref{thm:two-radii-stability}, one only imposes \(\mu_k>0\).  In both cases the weighted energy is $\mathfrak E_k^2
=
\mathfrak E_0^2+
\mu_k\mathfrak H_k^2$. \\

\noindent Choose an initial admissible candidate \((\wgray{\vartheta},\mu_k)\).

\hfill \\ 

\noindent\textbf{Step 2. Compute the linear damping for the current candidate.} \\

\noindent For the current values of \((\wgray{\vartheta},\mu_k)\), compute
\begin{equation}
\DLamLlow,
\qquad
\DLamLhigh,
\qquad
\DLamDLlow,
\qquad
\DLamDLhigh.
\end{equation}

\noindent The damping margin available for the combined energy is
\begin{equation}
\DLamStab
:=
\min\left\{
\DLamLlow-\mu_k\DLamDLlow,
\ \DLamDLhigh-\frac{\DLamLhigh}{\mu_k}
\right\}.
\label{eq:strategy-stab-margin}
\end{equation}
If \(\DLamStab\le0\), the candidate cannot close the theorem, so we penalize this term and move to a new candidate.

\clearpage

\noindent\textbf{Step 3. Compute the nonlinear and residual constants.}\\

\noindent If \(\DLamStab>0\), compute
\begin{equation}
\DLamThree
=
\DLamNthree+\DLamDNthree.
\end{equation}

\hfill 

\noindent For this candidate, the error that must be dominated by the damping is
\begin{equation}
\DLamThree\delta
+
\frac{\DLamR}{\delta}.
\label{eq:strategy-error-curve}
\end{equation}

\hfill \\

\noindent\textbf{Step 4. Optimize the radius and evaluate the gap.}\\

\noindent For the current candidate, let \(\mathcal A_\delta(\wgray{\vartheta},\mu_k)\subset(0,\infty)\) denote the nonempty set of radii for which every radius-dependent hypothesis used in the proof is valid. \\

\hfill

\noindent In the single-radius case, the admissible condition is
\begin{equation}
\delta\in\mathcal A_\delta(\wgray{\vartheta},\mu_k).
\label{eq:strategy-single-radius-range}
\end{equation}

\hfill 

\noindent In the two-radii case with prescribed target bounds \(\delta_0>0\) and \(\delta_k>0\), the admissible conditions are
\begin{equation}
\delta\in\mathcal A_\delta(\wgray{\vartheta},\mu_k),
\qquad
0<\delta\le \min\{\delta_0,\sqrt{\mu_k}\,\delta_k\}.
\label{eq:strategy-two-radii-range}
\end{equation}

\hfill

\noindent Define the optimized gap
\begin{equation}
\mathcal J(\wgray{\vartheta},\mu_k)
:=
\inf_{\delta\ \mathrm{admissible}}
\left(
\DLamThree\delta
+
\frac{\DLamR}{\delta}
\right)
-
\DLamStab.
\label{eq:strategy-weight-objective}
\end{equation}

\hfill \\

\noindent The candidate succeeds exactly when \(\mathcal J(\wgray{\vartheta},\mu_k)<0\).  In this case, choose an admissible radius \(\delta_*\) such that
\begin{equation}
\DLamThree\delta_*
+
\frac{\DLamR}{\delta_*}
<
\DLamStab.
\label{eq:strategy-optimized-closing}
\end{equation}

\hfill \\

\noindent\textbf{Step 5. Iterate the outer optimization.} \\ 

\noindent The outer objective is the gap \(\mathcal J(\wgray{\vartheta},\mu_k)\). \\

\noindent One evaluation of this objective consists of Steps 2--4: compute the damping margin, compute the error constants, optimize \(\delta\), and subtract the margin. \\

\noindent If \(\mathcal J(\wgray{\vartheta},\mu_k)<0\), the candidate satisfies the numerical closing condition and is retained for rigorous certification.  \\

\noindent If \(\mathcal J(\wgray{\vartheta},\mu_k)\ge0\), update \((\wgray{\vartheta},\mu_k)\) to decrease \(\mathcal J\), and repeat Steps 2--4.

\clearpage

\noindent\textbf{Step 6. Record the tentative certificate.}\\

\noindent Once \(\mathcal J<0\), record the weights and a selected admissible radius \(\delta_*\) satisfying \eqref{eq:strategy-optimized-closing}. The single-radius certificate gives
\begin{equation}
\mathfrak E_k(t)<\delta_*.
\end{equation}

\noindent Provided that $    0 <  \delta_*\le \min\{\delta_0,\sqrt{\mu_k}\delta_k\}$, the two-radii certificate gives
\begin{equation}
\mathfrak E_0(t)<\delta_0,
\qquad
\mathfrak H_k(t)<\delta_k.
\end{equation}

\hfill  

\noindent\textbf{Step 7. Validate rigorously the tentative certificate.}\\

\noindent Validate rigorously and analytically that the candidate weights satisfy the closing condition and all required estimates.  The final certificate also verifies every radius-dependent admissibility guard at the selected radius \(\delta_*\).

\hfill  

\noindent We summarize the numerical strategy used to identify and rigorously certify the conditions required by the stability theorems below

\hfill

\begin{certalgorithm}{Loop for optimizing the nonlinear stability certificate}

\item Initialize an admissible candidate \((\wgray{\vartheta},\mu_k)\). 

\vspace{2mm}

\item Repeat the following loop.
\begin{enumerate}[label=\textbf{\alph*.},leftmargin=2.5em,itemsep=1.8em,topsep=0.5em]
\item Compute
\begin{equation}
    \DLamStab
=
\min\left\{
\DLamLlow-\mu_k\DLamDLlow,
\ \DLamDLhigh-\frac{\DLamLhigh}{\mu_k}
\right\}.
\end{equation}
If \(\DLamStab\le0\), penalize the candidate and update \((\wgray{\vartheta},\mu_k)\).
\item If \(\DLamStab>0\), compute \(\DLamR\) and  $
    \DLamThree=\DLamNthree+\DLamDNthree$.
\item Determine the nonempty admissible set of radii for the current candidate, and compute
\begin{equation}
\mathcal J(\wgray{\vartheta},\mu_k)
=
\inf_{\delta\ \mathrm{admissible}}
\left(
\DLamThree\delta+
\frac{\DLamR}{\delta}
\right)
-
\DLamStab.
\end{equation}
\item If \(\mathcal J<0\), select an admissible radius \(\delta_*\) satisfying \eqref{eq:strategy-optimized-closing}, record the resulting tentative certificate, and proceed to rigorous validation.
\item If \(\mathcal J\ge0\), update \((\wgray{\vartheta},\mu_k)\) to decrease \(\mathcal J\), and repeat the loop.
\end{enumerate}

\vspace{2mm}

\item Record the final weights, the selected admissible radius \(\delta_*\), and the resulting bounds.

\vspace{2mm}

\item Validate rigorously the tentative proof using analytic certificates, including every radius-dependent admissibility guard at \(\delta_*\).

\end{certalgorithm}

\clearpage

\section{Partial Results: Linear Damping}
\label{sec:partial-linear-damping}

\subsection{Objective and Scope}

\noindent We now examine the signed linear damping mechanism that drives the stability argument. The objective is to obtain a coercive linear contribution to the full weighted energy after the low-order and high-order estimates are combined and the remaining linear terms are included.\\

\noindent The complete linear estimates we aim to prove have the two-level form
\begin{align}
\sum_{i\in \{\omega,r,z\}} \langle \mathcal L_i,\bvar{\delta v_i}\rangle_{\wgray{\Phi_i}}
&\le -\DLamLlow\mathfrak E_0^2+\DLamLhigh\mathfrak H_k^2,\label{eq:objective-low-linear-estimate}\\
\sum_{i\in \{\omega,r,z\}}\sum_{|\alpha|=k}\langle D^\alpha\mathcal L_i,D^\alpha\bvar{\delta v_i}\rangle_{\wgray{\Phi_{i,k}}}
&\le \DLamDLlow\mathfrak E_0^2-\DLamDLhigh\mathfrak H_k^2.\label{eq:objective-high-linear-estimate}
\end{align}
The four coefficients have complementary roles. The constant \(\DLamLlow\) is the damping margin supplied by the low-order estimate, while \(\DLamLhigh\) measures the accompanying loss at top order. Conversely, \(\DLamDLhigh\) is the damping margin supplied by the high-order estimate, while \(\DLamDLlow\) measures its accompanying low-order loss. The purpose of the remainder of this section is to construct and certify these four constants.\\

\noindent The reason the two estimates must be considered together is already visible from their form. Multiplying equation~\eqref{eq:objective-high-linear-estimate} by \(\mu_k\) and adding it to equation~\eqref{eq:objective-low-linear-estimate} gives
\begin{align}
&\sum_{i\in \{\omega,r,z\}} \langle \mathcal L_i,\bvar{\delta v_i}\rangle_{\wgray{\Phi_i}}
+\mu_k\sum_{i\in \{\omega,r,z\}}\sum_{|\alpha|=k}\langle D^\alpha\mathcal L_i,D^\alpha\bvar{\delta v_i}\rangle_{\wgray{\Phi_{i,k}}}\nonumber\\
&\qquad\qquad\qquad\qquad\qquad\qquad\le
-\left(\DLamLlow-\mu_k\DLamDLlow\right)\mathfrak E_0^2
-\mu_k\left(\DLamDLhigh-\frac{\DLamLhigh}{\mu_k}\right)\mathfrak H_k^2.
\label{eq:objective-combined-linear-estimate}
\end{align}

\hfill

\noindent Thus neither estimate needs to be coercive in both energy levels by itself. A positive \(\mathfrak H_k^2\) loss in the low-order estimate can be absorbed by the negative high-order contribution, while a positive \(\mathfrak E_0^2\) loss in the high-order estimate can be absorbed by the low-order damping. This is the basic role of the high-order estimate in the argument.\\

\noindent It is therefore natural to define the net linear damping margin by
\begin{equation}
\DLamStab(\mu_k)
=
\min\left\{
\DLamLlow-\mu_k\DLamDLlow,
\ \DLamDLhigh-\frac{\DLamLhigh}{\mu_k}
\right\}>0.
\label{eq:numerical-signed-linear-stability-margin}
\end{equation}
The two quantities inside the minimum are precisely the damping margins for \(\mathfrak E_0^2\) and \(\mu_k\mathfrak H_k^2\) in \eqref{eq:objective-combined-linear-estimate}. A positive value of \(\DLamStab\) gives a genuinely coercive linear estimate before the nonlinear and residual contributions are added. Unlike those later contributions, a missing linear sign cannot in general be repaired simply by restricting to a smaller perturbation neighborhood.\\

\noindent The rest of the section is organized as follows. First, we construct the low-order signed matrix and identify the pointwise criterion that yields low-order damping. Second, we choose admissible weights and certify the low-order matrix on the finite computational box and in the analytic far field. Third, once this direct low-order bound is established, we analyze any exceptional region where the same pointwise margin is not available. Fourth, we construct the complete high-order signed matrix and combine its damping with that low-order estimate. Finally, we incorporate the linear terms kept outside the two matrices, together with the nonlinear and residual contributions, to obtain the closing stability condition.

\hfill 

\subsection{Low-Order Signed Damping Matrix}
\label{sssec:numerical-low-order-signed-damping}

\subsubsection{Transport Field and Retained Linear Terms}

\noindent At low order, the question is whether the transport terms and a selected set of local linear couplings form a negative quadratic form at almost every point. The selection is made so that every retained term acts pointwise on the three variables already present in the low-order weighted energy and has a controlled coefficient after the corresponding weight conversion. Terms with a different structure are kept outside the pointwise matrix and estimated separately by specialized mechanisms adapted to them. \\

\noindent Recall that the transport field is
\begin{equation}
G_r=\lambda r+\bar u^r,
\qquad
G_z=\lambda z+\bar C+\bar u^z,
\label{eq:low-linear-profile-transport-coeffs}
\end{equation}
where
\begin{equation}
\bar u^r=-r\partial_z\bar\psi,
\qquad
\bar u^z=2\bar\psi+r\partial_r\bar\psi.
\end{equation}

\hfill 

\noindent Besides transport, the signed matrix retains the local diagonal coefficients multiplying \(\bvar{\delta\omega}\), \(\partial_r\bvar{\delta u}\), and \(\partial_z\bvar{\delta u}\), together with the two local \(r\)--\(z\) couplings between \(\partial_r\bvar{\delta u}\) and \(\partial_z\bvar{\delta u}\). These are exactly the local terms that can be represented in the same three-component pointwise quadratic form. \\

\noindent Note that the coupling term $2\bar u\, \partial_z\bvar{\delta u}$ is kept outside the low-order signed matrix, because after passing to the weighted variables its coefficient would contain
\(
2\bar u \sqrt{\Phi_\omega/\Phi_z}
\),
and the independent low-order singularity rates of \(\Phi_\omega\) and \(\Phi_z\) do not by themselves guarantee that this multiplier is bounded. The complete low-order term list and the precise partition between matrix and non-matrix terms are recorded in Appendix~\ref{app:low-order-linear-setup}. \\ 

\subsubsection{Weighted Transport Contribution}

\noindent The role of the weights is most transparent in the transport contribution. Assuming the boundary terms vanish under the boundary and decay conditions used in the weighted energy estimate, integration by parts gives
\begin{equation}
\begin{aligned}
\left\langle
-G_r\partial_r f-G_z\partial_z f,
 f
\right\rangle_{\wgray{\Phi_i}}
&=-\frac12\int_{\mathbb D}
\left(G_r\partial_r(f^2)+G_z\partial_z(f^2)\right)\wgray{\Phi_i}\,dr\,dz\\
&=\frac12\int_{\mathbb D}|f|^2
\left[\partial_r(G_r\wgray{\Phi_i})+\partial_z(G_z\wgray{\Phi_i})\right]dr\,dz\\
&=\frac12\int_{\mathbb D} |f|^2\wgray{\Phi_i}
\left(
\partial_rG_r+\partial_zG_z
+G_r\partial_r\log\wgray{\Phi_i}
+G_z\partial_z\log\wgray{\Phi_i}
\right)dr\,dz.
\end{aligned}
\label{eq:low-linear-transport-signed}
\end{equation}

\hfill 

\noindent Thus the coefficient multiplying \(|f|^2\wgray{\Phi_i}\) is
\begin{equation}
\frac12\left(\partial_rG_r+\partial_zG_z\right)
+\frac12\left(G_r\partial_r\log\wgray{\Phi_i}+G_z\partial_z\log\wgray{\Phi_i}\right).
\label{eq:numerical-transport-two-mechanisms}
\end{equation}
The first term is fixed by the divergence of the transport field. The second depends on the variation of the weight along the transport direction and is thus the part that can be adjusted through the choice of~\(\wgray{\Phi_i}\). This is the basic mechanism by which the weights can improve the signed transport contribution.  \\ 

\subsubsection{Explicit Low-Order Matrix}
\label{sssec:low-order-linear-matrix-construction}

\noindent We introduce the weighted low-order state
\begin{equation}
\bvar{d\xi}(r,z)
=
\left(
\wgray{\Phi_\omega}^{1/2}\bvar{\delta\omega},
\wgray{\Phi_r}^{1/2}\partial_r\bvar{\delta u},
\wgray{\Phi_z}^{1/2}\partial_z\bvar{\delta u}
\right)^{\!\top},
\qquad \quad 
\int_{\mathbb D}|\bvar{d\xi}|^2\,dr\,dz=\mathfrak E_0^2.
\label{eq:low-linear-xi-column}
\end{equation}
With the component order \((\omega,r,z)\), the retained low-order contribution is represented by a symmetric \(3\times3\) matrix \(\mathbb M_{\mathcal L}(r,z)\). We now define its entries from the transport contribution and the retained local coefficients. \\

\paragraph{Transport coefficients.}
\noindent The derivatives of the transport field that enter the matrix are
\begin{align}
\partial_r\bar u^r&=-\partial_z\bar\psi-r\partial_{rz}\bar\psi,
&
\partial_z\bar u^r&=-r\partial_{zz}\bar\psi,\\
\partial_r\bar u^z&=3\partial_r\bar\psi+r\partial_{rr}\bar\psi,
&
\partial_z\bar u^z&=2\partial_z\bar\psi+r\partial_{rz}\bar\psi,
\end{align}
so that
\begin{equation}
\partial_rG_r=\lambda+\partial_r\bar u^r,
\qquad
\partial_zG_z=\lambda+\partial_z\bar u^z.
\end{equation}

\hfill 

\noindent For each \(i\in\{\omega,r,z\}\), we denote the transport contribution to the corresponding diagonal entry by
\begin{equation}
D_i
:=\frac12\left(
\partial_rG_r+\partial_zG_z
+G_r\partial_r\log\wgray{\Phi_i}
+G_z\partial_z\log\wgray{\Phi_i}
\right).
\end{equation}

\hfill 

\paragraph{Diagonal entries.}
\noindent The constant modulation offset \(\cures\) is absorbed into the reference amplitude coefficient and therefore contributes to the zeroth-order coefficients added to the transport contribution, i.e.,
\begin{align}
(\mathbb M_{\mathcal L})_{\omega\omega}
&=D_\omega+\bar c_u+\cures-\lambda,
\label{eq:low-linear-matrix-diag-omega}\\
(\mathbb M_{\mathcal L})_{rr}
&=D_r+2\partial_z\bar\psi+\bar c_u+\cures-\lambda-\partial_r\bar u^r,
\label{eq:low-linear-matrix-diag-r}\\
(\mathbb M_{\mathcal L})_{zz}
&=D_z+2\partial_z\bar\psi+\bar c_u+\cures-\lambda-\partial_z\bar u^z.
\label{eq:low-linear-matrix-diag-z}
\end{align}

\hfill 

\paragraph{Cross entry.}
\noindent The two retained \(r\)--\(z\) couplings contribute to the same symmetric off-diagonal entry. After converting each term to the weighted variables and symmetrizing the quadratic form, we obtain
\begin{equation}
\begin{aligned}
(\mathbb M_{\mathcal L})_{rz}
=(\mathbb M_{\mathcal L})_{zr}
=-\frac12\left[
\partial_r\bar u^z\left(\frac{\Phi_r}{\Phi_z}\right)^{1/2}
+\partial_z\bar u^r\left(\frac{\Phi_z}{\Phi_r}\right)^{1/2}
\right].
\end{aligned}
\label{eq:low-linear-matrix-rz}
\end{equation}
There are no retained \(\omega\)--\(r\) or \(\omega\)--\(z\) entries at low order. In particular, the \(\omega\)--\(z\) coupling term $2\bar u\, \partial_z\bvar{\delta u}$ discussed earlier remains outside the matrix because its weight conversion is not controlled by the low-order weight assumptions. \\

\hfill \\

\noindent Altogether, these entries form the symmetric matrix \(\mathbb M_{\mathcal L}\). By construction, the retained transport and selected couplings contribute to the weighted energy through
\begin{equation}
\int_{\mathbb D}
\bvar{d\xi}^{\top}\mathbb M_{\mathcal L}(r,z)\bvar{d\xi}\,dr\,dz.
\label{eq:low-linear-sign-matrix-form}
\end{equation}

\hfill 

\noindent The pointwise sign of this quadratic form is determined by the largest eigenvalue of \(\mathbb M_{\mathcal L}\):
\begin{equation}
\lambda_{\max}\!\left(\mathbb M_{\mathcal L}(r,z)\right)<0
\qquad \Longleftrightarrow \qquad
\bvar{\eta}^{\top}\mathbb M_{\mathcal L}(r,z)\bvar{\eta}<0
\quad\text{for every }\bvar{\eta}\neq0.
\label{eq:numerical-low-local-margin-equivalence}
\end{equation}

\hfill 

\noindent Thus \(-\lambda_{\max}(\mathbb M_{\mathcal L})\) is the local raw damping margin whenever the largest eigenvalue is negative.

\subsection{Low-order Damping Certification}

\subsubsection{Weight Optimization and Finite-Box Evaluation}
\label{sssec:numerical-weight-parametrization}

\paragraph{Weight optimization.}
\noindent In the transport diagonal, the weight enters only through
\begin{equation}
G_r\partial_r\log\wgray{\Phi_i}
+G_z\partial_z\log\wgray{\Phi_i},
\qquad i\in\{\omega,r,z\}.
\label{eq:numerical-weight-directional-derivative}
\end{equation}
The free parameters in the admissible low-order weights are therefore chosen to improve the pointwise sign of the complete matrix while preserving the admissibility conditions imposed on the low-order weights. In our numerical experiments, we represent the weight functions as a combination of locally supported Gaussian basis functions, a singular component near the origin, asymptotically growing neural networks in the far field, and a global neural-network correction, with all components parametrized and jointly optimized.\\

\noindent Because the \(r\)--\(z\) cross entry also depends on weight ratios, it is not sufficient to optimize the transport diagonal component by component. The numerical objective is instead based on the largest eigenvalue of the full symmetric matrix \(\mathbb M_{\mathcal L}\), with the aim of making its worst value as negative as possible and reducing the size of any region where strict negativity fails. \\

\noindent For the optimized candidate, we evaluate \(\lambda_{\max}(\mathbb M_{\mathcal L}(r,z))\) throughout the finite half-plane box on which the reference profile is represented by splines. The rigorous finite-box step is then to enclose the largest eigenvalue with interval arithmetic, requiring a uniform negative upper bound on the region treated directly and only a finite upper bound on any localized pieces introduced below. Numerically, the largest eigenvalue is negative on almost all of the box, and the weakest margin is confined to small neighborhoods where the meridional transport field is small. These locations are analyzed in Subsection~\ref{sssec:numerical-low-transport-points}. The numerical evaluation identifies two low-transport neighborhoods: one around the normalized origin and one around a single off-axis fixed point. Their treatment is determined by the sign of the full low-order matrix, not by the smallness of the transport alone. Around the origin the full matrix remains in the directly damped regime, whereas the off-axis neighborhood is the single region on which the uniform negative matrix margin is not imposed and is therefore treated by localization.\\

\noindent The optimization step is used to identify a promising admissible set of weights. It is not itself the rigorous damping certificate. Once a candidate set of weights is fixed, the largest eigenvalue of the resulting matrix is evaluated on the computational box and then enclosed with interval arithmetic. The exterior region is treated separately by the analytic tail argument discussed in the subsequent section. 

\hfill 

\subsubsection{Analytic Far-Field Closure}
\label{sssec:low-order-linear-tail-closure}

\noindent Define the numerical enclosure box by
\begin{equation}
\mathbb{D}_{\rm box}:=[0,R_{\rm box}]\times[-R_{\rm box},R_{\rm box}].
\label{eq:low-linear-tail-box}
\end{equation}
On \(\mathbb D_{\rm box}\), the finite-box certification step consists of interval upper bounds for the largest eigenvalue of the pointwise matrix as described in the previous subsection. To extend these bounds to the full half-plane, it remains to control the exterior region \(\mathbb D\setminus\mathbb D_{\rm box}\) analytically. \\

\noindent Throughout the tail calculation,
\begin{equation}
\lambda=\frac12,
\qquad
\bar c_u=-1.
\label{eq:low-linear-tail-normalization}
\end{equation}
The profile-dependent constant \(\bar C\) remains symbolic until the numerical certificate is evaluated. \\

\noindent In our spline representation, the steady profile and every profile derivative entering \(\mathbb M_{\mathcal L}\) vanish exactly on \(\mathbb D\setminus\mathbb D_{\rm box}\). Then
\begin{equation}
\bar u=\bar\omega=\bar\psi=0,
\qquad
\bar u^r=\bar u^z=0,
\qquad
G_r=\frac r2,
\qquad
G_z=\frac z2+\bar C.
\label{eq:numerical-tail-effective-transport}
\end{equation}
All retained cross couplings vanish in this region. As a result, the tail damping matrix is diagonal, and its \(i\)-th diagonal entry is
\begin{equation}
-1+\cures
+\frac14(r\partial_r+z\partial_z)\log\wgray{\Phi_i}
+\frac12\bar C\,\partial_z\log\wgray{\Phi_i},
\qquad i\in\{\omega,r,z\}.
\label{eq:low-linear-tail-diagonal-entry}
\end{equation}

\hfill 

\paragraph{Scalar tail condition.}
\noindent It is sufficient to find \(\gamma_{\rm tail}>0\) such that
\begin{equation}
\sup_{\mathbb{D}\setminus\mathbb{D}_{\rm box}}
\left[
(r\partial_r+z\partial_z)\log\wgray{\Phi_i}
+2\bar C\,\partial_z\log\wgray{\Phi_i}
\right]
\le 4(1-\cures-\gamma_{\rm tail})
\label{eq:low-linear-tail-scalar-condition}
\end{equation}
for every \(i\in\{\omega,r,z\}\). Substituting this bound into \eqref{eq:low-linear-tail-diagonal-entry} yields
\begin{equation}
\lambda_{\max}\bigl(\mathbb M_{\mathcal L}(r,z)\bigr)
\le-\gamma_{\rm tail}
\qquad\text{on }\mathbb{D}\setminus\mathbb{D}_{\rm box}.
\label{eq:low-linear-tail-eigenvalue-bound}
\end{equation}

\noindent The weights can be continued analytically in the tail. For the polynomial-growth continuation,
\begin{equation}
\wgray{\Phi_i}(r,z)
=C_{\rm floor}+a_i\rho,
\qquad
\rho=r^2+z^2,
\qquad
a_i>0,
\label{eq:numerical-tail-weight-family}
\end{equation}
with \(p=1\), the calculation in Appendix~\ref{app:tail-weight-models} gives the sufficient condition
\begin{equation}
\left(1+\frac{2|\bar C|}{R_{\rm box}}\right)
<2(1-\cures).
\label{eq:numerical-tail-algebraic-condition}
\end{equation}
Whenever this inequality holds, the tail matrix has a uniform strictly negative largest eigenvalue. Once the finite-box interval enclosure is certified, combining it with this tail bound gives a pointwise low-order matrix estimate on the entire half-plane except for any problematic small region of low transport on which the uniform negative margin is not obtained. In our case, $\bar C \approx 169 $, $R_{\rm box} \approx 5.34 \times 10^{12}$ and $\cures \approx 2 \times 10^{-15}$, so \eqref{eq:numerical-tail-algebraic-condition} holds true and the scalar tail condition is satisfied with the $p=1$ polynomial growth continuation. \\

\subsection{Meridional Fixed Points and Localized Low-Order Control}
\label{sssec:numerical-low-transport-points}

\subsubsection{Fixed Points and Low-Transport Regions}

\noindent The low-order damping is weakest near locations where the meridional transport field becomes small. We refer to \((r_j^\odot,z_j^\odot)\) as a \emph{meridional fixed point} if the meridional transport field vanishes there, i.e.,
\begin{equation}
G_r(r_j^\odot,z_j^\odot)=0,
\qquad
G_z(r_j^\odot,z_j^\odot)=0.
\label{eq:numerical-fixed-point-definition}
\end{equation}

\hfill 

\noindent We also refer to a neighborhood as a \emph{low-transport region} if the meridional transport field is weak there, i.e.
\begin{equation}
\sqrt{|G_r(r,z)|^2+|G_z(r,z)|^2} \ll 1.
\label{eq:numerical-low-transport-region}
\end{equation}

\hfill 

\noindent The weight-dependent part of the low-order transport coefficient is
\begin{equation}
    G_r\partial_r\log\wgray{\Phi_i}+G_z\partial_z\log\wgray{\Phi_i}.
\end{equation}
At a meridional fixed point this contribution vanishes, and it becomes small in a low-transport neighborhood when the logarithmic derivatives of the admissible weights remain controlled. This weakens the main tunable contributions to the low-order diagonal entries, although it does not by itself determine the sign of the full damping matrix, because the divergence term and the retained local coefficients remain present. \\

\begin{figure}[t] \centering \includegraphics[width=\linewidth]{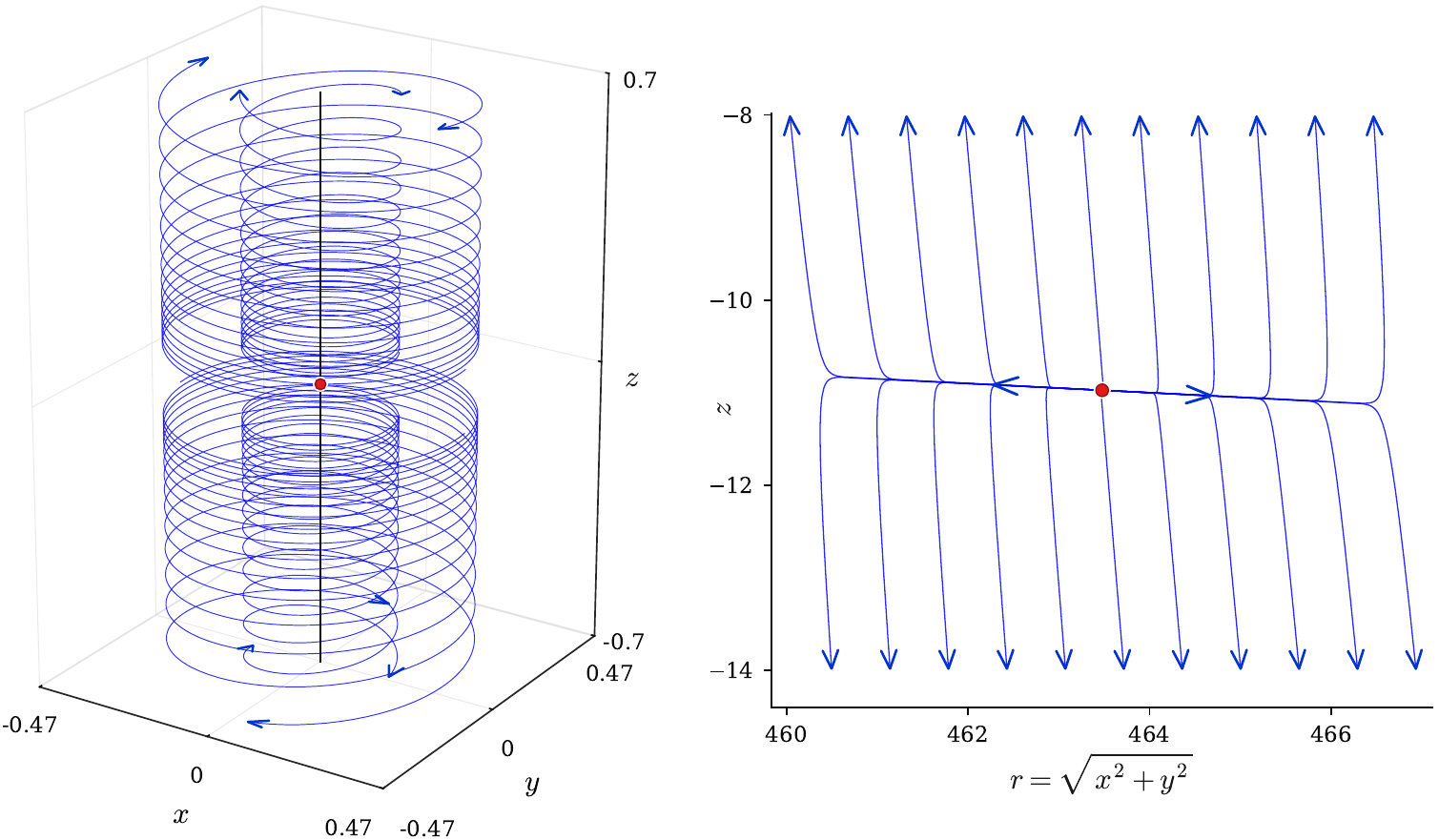} \caption{Local transport field near the two meridional fixed points of the rescaled Euler profile. The red markers denote the fixed points. \textbf{(Left)}  Transport field near the fixed point on the symmetry axis, illustrating its nearly degenerate locally outgoing structure. \textbf{(Right)} Meridional transport field $(G_r,G_z)$ near the off-axis fixed point, illustrating its strictly locally outgoing structure.} \label{fig:transport-fixed-points} \end{figure}

\noindent Numerically, the meridional transport field has two fixed points in the computational domain. One lies on the symmetry axis and is shifted to the origin by the normalization condition~\eqref{eq: perturbation normalization}, while the second lies off axis. Their local transport geometry is illustrated in Figure~\ref{fig:transport-fixed-points}: the on-axis fixed point exhibits a nearly degenerate locally outgoing structure, whereas the off-axis fixed point is strictly locally outgoing.

\noindent Although the transport field vanishes at both fixed points, the condition \(G_r=G_z=0\) alone does not determine the sign of the full low-order damping matrix, since the divergence contribution and the retained local coefficients remain nonzero in general. For the optimized candidate, the numerical evaluation indicates that these remaining terms are sufficient to keep the largest eigenvalue negative in a neighborhood of the normalized origin, so the on-axis fixed point remains within the region controlled by the direct finite-box bound. In contrast, the same uniform negative margin is not imposed in a small neighborhood of the off-axis fixed point; this region is instead controlled by the localized estimate developed below. Moreover, the first derivatives of \(G_r\) and \(G_z\) need not vanish at either fixed point. The high-order transport commutators can therefore retain a nontrivial signed contribution even where the low-order weight-dependent transport contribution vanishes, providing the complementary mechanism used in the high-order estimate.

\hfill \\ 

\subsubsection{Localized Matrix Estimate}
\label{sssec:localized-matrix-certificate}
\label{app:localized-matrix-certificate}

\noindent To quantify the off-axis low-transport region without imposing the same pointwise margin there as on the rest of the domain, we use a localized estimate. A uniform negative upper bound for \(\lambda_{\max}(\mathbb M_{\mathcal L})\) on the entire half-plane would be sufficient for the integrated energy estimate, but it is stronger than necessary. A failure of the pointwise margin at an isolated point has no effect on the integral, because a set of measure zero carries no \(L^2\) mass. What must be controlled is a positive-measure neighborhood on which the uniform negative margin is unavailable. The localized estimate permits such neighborhoods provided that the weighted energy that can concentrate in them is quantitatively controlled. \\

\noindent Small measure of the neighborhood alone is not sufficient, since an \(L^2\) perturbation can concentrate most of its mass in a small region. We thus combine two inputs: a uniform negative matrix bound on the complement of the localized regions and a concentration estimate that controls the weighted energy inside each localized piece by the global low-order and top-order energies. The low-order part of this concentration estimate reduces the effective coefficient of \(\mathfrak E_0^2\), while the top-order part produces an additional \(\mathfrak H_k^2\) term that must be absorbed by the negative high-order matrix contribution after the high-order estimate is weighted by \(\mu_k\).

\hfill  \\

\paragraph{Localized regions.}
Let the region on which the uniform low-order matrix margin is not obtained be a finite union of pairwise disjoint measurable pieces
\begin{equation}
\mathbb D_{\mathcal L}^{\rm loc}
=
\bigcup_{\ell\in\mathcal I_{\mathcal L}^{\rm loc}}B_\ell,
\qquad
\mathbb D
=
\mathbb D_{\mathcal L}^{\rm damp}\cup\mathbb D_{\mathcal L}^{\rm loc},
\qquad
\mathbb D_{\mathcal L}^{\rm damp}\cap\mathbb D_{\mathcal L}^{\rm loc}=\emptyset.
\label{eq:localized-matrix-region-split}
\end{equation}
We keep the notation above in a general finite-union form because the same argument would apply if several disjoint low-transport regions had to be localized. In the present computation, however, \(\mathcal I_{\mathcal L}^{\rm loc}=\{\ell_{\rm off}\}\) and \(B_{\ell_{\rm off}}\) is the single neighborhood of the unique off-axis fixed point. In particular, \(B_{\ell_{\rm off}}\) stays a positive distance from the origin, which is useful for the concentration estimate because the low-order weights are singular only at the origin and are therefore bounded on \(B_{\ell_{\rm off}}\).

\hfill \\

\paragraph{Inputs for the localized certificate.}
Once the finite-box enclosure is validated, the direct-region constant \(\Lambda_{\mathcal L}^{\rm damp}\) can be chosen as the minimum of the certified finite-box margin on the directly treated portion of \(\mathbb D_{\rm box}\) and the analytic tail margin \(\gamma_{\rm tail}\). On \(\mathbb D_{\mathcal L}^{\rm damp}\), we then require
\begin{equation}
\lambda_{\max}\bigl(\mathbb M_{\mathcal L}(r,z)\bigr)
\le -\Lambda_{\mathcal L}^{\rm damp}
\qquad
(r,z)\in\mathbb D_{\mathcal L}^{\rm damp},
\qquad
\Lambda_{\mathcal L}^{\rm damp}>0.
\label{eq:localized-matrix-damp-bound}
\end{equation}
On each localized piece, we require the finite upper bound
\begin{equation}
\lambda_{\max}\bigl(\mathbb M_{\mathcal L}(r,z)\bigr)
\le \Lambda_{\mathcal L}^{{\rm loc},\ell}
\qquad
(r,z)\in B_\ell,
\qquad
\Lambda_{\mathcal L}^{{\rm loc},\ell}\ge0.
\label{eq:localized-matrix-loc-bound}
\end{equation}
Thus the largest eigenvalue is uniformly negative on \(\mathbb D_{\mathcal L}^{\rm damp}\), while on each \(B_\ell\) only a finite upper bound is required. \\

\noindent The second input is a concentration estimate of the form
\begin{equation}
\int_{B_\ell}|\bvar{d\xi}|^2\,dr\,dz
\le
m_{\mathcal L}^{0,{\rm loc},\ell}\mathfrak E_0^2
+
m_{\mathcal L}^{k,{\rm loc},\ell}\mathfrak H_k^2,
\qquad \ell\in\mathcal I_{\mathcal L}^{\rm loc},
\label{eq:localized-matrix-loc-concentration}
\end{equation}
with nonnegative constants \(m_{\mathcal L}^{0,{\rm loc},\ell}\) and \(m_{\mathcal L}^{k,{\rm loc},\ell}\). The first coefficient measures the part of the localized weighted mass controlled at low order. The second measures the part that must be controlled through the top-order energy.

\hfill

\paragraph{Effective matrix estimate.}
Define
\begin{equation}
\gamma_\ell
:=
\Lambda_{\mathcal L}^{\rm damp}+\Lambda_{\mathcal L}^{{\rm loc},\ell}.
\label{eq:localized-matrix-piece-cost}
\end{equation}

\hfill 

\noindent Using the domain split and \(\int_{\mathbb D}|\bvar{d\xi}|^2=\mathfrak E_0^2\),
\begin{align}
\int_{\mathbb D}\bvar{d\xi}^{\top}\mathbb M_{\mathcal L}\bvar{d\xi}\,dr\,dz
&\le
-\Lambda_{\mathcal L}^{\rm damp}\mathfrak E_0^2
+
\sum_{\ell\in\mathcal I_{\mathcal L}^{\rm loc}}\gamma_\ell
\int_{B_\ell}|\bvar{d\xi}|^2\,dr\,dz.
\label{eq:localized-matrix-split-derivation}
\end{align}
Define
\begin{align}
\Lambda_{\mathcal L}^{\rm mat,loc}
&:=
\Lambda_{\mathcal L}^{\rm damp}
-
\sum_{\ell\in\mathcal I_{\mathcal L}^{\rm loc}}\gamma_\ell m_{\mathcal L}^{0,{\rm loc},\ell},
\label{eq:localized-matrix-loc-Lambda-choice}\\
A_{\mathcal L}^{k,\rm loc}
&:=
\sum_{\ell\in\mathcal I_{\mathcal L}^{\rm loc}}\gamma_\ell m_{\mathcal L}^{k,{\rm loc},\ell}.
\label{eq:localized-matrix-loc-Ak-choice}
\end{align}
Then the concentration estimate gives
\begin{equation}
\int_{\mathbb D}\bvar{d\xi}^{\top}\mathbb M_{\mathcal L}\bvar{d\xi}\,dr\,dz
\le
-\Lambda_{\mathcal L}^{\rm mat,loc}\mathfrak E_0^2
+
A_{\mathcal L}^{k,\rm loc}\mathfrak H_k^2.
\label{eq:localized-matrix-estimate}
\end{equation}
Thus localization reduces the effective low-order matrix margin from \(\Lambda_{\mathcal L}^{\rm damp}\) to \(\Lambda_{\mathcal L}^{\rm mat,loc}\) and introduces the top-order term \(A_{\mathcal L}^{k,\rm loc}\mathfrak H_k^2\). To retain a genuine low-order damping term, the localized certificate must in particular satisfy \(\Lambda_{\mathcal L}^{\rm mat,loc}>0\). The next subsection constructs the high-order matrix and provides the negative \(\mathfrak H_k^2\) contribution that is compared directly with \(A_{\mathcal L}^{k,\rm loc}\mathfrak H_k^2\).

\hfill \\

\subsection{High-Order Signed Damping Matrix}
\label{sssec:numerical-higher-order-signed-damping}

\subsubsection{High-Order Transport Commutators and Meridional Fixed Points}

\noindent The first transport commutators depend on the full Jacobian of the meridional transport field,
\begin{equation}
\partial_rG_r=\lambda+\partial_r\bar u^r,
\qquad
\partial_zG_r=\partial_z\bar u^r,
\qquad
\partial_rG_z=\partial_r\bar u^z,
\qquad
\partial_zG_z=\lambda+\partial_z\bar u^z.
\label{eq:high-linear-explicit-G-gradient-add}
\end{equation}

\hfill 

\noindent We consider \(\alpha=(q,k-q)\) with \(0\le q\le k\) and apply \(D^\alpha\) to \(G_r\partial_r f+G_z\partial_z f\). The top-order commutator terms are exactly those in which one derivative lands on \(G_r\) or \(G_z\). Each such term still contains \(k\) derivatives of \(f\), so it belongs to the top-order signed matrix rather than to a lower-order remainder. Terms in which at least two derivatives land on the transport coefficients contain at most \(k-1\) derivatives of the perturbation. \\

\paragraph{Top-order derivative variables.}
For \(0\le q\le k\), set
\begin{equation}
V_{i,q}:=\partial_r^q\partial_z^{k-q}\bvar{\delta v_i},
\qquad i\in\{\omega,r,z\}.
\label{eq:high-linear-top-derivative-array-add}
\end{equation}
 Here \(q\) is the number of \(r\)-derivatives and \(k-q\) is the number of \(z\)-derivatives. The neighboring terms in the formula below occur only when the corresponding derivative is available: the \(V_{i,q-1}\) term carries a factor \(q\) and therefore disappears when \(q=0\), while the \(V_{i,q+1}\) term carries a factor \(k-q\) and disappears when \(q=k\).\\
 
 \noindent For each \((i,q)\), the complete top-order contribution coming from  high-order transport and its first commutators is
\begin{equation}
\begin{aligned}
&-G_r\partial_rV_{i,q}-G_z\partial_zV_{i,q}
-\bigl(q\,\partial_rG_r+(k-q)\,\partial_zG_z\bigr)V_{i,q}
-q\,\partial_rG_z\,V_{i,q-1}
-(k-q)\,\partial_zG_r\,V_{i,q+1}.
\end{aligned}
\label{eq:high-linear-complete-transport-add}
\end{equation}
The first two terms transport the top derivative itself. The coefficient of \(V_{i,q}\) contains the two diagonal commutator contributions \(q\partial_rG_r\) and \((k-q)\partial_zG_z\).  The index shifts can be read directly from the Leibniz rule. In the term \(G_z\partial_z f\), if one of the \(q\) \(r\)-derivatives falls on \(G_z\), then only \(q-1\) \(r\)-derivatives remain on \(f\), while the original \(\partial_z\) supplies one additional \(z\)-derivative. This gives \(V_{i,q-1}\) with coefficient \(q\,\partial_rG_z\). Similarly, in \(G_r\partial_r f\), if one of the \(k-q\) \(z\)-derivatives falls on \(G_r\), the perturbation factor gains one \(r\)-derivative and loses one \(z\)-derivative, giving \(V_{i,q+1}\) with coefficient \((k-q)\,\partial_zG_r\). \\

\noindent These are all of the transport contributions that remain at top perturbation order. \\ 

\paragraph{Behavior at a meridional fixed point.}
\noindent  This is precisely where the high-order estimate differs from the low-order one. At low order, the tunable weight contribution is proportional to \(G_r\) and \(G_z\) and therefore vanishes at a meridional fixed point. At high order, the first commutators contain derivatives of the transport field, such as \(\partial_rG_r\), \(\partial_zG_z\), \(\partial_rG_z\), and \(\partial_zG_r\). These quantities need not vanish when \(G_r=G_z=0\), so the high-order transport matrix can retain a nontrivial signed contribution at the fixed point. \\

\noindent After weighted integration by parts, the transport contribution to the diagonal \((i,q)\) entry contains
\begin{equation}
\begin{aligned}
D_{i,q}^{(k),\rm tr}
={}&
\frac12\left(
\partial_rG_r+\partial_zG_z
+G_r\partial_r\log\wgray{\Phi_{i,k}}
+G_z\partial_z\log\wgray{\Phi_{i,k}}
\right)
-q\,\partial_rG_r-(k-q)\,\partial_zG_z.
\end{aligned}
\label{eq:high-linear-complete-transport-add-diagonal}
\end{equation}

\hfill 

\noindent At a meridional fixed point, the terms multiplied by \(G_r\) and \(G_z\) vanish, while the derivatives of the transport field remain, so
\begin{equation}
\begin{aligned}
D_{i,q}^{(k),\rm tr}(r_j^\odot,z_j^\odot)
={}&
\left(\frac12-q\right)\partial_rG_r(r_j^\odot,z_j^\odot)
+\left(\frac12-(k-q)\right)\partial_zG_z(r_j^\odot,z_j^\odot).
\end{aligned}
\label{eq:numerical-high-order-fixed-point-diagonal}
\end{equation}
The neighboring derivative couplings also retain the cross derivatives \(\partial_zG_r\) and \(\partial_rG_z\). Thus the  high-order transport contribution can remain nontrivial even where the low-order weight-dependent transport contribution vanishes or becomes weak. This observation explains why high-order damping is available as a complementary mechanism near the fixed points, although its sign is determined only after the full high-order matrix is assembled.

\hfill 

\subsubsection{Complete High-Order Signed Matrix}
\label{sssec:complete-high-order-signed-matrix}

\noindent The  high-order transport contribution described in the previous subsection determines the transport part of the top-order matrix. It remains to add the linear terms for which all \(k\) derivatives stay on the perturbation. These terms are, component by component,
\begin{align}
\omega\text{-term}:\quad
& (\bar c_u-\lambda+\cures)V_{\omega,q}+2\bar u\,V_{z,q},
\label{eq:high-linear-explicit-top-local-matrix-add}\\
r\text{-term}:\quad
& \left(2\partial_z\bar\psi+\bar c_u-\lambda-\partial_r\bar u^r+\cures\right)V_{r,q}
-\partial_r\bar u^z\,V_{z,q},
\label{eq:high-linear-explicit-top-local-r-add}\\
z\text{-term}:\quad
& -\partial_z\bar u^r\,V_{r,q}
+\left(2\partial_z\bar\psi+\bar c_u-\lambda-\partial_z\bar u^z+\cures\right)V_{z,q}.
\label{eq:high-linear-explicit-top-local-z-add}
\end{align}
\noindent These are exactly the local terms for which all \(k\) derivatives remain on the perturbation.  If a derivative lands on one of their spatially varying coefficients, the resulting perturbation factor has order at most \(k-1\) and is treated as a remainder.

\hfill \\
\noindent As in the low-order construction, it is convenient to write the signed quadratic form directly as a symmetric matrix.  Order the weighted top derivatives by component and derivative index:
\begin{equation}
\mathbf V_k:=
\left(
\Phi_{\omega,k}^{1/2}V_{\omega,0},\ldots,\Phi_{\omega,k}^{1/2}V_{\omega,k},
\Phi_{r,k}^{1/2}V_{r,0},\ldots,\Phi_{r,k}^{1/2}V_{r,k},
\Phi_{z,k}^{1/2}V_{z,0},\ldots,\Phi_{z,k}^{1/2}V_{z,k}
\right)^{\!\top}.
\label{eq:high-linear-weighted-vector-add}
\end{equation}
\noindent With this ordering, the matrix separates into diagonal \((i,q)\) terms, neighboring \(q\)-couplings within each component, and same-\(q\) cross-component couplings.\\

\noindent The matrix \(\mathbb M_{\mathcal L}^{(k),\rm full}\) is indexed by the pairs \((i,q)\), with \(i\in\{\omega,r,z\}\) and \(0\le q\le k\):
\begin{equation}
\mathbb M_{\mathcal L}^{(k),\rm full}
=
\left[
\left(\mathbb M_{\mathcal L}^{(k),\rm full}\right)_{(i,q),(j,p)}
\right]_{\substack{i,j\in\{\omega,r,z\}\\0\le q,p\le k}},
\qquad
\mathbb M_{\mathcal L}^{(k),\rm full}
=\left(\mathbb M_{\mathcal L}^{(k),\rm full}\right)^{\!\top}.
\label{eq:high-linear-explicit-block-matrix-add}
\end{equation}

\hfill \\

\paragraph{Diagonal entries.} The diagonal entries are obtained exactly as in the low-order matrix: weighted integration by parts contributes the transport divergence and weight-gradient terms, the first commutators contribute the terms proportional to \(q\) and \(k-q\), and equations \eqref{eq:high-linear-explicit-top-local-matrix-add}--\eqref{eq:high-linear-explicit-top-local-z-add} contribute the local diagonal coefficients.  Explicitly, for \(0\le q\le k\),
\begin{align}
\left(\mathbb M_{\mathcal L}^{(k),\rm full}\right)_{(\omega,q),(\omega,q)}
={}&
\frac{1}{2}\left(
\partial_rG_r+\partial_zG_z
+G_r\partial_r\log\Phi_{\omega,k}
+G_z\partial_z\log\Phi_{\omega,k}
\right)\nonumber\\
& \quad +\bar c_u-\lambda+\cures
-q\,\partial_rG_r-(k-q)\,\partial_zG_z,
\label{eq:high-linear-explicit-diagonal-omega-add}\\
\left(\mathbb M_{\mathcal L}^{(k),\rm full}\right)_{(r,q),(r,q)}
={}&
\frac{1}{2}\left(
\partial_rG_r+\partial_zG_z
+G_r\partial_r\log\Phi_{r,k}
+G_z\partial_z\log\Phi_{r,k}
\right)\nonumber\\
& \quad  +2\partial_z\bar\psi+\bar c_u-\lambda-\partial_r\bar u^r+\cures
-q\,\partial_rG_r-(k-q)\,\partial_zG_z,
\label{eq:high-linear-explicit-diagonal-r-add}\\
\left(\mathbb M_{\mathcal L}^{(k),\rm full}\right)_{(z,q),(z,q)}
={}&
\frac{1}{2}\left(
\partial_rG_r+\partial_zG_z
+G_r\partial_r\log\Phi_{z,k}
+G_z\partial_z\log\Phi_{z,k}
\right)\nonumber\\
& \quad  +2\partial_z\bar\psi+\bar c_u-\lambda-\partial_z\bar u^z+\cures
-q\,\partial_rG_r-(k-q)\,\partial_zG_z.
\label{eq:high-linear-explicit-diagonal-z-add}
\end{align}

\hfill \\

\paragraph{Same-component off-diagonal entries.}  For a fixed component \(i\), equation~\eqref{eq:high-linear-complete-transport-add} couples \(V_{i,q}\) only to the neighboring derivative splits \(V_{i,q-1}\) and \(V_{i,q+1}\). After symmetrization, each same-component block is therefore tridiagonal in the derivative index \(q\): its only off-diagonal entries connect adjacent indices \((i,q)\) and \((i,q+1)\). For \(0\le q<k\),
\begin{equation}
\left(\mathbb M_{\mathcal L}^{(k),\rm full}\right)_{(i,q),(i,q+1)}
=
\left(\mathbb M_{\mathcal L}^{(k),\rm full}\right)_{(i,q+1),(i,q)}
=
-\frac{1}{2}\left[(k-q)\,\partial_z\bar u^r+(q+1)\,\partial_r\bar u^z\right].
\label{eq:high-linear-explicit-neighbor-entry-add}
\end{equation}

\hfill \\ 

\paragraph{Cross-component entries.} The cross-component entries are obtained as in the low-order matrix by keeping each coefficient together with its required weight conversion and then symmetrizing.  For \(0\le q\le k\),
\begin{align}
\left(\mathbb M_{\mathcal L}^{(k),\rm full}\right)_{(\omega,q),(z,q)}
&=
\left(\mathbb M_{\mathcal L}^{(k),\rm full}\right)_{(z,q),(\omega,q)}
=
\bar u\left(\frac{\Phi_{\omega,k}}{\Phi_{z,k}}\right)^{1/2},
\label{eq:high-linear-explicit-cross-blocks-add}\\
\left(\mathbb M_{\mathcal L}^{(k),\rm full}\right)_{(r,q),(z,q)}
&=
\left(\mathbb M_{\mathcal L}^{(k),\rm full}\right)_{(z,q),(r,q)}=
-\frac{1}{2}\left[
\partial_r\bar u^z\left(\frac{\Phi_{r,k}}{\Phi_{z,k}}\right)^{1/2}
+\partial_z\bar u^r\left(\frac{\Phi_{z,k}}{\Phi_{r,k}}\right)^{1/2}
\right].
\label{eq:high-linear-explicit-cross-rz-add}
\end{align}

\hfill 

\noindent For the $\omega$--$z$ entry, the quantity to certify is the complete multiplier $
\bar u \sqrt{\Phi_{\omega,k}/ \Phi_{z,k}}$. Near the origin, $\wgray{\Phi_{i,k}}=1+\rho^k\wgray{\Phi_i}$ and $\rho^k\wgray{\Phi_i}\to0$, so both high-order weights tend to $1$ and this multiplier has no singular weight-conversion obstruction there. On the rest of the finite box and on the tail, the complete product is enclosed directly. \\

\noindent There are no \(\omega\)-\(r\) entries.  There are also no cross-component entries with different derivative indices \(q\ne p\).  \\

\hfill 

\noindent Together with symmetry, equations \eqref{eq:high-linear-explicit-diagonal-omega-add}--\eqref{eq:high-linear-explicit-cross-rz-add} therefore specify every entry of \(\mathbb M_{\mathcal L}^{(k),\rm full}\).

\hfill \\ 

\paragraph{High-order matrix certificate.}
\noindent The high-order sign is tested on the complete symmetric matrix, because the off-diagonal derivative and component couplings can affect the largest eigenvalue even when all diagonal entries are negative. The required certificate is
\begin{equation}
\operatorname*{ess\,sup}_{(r,z)\in\mathbb D}
\lambda_{\max}\!\left(\mathbb M_{\mathcal L}^{(k),\rm full}(r,z)\right)
\le-\Lambda_{D^\alpha L}^{\rm full}<0.
\label{eq:high-linear-full-matrix-certificate-add}
\end{equation}
The current numerical evaluation indicates a favorable high-order margin in the same neighborhoods where the low-order margin is weakest. A rigorous certificate requires a validated enclosure of the largest eigenvalue over the full domain. A successful full-domain enclosure defines the high-order matrix margin \(\Lambda_{D^\alpha L}^{\rm full}\) carried into the compatibility condition below. 

\hfill \\

\subsection{Compatibility of the Low-Order and High-Order Damping}

\noindent Assuming the localized low-order bound and the high-order matrix bound above, the two estimates have complementary roles. The low-order estimate supplies the negative \(\mathfrak E_0^2\) term, but localization introduces the positive contribution \(A_{\mathcal L}^{k,\rm loc}\mathfrak H_k^2\).  The high-order matrix supplies a negative \(\mathfrak H_k^2\) contribution with margin \(\Lambda_{D^\alpha L}^{\rm full}\). Since the full energy uses the high-order estimate with weight \(\mu_k\), these two top-order contributions must be compared after multiplying the high-order matrix estimate by \(\mu_k\). \\

\noindent Using \eqref{eq:localized-matrix-estimate} and \eqref{eq:high-linear-full-matrix-certificate-add} gives
\begin{equation}
\begin{aligned}
&\int_{\mathbb D}
\bvar{d\xi}^{\top}\mathbb M_{\mathcal L}\bvar{d\xi}\,dr\,dz
+
\mu_k\int_{\mathbb D}
\mathbf V_k^{\top}\mathbb M_{\mathcal L}^{(k),\rm full}\mathbf V_k\,dr\,dz\le
-\Lambda_{\mathcal L}^{\rm mat,loc}\mathfrak E_0^2
-\left(\mu_k\Lambda_{D^\alpha L}^{\rm full}-A_{\mathcal L}^{k,\rm loc}\right)\mathfrak H_k^2.
\end{aligned}
\label{eq:numerical-localized-matrix-combined-estimate}
\end{equation}

\hfill 

\noindent This inequality is the precise sense in which high-order damping compensates for a localized failure of low-order pointwise damping. The high-order estimate does not change \(\lambda_{\max}(\mathbb M_{\mathcal L})\) inside the localized region. Instead, the concentration estimate converts the low-order contribution from that region into the top-order cost \(A_{\mathcal L}^{k,\rm loc}\mathfrak H_k^2\), and the negative high-order term absorbs that cost after multiplication by \(\mu_k\).

\hfill 

\noindent Therefore the localized low-order loss is compatible with the high-order damping provided
\begin{equation}
A_{\mathcal L}^{k,\rm loc}
<
\mu_k\Lambda_{D^\alpha L}^{\rm full}.
\label{eq:numerical-localized-matrix-compatibility-reorg}
\end{equation}
Together with \(\Lambda_{\mathcal L}^{\rm mat,loc}>0\), this condition gives strict negativity of both matrix-level energy coefficients. This comparison is still at the matrix level. The linear terms kept outside the two signed matrices are incorporated in Subsection~\ref{sssec:numerical-remaining-closure} when the full stability margin is assembled.

\hfill \\

\subsection{Weighted PDE Residuals}
\label{sssec:numerical-residuals-damping-weights}

\noindent The same weights used in the damping estimate also enter the norm in which the PDE residual is measured. Improving the signed linear matrix could therefore, in principle, amplify the weighted residual. It is important to check this after the damping weights have been selected rather than treating the two calculations independently. \\

\noindent For the present optimized weights, the current numerical evaluation indicates that the corresponding value of \(\DLamR\) remains small. Thus the weight choice that improves the linear damping does not appear to introduce a competing residual penalty. \\

\noindent The size of this residual matters differently from the nonlinear constants. In the final closing inequality, the residual enters through \(\DLamR/\delta_*\), whereas the nonlinear losses decrease as the stability radius \(\delta_*\) is reduced. Consequently, a large residual cannot be compensated for simply by taking a smaller perturbation neighborhood. The weighted residual must therefore be evaluated after the damping weights have been fixed and carried explicitly into the closing inequality. If it is too large relative to the available margin, the profile itself must be refined in the weighted norms selected by the damping calculation.  

\hfill \\

\subsection{Remaining Certification, Practical Next Steps, and Formalization}
\label{sssec:numerical-remaining-closure}

\noindent The preceding subsections reduce the required linear sign to explicit low-order and high-order matrix conditions. The low-order matrix is responsible for the principal \(\mathfrak E_0^2\) damping, the localized estimate accounts for the small off-axis region where the same uniform pointwise margin is not imposed, and the  high-order matrix provides the complementary top-order mechanism. The remaining task is to pass from these matrix-level conditions to the full stability inequality by incorporating the linear terms kept outside the matrices, the nonlinear losses, and the weighted residual.

\hfill \\ 

\noindent\paragraph{From matrix damping to the full stability margin.}
To pass from the matrix-level bounds to the linear coefficients in the full energy estimate, we now include the linear terms that were intentionally kept outside the signed matrices. In the localized low-order estimate, the resulting constants are
\begin{equation}
\DLamLlow
=
\Lambda_{\mathcal L}^{\rm mat,loc}-A_{\mathcal L}^{0},
\qquad
\DLamLhigh
=
A_{\mathcal L}^{k}+A_{\mathcal L}^{k,\rm loc}.
\label{eq:numerical-low-net-margin}
\end{equation}
At high order, the corresponding constants are
\begin{equation}
\DLamDLhigh
=
\Lambda_{D^\alpha L}^{\rm full}-A_{\rm lin}^{k,\rm rem},
\qquad
\DLamDLlow
=
A_{\rm lin}^{0,\rm rem}.
\label{eq:numerical-high-net-margin}
\end{equation}
After choosing \(\mu_k\), the actual linear coercive margin is
\begin{equation}
\DLamStab(\mu_k)
=
\min\left\{
\DLamLlow-\mu_k\DLamDLlow,
\ \DLamDLhigh-\frac{\DLamLhigh}{\mu_k}
\right\}>0.
\label{eq:numerical-net-stability-margin}
\end{equation}
A positive value of \(\DLamStab\) establishes coercivity before the nonlinear terms are invoked. On the boundary of a ball of radius \(\delta_*\), the cubic contribution is \(\DLamThree\delta_*\). These nonlinear contributions decrease when \(\delta_*\) is reduced. Consequently, large but finite nonlinear constants usually shrink the stability neighborhood that can be certified rather than destroy an already positive linear sign.  This makes the damping calculation particularly informative, because it isolates the linear sign condition that cannot be recovered merely by shrinking the perturbation neighborhood. The full closing condition is
\begin{equation}
\DLamThree\delta_*
+\frac{\DLamR}{\delta_*}
<
\DLamStab.
\label{eq:numerical-exact-euler-smallness}
\end{equation}

\hfill \\

\paragraph{Remaining quantitative certification.}

Closing the proof requires rigorous certification of the explicit constants entering the energy argument, which are presented in full detail in a separate stability paper~\citep{EulerBlowupStability}. These include the direct linear remainders, elliptic and interpolation constants, point-evaluation and modulation bounds, weighted multiplier constants, nonlinear product constants, and the parameters introduced by Cauchy--Schwarz and Young inequalities. This leads to a large but finite collection of scalar inequalities. The task is computationally extensive, but each required quantity is explicitly defined and has a prescribed route to rigorous certification.\\

\noindent The spline representation, interval arithmetic, matrix eigenvalue enclosures, and finite-dimensional optimization provide a direct route to this remaining verification with rigorous error control. If a particular contribution limits the final margin, the corresponding part of the construction can be refined without changing the overall proof mechanism. For example, the profile can be refined preferentially where the weighted residual is largest, while the weights, spline resolution, localization parameters, or analytic estimates can be adjusted when they are the limiting factors. The modular organization of the proof is intended precisely to make these quantitative refinements local to the relevant estimate. \\

\noindent If the complete collection of remaining constants can be rigorously certified with a positive closing margin, then the nonlinear stability argument developed here is fully closed. Together with the finite-time reconstruction from the dynamically rescaled variables, this yields an admissible solution that remains asymptotic to the candidate self-similar profile and becomes singular in finite physical time. Thus the principal immediate task along this route is to determine whether the required constants and inequalities can be rigorously certified with a positive closing margin within the proposed framework. If sufficient margin cannot be obtained, further refinement of the estimates or stability mechanism may be required.\\

\clearpage

\appendix

\addcontentsline{toc}{part}{Appendix}

\section{Convection--Varied Euler System} \label{appx: convection-varied}

\hfill

\subsection{Convection-varied equations}

\noindent Starting from the axisymmetric equations, consider the following variation on the convection terms:
\begin{align}
\label{eq:transformed_1}
\partial_t u_{\bullet}
+ \varepsilon\, u^r\,\partial_{\mathring{r}} u_{\bullet}
+ \varepsilon\, u^z\,\partial_{\mathring{z}} u_{\bullet}
&= 2u_{\bullet}\,\partial_{\mathring{z}} \psi_{\bullet} , \\
\label{eq:transformed_2}
\partial_t \omega_{\bullet}
+ \varepsilon\, u^r\,\partial_{\mathring{r}} \omega_{\bullet}
+ \varepsilon\, u^z\,\partial_{\mathring{z}} \omega_{\bullet}
&= 2u_{\bullet}\,\partial_{\mathring{z}} u_{\bullet} , \\
\label{eq:transformed_3}
-\Bigl[\partial_{\mathring{r}}^2 + \frac{3}{\mathring{r}}\partial_{\mathring{r}} + \partial_{\mathring{z}}^2\Bigr]\psi_{\bullet}
&= \omega_{\bullet}, \\
u^{r} = -\,\mathring{r}\,\partial_{\mathring{z}} \psi_{\bullet}, \qquad \quad \ &  
u^{z} = 2\psi_{\bullet} + \mathring{r}\,\partial_{\mathring{r}} \psi_{\bullet},
\end{align}
with amplification factor $\varepsilon$, which is used to vary the strength of convection that exist within the Euler equations, with $\varepsilon = 1$ corresponding to the full Euler system. \\

\noindent When $\varepsilon = 0$, this corresponds to the non-convective blowup scenario proposed by \citet{HouLei2009}. However, when convection is included, this finite-time blowup scenario can be destroyed. \citet{pengfeiPaper} studied the amplification of this $\varepsilon$ term, and observed a traveling self-similar singularity scenario at the tested values \(\varepsilon \in \{0,0.1,0.2,0.3\} \), whereas the evolution at \(\varepsilon=0.4\) did not exhibit the same behavior. Their dynamic-rescaling computations further indicated that the computed self-similar state became progressively less stable as \(\varepsilon\) increased. These results left open whether self-similar singular profiles persist for stronger convection but are not reached by the forward dynamics.

\hfill \\

\subsection{Traveling-wave profile equations}

\hfill

\begin{proposition}[Traveling-wave profile equations for the convection-varied system] \label{thm: Traveling-wave profile equations for the convection-varied system}
Assume that $(u_\bullet,\omega_\bullet,\psi_\bullet)$ is given by the ansatz \eqref{ansantz1}--\eqref{ansantz3}, that the center satisfies \eqref{eq:center_speed}, and that the exponents are chosen via \eqref{eq:scaling_exponents}. Then the transformed profile $(U,\Omega,\Psi)$ satisfies
\begin{align}
U + \bigl(\lambda r+\varepsilon U^r\bigr)\partial_r U
+ \bigl(C+\lambda z+\varepsilon U^z\bigr)\partial_z U
&=2U\partial_z\Psi,
\label{eq:ProfileEq1_appx}\\
(1+\lambda)\Omega + \bigl(\lambda r+\varepsilon U^r\bigr)\partial_r\Omega
+ \bigl(C+\lambda z+\varepsilon U^z\bigr)\partial_z\Omega
&=2U\partial_z U,
\label{eq:ProfileEq2_appx}\\
-\, \mathcal E \, \Psi&=\Omega,
\label{eq:ProfileEq3_appx}
\end{align}
where
\begin{equation}
U^r=-r\partial_z\Psi,
\qquad
U^z=2\Psi+r\partial_r\Psi.
\label{eq:profile_velocity_recovery_appx}
\end{equation}

\hfill 

\noindent Equivalently, the system can be written in terms of $U,\Omega,\Psi$ only as
\begin{align}
U+r\bigl(\lambda-\varepsilon\partial_z\Psi\bigr)\partial_r U
+\bigl(C+\lambda z+\varepsilon(2\Psi+r\partial_r\Psi)\bigr)\partial_z U
&=2U\partial_z\Psi,
\label{eq:ProfileEq1_expanded_appx}\\
(1+\lambda)\Omega+r\bigl(\lambda-\varepsilon\partial_z\Psi\bigr)\partial_r\Omega
+\bigl(C+\lambda z+\varepsilon(2\Psi+r\partial_r\Psi)\bigr)\partial_z\Omega
&=2U\partial_z U,
\label{eq:ProfileEq2_expanded_appx}\\
-\, \mathcal E \,\Psi&=\Omega.
\label{eq:ProfileEq3_expanded_appx}
\end{align}
\end{proposition}

\begin{grayproof} \Lean \\

\noindent \textbf{Roadmap.} The factors produced by time differentiation, spatial differentiation, convection, stretching, and the elliptic streamfunction relation are balanced to determine the amplitude exponents in \eqref{eq:scaling_exponents}. After the common power of $T-t$ is removed, the Euler equations yield a time-independent profile problem. 

\hfill \\ 

\noindent Set
\begin{equation}
\tau=T-t, \qquad
r=\frac{\mathring r}{\tau^\lambda}, \qquad z=\frac{\mathring z-\mathring z_c(t)}{\tau^\lambda}.
\label{eq:profile_derivation_coordinates_appx}
\end{equation}

\hfill 

\noindent The traveling-center condition \eqref{eq:center_speed} is
\begin{equation}
\frac{d}{dt}\mathring z_c(t)=-C\tau^{\lambda-1}.
\label{eq:profile_derivation_center_speed_appx}
\end{equation}

\hfill 

\noindent Differentiating at fixed physical coordinates $(\mathring r,\mathring z)$ gives
\begin{equation}
\partial_t r=\frac{\lambda r}{\tau}, \qquad \partial_t z=\frac{\lambda z+C}{\tau}.
\label{eq:profile_derivation_coordinate_time_derivatives_appx}
\end{equation}

\noindent  Moreover,
\begin{equation}
\partial_{\mathring r}=\tau^{-\lambda}\partial_r, \qquad \partial_{\mathring z}=\tau^{-\lambda}\partial_z.
\label{eq:profile_derivation_spatial_derivatives_appx}
\end{equation}

\hfill 

\noindent From the streamfunction relations in the physical variables,
\begin{equation}
u^r_{\rm phys}=-\mathring r\,\partial_{\mathring z}\psi_\bullet, \quad \qquad u^z_{\rm phys}=2\psi_\bullet+\mathring r\,\partial_{\mathring r}\psi_\bullet.
\label{eq:profile_physical_velocity_recovery_appx}
\end{equation}

\hfill 

\noindent Therefore,
\begin{equation}
\psi_\bullet(\mathring r,\mathring z,t)=\tau^{c_\psi}\Psi(r,z),
\end{equation}
gives
\begin{align}
u^r_{\rm phys}
&=-\mathring r\,\tau^{c_\psi-\lambda}\partial_z\Psi
=\tau^{c_\psi}\bigl(-r\partial_z\Psi\bigr),\\
u^z_{\rm phys}
&=2\tau^{c_\psi}\Psi
+\mathring r\,\tau^{c_\psi-\lambda}\partial_r\Psi
=\tau^{c_\psi}\bigl(2\Psi+r\partial_r\Psi\bigr).
\end{align}

\hfill 

\noindent Thus
\begin{equation}
u^r_{\rm phys}=\tau^{c_\psi}U^r,
\qquad
u^z_{\rm phys}=\tau^{c_\psi}U^z,
\qquad
U^r=-r\partial_z\Psi,
\qquad
U^z=2\Psi+r\partial_r\Psi.
\label{eq:profile_derivation_velocity_scaling_appx}
\end{equation}

\hfill \\ 

\noindent We first derive the equation for $U$. The ansatz
\begin{equation}
u_\bullet(\mathring r,\mathring z,t)=\tau^{c_u}U(r,z)
\end{equation}
and \eqref{eq:profile_derivation_coordinate_time_derivatives_appx} give
\begin{align}
\partial_t u_\bullet
&=-c_u\tau^{c_u-1}U
+\tau^{c_u}\left(\partial_rU\,\partial_t r
+\partial_zU\,\partial_t z\right) =\tau^{c_u-1}
\left[-c_uU+\lambda r\partial_rU+(\lambda z+C)\partial_zU\right].
\label{eq:profile_derivation_ut_appx}
\end{align}

\hfill 

\noindent The spatial derivatives are
\begin{equation}
\partial_{\mathring r}u_\bullet=\tau^{c_u-\lambda} \, \partial_rU,
\qquad
\partial_{\mathring z}u_\bullet=\tau^{c_u-\lambda} \,\partial_zU.
\label{eq:profile_derivation_u_spatial_appx}
\end{equation}

\hfill 

\noindent Using \eqref{eq:profile_derivation_velocity_scaling_appx}, the convection term is
\begin{align}
\varepsilon \, u^r_{\rm phys} \, \partial_{\mathring r}u_\bullet
+\varepsilon \, u^z_{\rm phys} \, \partial_{\mathring z}u_\bullet
&=\tau^{c_u+c_\psi-\lambda} \, \varepsilon
\left(U^r\partial_rU+U^z\partial_zU\right).
\label{eq:profile_derivation_u_convection_appx}
\end{align}

\hfill 

\noindent The stretching term is
\begin{equation}
2u_\bullet \, \partial_{\mathring z}\psi_\bullet
=2\tau^{c_u+c_\psi-\lambda} \, U \, \partial_z\Psi.
\label{eq:profile_derivation_u_rhs_appx}
\end{equation}

\hfill  \\ 

\noindent Substitution into the physical $u_\bullet$ equation gives
\begin{align}
&\tau^{c_u-1}
\left[-c_uU+\lambda r\partial_rU+(\lambda z+C)\partial_zU\right]
+\tau^{c_u+c_\psi-\lambda} \, \varepsilon \,
\left(U^r\partial_rU+U^z\partial_zU\right) 
=2\, \tau^{c_u+c_\psi-\lambda} \,U\partial_z\Psi.
\label{eq:profile_derivation_u_transformed_appx}
\end{align}

\hfill 

\noindent The time derivative, convection, and stretching terms have the same power of
$\tau$ when
\begin{equation}
c_u-1=c_u+c_\psi-\lambda,
\qquad\text{hence}\qquad
c_\psi=\lambda-1.
\label{eq:profile_derivation_cpsi_appx}
\end{equation}

\hfill 

\noindent After division by $\tau^{c_u-1}$, the equation becomes
\begin{equation}
-c_uU+\bigl(\lambda r+\varepsilon U^r\bigr)\partial_rU
+\bigl(C+\lambda z+\varepsilon U^z\bigr)\partial_zU
=2U\partial_z\Psi.
\label{eq:profile_derivation_u_profile_general_cu_appx}
\end{equation}

\hfill 

\noindent Using $c_u=-1$ gives \eqref{eq:ProfileEq1_appx}.

\hfill \\

\noindent We next derive the equation for $\Omega$. \\

\hfill 

\noindent From
\begin{equation}
\omega_\bullet(\mathring r,\mathring z,t)=\tau^{c_\omega}\Omega(r,z)
\end{equation}
we obtain
\begin{equation}
\partial_t\omega_\bullet
=\tau^{c_\omega-1}
\left[-c_\omega\Omega+\lambda r\partial_r\Omega
+(\lambda z+C)\partial_z\Omega\right],
\label{eq:profile_derivation_omegat_appx}
\end{equation}
and
\begin{equation}
\partial_{\mathring r}\omega_\bullet
=\tau^{c_\omega-\lambda} \,\partial_r\Omega,
\qquad\quad 
\partial_{\mathring z}\omega_\bullet
=\tau^{c_\omega-\lambda} \, \partial_z\Omega.
\label{eq:profile_derivation_omega_spatial_appx}
\end{equation}

\hfill 

\noindent The convection term is
\begin{equation}
\varepsilon \, u^r_{\rm phys} \, \partial_{\mathring r}\omega_\bullet
+\varepsilon \, u^z_{\rm phys} \, \partial_{\mathring z}\omega_\bullet
=\tau^{c_\omega+c_\psi-\lambda} \, \varepsilon \, 
\left(U^r\partial_r\Omega+U^z\partial_z\Omega\right).
\label{eq:profile_derivation_omega_convection_appx}
\end{equation}

\hfill 

\noindent Using \eqref{eq:profile_derivation_cpsi_appx}, this exponent is
$c_\omega-1$.  \\

\hfill 

\noindent The stretching term is
\begin{equation}
2u_\bullet \, \partial_{\mathring z}u_\bullet
=2 \, \tau^{2c_u-\lambda} \, U\partial_zU.
\label{eq:profile_derivation_omega_rhs_appx}
\end{equation}

\hfill

\noindent Consequently,
\begin{align}
&\tau^{c_\omega-1}
\left[-c_\omega\Omega+\lambda r\partial_r\Omega
+(\lambda z+C)\partial_z\Omega\right]
+\tau^{c_\omega-1} \, \varepsilon
\left(U^r\partial_r\Omega+U^z\partial_z\Omega\right) 
=2 \, \tau^{2c_u-\lambda} \, U\partial_zU.
\label{eq:profile_derivation_omega_transformed_appx}
\end{align}

\hfill 

\noindent Matching the remaining powers gives
\begin{equation}
c_\omega-1=2c_u-\lambda.
\label{eq:profile_derivation_omega_scaling_relation_appx}
\end{equation}

\hfill 

\noindent After division by $\tau^{c_\omega-1}$, we obtain
\begin{equation}
-c_\omega\Omega
+\bigl(\lambda r+\varepsilon U^r\bigr)\partial_r\Omega
+\bigl(C+\lambda z+\varepsilon U^z\bigr)\partial_z\Omega
=2U\partial_zU.
\label{eq:profile_derivation_omega_profile_general_comega_appx}
\end{equation}

\hfill 

\noindent Using $c_\omega=-(1+\lambda)$ gives \eqref{eq:ProfileEq2_appx}. 

\hfill \\

\noindent It remains to check the elliptic streamfunction equation. \\

\noindent Since
\begin{equation}
\partial_{\mathring r\mathring r}\psi_\bullet
=\tau^{c_\psi-2\lambda} \, \partial_{rr}\Psi,
\qquad
\partial_{\mathring z\mathring z}\psi_\bullet
=\tau^{c_\psi-2\lambda} \, \partial_{zz}\Psi,
\label{eq:profile_derivation_elliptic_second_derivatives_appx}
\end{equation}
and since $\mathring r=r\tau^\lambda$, we also have
\begin{equation}
\frac{3}{\mathring r} \, \partial_{\mathring r}\psi_\bullet
=\frac{3}{r\tau^\lambda} \, \tau^{c_\psi-\lambda} \, \partial_r\Psi
=\tau^{c_\psi-2\lambda} \, \frac{3}{r} \, \partial_r\Psi.
\label{eq:profile_derivation_elliptic_radial_weight_appx}
\end{equation}

\hfill 

\noindent Therefore,
\begin{equation}
-\mathcal E_{\mathring r,\mathring z}\psi_\bullet
=-\tau^{c_\psi-2\lambda}\mathcal E\Psi.
\label{eq:profile_derivation_elliptic_operator_appx}
\end{equation}

\hfill \\

\noindent The right-hand side of the physical elliptic equation is
\begin{equation}
\omega_\bullet=\tau^{c_\omega}\Omega.
\label{eq:profile_derivation_elliptic_vorticity_appx}
\end{equation}

\hfill 

\noindent The common power cancels when
\begin{equation}
c_\omega=c_\psi-2\lambda.
\label{eq:profile_derivation_elliptic_scaling_relation_appx}
\end{equation}

\hfill 

\noindent Together with \eqref{eq:profile_derivation_cpsi_appx}, this gives
$c_\omega=-1-\lambda$. Combining this identity with
\eqref{eq:profile_derivation_omega_scaling_relation_appx} gives $c_u=-1$. \\

\noindent After cancelling the common power of $\tau$, the elliptic equation is
\begin{equation}
-\mathcal E\Psi=\Omega.
\label{eq:profile_derivation_elliptic_final_appx}
\end{equation}
\end{grayproof}

\noindent In these profile equations,
\begin{itemize}
    \item The terms $\lambda r\partial_r$ and $\lambda z\partial_z$ come from observing the solution in a shrinking coordinate system.
    \item The term $C\partial_z$ comes from following the center of the profile as it travels along the symmetry axis. 
    \item The terms involving $\varepsilon U^r$ and $\varepsilon U^z$ are the convection terms from the convection-varied Euler system.
    \item  The terms $2U\partial_z\Psi$ and $2U\partial_z U$ on the right-hand side are the stretching and coupling terms inherited from the axisymmetric formulation.
\end{itemize}

\noindent For the convection-varied Euler equations, these profile equations are the steady-state equations for the traveling self-similar profile. Thus, for the convection-varied Euler equations, the traveling-wave ansatz reduces the study of a possible finite-time singularity to an autonomous nonlinear profile problem for $(U,\Omega,\Psi)$ and the parameters $\lambda$ and $C$. The parameter $\lambda$ determines the rate at which the spatial scale collapses, while $C$ determines the speed of the axial drift in the rescaled frame. 

\clearpage

\subsection{Dynamic rescaling equations for the convection-varied systems}

As before, define the elliptic operator
\begin{equation}
    \mathcal{E} \coloneqq \partial_r^2 + \frac{3}{r}\partial_r + \partial_z^2.
\end{equation}

\hfill \\ 

\noindent The proposition below provides the convection-varied axisymmetric Euler equations in the dynamically rescaled variables. 

\hfill \\

\begin{proposition}[Dynamic rescaling equations]
\label{thm:dynamic_rescaling_equations}
In the dynamically rescaled variables, the convection-varied axisymmetric system takes the form
\begin{align}
\partial_t u
+\bigl(\lambda(t)r+\varepsilon u^r\bigr)\partial_r u
+\bigl(\lambda(t)z+C(t)+\varepsilon u^z\bigr)\partial_z u
&=2u\partial_z\psi+c_u(t)u ,
\label{eq:dynamic_rescaled_u}\\
\partial_t\omega
+\bigl(\lambda(t)r+\varepsilon u^r\bigr)\partial_r\omega
+\bigl(\lambda(t)z+C(t)+\varepsilon u^z\bigr)\partial_z\omega
&=2u\partial_z u+c_\omega(t)\omega,
\label{eq:dynamic_rescaled_omega}\\
- \mathcal{E}  \psi
&=\omega,
\label{eq:dynamic_rescaled_psi}
\end{align}
with
\begin{equation}
    u^r=-r\partial_z\psi,
    \qquad
    u^z=2\psi+r\partial_r\psi,
    \label{eq:dynamic_rescaled_stream_intro}
\end{equation}
and
\begin{equation}
c_\omega(t)=c_u(t)-\lambda(t).
\label{eq:dynamic_rescaled_rate_relation}
\end{equation}
\end{proposition}

\begin{grayproof}
The proof is provided in Appendix~\ref{appx: dynamic rescaling equations derivation}.
\end{grayproof}

\hfill  \\

\paragraph{Normalization conditions.}
\noindent
The independent \emph{modulation parameters} are $\lambda(t)$, $C(t)$, and $c_u(t)$. The vorticity amplitude rate is determined by the scaling relation
\begin{equation}
    c_\omega(t)=c_u(t)-\lambda(t).
    \label{eq:dynamic_amplitude_rate_relation}
\end{equation}
Choose the axial translation of the reference profile so that
\begin{equation}
\bar C+\varepsilon\bar u^z_0=\bar C+2\varepsilon\bar\psi_0=0.
\label{eq:convection-varied-stagnation-condition}
\end{equation}

\noindent We impose the normalization conditions
\begin{equation}
\begin{aligned}
&
C(t)+\varepsilon u^z(0,0,t)=0,\\
&\partial_t\bigl(\partial_z\omega\bigr)(0,0,t)=0,
\qquad u(0,0,t)=\bar u_0.
\end{aligned}
\label{eq:normalization_conditions_intro}
\end{equation}

\noindent The first condition fixes the meridional fixed point at the origin. The second preserves the initially prescribed value of $\partial_z\omega(0,0,t)$; evaluating the $z$-differentiated vorticity equation at the origin determines $\lambda(t)$. The third fixes the amplitude at the origin and therefore implies $\partial_tu(0,0,t)=0$. Together these conditions determine $\lambda(t)$, $C(t)$, and $c_u(t)$. In perturbation variables, the stagnation and amplitude conditions are
\begin{equation}
\bvar{\delta C}+2\varepsilon\bvar{\delta\psi}_0=0,
\qquad
\bvar{\delta u}_0=0,
\label{eq:convection-varied-perturbation-normalization}
\end{equation}
so that $\bvar{\delta C}=-2\varepsilon\bvar{\delta\psi}_0$. \\

\clearpage

\paragraph{Limiting behaviour.}
\noindent
If the dynamic rescaling captures a stable traveling self-similar blowup, the rescaled solution should converge as rescaled time advances to a time-independent profile,
\begin{equation}
    u(r,z,t)\to u_\infty(r,z),
    \qquad
    \omega(r,z,t)\to \omega_\infty(r,z),
    \qquad
    \psi(r,z,t)\to \psi_\infty(r,z),\label{eq:profile_convergence_intro_appx}
\end{equation}
\noindent
and the modulation parameters should converge to constants,
\begin{equation}
    \lambda(t)\to \lambda_\infty,
    \qquad
    C(t)\to C_\infty,
    \qquad
    c_u(t)\to c_{u,\infty},
    \qquad
    c_\omega(t)\to c_{\omega,\infty},.
\label{eq:modulation_convergence_intro_appx}
\end{equation}
\noindent
In that limit, the time derivatives in \eqref{eq:dynamic_rescaled_u}--\eqref{eq:dynamic_rescaled_omega} disappear. The steady dynamic equations are
\begin{align}
\Bigl(\lambda_\infty \, r+\varepsilon u^r_\infty\Bigr)\partial_r u_\infty
+\Bigl(\lambda_\infty \, z+ C_\infty +\varepsilon u^z_\infty \Bigr)\partial_z u_\infty
&=2 u_\infty \, \partial_z \psi_\infty+ c_{u,\infty} \, u_\infty , 
\label{eq:steady_dynamic_u_appx}\\
\Bigl(\lambda_\infty \, r+\varepsilon u^r_\infty\Bigr)\partial_r\omega_\infty
+\Bigl(\lambda_\infty \, z+C_\infty+\varepsilon u^z_\infty\Bigr)\partial_z\omega_\infty
&=2 u_\infty \, \partial_z u_\infty+ c_{\omega,\infty} \, \omega_\infty,
\label{eq:steady_dynamic_w_appx}
\end{align}
\noindent
with $u^r_\infty =-r\partial_z \psi_\infty$ and $u^z_\infty=2 \psi_\infty+r\partial_r\psi_\infty$. \\ 

\hfill 

\noindent
For the traveling-wave ansatz in \eqref{ansantz1}--\eqref{ansantz3}, the limiting amplitude rates are
\begin{equation}
    c_{u,\infty}=-1,
    \qquad
    c_{\omega,\infty} = -\bigl(1+\lambda_\infty\bigr).
    \label{eq:limiting_amplitude_rates_appx}
\end{equation}

\noindent In the Euler setting, substituting \eqref{eq:limiting_amplitude_rates_appx} into \eqref{eq:steady_dynamic_u_appx}--\eqref{eq:steady_dynamic_w_appx} gives the traveling-wave profile equations for general $\lambda_\infty$. \\

\noindent Therefore, the steady states of the Euler dynamic rescaling system reproduce the Euler traveling-wave profile equations for general $\lambda_\infty$. After writing
\begin{equation}
\bigl(u_\infty,\omega_\infty,\psi_\infty,C_\infty\bigr) = \bigl(U,\Omega,\Psi,C\bigr),
\end{equation}
the limiting dynamic equations agree with the corresponding profile equations. In this sense, the traveling-wave profiles are fixed points of the appropriate dynamically rescaled evolution.

\newpage

\subsection{Numerical results for the convection-varied equations} \label{appx: convection-varied numerics}

We next investigate the family of self-similar profiles generated by varying the convection amplitude $\varepsilon\in[0,1]$. Concretely, we optimize separate PINN models using \sssoap for each
\begin{equation}
    \varepsilon \in \{0, 0.1,0.2,0.3,0.4,0.5,0.6,0.7,0.8,0.9,1\}.
\end{equation}
In this setting the scaling parameter $\lambda$ and $C$ are treated as free parameters and optimized jointly with the PINN models.\\

\begin{figure}[h]
  \centering
  \includegraphics[width=1\linewidth]{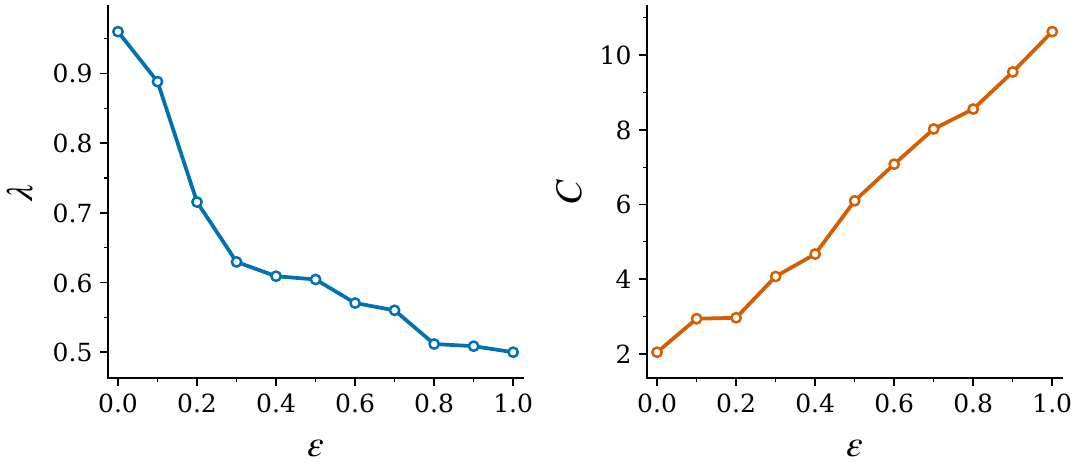} \vspace{-6mm}
  \caption{Scaling parameter $\lambda$ as a function of $\varepsilon$ for the convection-varied Euler system.}
  \label{fig:lambdaC_vs_eps}
\end{figure}

\noindent  We obtain approximate profiles for the different values of $\varepsilon$ considered, and \Cref{fig:lambdaC_vs_eps} displays the converged values of $\lambda$ and $C$ as functions of the convection parameter $\varepsilon$. We see that $\lambda$ decreases steadily toward the critical value $0.5$ as $\varepsilon$ approaches $1$, while $C$ increases approximately linearly.

\hfill \\

\clearpage

\section{PDE in the Computational Domain} \label{app:sinhDerivation}

\hfill \\

\noindent
We use the coordinate transformation
\begin{equation}
    r=\sinh\!\left(r_{\!\mathfrak{comp}}\right),
    \qquad
    z=\sinh\!\left(z_{\mathfrak{comp}}\right),
\end{equation}
and write the corresponding profiles in computational coordinates as
\begin{equation}
    (u,\omega,\psi)\bigl(r_{\!\mathfrak{comp}},z_{\mathfrak{comp}}\bigr)
    =
    (U,\Omega,\Psi)\!\left(
    \sinh(r_{\!\mathfrak{comp}}),
    \sinh(z_{\mathfrak{comp}})
    \right).
\end{equation}

\hfill 

\noindent Since
\begin{equation}
    \frac{dr}{dr_{\!\mathfrak{comp}}}
    =\cosh\!\left(r_{\!\mathfrak{comp}}\right),
    \qquad
    \frac{dz}{dz_{\mathfrak{comp}}}
    =\cosh\!\left(z_{\mathfrak{comp}}\right),
\end{equation}
the chain rule gives
\begin{equation}
    \partial_r
    =\frac{1}{\cosh(r_{\!\mathfrak{comp}})}
      \partial_{r_{\!\mathfrak{comp}}},
    \qquad
    \partial_z
    =\frac{1}{\cosh(z_{\mathfrak{comp}})}
      \partial_{z_{\mathfrak{comp}}}.
\end{equation}

\hfill \\

\noindent The derivatives are related by
\begin{equation}
    f_r
    =\frac{f_{r_{\!\mathfrak{comp}}}}
    {\cosh(r_{\!\mathfrak{comp}})},
    \qquad
    f_z
    =\frac{f_{z_{\mathfrak{comp}}}}
    {\cosh(z_{\mathfrak{comp}})}.
    \label{eq:firstDerivativesComp}
\end{equation}
Applying the chain rule once more yields
\begin{align}
    f_{rr}
    &=
      \frac{f_{r_{\!\mathfrak{comp}}r_{\!\mathfrak{comp}}}}
      {\cosh^2(r_{\!\mathfrak{comp}})}
      -
      \frac{\sinh(r_{\!\mathfrak{comp}})}
      {\cosh^3(r_{\!\mathfrak{comp}})}
      f_{r_{\!\mathfrak{comp}}}, \qquad 
    f_{zz}
    =
\frac{f_{z_{\mathfrak{comp}}z_{\mathfrak{comp}}}}
      {\cosh^2(z_{\mathfrak{comp}})}
      -
      \frac{\sinh(z_{\mathfrak{comp}})}
      {\cosh^3(z_{\mathfrak{comp}})}
      f_{z_{\mathfrak{comp}}}. \label{eq:secondDerivativeZComp}
\end{align}

\hfill \\

\noindent
The elliptic equation involves
$f_{rr}+\frac{3}{r}f_r+f_{zz}$. In computational coordinates, we denote this expression by $\mathcal E_{\!\mathfrak{comp}}f$. Using \eqref{eq:firstDerivativesComp}--\eqref{eq:secondDerivativeZComp},
\begin{align}
    \mathcal E_{\!\mathfrak{comp}}f
    :={}&
    \frac{f_{r_{\!\mathfrak{comp}}r_{\!\mathfrak{comp}}}}
    {\cosh^2(r_{\!\mathfrak{comp}})}
    -\frac{\sinh(r_{\!\mathfrak{comp}})}
    {\cosh^3(r_{\!\mathfrak{comp}})}
    f_{r_{\!\mathfrak{comp}}}
    +\frac{3f_{r_{\!\mathfrak{comp}}}}
    {\sinh(r_{\!\mathfrak{comp}})\cosh(r_{\!\mathfrak{comp}})}
    \notag\\
    &+
    \frac{f_{z_{\mathfrak{comp}}z_{\mathfrak{comp}}}}
    {\cosh^2(z_{\mathfrak{comp}})}
    -\frac{\sinh(z_{\mathfrak{comp}})}
    {\cosh^3(z_{\mathfrak{comp}})}
    f_{z_{\mathfrak{comp}}}.
    \label{eq:Ecomp}
\end{align}

\noindent
The factor $1/\sinh(r_{\!\mathfrak{comp}})$ in \eqref{eq:Ecomp} appears singular at the axis. If the profile is smooth and even in $r$, then it is also even in $r_{\!\mathfrak{comp}}$. Near the axis,
\begin{equation}
    f_{r_{\!\mathfrak{comp}}}
    =f_{r_{\!\mathfrak{comp}}r_{\!\mathfrak{comp}}}(0,z_{\mathfrak{comp}})
    \,r_{\!\mathfrak{comp}}+\mathcal{O}(r_{\!\mathfrak{comp}}),
    \qquad
    \sinh(r_{\!\mathfrak{comp}})\cosh(r_{\!\mathfrak{comp}})
    =r_{\!\mathfrak{comp}}+\mathcal{O}(r_{\!\mathfrak{comp}}).
\end{equation}

\hfill 

\noindent Therefore,
\begin{equation}
    \lim_{r_{\!\mathfrak{comp}}\to0} \ 
    \frac{f_{r_{\!\mathfrak{comp}}}}
    {\sinh(r_{\!\mathfrak{comp}})\cosh(r_{\!\mathfrak{comp}})}
    =f_{r_{\!\mathfrak{comp}}r_{\!\mathfrak{comp}}}(0,z_{\mathfrak{comp}}),
\end{equation}
The apparent singularity is therefore removable. At the axis, the radial part of $\mathcal E_{\!\mathfrak{comp}}f$ extends continuously to $4f_{r_{\!\mathfrak{comp}}r_{\!\mathfrak{comp}}}(0,z_{\mathfrak{comp}})$. \\

\noindent
The remaining factors in the transport terms are
\begin{equation}
    r f_r
    =\frac{\sinh(r_{\!\mathfrak{comp}})}
    {\cosh(r_{\!\mathfrak{comp}})}
    f_{r_{\!\mathfrak{comp}}},
    \qquad
    z=\sinh(z_{\mathfrak{comp}}).
\end{equation}

\hfill 

\noindent Substituting these identities into the profile equations gives the Euler system in computational coordinates
\vspace{4mm}
\begin{align}
0={}&
 u
-\frac{2u\,\psi_{z_{\mathfrak{comp}}}}
{\cosh(z_{\mathfrak{comp}})}
+\frac{\sinh(r_{\!\mathfrak{comp}})}
{\cosh(r_{\!\mathfrak{comp}})}
\,u_{r_{\!\mathfrak{comp}}}
\left(
\lambda
-\varepsilon\,\frac{\psi_{z_{\mathfrak{comp}}}}
{\cosh(z_{\mathfrak{comp}})}
\right)
\notag\\
&\qquad
+\frac{u_{z_{\mathfrak{comp}}}}
{\cosh(z_{\mathfrak{comp}})}
\left(
C+\lambda\sinh(z_{\mathfrak{comp}})
+\varepsilon\left(
2\psi
+\frac{\sinh(r_{\!\mathfrak{comp}})}
{\cosh(r_{\!\mathfrak{comp}})}
\psi_{r_{\!\mathfrak{comp}}}
\right)
\right),
\label{eq:NSCompU}
\\[24pt]
0={}&
(1+\lambda)\omega
-\frac{2u\,u_{z_{\mathfrak{comp}}}}
{\cosh(z_{\mathfrak{comp}})}
+\frac{\sinh(r_{\!\mathfrak{comp}})}
{\cosh(r_{\!\mathfrak{comp}})}
\,\omega_{r_{\!\mathfrak{comp}}}
\left(
\lambda
-\varepsilon\,\frac{\psi_{z_{\mathfrak{comp}}}}
{\cosh(z_{\mathfrak{comp}})}
\right)
\notag\\
&\qquad
+\frac{\omega_{z_{\mathfrak{comp}}}}
{\cosh(z_{\mathfrak{comp}})}
\left(
C+\lambda\sinh(z_{\mathfrak{comp}})
+\varepsilon\left(
2\psi
+\frac{\sinh(r_{\!\mathfrak{comp}})}
{\cosh(r_{\!\mathfrak{comp}})}
\psi_{r_{\!\mathfrak{comp}}}
\right)
\right),
\label{eq:NSCompOmega}
\\[24pt]
0={}&
\omega+\mathcal E_{\!\mathfrak{comp}} \, \psi,
\label{eq:NSCompPsi}
\end{align}
\vspace{2mm}

\noindent where $\varepsilon$, as before, denotes the convection amplification factor, and
\vspace{2mm}
\begin{align}
    \mathcal E_{\!\mathfrak{comp}}f
    :={}&
    \frac{f_{r_{\!\mathfrak{comp}}r_{\!\mathfrak{comp}}}}
    {\cosh^2(r_{\!\mathfrak{comp}})}
    -\frac{\sinh(r_{\!\mathfrak{comp}})}
    {\cosh^3(r_{\!\mathfrak{comp}})}
    f_{r_{\!\mathfrak{comp}}}
    +\frac{3f_{r_{\!\mathfrak{comp}}}}
    {\sinh(r_{\!\mathfrak{comp}})\cosh(r_{\!\mathfrak{comp}})}
    \notag\\
    &+
    \frac{f_{z_{\mathfrak{comp}}z_{\mathfrak{comp}}}}
    {\cosh^2(z_{\mathfrak{comp}})}
    -\frac{\sinh(z_{\mathfrak{comp}})}
    {\cosh^3(z_{\mathfrak{comp}})}
    f_{z_{\mathfrak{comp}}}.
\end{align}

\hfill \\

\noindent As before, the full Euler system is obtained by setting $\varepsilon=1$.

\hfill \\ 

\noindent \emph{ \checkmark \ These formulas have been formalized and verified in Lean.}

\clearpage

\section{Derivations of the Profile Equations} \label{appx: Traveling-wave profile equations derivation}

\hfill  

\begin{proposition}[Traveling-wave profile equations]
Assume that $(u_\bullet,\omega_\bullet,\psi_\bullet)$ is given by the ansatz \eqref{ansantz1}--\eqref{ansantz3}, that the center satisfies \eqref{eq:center_speed}, and that the exponents are chosen as in \eqref{eq:scaling_exponents}.  Set
\begin{equation}
 \lambda=\frac12. \label{eq:profile_eps1_fixed_lambda_appx}
\end{equation}
 Then $(U,\Omega,\Psi)$ satisfies
\begin{align}
U + \bigl(\lambda r+U^r\bigr)\partial_rU + \bigl(C+\lambda z+U^z\bigr)\partial_zU
&=2U\partial_z\Psi,
\label{eq:ProfileEq1_eps1_appx}\\
(1+\lambda)\Omega + \bigl(\lambda r+U^r\bigr)\partial_r\Omega
+ \bigl(C+\lambda z+U^z\bigr)\partial_z\Omega
&=2U\partial_zU,
\label{eq:ProfileEq2_eps1_appx}\\
-\, \mathcal E \,\Psi&=\Omega,
\label{eq:ProfileEq3_eps1_appx}
\end{align}
where
\begin{equation}
U^r=-r\partial_z\Psi,
\qquad
U^z=2\Psi+r\partial_r\Psi.
\label{eq:profile_velocity_recovery_eps1_appx}
\end{equation}
Equivalently,
\begin{align}
U+r\bigl(\lambda-\partial_z\Psi\bigr)\partial_rU
+\bigl(C+\lambda z+2\Psi+r\partial_r\Psi\bigr)\partial_zU
&=2U\partial_z\Psi,
\label{eq:ProfileEq1_expanded_eps1_appx}\\
(1+\lambda)\Omega+r\bigl(\lambda-\partial_z\Psi\bigr)\partial_r\Omega
+\bigl(C+\lambda z+2\Psi+r\partial_r\Psi\bigr)\partial_z\Omega
&=2U\partial_zU,
\label{eq:ProfileEq2_expanded_eps1_appx}\\
-\, \mathcal E \,\Psi&=\Omega.
\label{eq:ProfileEq3_expanded_eps1_appx}
\end{align}
\end{proposition}

\begin{grayproof} \Lean \\

\noindent  We specialize the convection-varied calculation in the proof of \Cref{thm: Traveling-wave profile equations for the convection-varied system}  to $\varepsilon=1$ and $\lambda=1/2$.  \\
 
 \noindent The scaling exponents in \eqref{eq:scaling_exponents} become
\begin{equation}
c_u=-1,
\qquad
c_\omega=-\frac32,
\qquad
c_\psi=-\frac12.
\end{equation}

\hfill 

\noindent The coordinate identities in
\eqref{eq:profile_derivation_coordinate_time_derivatives_appx} and
\eqref{eq:profile_derivation_spatial_derivatives_appx} therefore reduce to
\begin{equation}
\partial_t r=\frac{r}{2\tau},
\qquad
\partial_t z=\frac{z/2+C}{\tau},
\qquad
\partial_{\mathring r}=\tau^{-1/2}\partial_r,
\qquad
\partial_{\mathring z}=\tau^{-1/2}\partial_z.
\end{equation}

\hfill \\

\noindent For the $u_\bullet$ equation, the time derivative is
\begin{equation}
\partial_tu_\bullet
=\tau^{-2}\left[
U+\frac12 r\partial_rU+\left(C+\frac12 z\right)\partial_zU
\right].
\end{equation}
\hfill 

\noindent Because $c_\psi=-\frac12$, the physical meridian velocity components satisfy
\begin{equation}
u^r_{\rm phys}=\tau^{-1/2} \, U^r,
\qquad
u^z_{\rm phys}=\tau^{-1/2} \, U^z,
\qquad
U^r=-r\partial_z\Psi,
\qquad
U^z=2\Psi+r\partial_r\Psi.
\end{equation}
\hfill 

\noindent Also,
\begin{equation} \partial_{\mathring r}u_\bullet=\tau^{-3/2}\partial_rU, \qquad \partial_{\mathring z}u_\bullet=\tau^{-3/2}\partial_zU.\end{equation}

\hfill 

\noindent Hence the convection terms and the stretching term
$2u_\bullet \, \partial_{\mathring z}\psi_\bullet$ all carry the same factor
$\tau^{-2}$ as the time derivative.  

\clearpage

\hfill 

\noindent For the vorticity equation,
\begin{equation}
\partial_t\omega_\bullet
=\tau^{-5/2}\left[
\frac32\Omega+\frac12 r\partial_r\Omega
+\left(C+\frac12 z\right)\partial_z\Omega
\right].
\end{equation}

\hfill 

\noindent The convection terms scale like
$\tau^{-1/2}\tau^{-2}=\tau^{-5/2}$. \\

\noindent The stretching term
$2u_\bullet \, \partial_{\mathring z}u_\bullet$ scales like
$\tau^{-1}\tau^{-3/2}=\tau^{-5/2}$.\\

\noindent Canceling the common factor gives the desired equation. \\

\hfill \\

\noindent Finally,
\begin{equation}
    c_\psi-2\lambda=-\frac12-1=-\frac32=c_\omega,
\end{equation}
so both sides of the
physical elliptic equation have the factor $\tau^{-3/2}$.  Canceling it gives the desired equation. \\
\end{grayproof}

\clearpage 

\section{Derivation of the Dynamic Rescaling Equations} \label{appx: dynamic rescaling equations derivation}

\hfill

\noindent In this Appendix, we derive the dynamic rescaling equations for the convection-varied system in the more general setting where the spatial rescaling rate $\lambda$ is time-dependent, since in the convection-varied case one does not necessarily have $\lambda=1/2$. The corresponding equations for the Euler system are recovered as the special case $\varepsilon=1$ and $\lambda=1/2$.

\hfill \\

\begin{proposition}[Dynamic rescaling equations]
In the dynamically rescaled variables, the convection-varied axisymmetric system takes the form
\begin{align}
\partial_t u
+\bigl(\lambda(t)r+\varepsilon u^r\bigr)\partial_r u
+\bigl(\lambda(t)z+C(t)+\varepsilon u^z\bigr)\partial_z u
&=2u\partial_z\psi+c_u(t)u  \\
\partial_t\omega
+\bigl(\lambda(t)r+\varepsilon u^r\bigr)\partial_r\omega
+\bigl(\lambda(t)z+C(t)+\varepsilon u^z\bigr)\partial_z\omega
&=2u\partial_z u+c_\omega(t)\omega \\
- \mathcal{E} \psi
&=\omega,
\end{align}
with
\begin{equation}
    u^r=-r\partial_z\psi,
    \qquad
    u^z=2\psi+r\partial_r\psi,
    \label{eq:dynamic_rescaled_stream_intro_appx}
\end{equation}
and
\begin{equation}
c_\omega(t)=c_u(t)-\lambda(t).
\label{eq:dynamic_rescaled_rate_relation_appx}
\end{equation}
\end{proposition}

\begin{grayproof} \Lean \\

\noindent\textit{Coordinate and amplitude scaling.}
At fixed physical coordinates, differentiating equation \eqref{eq:dynamic_spatial_rescaling_setup} and using \eqref{eq:dynamic_rate_definitions_setup}--\eqref{eq:dynamic_axial_drift_setup} gives
\begin{equation}
\left(\partial_t r\right)_{\mathring r,\mathring z}=\lambda(t)r,
\qquad
\left(\partial_t z\right)_{\mathring r,\mathring z}
=\lambda(t)z+C(t).
\label{eq:dynamic_coordinate_drift_appx}
\end{equation}

\hfill 

\noindent For any rescaled scalar field $f(r,z,t)$,
\begin{equation}
\left(\partial_t f\right)_{\mathring r,\mathring z}
=\partial_tf+\lambda(t)  r \, \partial_rf
+\bigl(\lambda(t)z+C(t)\bigr) \, \partial_zf.
\label{eq:dynamic_chain_rule_appx}
\end{equation}

\hfill 

\noindent Applying this identity to the amplitudes in \eqref{eq:dynamic_u_omega_rescaling_setup}--\eqref{eq:dynamic_psi_rescaling_setup} gives
\begin{align}
\partial_{\mathring t}u_\bullet
&=\mathrm{s}_u \, \frac{dt}{d\mathring t} \, 
\left[
\partial_t u+\lambda r \, \partial_r u+(\lambda z+C)\partial_z u-c_u u
\right],
\label{eq:dynamic_physical_time_u_appx}\\
\partial_{\mathring t}\omega_\bullet
&=\mathrm{s}_\omega \, \frac{dt}{d\mathring t} \, 
\left[
\partial_t\omega+\lambda r \, \partial_r\omega
+(\lambda z+C)\partial_z\omega-c_\omega\omega
\right].
\label{eq:dynamic_physical_time_omega_appx}
\end{align}

\hfill 

\noindent The streamfunction scaling gives the physical meridian velocities
\begin{equation}
u^r_{\rm phys}
=\mathrm{s}_\omega \, \mathrm{s}_r^2 \, u^r,
\qquad \quad 
u^z_{\rm phys}
=\mathrm{s}_\omega \, \mathrm{s}_r^2 \, u^z,
\label{eq:dynamic_physical_velocity_scaling_appx}
\end{equation}
where
\begin{equation}
u^r=-r\partial_z\psi,
\qquad
u^z=2\psi+r\partial_r\psi.
\label{eq:dynamic_velocity_recovery_proof_appx}
\end{equation}

\hfill 

\noindent Since
$\partial_{\mathring r}=\mathrm{s}_r^{-1}\partial_r$ and
$\partial_{\mathring z}=\mathrm{s}_r^{-1}\partial_z$, the physical
convection and stretching terms in the $u_\bullet$ equation have the common
coefficient
\begin{equation}
\mathrm{s}_\omega \, \mathrm{s}_r^2 \, 
\frac{\mathrm{s}_u}{\mathrm{s}_r}
=\mathrm{s}_u \, \mathrm{s}_\omega \, \mathrm{s}_r
=\mathrm{s}_u \, \frac{dt}{d\mathring t}.
\label{eq:dynamic_u_common_coefficient_appx}
\end{equation}

\hfill 

\noindent In the vorticity equation, the convection coefficient is $\mathrm{s}_\omega^2 \, \mathrm{s}_r$, while the stretching coefficient is
\begin{equation}
\frac{\mathrm{s}_u^2}{\mathrm{s}_r}
=\mathrm{s}_\omega^2 \, \mathrm{s}_r
=\mathrm{s}_\omega \, \frac{dt}{d\mathring t}.
\label{eq:dynamic_omega_common_coefficient_appx}
\end{equation}

\hfill 

\noindent Thus the relations in
\eqref{eq:dynamic_amplitude_time_rescaling_setup} normalize all terms in both equations.  Taking the logarithmic derivative of
$\mathrm{s}_u=\mathrm{s}_\omega \, \mathrm{s}_r$ gives
\begin{equation}
c_\omega=c_u-\lambda.
\label{eq:dynamic_rate_relation_proof_appx}
\end{equation}

\hfill \\

\noindent\textit{The rescaled $u$ equation.}
Applying the coordinate chain rule \eqref{eq:dynamic_chain_rule_appx}, adding the physical convection terms, and including the amplitude normalization gives
\begin{align}
\partial_t u
&+\lambda(t)r\partial_r u
+\bigl(\lambda(t)z+C(t)\bigr)\partial_z u
+\varepsilon u^r\partial_r u
+\varepsilon u^z\partial_z u
-c_u(t)u =2u\partial_z\psi.
\label{eq:dynamic_u_before_rearranging_appx}
\end{align}

\hfill 

\noindent The first two transport terms come from the moving rescaled frame.  The two terms multiplied by $\varepsilon$ are the physical convection terms. \\

\noindent Grouping the radial and axial derivatives gives
\begin{align}
\partial_t u
&+\bigl(\lambda(t)r+\varepsilon u^r\bigr)\partial_r u
+\bigl(\lambda(t)z+C(t)+\varepsilon u^z\bigr)\partial_z u
-c_u(t)u =2u\partial_z\psi.
\label{eq:dynamic_u_grouped_appx}
\end{align}

\hfill 

\noindent Moving the amplitude term to the right gives
\begin{align}
\partial_t u
&+\bigl(\lambda(t)r+\varepsilon u^r\bigr)\partial_r u
+\bigl(\lambda(t)z+C(t)+\varepsilon u^z\bigr)\partial_z u =2u\partial_z\psi+c_u(t)u.
\label{eq:dynamic_u_final_derivation_appx}
\end{align}

\hfill \\

\noindent\textit{The rescaled vorticity equation.}
The moving frame contributes
\begin{equation}
\partial_t\omega
+\lambda(t)r\partial_r\omega
+\bigl(\lambda(t)z+C(t)\bigr)\partial_z\omega,
\label{eq:dynamic_omega_frame_transport_appx}
\end{equation}
while the physical convection contributes
\begin{equation}
\varepsilon u^r\partial_r\omega
+\varepsilon u^z\partial_z\omega.
\label{eq:dynamic_omega_physical_convection_appx}
\end{equation}
The stretching term is $2u\partial_z u$, and the amplitude normalization contributes $-c_\omega(t)\omega$ on the left. \\

\noindent Thus
\begin{align}
\partial_t\omega
&+\lambda(t)r\partial_r\omega
+\bigl(\lambda(t)z+C(t)\bigr)\partial_z\omega
+\varepsilon u^r\partial_r\omega
+\varepsilon u^z\partial_z\omega
-c_\omega(t)\omega =2u\partial_z u.
\label{eq:dynamic_omega_before_rearranging_appx}
\end{align}

\hfill 

\noindent Combining the frame transport with the physical convection gives
\begin{align}
\partial_t\omega
&+\bigl(\lambda(t)r+\varepsilon u^r\bigr)\partial_r\omega
+\bigl(\lambda(t)z+C(t)+\varepsilon u^z\bigr)\partial_z\omega
-c_\omega(t)\omega =2u\partial_z u.
\label{eq:dynamic_omega_grouped_appx}
\end{align}

\hfill 

\noindent Moving the amplitude term to the right gives
\begin{align}
\partial_t\omega
&+\bigl(\lambda(t)r+\varepsilon u^r\bigr)\partial_r\omega
+\bigl(\lambda(t)z+C(t)+\varepsilon u^z\bigr)\partial_z\omega =2u\partial_z u+c_\omega(t)\omega.
\label{eq:dynamic_omega_final_derivation_appx}
\end{align}

\hfill 

\noindent\textit{The elliptic equation and velocity recovery.}
The streamfunction scaling gives
\begin{equation}
-\mathcal E_{\mathring r,\mathring z}\psi_\bullet
=-\mathrm{s}_\omega(t) \, \mathcal E \, \psi,
\qquad
\omega_\bullet=\mathrm{s}_\omega(t) \, \omega.
\label{eq:dynamic_elliptic_scaling_appx}
\end{equation}

\hfill 

\noindent Canceling the common factor $\mathrm{s}_\omega(t)$ gives
\begin{equation}
-\mathcal E\psi=\omega.
\label{eq:dynamic_elliptic_final_proof_appx}
\end{equation}

\hfill 

\noindent Starting from the
physical recovery formulas
\begin{equation}
u^r_{\rm phys}
=-\mathring r\partial_{\mathring z}\psi_\bullet,
\qquad
u^z_{\rm phys}
=2\psi_\bullet+\mathring r\partial_{\mathring r}\psi_\bullet,
\label{eq:dynamic_physical_velocity_recovery_appx}
\end{equation}
and using \eqref{eq:dynamic_spatial_rescaling_setup}--\eqref{eq:dynamic_psi_rescaling_setup}, the common factor $\mathrm{s}_\omega(t) \, \mathrm{s}_r(t)^2$ cancels.  This gives
\begin{equation}
u^r=-r\partial_z\psi,
\qquad
u^z=2\psi+r\partial_r\psi.
\label{eq:dynamic_velocity_recovery_final_proof_appx}
\end{equation}
\end{grayproof}

\clearpage

\section{Energy Identities}
\label{Energy Derivations_appx}

\noindent In this appendix, we consider the energy functional for the convection-varied system,
\begin{equation}
E_{\varepsilon}
=\int_{-\infty}^{\infty}\int_{0}^{\infty}
\left(\lvert u_\bullet\rvert^2
+ (2-\varepsilon)\lvert \nabla_{\mathring r,\mathring z} \psi_\bullet\rvert^2\right)
\, \mathring r^3\,d\mathring r\,d\mathring z .
\label{eq:energy_functional_appx}
\end{equation}
As before, the energy functional for the Euler system can be obtained as the special case where $\varepsilon=1$, i.e.
\begin{equation}
E
=\int_{-\infty}^{\infty}\int_{0}^{\infty}
\left(\lvert u_\bullet\rvert^2
+ \lvert \nabla_{\mathring r,\mathring z} \psi_\bullet\rvert^2\right)
\, \mathring r^3\,d\mathring r\,d\mathring z .
\end{equation}

\hfill   \\ 

\begin{proposition}[Energy scaling under the self-similar ansatz]
\label{thm:energy-scaling_appx}
Let $\tau=T-t$, and assume that
\begin{equation}
u_\bullet(\mathring r,\mathring z,t)
=
\tau^{-1}
u\!\left(
\frac{\mathring r}{\tau^\lambda},
\frac{\mathring z-\mathring z(t)}{\tau^\lambda}
\right),
\qquad
\psi_\bullet(\mathring r,\mathring z,t)
=
\tau^{\lambda-1}
\psi\!\left(
\frac{\mathring r}{\tau^\lambda},
\frac{\mathring z-\mathring z(t)}{\tau^\lambda}
\right).
\label{eq:energy_scaling_ansatz_appx}
\end{equation}
Equivalently, introduce the rescaled variables
\begin{equation}
r = \frac{\mathring r}{\tau^\lambda},
\qquad
z = \frac{\mathring z - \mathring z(t)}{\tau^\lambda}.
\label{eq:energy_scaling_coordinates_appx}
\end{equation}
Then
\begin{equation}
E_{\varepsilon}(t)
=
\tau^{5\lambda-2}
\int_{-\infty}^{\infty}\int_{0}^{\infty}
\left(
\lvert u\rvert^2
+
(2-\varepsilon)\lvert \nabla_{r,z}\psi\rvert^2
\right)
\,r^3\,dr\,dz .
\label{eq:energy_scaling_relation_appx}
\end{equation}
If the rescaled profile energy is finite and nonzero, then boundedness
of $E_\varepsilon(t)$ as $t\to T^-$ requires $\lambda\ge 2/5$.
\end{proposition}

\begin{grayproof} \Lean \\

\noindent We have
\begin{equation}
\mathring r = r\tau^\lambda,
\qquad
\mathring z = z\tau^\lambda + \mathring z(t).
\label{eq:energy_inverse_coordinates_appx}
\end{equation}
Therefore
\begin{equation}
d\mathring r\,d\mathring z
=
\tau^{2\lambda}\,dr\,dz,
\qquad
\mathring r^3
=
r^3\tau^{3\lambda}.
\label{eq:energy_measure_components_appx}
\end{equation}

\hfill 

\noindent Hence the weighted measure transforms as
\begin{equation}
\mathring r^3\,d\mathring r\,d\mathring z
=
r^3\tau^{5\lambda}\,dr\,dz.
\label{eq:energy_weighted_measure_appx}
\end{equation}

\hfill 

\noindent Next, since $\mathring z(t)$ is independent of $z$, we have
\begin{equation}
\partial_{\mathring r}
=
\tau^{-\lambda}\partial_r,
\qquad
\partial_{\mathring z}
=
\tau^{-\lambda}\partial_z.
\label{eq:energy_derivative_scaling_appx}
\end{equation}

\hfill 

\noindent Thus
\begin{equation}
\nabla_{\mathring r,\mathring z}
=
\tau^{-\lambda}\nabla_{r,z}.
\label{eq:energy_gradient_scaling_appx}
\end{equation}
Using $\psi_\bullet=\tau^{\lambda-1}\psi$, it follows that
\begin{equation}
\nabla_{\mathring r,\mathring z}\psi_\bullet
=
\tau^{\lambda-1}
\nabla_{\mathring r,\mathring z}\psi
=
\tau^{\lambda-1}
\tau^{-\lambda}
\nabla_{r,z}\psi
=
\tau^{-1}\nabla_{r,z}\psi.
\label{eq:energy_stream_gradient_scaling_appx}
\end{equation}

\hfill 

\noindent Also, by the ansatz,
\begin{equation}
u_\bullet=\tau^{-1}u.
\label{eq:energy_u_scaling_appx}
\end{equation}

\hfill 

\noindent Substituting these identities into the definition of $E_\varepsilon$ gives
\begin{align}
E_{\varepsilon}(t)
&=
\int_{-\infty}^{\infty}\int_{0}^{\infty}
\left(
\left\lvert \tau^{-1}u\right\rvert^2
+
(2-\varepsilon)
\left\lvert \tau^{-1}\nabla_{r,z}\psi\right\rvert^2
\right)
r^3\tau^{5\lambda}\,dr\,dz
\label{eq:energy_substitution_step_appx}\\
&=
\tau^{5\lambda-2}
\int_{-\infty}^{\infty}\int_{0}^{\infty}
\left(
\lvert u\rvert^2
+
(2-\varepsilon)\lvert \nabla_{r,z}\psi\rvert^2
\right)
r^3\,dr\,dz .
\label{eq:energy_final_scaling_step_appx}
\end{align}
This proves the claimed scaling relation.\\

\noindent Assume now that the rescaled profile energy
\begin{equation}
\int_{-\infty}^{\infty}\int_{0}^{\infty}
\left(
\lvert u\rvert^2
+
(2-\varepsilon)\lvert \nabla_{r,z}\psi\rvert^2
\right)
r^3\,dr\,dz
\label{eq:rescaled_profile_energy_appx}
\end{equation}
is finite and nonzero. Since $\tau\to0^+$ as $t\to T^-$, boundedness of
$E_\varepsilon(t)$ requires
\begin{equation}
5\lambda-2\ge 0.
\label{eq:energy_exponent_nonnegative_appx}
\end{equation}
Therefore
\begin{equation}
\lambda\ge \frac{2}{5}.
\label{eq:energy_lambda_bound_proof_appx}
\end{equation}
\end{grayproof}

\clearpage

\section{Scaling Symmetry and the Choice of Variables for Linear Stability}
\label{appx: linear stability variables}

\hfill 

\noindent We consider the three-dimensional incompressible Euler equations
\begin{equation}
\partial_t u + (u\cdot \nabla)u + \nabla p = 0,
\qquad
\nabla\cdot u = 0,
\label{eq:euler}
\end{equation}
where $u=u(x,t)\in \mathbb{R}^3$ is the velocity field and $p=p(x,t)\in\mathbb{R}$ is the pressure.

\hfill \\

\noindent\textbf{Scaling of the velocity.}
\emph{Let $(u,p)$ be a solution of \eqref{eq:euler}. For any $\eta>0$ and any $\alpha\in\mathbb{R}$, define}
\begin{equation}
u_\eta(x,t):=\eta^\alpha u(\eta x,\eta^{\alpha+1}t),
\qquad
p_\eta(x,t):=\eta^{2\alpha}p(\eta x,\eta^{\alpha+1}t).
\label{eq:euler-scaling}
\end{equation}
\emph{Then $(u_\eta,p_\eta)$ is again a solution of \eqref{eq:euler}.}

\hfill 

\begin{grayproof} \Lean \\

\noindent Set
\begin{equation}
X=\eta x,
\qquad
T=\eta^{\alpha+1}t.
\end{equation}
Then
\begin{equation}
u_\eta(x,t)=\eta^\alpha u(X,T),
\qquad
p_\eta(x,t)=\eta^{2\alpha}p(X,T).
\end{equation}

\hfill

\noindent We first compute the time derivative. By the chain rule,
\begin{equation}
\partial_t u_\eta(x,t)
=
\partial_t\bigl(\eta^\alpha u(X,T)\bigr)
=
\eta^\alpha(\partial_t u)(X,T)\,\partial_tT.
\end{equation}
Since $\partial_tT=\eta^{\alpha+1}$, this gives
\begin{equation}
\partial_tu_\eta(x,t)
=
\eta^\alpha\eta^{\alpha+1}(\partial_tu)(X,T)
=
\eta^{2\alpha+1}(\partial_tu)(X,T).
\label{eq:proof-time}
\end{equation}

\hfill

\noindent Next we compute the spatial derivatives. Since $X=\eta x$, each derivative with respect to $x_j$ contributes one factor of $\eta$. More precisely,
\begin{equation}
\partial_{x_j}u_\eta(x,t)
=
\partial_{x_j}\bigl(\eta^\alpha u(X,T)\bigr)
=
\eta^\alpha\sum_{k=1}^3(\partial_{X_k}u)(X,T)\,\partial_{x_j}X_k.
\end{equation}
Because $\partial_{x_j}X_k=\eta\delta_{jk}$, we obtain
\begin{equation}
\partial_{x_j}u_\eta(x,t)
=
\eta^\alpha\eta(\partial_{X_j}u)(X,T).
\end{equation}
Therefore
\begin{equation}
\nabla u_\eta(x,t)
=
\eta^{\alpha+1}(\nabla u)(X,T).
\label{eq:proof-grad}
\end{equation}

\hfill

\noindent Using \eqref{eq:proof-grad}, we can now compute the transport term:
\begin{align}
(u_\eta\cdot\nabla)u_\eta
&=
\bigl(\eta^\alpha u(X,T)\bigr)\cdot\nabla\bigl(\eta^\alpha u(X,T)\bigr)
\nonumber \\
&=
\bigl(\eta^\alpha u(X,T)\bigr)\cdot\bigl(\eta^{\alpha+1}(\nabla u)(X,T)\bigr)
\nonumber \\
&=
\eta^{2\alpha+1}\bigl((u\cdot\nabla)u\bigr)(X,T).
\label{eq:proof-transport}
\end{align}

\hfill

\noindent We now compute the pressure gradient. Again by the chain rule,
\begin{equation}
\partial_{x_j}p_\eta(x,t)
=
\partial_{x_j}\bigl(\eta^{2\alpha}p(X,T)\bigr)
=
\eta^{2\alpha}\sum_{k=1}^3(\partial_{X_k}p)(X,T)\,\partial_{x_j}X_k.
\end{equation}
Since $\partial_{x_j}X_k=\eta\delta_{jk}$, it follows that
\begin{equation}
\partial_{x_j}p_\eta(x,t)
=
\eta^{2\alpha}\eta(\partial_{X_j}p)(X,T),
\end{equation}
and therefore
\begin{equation}
\nabla p_\eta(x,t)
=
\eta^{2\alpha+1}(\nabla p)(X,T).
\label{eq:proof-pressure}
\end{equation}

\hfill 

\noindent Combining \eqref{eq:proof-time}, \eqref{eq:proof-transport}, and \eqref{eq:proof-pressure}, we obtain
\begin{align}
\partial_tu_\eta+(u_\eta\cdot\nabla)u_\eta+\nabla p_\eta
&=
\eta^{2\alpha+1}(\partial_tu)(X,T)
+\eta^{2\alpha+1}\bigl((u\cdot\nabla)u\bigr)(X,T)
+\eta^{2\alpha+1}(\nabla p)(X,T)
\nonumber \\
&=
\eta^{2\alpha+1}
\Bigl[\partial_tu+(u\cdot\nabla)u+\nabla p\Bigr](X,T).
\end{align}

\hfill

\noindent Since $(u,p)$ solves \eqref{eq:euler}, the bracketed term vanishes. Hence
\begin{equation}
\partial_tu_\eta+(u_\eta\cdot\nabla)u_\eta+\nabla p_\eta=0.
\end{equation}

\hfill 

\noindent It remains to verify the divergence-free condition. Using the same computation as above,
\begin{align}
\nabla\cdot u_\eta(x,t)
&=
\sum_{j=1}^3\partial_{x_j}\bigl(\eta^\alpha u_j(X,T)\bigr)
=
\eta^\alpha\sum_{j=1}^3\sum_{k=1}^3(\partial_{X_k}u_j)(X,T)\,\partial_{x_j}X_k
\nonumber \\
&=
\eta^\alpha\eta\sum_{j=1}^3(\partial_{X_j}u_j)(X,T)
=
\eta^{\alpha+1}(\nabla\cdot u)(X,T).
\end{align}
Since $\nabla\cdot u=0$, it follows that $\nabla\cdot u_\eta=0$ as well. Thus $(u_\eta,p_\eta)$ satisfies \eqref{eq:euler}.
\end{grayproof}

\hfill \\

\hfill

\noindent We now turn to the vorticity $\omega=\nabla\times u$ and the velocity gradient $\nabla u$.

\hfill

\noindent\textbf{Scaling of the vorticity.}
\emph{Let $\omega=\nabla\times u$. If $u_\eta$ is defined by \eqref{eq:euler-scaling}, then}
\begin{equation}
\omega_\eta(x,t):=\nabla\times u_\eta(x,t)
=
\eta^{\alpha+1}\omega(\eta x,\eta^{\alpha+1}t).
\label{eq:omega-scaling}
\end{equation}

\medskip

\begin{grayproof} \Lean \\ 

\noindent As before, we write
\begin{equation}
X=\eta x,
\qquad
T=\eta^{\alpha+1}t.
\end{equation}
Then
\begin{equation}
u_\eta(x,t)=\eta^\alpha u(X,T).
\end{equation}

\hfill 

\noindent We compute the curl componentwise. For example, the first component is
\begin{align}
(\omega_\eta)_1
&=
\partial_{x_2}(u_\eta)_3-\partial_{x_3}(u_\eta)_2
=
\partial_{x_2}\bigl(\eta^\alpha u_3(X,T)\bigr)
-
\partial_{x_3}\bigl(\eta^\alpha u_2(X,T)\bigr).
\end{align}

\hfill 

\noindent Using the chain rule,
\begin{equation}
\partial_{x_2}\bigl(\eta^\alpha u_3(X,T)\bigr)
=
\eta^\alpha\eta(\partial_{X_2}u_3)(X,T),
\qquad
\partial_{x_3}\bigl(\eta^\alpha u_2(X,T)\bigr)
=
\eta^\alpha\eta(\partial_{X_3}u_2)(X,T).
\end{equation}
Therefore
\begin{equation}
(\omega_\eta)_1
=
\eta^{\alpha+1}\bigl(\partial_{X_2}u_3-\partial_{X_3}u_2\bigr)(X,T)
=
\eta^{\alpha+1}\omega_1(X,T).
\end{equation}

\hfill 

\noindent The same argument gives
\begin{equation}
(\omega_\eta)_2=\eta^{\alpha+1}\omega_2(X,T),
\qquad
(\omega_\eta)_3=\eta^{\alpha+1}\omega_3(X,T).
\end{equation}
Hence
\begin{equation}
\omega_\eta(x,t)=\eta^{\alpha+1}\omega(X,T),
\end{equation}
that is,
\begin{equation}
\omega_\eta(x,t)
=
\eta^{\alpha+1}\omega(\eta x,\eta^{\alpha+1}t).
\end{equation}
\end{grayproof}

\hfill \\

\noindent\textbf{Scaling of the velocity gradient.}
\emph{If $u_\eta$ is defined by \eqref{eq:euler-scaling}, then}
\begin{equation}
\nabla u_\eta(x,t)
=
\eta^{\alpha+1}(\nabla u)(\eta x,\eta^{\alpha+1}t).
\label{eq:gradu-scaling}
\end{equation}

\medskip

\begin{grayproof} \Lean \\

\noindent As before, let $X=\eta x$ and $T=\eta^{\alpha+1}t$. For each pair of indices $i,j$,
\begin{equation}
\partial_{x_j}(u_\eta)_i(x,t)
=
\partial_{x_j}\bigl(\eta^\alpha u_i(X,T)\bigr).
\end{equation}
Applying the chain rule gives
\begin{equation}
\partial_{x_j}(u_\eta)_i(x,t)
=
\eta^\alpha\sum_{k=1}^3(\partial_{X_k}u_i)(X,T)\,\partial_{x_j}X_k.
\end{equation}
Since $\partial_{x_j}X_k=\eta\delta_{jk}$, we obtain
\begin{equation}
\partial_{x_j}(u_\eta)_i(x,t)
=
\eta^\alpha\eta(\partial_{X_j}u_i)(X,T)
=
\eta^{\alpha+1}(\partial_{X_j}u_i)(X,T).
\end{equation}

\hfill

\noindent Because this holds for every $i$ and $j$, it follows that
\begin{equation}
\nabla u_\eta(x,t)
=
\eta^{\alpha+1}(\nabla u)(X,T),
\end{equation}
which is exactly \eqref{eq:gradu-scaling}.
\end{grayproof}

\hfill \\

\noindent The variables $u$, $\omega$, and $\nabla u$ do not all transform in the same way. Indeed, from \eqref{eq:euler-scaling} and \eqref{eq:omega-scaling},
\begin{equation}
u_\eta(x,t)=\eta^\alpha u(\eta x,\eta^{\alpha+1}t),
\qquad
\omega_\eta(x,t)=\eta^{\alpha+1}\omega(\eta x,\eta^{\alpha+1}t).
\end{equation}
Thus the full vorticity carries one extra factor of $\eta$ relative to the full velocity. This is expected since $\omega=\nabla\times u$ contains one spatial derivative of $u$. By contrast, \eqref{eq:omega-scaling} and \eqref{eq:gradu-scaling} show that
\begin{equation}
\omega_\eta(x,t)=\eta^{\alpha+1}\omega(\eta x,\eta^{\alpha+1}t),
\qquad
\nabla u_\eta(x,t)=\eta^{\alpha+1}(\nabla u)(\eta x,\eta^{\alpha+1}t).
\end{equation}

\hfill

\noindent The preceding scaling laws concern the full three-dimensional velocity and vorticity fields. In the stability analysis, however, $u$ and $\omega$ denote the reduced variables
\begin{equation}
u=\frac{u^\theta}{r},
\qquad
\omega=\frac{\omega^\theta}{r}.
\end{equation}
Since division by $r$ contributes one additional factor of $\eta$ under the rescaling $(r,z)\mapsto(\eta r,\eta z)$, the reduced variables satisfy
\begin{equation}
u_\eta(r,z,t)
=
\eta^{\alpha+1}u(\eta r,\eta z,\eta^{\alpha+1}t),
\qquad
\omega_\eta(r,z,t)
=
\eta^{\alpha+2}\omega(\eta r,\eta z,\eta^{\alpha+1}t).
\end{equation}

\hfill 

\noindent Consequently,
\begin{equation}
\nabla u_\eta(r,z,t)
=
\eta^{\alpha+2}(\nabla u)(\eta r,\eta z,\eta^{\alpha+1}t),
\end{equation}
where $\nabla$ acts on the reduced spatial variables $(r,z)$.

\hfill

\noindent Thus the reduced variables $\omega$ and $\nabla u$ transform in the same way. This is the reason why they are the natural pair of quantities in a scale-consistent linear stability analysis. Quantities with different scaling correspond to different differential orders, so comparing them directly mixes objects of different analytic strength. \\

 \clearpage

\section{Linearization and Modulation Equations}
\label{appx: linearization derivation}

\subsection{Setup}

We now derive the linearized Euler perturbation equations, together with the modulation equations used to close the perturbation system. We use \textcolor{blue}{blue} for perturbation quantities, \hlred{\text{red highlighting}} for nonlinear perturbation products. A subscript \(_0\) denotes evaluation at \((r,z)=(0,0)\). \\

\noindent \emph{\checkmark \ All the derivations in Appendix~\ref{appx: linearization derivation} have been formalized and verified in Lean.}

\hfill \\ 
 
 \noindent We start from the dynamically rescaled equations
\begin{align}
\partial_t u
+\bigl(\lambda r+u^r\bigr)\partial_r u
+\bigl(\lambda z+C(t)+u^z\bigr)\partial_z u
&=2u\partial_z\psi+c_u(t)u,\\
\partial_t\omega
+\bigl(\lambda r+u^r\bigr)\partial_r\omega
+\bigl(\lambda z+C(t)+u^z\bigr)\partial_z\omega
&=2u\partial_z u+c_\omega(t)\omega,
\label{eq:dynamic_rescaled_omega_appx}\\
-\mathcal E\psi&=\omega,
\end{align}
where
\begin{align}
\label{eq:radial_velocity_definition_appx}u^r&=-r\partial_z\psi,\\
\label{eq:vertical_velocity_definition_appx}u^z&=2\psi+r\partial_r\psi,\\
\label{eq:comega_cu_relation_appx}c_\omega&=c_u-\lambda,
\end{align}
and
\begin{equation}
\mathcal E f:=\partial_r^2 f+\frac{3}{r}\partial_r f+\partial_z^2 f.
\label{eq:elliptic_operator_definition_appx}
\end{equation}

\hfill \\

\noindent We denote steady-state values using an overhead bar.  Thus
\(\bar u,\bar\omega,\bar\psi,\bar C,\bar c_u,\bar c_\omega\)
form an Euler steady state. The profile may be approximate, so the steady state equations are
\begin{align}
\partial_t\bar u
+(\lambda r+\bar u^r)\partial_r\bar u
+(\lambda z+\bar C+\bar u^z)\partial_z\bar u
&=2\bar u\partial_z\bar\psi+\bar c_u\bar u
\ -\ \mathcal R_u^{\mathsf{raw}},
\label{eq:profile_u_residual_definition_appx}\\
\partial_t\bar\omega
+(\lambda r+\bar u^r)\partial_r\bar\omega
+(\lambda z+\bar C+\bar u^z)\partial_z\bar\omega
&=2\bar u\partial_z\bar u+\bar c_\omega\bar\omega
\ -\ \mathcal R_\omega^{\mathsf{raw}}.
\label{eq:profile_omega_residual_definition_appx}
\end{align}

\hfill 

\noindent We write the perturbation ansatz as
\begin{equation}
u=\bar u+\textcolor{blue}{\delta u},\qquad
\omega=\bar\omega+\textcolor{blue}{\delta\omega},\qquad
\psi=\bar\psi+\textcolor{blue}{\delta\psi},
\label{eq:field_perturbation_ansatz_appx}
\end{equation}
\begin{equation}
C=\bar C+\textcolor{blue}{\delta C},\qquad
c_u=\bar c_u+\cures+\textcolor{blue}{\delta c_u},\qquad
c_\omega=\bar c_\omega+\cures+\textcolor{blue}{\delta c_\omega}.
\end{equation}

\hfill \\

\noindent The induced meridional velocity perturbations are
\begin{equation}
\label{eq:meridional_velocity_decomposition_appx}u^r=\bar u^r+\textcolor{blue}{\delta u^r},
\qquad
u^z=\bar u^z+\textcolor{blue}{\delta u^z},
\end{equation}
with
\begin{equation}
\label{eq:meridional_velocity_perturbations_appx}
\textcolor{blue}{\delta u^r}=-r\partial_z\textcolor{blue}{\delta\psi},
\qquad
\textcolor{blue}{\delta u^z}=2\textcolor{blue}{\delta\psi}+r\partial_r\textcolor{blue}{\delta\psi},
\qquad
\bar c_\omega=\bar c_u-\lambda,
\qquad
\textcolor{blue}{\delta c_\omega}=\textcolor{blue}{\delta c_u}.
\end{equation}

\hfill \\

\noindent \noindent Before linearizing, write the radial and axial transport coefficients of the numerically obtained profile as
\begin{equation}
G_r(r,z)=\lambda r+u^r_{\mathsf{num}}(r,z),
\qquad
G_z(r,z)=\lambda z+C_{\mathsf{num}}+u^z_{\mathsf{num}}(r,z).
\label{eq:numerical-transport-appx-add}
\end{equation}
We choose $z_{\mathsf{shift}}$ so that $(0,z_{\mathsf{shift}})$ becomes the analysis origin. For this point to be a meridional fixed point, we require
\begin{equation}
G_r(0,z_{\mathsf{shift}})=0,
\qquad
G_z(0,z_{\mathsf{shift}})=0.
\label{eq:numerical-transport-stagnation-appx-add}
\end{equation}

\hfill 

\noindent The radial condition is automatic because $u^r_{\mathsf{num}}(0,z)=0$ on the axis. Using $u^z_{\mathsf{num}}(0,z)=2\psi_{\mathsf{num}}(0,z)$, the axial condition $G_z(0,z_{\mathsf{shift}})=0$ becomes
\begin{equation}
C_{\mathsf{num}}+\lambda z_{\mathsf{shift}}+2\psi_{\mathsf{num}}(0,z_{\mathsf{shift}})=0.
\label{eq:shift-root-appx-add}
\end{equation}
We therefore impose this scalar condition and choose $z_{\mathsf{shift}}$ as a numerically certified root. We then recenter each numerical profile field by
\begin{equation}
\bar f(r,z):=f_{\mathsf{num}}(r,z+z_{\mathsf{shift}}).
\label{eq:profile-recentering-appx-add}
\end{equation}

\hfill 

\noindent Denoting the constant term in the recentered axial transport by $\bar C$, we continue to write the reference transport coefficients as
\begin{equation}
G_r(r,z)=\lambda r+\bar u^r(r,z),
\qquad
G_z(r,z)=\lambda z+\bar C+\bar u^z(r,z).
\label{eq:recentered-transport-appx-add}
\end{equation}
With this notation, the imposed root condition gives
\begin{equation}
G_r(0,0)=0,
\qquad
G_z(0,0)=\bar C+2\bar\psi_0=0.
\label{eq:recentered-stagnation-appx-add}
\end{equation}

\noindent Having fixed the shift, set
\begin{equation}
\bar u_0:=\bar u(0,0),
\qquad
\bar u_0\ne0.
\label{eq:recentered-origin-amplitude-appx-add}
\end{equation}
and certify a positive lower bound for $|\bar u_0|$. Thus the recentered numerical profile used in the analysis has a fixed point of the meridional flow at the origin.

\hfill 

\noindent The two scalar modulation conditions used below preserve this choice:
\begin{equation}
C(t)+u^z_0(t)=0,
\qquad
u_0(t)=\bar u_0.
\end{equation}

\hfill 

\noindent Since $u^r_0(t)=0$ identically, these conditions imply
\begin{equation}
\bigl(\lambda r+u^r\bigr)_0=0,
\qquad
\bigl(\lambda z+C+u^z\bigr)_0=0
\end{equation}
for every rescaled time.

\hfill 

\noindent Equivalently, relative to the reference profile,
\begin{equation}
\textcolor{blue}{\delta C}+2\textcolor{blue}{\delta\psi_0}=0,
\qquad
\textcolor{blue}{\delta u_0}=0,
\end{equation}
and the time-differentiated normalization condition used in the second modulation equation is
\begin{equation}
\partial_t\textcolor{blue}{\delta u_0}=0.
\end{equation}

\hfill \\

\subsection{Linearization of the Euler equations}

\subsubsection{Linearization of the first equation}
\hfill

\noindent We consider the $u$-equation
\begin{align}
\partial_t u
+ (\lambda r + u^r)\,\partial_r u
+ (\lambda z+C+ u^z)\,\partial_z u
&= 2u\,\partial_z \psi + c_u(t) u.
\label{eq:full_u_equation_for_linearization_appx}
\end{align}

\hfill 

\noindent Substituting \eqref{eq:field_perturbation_ansatz_appx}--\eqref{eq:meridional_velocity_perturbations_appx} into \eqref{eq:full_u_equation_for_linearization_appx} gives
\begin{align}
&\partial_t(\bar u+\textcolor{blue}{\delta u})
+\bigl(\lambda r+(\bar u^r+\textcolor{blue}{\delta u^r})\bigr)\,
\partial_r(\bar u+\textcolor{blue}{\delta u})
+\bigl(\lambda z+(\bar C+\textcolor{blue}{\delta C})+(\bar u^z+\textcolor{blue}{\delta u^z})\bigr)\,
\partial_z(\bar u+\textcolor{blue}{\delta u})
\nonumber\\
& \qquad \quad =
2(\bar u+\textcolor{blue}{\delta u})\,\partial_z(\bar\psi+\textcolor{blue}{\delta\psi})
+(\bar c_u+\cures+\textcolor{blue}{\delta c_u})(\bar u+\textcolor{blue}{\delta u}).
\label{eq:u_equation_substituted_perturbations_appx}
\end{align}

\hfill 

\noindent Since $\bar u,\bar\psi,\bar C,\bar c_u$ satisfy the steady-state equation, we have
\begin{align}
\partial_t\bar u
+(\lambda r+\bar u^r)\,\partial_r\bar u
+(\lambda z+\bar C+\bar u^z)\,\partial_z\bar u
=
2\bar u\,\partial_z\bar\psi+\bar c_u\,\bar u 
\ -\ \mathcal R_u^{\mathsf{raw}}.
\label{eq:u_steady_state_equation_appx}
\end{align}

\hfill 

\noindent We subtract \eqref{eq:u_steady_state_equation_appx} from \eqref{eq:u_equation_substituted_perturbations_appx} and group the nonlinear perturbation products in red to obtain
\begin{align}
\partial_t \textcolor{blue}{\delta u}
&+(\lambda r+\bar u^r)\,\partial_r \textcolor{blue}{\delta u}
+(\lambda z + \bar C + \bar u^z)\,\partial_z \textcolor{blue}{\delta u}
+(\textcolor{blue}{\delta u^r})\,\partial_r \bar u
+(\textcolor{blue}{\delta C}+\textcolor{blue}{\delta u^z})\,\partial_z \bar u
\label{eq:u_linearization_raw_grouped_appx} \\
&=(2\,\partial_z\bar\psi+\bar c_u+\cures)\,\textcolor{blue}{\delta u}
+2\,\bar u\,\partial_z\textcolor{blue}{\delta\psi}
+\textcolor{blue}{\delta c_u}\,\bar u \ + \ \cures\,\bar u\ + \ \mathcal R_u^{\mathsf{raw}}
\nonumber\\
& \qquad 
\hlred{
-\Bigl(\textcolor{blue}{\delta u^r}\Bigr)\,\partial_r \textcolor{blue}{\delta u}
-\Bigl(\textcolor{blue}{\delta C}+\textcolor{blue}{\delta u^z}\Bigr)\,\partial_z \textcolor{blue}{\delta u}
+2\,\textcolor{blue}{\delta u}\,\partial_z\textcolor{blue}{\delta\psi}
+\textcolor{blue}{\delta c_u}\,\textcolor{blue}{\delta u}}
\end{align}

\hfill 

\noindent Using \eqref{eq:meridional_velocity_perturbations_appx}, we rewrite \eqref{eq:u_linearization_raw_grouped_appx} as
\begin{align}
& \partial_t \textcolor{blue}{\delta u}
+(\lambda r+\bar u^r)\,\partial_r \textcolor{blue}{\delta u}
+(\lambda z + \bar C + \bar u^z)\,\partial_z \textcolor{blue}{\delta u}
- r\,\partial_z\textcolor{blue}{\delta\psi}\,\partial_r\bar u  
+\bigl(\textcolor{blue}{\delta C}+(2\textcolor{blue}{\delta\psi}+r\partial_r\textcolor{blue}{\delta\psi})\bigr)\,\partial_z\bar u
\nonumber\\
& \quad   =
(2\,\partial_z\bar\psi+\bar c_u+\cures)\,\textcolor{blue}{\delta u}
+2\,\bar u\,\partial_z\textcolor{blue}{\delta\psi}
+\textcolor{blue}{\delta c_u}\,\bar u  \ + \ \cures\,\bar u\ + \ \mathcal R_u^{\mathsf{raw}}
\nonumber\\
& \qquad \quad  \ 
\hlred{
+\  r\,\partial_z\textcolor{blue}{\delta\psi}\,\partial_r \textcolor{blue}{\delta u}
-\Bigl(\textcolor{blue}{\delta C}+(2\textcolor{blue}{\delta\psi}+r\partial_r\textcolor{blue}{\delta\psi})\Bigr)\,\partial_z \textcolor{blue}{\delta u}
+2\,\textcolor{blue}{\delta u}\,\partial_z\textcolor{blue}{\delta\psi}
+\textcolor{blue}{\delta c_u}\,\textcolor{blue}{\delta u}
} .
\label{eq:u_linearization_velocity_perturbations_substituted_appx}
\end{align}

\hfill  \\

\noindent Collecting and reorganizing terms gives
\begin{align}
\partial_t \textcolor{blue}{\delta u}
&=
-(\lambda r+\bar u^r)\,\partial_r \textcolor{blue}{\delta u}
-(\lambda z + \bar C + \bar u^z)\,\partial_z \textcolor{blue}{\delta u}
+(2\,\partial_z\bar\psi+\bar c_u+\cures)\,\textcolor{blue}{\delta u}
+2\,\bar u\,\partial_z\textcolor{blue}{\delta\psi}
+\bar u\,\textcolor{blue}{\delta c_u}
\nonumber\\
&\qquad
+ r\,(\partial_z\textcolor{blue}{\delta\psi})\,\partial_r\bar u
-\bigl(2\textcolor{blue}{\delta\psi}+r\partial_r\textcolor{blue}{\delta\psi}\bigr)\,\partial_z\bar u
-(\partial_z\bar u)\,\textcolor{blue}{\delta C} \ + \ \cures\,\bar u\ + \ \mathcal R_u^{\mathsf{raw}}
\nonumber\\
&\qquad \hlred{ \ 
+\  r\,\partial_z\textcolor{blue}{\delta\psi}\,\partial_r \textcolor{blue}{\delta u}
-\Bigl(\textcolor{blue}{\delta C}+(2\textcolor{blue}{\delta\psi}+r\partial_r\textcolor{blue}{\delta\psi})\Bigr)\,\partial_z \textcolor{blue}{\delta u}
+2\,\textcolor{blue}{\delta u}\,\partial_z\textcolor{blue}{\delta\psi}
+\textcolor{blue}{\delta c_u}\,\textcolor{blue}{\delta u}
}.
\label{eq:u_linearization_collected_appx}
\end{align}

\hfill \\

\noindent Define the residual amplitude offset by
\begin{equation}
\cures:=-\frac{\mathcal R_u^{\mathsf{raw}}(0)}{\bar u_0}
\label{eq:residual-phase-offsets-add}
\end{equation}
under the standing nondegeneracy condition $\bar u_0\ne0$. \\

\noindent We absorb the residual amplitude offset into the reference amplitude rate and use the modulation coefficients directly:
\begin{equation}
C=\bar C+\bvar{\delta C},
\qquad
c_u=\bar c_u+\cures+\bvar{\delta c_u},
\qquad
c_\omega=\bar c_\omega+\cures+\bvar{\delta c_u}.
\label{eq:effective-reference-coefficients-appx}
\end{equation}
The residuals associated with this decomposition are
\begin{align}
\mathcal R_u
&:=\mathcal R_u^{\mathsf{raw}}+\cures\,\bar u,
\label{eq:effective-u-residual-appx}\\
\mathcal R_\omega
&:=\mathcal R_\omega^{\mathsf{raw}}+\cures\,\bar\omega.
\label{eq:effective-omega-residual-appx}
\end{align}

\hfill

\noindent The linear operator is
\begin{align}
\mathcal L_u(\textcolor{blue}{\delta})
:={}&
\bigl(2\partial_z\bar\psi+\bar c_u+\cures\bigr)\,\textcolor{blue}{\delta u}
-\bigl(\lambda r+\bar u^r\bigr)\,\partial_r \textcolor{blue}{\delta u}
-\bigl(\lambda z + \bar C + \bar u^z\bigr)\,\partial_z \textcolor{blue}{\delta u}
+\bigl(2\bar u+ r\,\partial_r\bar u\bigr)\partial_z\textcolor{blue}{\delta\psi} \nonumber\\
&\quad
-\bigl( r\,\partial_z\bar u\bigr)\partial_r\textcolor{blue}{\delta\psi}
-\bigl(2\partial_z\bar u\bigr)\textcolor{blue}{\delta\psi}
-(\partial_z\bar u)\,\,\bvar{\delta C}
+\bar u\,\,\bvar{\delta c_u}.
\label{eq:u_euler_linear_operator_definition_appx}
\end{align}

\hfill

\noindent The nonlinear operator is
\begin{align}
\mathcal L_\omega(\textcolor{blue}{\delta})
={}&\hlred{ r\,\partial_z\textcolor{blue}{\delta\psi}\,\partial_r\textcolor{blue}{\delta u}-\left(\,\bvar{\delta C}+\left(2\textcolor{blue}{\delta\psi}+r\partial_r\textcolor{blue}{\delta\psi}\right)\right)\partial_z\textcolor{blue}{\delta u}+2\,\textcolor{blue}{\delta u}\,\partial_z\textcolor{blue}{\delta\psi}+\,\bvar{\delta c_u}\,\textcolor{blue}{\delta u}.}
\label{eq:u_euler_nonlinear_remainder_definition_appx}
\end{align}

\hfill \\

\subsubsection{Linearization of the second equation}

\noindent We consider the $\omega$-equation
\begin{align}
\partial_t \omega
+ (\lambda r + u^r)\,\partial_r \omega
+ (\lambda z+C+ u^z)\,\partial_z \omega
&= 2u\,\partial_z u + c_\omega(t)\,\omega.
\label{eq:full_omega_equation_for_linearization_appx}
\end{align}

\hfill 

\noindent Substituting \eqref{eq:field_perturbation_ansatz_appx}--\eqref{eq:meridional_velocity_perturbations_appx} into \eqref{eq:full_omega_equation_for_linearization_appx} gives
\begin{align}
&\partial_t(\bar\omega+\textcolor{blue}{\delta\omega})
+\bigl(\lambda r+(\bar u^r+\textcolor{blue}{\delta u^r})\bigr)\,
\partial_r(\bar\omega+\textcolor{blue}{\delta\omega})
+\bigl(\lambda z+(\bar C+\textcolor{blue}{\delta C})+(\bar u^z+\textcolor{blue}{\delta u^z})\bigr)\,
\partial_z(\bar\omega+\textcolor{blue}{\delta\omega})
\nonumber\\
& \qquad  \quad = 
2(\bar u+\textcolor{blue}{\delta u})\,\partial_z(\bar u+\textcolor{blue}{\delta u})
+(\bar c_\omega+\cures+\textcolor{blue}{\delta c_\omega})(\bar\omega+\textcolor{blue}{\delta\omega}).
\label{eq:omega_equation_substituted_perturbations_appx}
\end{align}

\hfill

\noindent Now, $\bar u,\bar\omega,\bar C,\bar c_\omega$ satisfy the steady-state equation, so
\begin{align}
\partial_t \bar\omega
+ (\lambda r + \bar u^r)\,\partial_r \bar\omega
+ (\lambda z+\bar C+ \bar u^z)\,\partial_z \bar\omega
&= 2\bar u\,\partial_z \bar u + \bar c_\omega\,\bar\omega 
\ -\ \mathcal R_\omega^{\mathsf{raw}}.
\label{eq:omega_steady_state_equation_appx}
\end{align}

\hfill 

\noindent We subtract \eqref{eq:omega_steady_state_equation_appx} from \eqref{eq:omega_equation_substituted_perturbations_appx} and group the nonlinear perturbation products in red to obtain
\begin{align}
\partial_t \textcolor{blue}{\delta\omega}
+(\lambda r+ & \bar u^r) \,\partial_r \textcolor{blue}{\delta\omega} 
+(\lambda z+\bar C+\bar u^z)\,\partial_z \textcolor{blue}{\delta\omega}
+(\textcolor{blue}{\delta u^r})\,\partial_r\bar\omega 
+(\textcolor{blue}{\delta C}+\textcolor{blue}{\delta u^z})\,\partial_z\bar\omega
\nonumber\\
&=
2\bar u\,\partial_z\textcolor{blue}{\delta u}
+2(\partial_z\bar u)\,\textcolor{blue}{\delta u}
+(\bar c_\omega+\cures)\,\textcolor{blue}{\delta\omega}
+\bar\omega\,\textcolor{blue}{\delta c_\omega} 
\ + \ \cures\,\bar\omega\ + \ \mathcal R_\omega^{\mathsf{raw}} 
\nonumber\\
& \qquad
\hlred{
+ \ 2\, \textcolor{blue}{\delta u}\,\partial_z\textcolor{blue}{\delta u}
-\Bigl(\textcolor{blue}{\delta u^r}\Bigr)\,\partial_r \textcolor{blue}{\delta\omega}
-\Bigl(\textcolor{blue}{\delta C}+\textcolor{blue}{\delta u^z}\Bigr)\,\partial_z \textcolor{blue}{\delta\omega}
+\textcolor{blue}{\delta c_u}\,\textcolor{blue}{\delta\omega}
}.
\label{eq:omega_linearization_raw_grouped_appx}
\end{align}

\hfill \\

\noindent We collect and reorganize terms
\begin{align}
\partial_t \textcolor{blue}{\delta\omega}
&=
-(\lambda r+\bar u^r)\,\partial_r \textcolor{blue}{\delta\omega}
-(\lambda z+\bar C+\bar u^z)\,\partial_z \textcolor{blue}{\delta\omega}
+(\bar c_\omega+\cures)\,\textcolor{blue}{\delta\omega}
+2(\partial_z\bar u)\,\textcolor{blue}{\delta u}
+2\bar u\,\partial_z\textcolor{blue}{\delta u}
+\bar\omega\,\textcolor{blue}{\delta c_\omega}
\nonumber\\
&\qquad
+ r\,(\partial_z\textcolor{blue}{\delta\psi})\,\partial_r\bar\omega
-\bigl(2\textcolor{blue}{\delta\psi}+r\partial_r\textcolor{blue}{\delta\psi}\bigr)\,\partial_z\bar\omega
-(\partial_z\bar\omega)\,\textcolor{blue}{\delta C} \ + \ \cures\,\bar\omega\ + \ \mathcal R_\omega^{\mathsf{raw}}
\nonumber\\
& \qquad
\hlred{
-\Bigl(\textcolor{blue}{\delta u^r}\Bigr)\,\partial_r \textcolor{blue}{\delta\omega}
-\Bigl(\textcolor{blue}{\delta C}+\textcolor{blue}{\delta u^z}\Bigr)\,\partial_z \textcolor{blue}{\delta\omega}
+2\,\textcolor{blue}{\delta u}\,\partial_z\textcolor{blue}{\delta u}
+\textcolor{blue}{\delta c_u}\,\textcolor{blue}{\delta\omega}
}.
\label{eq:omega_linearization_collected_appx}
\end{align}

\hfill \\

\noindent Then, substituting $\textcolor{blue}{\delta c_\omega} = \textcolor{blue}{\delta c_u} $, we obtain
\begin{align}
\partial_t \textcolor{blue}{\delta\omega}
&=
-(\lambda r+\bar u^r)\,\partial_r \textcolor{blue}{\delta\omega}
-(\lambda z+\bar C+\bar u^z)\,\partial_z \textcolor{blue}{\delta\omega}
+(\bar c_\omega+\cures)\,\textcolor{blue}{\delta\omega}
+2(\partial_z\bar u)\,\textcolor{blue}{\delta u}
+2\bar u\,\partial_z\textcolor{blue}{\delta u}
\nonumber\\
& \quad  +\bar\omega\,(\textcolor{blue}{\delta c_u})
+ r\,(\partial_z\textcolor{blue}{\delta\psi})\,\partial_r\bar\omega
-\bigl(2\textcolor{blue}{\delta\psi}+r\partial_r\textcolor{blue}{\delta\psi}\bigr)\,\partial_z\bar\omega
-(\partial_z\bar\omega)\,\textcolor{blue}{\delta C} \ + \ \cures\,\bar\omega\ + \ \mathcal R_\omega^{\mathsf{raw}}
\nonumber\\
&\quad \hlred{
-\Bigl(\textcolor{blue}{\delta u^r}\Bigr)\,\partial_r \textcolor{blue}{\delta\omega}
-\Bigl(\textcolor{blue}{\delta C}+\textcolor{blue}{\delta u^z}\Bigr)\,\partial_z \textcolor{blue}{\delta\omega}
+2\,\textcolor{blue}{\delta u}\,\partial_z\textcolor{blue}{\delta u}
+\textcolor{blue}{\delta c_u}\,\textcolor{blue}{\delta\omega}
}.
\label{eq:omega_linearization_comega_substituted_appx}
\end{align}

\hfill \\

\noindent The linear operator is
\begin{align}
\mathcal L_\omega(\textcolor{blue}{\delta})
:={}&
(\bar c_u+\cures-\lambda)\,\textcolor{blue}{\delta\omega}
-\bigl(\lambda r + \bar u^r \bigr)\partial_r\textcolor{blue}{\delta\omega}
-\left(\lambda z+\bar C+ \bar u^z \right)\partial_z\textcolor{blue}{\delta\omega}
+2(\partial_z\bar u)\,\textcolor{blue}{\delta u}
+2\bar u\,\partial_z\textcolor{blue}{\delta u}
\nonumber\\
&\quad
-2(\partial_z\bar\omega)\,\textcolor{blue}{\delta\psi}
+ r\,(\partial_r\bar\omega)\,\partial_z\textcolor{blue}{\delta\psi}
- r\,(\partial_z\bar\omega)\,\partial_r\textcolor{blue}{\delta\psi}
-(\partial_z\bar\omega)\,\,\bvar{\delta C}
+\bar\omega\,\,\bvar{\delta c_u}.
\label{eq:omega_euler_linear_operator_definition_appx}
\end{align}

\hfill \\

\noindent The nonlinear operator is
\begin{align}
\mathcal N_\omega(\textcolor{blue}{\delta})
={}&\hlred{- \textcolor{blue}{\delta u^r} \partial_r\textcolor{blue}{\delta\omega}-\left(\,\bvar{\delta C}+\textcolor{blue}{\delta u^z}\right)\partial_z\textcolor{blue}{\delta\omega}+2\,\textcolor{blue}{\delta u}\,\partial_z\textcolor{blue}{\delta u}+\,\bvar{\delta c_u}\,\textcolor{blue}{\delta\omega}.}
\label{eq:omega_euler_nonlinear_remainder_definition_appx}
\end{align}

\hfill \\

\subsubsection{Linearization of the elliptic equation}

\hfill 

\noindent We consider
\begin{align}
-\Bigl[\partial_r^2 + \frac{3}{r}\partial_r + \partial_z^2\Bigr]\psi
&= \omega.
\label{eq:elliptic_equation_for_linearization_appx}
\end{align}

\noindent Then
\begin{align}
-\Bigl[\partial_r^2 + \frac{3}{r}\partial_r + \partial_z^2\Bigr]\bigl(\bar\psi+\textcolor{blue}{\delta\psi}\bigr)
&= \bar\omega+\textcolor{blue}{\delta\omega}.
\label{eq:elliptic_equation_substituted_perturbations_appx}
\end{align}

\hfill

\noindent Since $(\bar\psi,\bar\omega)$ form a steady-state, they satisfy
\begin{align}
-\Bigl[\partial_r^2 + \frac{3}{r}\partial_r + \partial_z^2\Bigr]\bar\psi
&= \bar\omega.
\label{eq:elliptic_steady_state_equation_appx}
\end{align}

\hfill 

\noindent There is no elliptic residual $\mathcal{R}_\psi$ here because the spline representation of \(\bar\omega\) is obtained by applying the elliptic operator to \(\bar\psi\) exactly, as discussed in \Cref{sec: analytic representation splines}. \\

\hfill

\noindent We subtract \eqref{eq:elliptic_steady_state_equation_appx} from \eqref{eq:elliptic_equation_substituted_perturbations_appx} and obtain only linear terms:
\begin{align}
-\Bigl[\partial_r^2 + \frac{3}{r}\partial_r + \partial_z^2\Bigr]\textcolor{blue}{\delta\psi}
&= \textcolor{blue}{\delta\omega}.
\label{eq:elliptic_linearized_relation_appx}
\end{align}

\newpage

\subsubsection{\texorpdfstring{$z$-differentiated form of the linearized first equation}{z-differentiated form of the linearized first equation}}

\noindent Taking one $z$-derivative of the full linearized first equation, we use the notation
\begin{align}
\mathcal L_z(\textcolor{blue}{\delta})&:=\partial_z\mathcal L_u(\textcolor{blue}{\delta}),
\label{eq:z_differentiated_euler_linear_definition_appx}\\
\mathcal N_z(\textcolor{blue}{\delta})&:=\partial_z\mathcal L_\omega(\textcolor{blue}{\delta}).
\label{eq:z_differentiated_euler_nonlinear_definition_appx}
\end{align}

\hfill

\noindent The main product-rule contributions are
\begin{align}
\partial_z\Bigl[-(\lambda r+\bar u^r)\partial_r\textcolor{blue}{\delta u}\Bigr]
&=-(\lambda r+\bar u^r)\partial_{rz}\textcolor{blue}{\delta u}-(\partial_z\bar u^r)\partial_r\textcolor{blue}{\delta u},
\label{eq:z_product_rule_radial_transport_appx}\\
\partial_z\Bigl[-(\lambda z+\bar C+\bar u^z)\partial_z\textcolor{blue}{\delta u}\Bigr]
&=-(\lambda z+\bar C+\bar u^z)\partial_{zz}\textcolor{blue}{\delta u}-(\lambda+\partial_z\bar u^z)\partial_z\textcolor{blue}{\delta u},
\label{eq:z_product_rule_vertical_transport_appx}\\
\partial_z\Bigl[(2\partial_z\bar\psi+\bar c_u)\textcolor{blue}{\delta u}\Bigr]
&=(2\partial_z\bar\psi+\bar c_u)\partial_z\textcolor{blue}{\delta u}+2\partial_{zz}\bar\psi\,\textcolor{blue}{\delta u}.
\label{eq:z_product_rule_source_appx}
\end{align}
These identities, together with term-by-term differentiation of the \(\textcolor{blue}{\delta\psi}\), \(\textcolor{blue}{\delta C}\), and \(\textcolor{blue}{\delta c_u}\) terms, give the following explicit formulas.

\hfill \\ 

\noindent The linear part is
\begin{align}
\mathcal L_z(\textcolor{blue}{\delta})
={}&
-(\lambda r+\bar u^r)\,\partial_{rz}\textcolor{blue}{\delta u}
-(\lambda z+\bar C+\bar u^z)\,\partial_{zz}\textcolor{blue}{\delta u}
\nonumber\\
&+\bigl(2\partial_z\bar\psi+\bar c_u+\cures-\lambda-\partial_z\bar u^z\bigr)
\partial_z\textcolor{blue}{\delta u}
-(\partial_z\bar u^r)\,\partial_r\textcolor{blue}{\delta u}
+2\partial_{zz}\bar\psi\,\textcolor{blue}{\delta u}
\nonumber\\
&+(2\bar u+ r\partial_r\bar u)\,\partial_{zz}\textcolor{blue}{\delta\psi}
+\bigl(r\partial_{rz}\bar u\bigr)
\partial_z\textcolor{blue}{\delta\psi}
\nonumber\\
&- r(\partial_z\bar u)\,\partial_{rz}\textcolor{blue}{\delta\psi}
- r(\partial_{zz}\bar u)\,\partial_r\textcolor{blue}{\delta\psi}
-2\partial_{zz}\bar u\,\textcolor{blue}{\delta\psi}
\nonumber\\
&-(\partial_{zz}\bar u)\,\,\bvar{\delta C}
+(\partial_z\bar u)\,\,\bvar{\delta c_u}.
\label{eq:z_differentiated_u_euler_linear_operator_appx}
\end{align}

\hfill

\noindent The nonlinear operator is
\begin{align}
\mathcal N_z(\textcolor{blue}{\delta})
={}&
 r(\partial_{zz}\textcolor{blue}{\delta\psi})\,\partial_r\textcolor{blue}{\delta u}
+ r(\partial_z\textcolor{blue}{\delta\psi})\,\partial_{rz}\textcolor{blue}{\delta u}
\nonumber\\
&-\left(2\partial_z\textcolor{blue}{\delta\psi}
+r\partial_{rz}\textcolor{blue}{\delta\psi}\right)\partial_z\textcolor{blue}{\delta u}
-\left(\,\bvar{\delta C}+\left(2\textcolor{blue}{\delta\psi}+r\partial_r\textcolor{blue}{\delta\psi}\right)\right)
\partial_{zz}\textcolor{blue}{\delta u}
\nonumber\\
&+2(\partial_z\textcolor{blue}{\delta u})(\partial_z\textcolor{blue}{\delta\psi})
+2\textcolor{blue}{\delta u}\,\partial_{zz}\textcolor{blue}{\delta\psi}
+\,\bvar{\delta c_u}\,\partial_z\textcolor{blue}{\delta u}.
\label{eq:z_differentiated_u_euler_nonlinear_remainder_appx}
\end{align}

\hfill

\clearpage

\subsubsection{\texorpdfstring{$r$-differentiated form of the linearized first equation}{r-differentiated form of the linearized first equation}}

\noindent We now differentiate the full linearized first equation with respect to $r$. We write
\begin{align}
\mathcal L_r(\textcolor{blue}{\delta})&:=\partial_r\mathcal L_u(\textcolor{blue}{\delta}),
\label{eq:r_differentiated_euler_linear_definition_appx}\\
\mathcal N_r(\textcolor{blue}{\delta})&:=\partial_r\mathcal L_\omega(\textcolor{blue}{\delta}),
\label{eq:r_differentiated_euler_nonlinear_definition_appx}.
\end{align}

\hfill

\noindent The main product-rule contributions are
\begin{align}
\partial_r\Bigl[-(\lambda r+\bar u^r)\partial_r\textcolor{blue}{\delta u}\Bigr]
&=-(\lambda r+\bar u^r)\partial_{rr}\textcolor{blue}{\delta u}-(\lambda+\partial_r\bar u^r)\partial_r\textcolor{blue}{\delta u},
\label{eq:r_product_rule_radial_transport_appx}\\
\partial_r\Bigl[-(\lambda z+\bar C+\bar u^z)\partial_z\textcolor{blue}{\delta u}\Bigr]
&=-(\lambda z+\bar C+\bar u^z)\partial_{zr}\textcolor{blue}{\delta u}-(\partial_r\bar u^z)\partial_z\textcolor{blue}{\delta u},
\label{eq:r_product_rule_vertical_transport_appx}\\
\partial_r\Bigl[(2\partial_z\bar\psi+\bar c_u)\textcolor{blue}{\delta u}\Bigr]
&=(2\partial_z\bar\psi+\bar c_u)\partial_r\textcolor{blue}{\delta u}+2\partial_{rz}\bar\psi\,\textcolor{blue}{\delta u}.
\label{eq:r_product_rule_source_appx}
\end{align}
These identities, together with term-by-term differentiation of the \(\textcolor{blue}{\delta\psi}\), \(\textcolor{blue}{\delta C}\), and \(\textcolor{blue}{\delta c_u}\) terms, give the following explicit formulas.

\hfill \\

\noindent The linear part is
\begin{align}
\mathcal L_r(\textcolor{blue}{\delta})
={}&
-(\lambda r+\bar u^r)\,\partial_{rr}\textcolor{blue}{\delta u}
-(\lambda z+\bar C+\bar u^z)\,\partial_{zr}\textcolor{blue}{\delta u}
\nonumber\\
&+\bigl(2\partial_z\bar\psi+\bar c_u+\cures-\lambda-\partial_r\bar u^r\bigr)
\partial_r\textcolor{blue}{\delta u}
-(\partial_r\bar u^z)\,\partial_z\textcolor{blue}{\delta u}
+2\partial_{rz}\bar\psi\,\textcolor{blue}{\delta u}
\nonumber\\
&+(2\bar u+ r\partial_r\bar u)\,\partial_{zr}\textcolor{blue}{\delta\psi}
+\bigl(3\partial_r\bar u+ r\partial_{rr}\bar u\bigr)
\partial_z\textcolor{blue}{\delta\psi}
\nonumber\\
&- r(\partial_z\bar u)\,\partial_{rr}\textcolor{blue}{\delta\psi}
-\bigl(3\partial_z\bar u+ r\partial_{rz}\bar u\bigr)
\partial_r\textcolor{blue}{\delta\psi}
-2\partial_{rz}\bar u\,\textcolor{blue}{\delta\psi}
\nonumber\\
&-(\partial_{rz}\bar u)\,\,\bvar{\delta C}
+(\partial_r\bar u)\,\,\bvar{\delta c_u}.
\label{eq:r_differentiated_u_euler_linear_operator_appx}
\end{align}

\hfill

\noindent The nonlinear part is
\begin{align}
\mathcal N_r(\textcolor{blue}{\delta})
={}&
\left(\partial_z\textcolor{blue}{\delta\psi}+r\partial_{rz}\textcolor{blue}{\delta\psi}\right)
\partial_r\textcolor{blue}{\delta u}
+ r\partial_z\textcolor{blue}{\delta\psi}\,\partial_{rr}\textcolor{blue}{\delta u}
\nonumber\\
&-\left(3\partial_r\textcolor{blue}{\delta\psi}+r\partial_{rr}\textcolor{blue}{\delta\psi}\right)
\partial_z\textcolor{blue}{\delta u}
-\left(\,\bvar{\delta C}+\left(2\textcolor{blue}{\delta\psi}+r\partial_r\textcolor{blue}{\delta\psi}\right)\right)
\partial_{zr}\textcolor{blue}{\delta u}
\nonumber\\
&+2(\partial_r\textcolor{blue}{\delta u})(\partial_z\textcolor{blue}{\delta\psi})
+2\textcolor{blue}{\delta u}\,\partial_{zr}\textcolor{blue}{\delta\psi}
+\,\bvar{\delta c_u}\,\partial_r\textcolor{blue}{\delta u}.
\label{eq:r_differentiated_u_euler_nonlinear_remainder_appx}
\end{align}

\clearpage
\subsection{Modulation equations}

The modulation conditions determine \(\textcolor{blue}{\delta C}\) and \(\textcolor{blue}{\delta c_u}\). 

\hfill

\subsubsection{Modulation equation for \texorpdfstring{\(\bvar{\delta C}\)}{delta C}}

\noindent The first normalization condition is
\begin{equation}
C+u^z_0=0.
\label{eq:euler_stagnation_condition_appx}
\end{equation}

\hfill 

\noindent Since $u^z_0=2\psi_0$, substituting $C=\bar C+\bvar{\delta C}$ and $\psi=\bar\psi+\bvar{\delta\psi}$ gives
\begin{equation}
0=(\bar C+2\bar\psi_0)+(\bvar{\delta C}+2\bvar{\delta\psi}_0).
\end{equation}

\hfill 

\noindent The reference stagnation condition $\bar C+2\bar\psi_0=0$ therefore yields $\bvar{\delta C}+2\bvar{\delta\psi}_0=0,$ and thus
\begin{equation}
\bvar{\delta C}=-2\bvar{\delta\psi}_0.
\label{eq:euler_deltaC_explicit_formula_appx}
\end{equation}

\hfill \\

\subsubsection{Modulation equation for \texorpdfstring{\(\bvar{\delta c_u}\)}{delta c u}}

\noindent The second normalization condition fixes $u_0(t)=\bar u_0$, hence $\partial_tu_0=0$. \\

\noindent Evaluating the $u$-equation at the origin gives
\begin{equation}
\begin{aligned}
0
={}&-(\lambda r+u^r)_0(\partial_r u)_0
-(\lambda z+C+u^z)_0(\partial_z u)_0
+2u_0(\partial_z\psi)_0+c_u u_0.
\end{aligned}
\label{eq:euler_u_origin_exact_relation_appx_current}
\end{equation}

\hfill

\noindent Since $u^r_0=0$ and $C+u^z_0=0$, both transport terms vanish. \\

\noindent Using $u_0=\bar u_0$, the origin equation becomes
\begin{equation}
0=2\bar u_0(\partial_z\psi)_0+c_u\bar u_0.
\end{equation}
\noindent We introduce the re-centered residual
\begin{equation}
\mathcal R_u:=\mathcal R_u^{\mathsf{raw}}+\cures\bar u,
\qquad
\cures=-\frac{\mathcal R_u^{\mathsf{raw}}(0)}{\bar u_0}.
\end{equation}
\noindent By construction,
\begin{equation}
\mathcal R_u(0)
=\mathcal R_u^{\mathsf{raw}}(0)+\cures\bar u_0=0.
\label{eq:euler-recentered-residual-origin-zero}
\end{equation}

\hfill 

\noindent Substituting $\psi=\bar\psi+\bvar{\delta\psi}$ and $c_u=\bar c_u+\cures+\bvar{\delta c_u}$ separates the zero reference contribution from the perturbative terms:
\begin{equation}
0=
\underbrace{\bigl(2(\partial_z\bar\psi)_0+\bar c_u+\cures\bigr)\bar u_0}_{\mathcal R_u(0)=0}
+2\bar u_0(\partial_z\bvar{\delta\psi})_0
+\bar u_0\bvar{\delta c_u}.
\end{equation}

\hfill 

\noindent The underbraced term vanishes by \eqref{eq:euler-recentered-residual-origin-zero}, leaving
\begin{equation}
\bar u_0\bvar{\delta c_u}
=-2\bar u_0(\partial_z\bvar{\delta\psi})_0.
\label{eq:euler_second_modulation_equation_appx}
\end{equation}

\hfill 

\noindent Since $\bar u_0\ne0$,
\begin{equation}
\bvar{\delta c_u}=-2(\partial_z\bvar{\delta\psi})_0.
\label{eq:euler_deltacu_explicit_formula_appx}
\end{equation}

\clearpage

\section{Additional Details for the Linear Damping Estimate}
\label{app:linear-damping}

\subsection{Additional details for the low-order linear estimate}
\label{app:low-order-damping}

\noindent Here we record the complete low-order operator decomposition and identifies which terms are retained in the signed matrix and which are kept outside it in Section~\ref{sec:partial-linear-damping}. The separate mechanisms used to estimate the terms outside the matrix are developed in a separate stability paper~\citep{EulerBlowupStability}. \\

\subsubsection{Complete low-order operator decomposition}
\label{app:low-order-linear-setup}

\noindent The low-order linear pairing is
\begin{equation}
\label{eq:low-linear-pairing-def}
\begin{aligned}
\mathcal Q_{\mathcal L}^{\rm low}(\bvar{\delta})
&:=\langle \mathcal L_\omega(\bvar{\delta}),\bvar{\delta v_\omega}\rangle_{\wgray{\Phi_\omega}}
+\langle \mathcal L_r(\bvar{\delta}),\bvar{\delta v_r}\rangle_{\wgray{\Phi_r}}
+\langle \mathcal L_z(\bvar{\delta}),\bvar{\delta v_z}\rangle_{\wgray{\Phi_z}}.
\end{aligned}
\end{equation}
\noindent Recall that \(\mathcal Q_{\mathcal L}^{\rm low}\) also contains the linear terms obtained when the residual amplitude offset from the modulation equations multiplies a perturbation. Since \(\bvar{\delta C}\) and \(\bvar{\delta c_u}\) are exactly linear and therefore enter the linear operator only.\\

\noindent For each \(i\in\{\omega,r,z\}\), separate the transport part and non-transport terms as
\begin{equation}
\label{eq:low-linear-operator-split}
\mathcal L_i(\bvar{\delta})
=-G_r\partial_r\bvar{\delta v_i}-G_z\partial_z\bvar{\delta v_i}
+\sum_m B_{i,m}(\bvar{\delta}).
\end{equation}

\hfill 

\noindent  The complete list of non-transport terms is
\begin{align}
\label{eq:low-linear-B-omega-list}
B_{\omega,1}&=(\bar c_u-\lambda+\cures)\bvar{\delta v_\omega},
&B_{\omega,2}&=2\partial_z\bar u\,\bvar{\delta u},\nonumber\\
B_{\omega,3}&=2\bar u\,\bvar{\delta v_z},
&B_{\omega,4}&=-2\partial_z\bar\omega\,\bvar{\delta\psi},\nonumber\\
B_{\omega,5}&=r\partial_r\bar\omega\,\partial_z\bvar{\delta\psi},
&B_{\omega,6}&=-r\partial_z\bar\omega\,\partial_r\bvar{\delta\psi},\nonumber\\
B_{\omega,7}&=-\partial_z\bar\omega\,\bvar{\delta C},
&B_{\omega,8}&=\bar\omega\,\bvar{\delta c_u},
\end{align}
\begin{align}
\label{eq:low-linear-B-r-list}
B_{r,1}&=(2\partial_z\bar\psi+\bar c_u-\lambda-\partial_r\bar u^r+\cures)\bvar{\delta v_r},
&B_{r,2}&=-\partial_r\bar u^z\,\bvar{\delta v_z},\nonumber\\
B_{r,3}&=2\partial_{rz}\bar\psi\,\bvar{\delta u},
&B_{r,4}&=(2\bar u+r\partial_r\bar u)\partial_{rz}\bvar{\delta\psi},\nonumber\\
B_{r,5}&=(3\partial_r\bar u+r\partial_{rr}\bar u)\partial_z\bvar{\delta\psi},
&B_{r,6}&=-r\partial_z\bar u\,\partial_{rr}\bvar{\delta\psi},\nonumber\\
B_{r,7}&=-(3\partial_z\bar u+r\partial_{rz}\bar u)\partial_r\bvar{\delta\psi},
&B_{r,8}&=-2\partial_{rz}\bar u\,\bvar{\delta\psi},\nonumber\\
B_{r,9}&=-\partial_{rz}\bar u\,\bvar{\delta C},
&B_{r,10}&=\partial_r\bar u\,\bvar{\delta c_u},
\end{align}
\begin{align}
\label{eq:low-linear-B-z-list}
B_{z,1}&=-\partial_z\bar u^r\,\bvar{\delta v_r},
&B_{z,2}&=(2\partial_z\bar\psi+\bar c_u-\lambda-\partial_z\bar u^z+\cures)\bvar{\delta v_z},\nonumber\\
B_{z,3}&=2\partial_{zz}\bar\psi\,\bvar{\delta u},
&B_{z,4}&=(2\bar u+r\partial_r\bar u)\partial_{zz}\bvar{\delta\psi},\nonumber\\
B_{z,5}&=(r\partial_{rz}\bar u)\partial_z\bvar{\delta\psi},
&B_{z,6}&=-r\partial_z\bar u\,\partial_{rz}\bvar{\delta\psi},\nonumber\\
B_{z,7}&=-r\partial_{zz}\bar u\,\partial_r\bvar{\delta\psi},
&B_{z,8}&=-2\partial_{zz}\bar u\,\bvar{\delta\psi},\nonumber\\
B_{z,9}&=-\partial_{zz}\bar u\,\bvar{\delta C},
&B_{z,10}&=\partial_z\bar u\,\bvar{\delta c_u}.
\end{align}

\hfill 

\noindent Among the non-transport terms above, the low-order signed matrix retains
\begin{equation}
\mathcal S_{\mathcal L}
=
\{B_{\omega,1},B_{r,1},B_{r,2},B_{z,1},B_{z,2}\}.
\label{eq:low-linear-selected-set-fixed}
\end{equation}
The three diagonal labels contribute to the diagonal entries, while \(B_{r,2}\) and \(B_{z,1}\) contribute to the symmetric \(r\)--\(z\) entry. All other labels in the lists above are kept outside the signed matrix. \\

\noindent The terms outside the signed matrix are precisely the terms in the complete lists above that are not contained in \(\mathcal S_{\mathcal L}\). We denote by \(A_{\mathcal L}^{0}\) and \(A_{\mathcal L}^{k}\) the aggregate low-order and top-order losses obtained after these terms are estimated by the separate mechanisms developed in the stability argument and presented in a separate stability paper~\citep{EulerBlowupStability}. \\

\subsubsection{Alternative weight continuations outside the computational box}
\label{app:tail-weight-models}

\noindent Section~\ref{sssec:low-order-linear-tail-closure} uses the tail continuation employed in the numerical construction. For completeness, other radial tail models can be analyzed from the same scalar condition. Suppose
\begin{equation}
\label{eq:low-linear-tail-radial-weight}
\wgray{\Phi_i}(r,z)=\phi_i(\rho),
\qquad \rho=r^2+z^2.
\end{equation}
\noindent Then
\begin{equation}
\label{eq:low-linear-tail-radial-chain-rule}
\partial_r\log\wgray{\Phi_i}
=2r\frac{d}{d\rho}\log\phi_i,
\qquad
\partial_z\log\wgray{\Phi_i}
=2z\frac{d}{d\rho}\log\phi_i,
\end{equation}
and the left side of \eqref{eq:low-linear-tail-scalar-condition} becomes
\begin{equation}
\label{eq:low-linear-tail-radial-expression}
\left(1+\frac{2z}{\rho}\bar C\right)
\left(2\rho\frac{d}{d\rho}\log\phi_i(\rho)\right).
\end{equation}

\hfill 

\noindent Since \(\sqrt\rho\ge R_{\rm box}\) and \(|z|\le\sqrt\rho\) on the tail,
\begin{equation}
\label{eq:low-linear-tail-finite-box-factor}
1-\frac{2|\bar C|}{R_{\rm box}}
\le
1+\frac{2z}{\rho}\bar C
\le
1+\frac{2|\bar C|}{R_{\rm box}}.
\end{equation}

\hfill 

\noindent The convenient sufficient assumptions
\begin{equation}
\label{eq:low-linear-tail-basic-phase-condition}
\cures<1,
\qquad
R_{\rm box}>2|\bar C|
\end{equation}
\noindent leave a positive zeroth-order margin and keep the correction factor positive throughout the tail. \\

\noindent The main asymptotic guide is
\begin{equation}
\limsup_{\rho\to\infty}
2\rho\frac{d}{d\rho}\log\phi_i(\rho)
<4(1-\cures).
\label{eq:low-linear-tail-asymptotic-rate-condition}
\end{equation}

\hfill 

\noindent We consider different continuation families below:\\

\begin{itemize}[leftmargin=1.3em,itemsep=2.5em,topsep=0.1em]
\item \textbf{Exponential decay.} If
\begin{equation}
\phi_i(\rho)=C_{\rm floor}+a\exp(-b\rho^p),
\qquad a,b,p>0,
\end{equation}
then \(d(\log\phi_i)/d\rho\le0\). Under \eqref{eq:low-linear-tail-basic-phase-condition}, the left side of \eqref{eq:low-linear-tail-scalar-condition} is nonpositive. Any \(0<\gamma_{\rm tail}<1-\cures\) therefore works.

\item \textbf{Algebraic decay.} If
\begin{equation}
\phi_i(\rho)=C_{\rm floor}+a\rho^{-p},
\qquad a,p>0,
\end{equation}
In this case,
\[
2\rho\frac{d}{d\rho}\log\phi_i(\rho)
=-2p\frac{a\rho^{-p}}{C_{\rm floor}+a\rho^{-p}}\le0.
\]
Under \eqref{eq:low-linear-tail-basic-phase-condition}, the factor \(1+2z\bar C/\rho\) in \eqref{eq:low-linear-tail-radial-expression} is positive throughout the tail. Hence the entire expression in \eqref{eq:low-linear-tail-radial-expression} is nonpositive, so the scalar condition \eqref{eq:low-linear-tail-scalar-condition} holds for any \(0<\gamma_{\rm tail}<1-\cures\).

\item \textbf{Algebraic growth.} If
\begin{equation}
\phi_i(\rho)=C_{\rm floor}+a\rho^p,
\qquad a,p>0,
\end{equation}
then
\begin{equation}
\label{eq:low-linear-tail-algebraic-growth-rate}
2\rho\frac{d}{d\rho}\log\phi_i(\rho)
=2p\frac{a\rho^p}{C_{\rm floor}+a\rho^p}
\le2p.
\end{equation}
\noindent A sufficient finite-box condition is
\begin{equation}
\label{eq:low-linear-tail-algebraic-growth-threshold}
p\left(1+\frac{2|\bar C|}{R_{\rm box}}\right)
<2(1-\cures).
\end{equation}

\item \textbf{Positive exponential growth.} If
\begin{equation}
\phi_i(\rho)=C_{\rm floor}+a\exp(b\rho^p),
\qquad a,b,p>0,
\end{equation}
then
\begin{equation}
2\rho\frac{d}{d\rho}\log\phi_i(\rho)
=2bp\rho^p\frac{a\exp(b\rho^p)}{C_{\rm floor}+a\exp(b\rho^p)},
\end{equation}
\noindent which is unbounded as \(\rho\to\infty\). Hence the scalar tail certificate cannot hold.
\end{itemize}

\hfill \\

\hfill

\noindent The four continuation families above are summarized in the following table, separating the asymptotic rate restriction from the finite-tail condition used in the certificate.

\begin{center}
\renewcommand{\arraystretch}{1.17}
\setlength{\tabcolsep}{4pt}
\textbf{Summary of radial model rates.}\vspace{0.5em}

\newcolumntype{Y}{>{\raggedright\arraybackslash}X}
\begin{tabularx}{\textwidth}{@{}lXYY@{}}
\toprule
\textbf{Rate}
& \textbf{Model \(\phi_i(\rho)\)}
& \textbf{Asymptotic condition}
& \textbf{Finite-tail condition} \\
\midrule
Exponential decay
& \(C_{\rm floor}+a\exp(-b\rho^p)\)
& \(a,b,p>0\), with \(\cures<1\)
& \eqref{eq:low-linear-tail-basic-phase-condition} \\
\addlinespace[0.15em]
Algebraic decay
& \(C_{\rm floor}+a\rho^{-p}\)
& \(a,p>0\), with \(\cures<1\)
& \eqref{eq:low-linear-tail-basic-phase-condition} \\
\addlinespace[0.15em]
Algebraic growth
& \(C_{\rm floor}+a\rho^p\)
& \(p<2(1-\cures)\)
& \eqref{eq:low-linear-tail-algebraic-growth-threshold} \\
\addlinespace[0.15em]
Exponential growth
& \(C_{\rm floor}+a\exp(b\rho^p)\)
& Not possible
& Not possible \\
\bottomrule
\end{tabularx} \vspace{12mm}
\end{center}

\noindent For the \(r,z\) weights, compatibility with the low-order weighted energy also requires \(p\ge1\). Thus an algebraically growing tail is compatible when one can choose
\begin{equation}
1\le p<
\frac{2(1-\cures)}{1+2|\bar C|/R_{\rm box}}.
\label{eq:low-linear-tail-compatible-growth-range}
\end{equation}
\noindent This interval is nonempty if and only if
\begin{equation}
1+\frac{2|\bar C|}{R_{\rm box}}
<2(1-\cures).
\label{eq:low-linear-tail-compatible-growth-nonempty}
\end{equation}

\clearpage

\subsubsection{Local concentration estimates}

\noindent To apply the localized certificate, it is enough to prove \eqref{eq:localized-matrix-loc-concentration} on each piece \(B_\ell\). Since every \(B_\ell\) is bounded and stays away from the origin, and the weights are singular only at the origin,
\(
\sup_{B_\ell}\wgray{\Phi_i}<\infty
\)
for every \(i\in\{\omega,r,z\}\). The weighted energy can thus be controlled by ordinary local estimates for the components \(\bvar{\delta v_i}\). \\ 

\paragraph{Pointwise estimate.}
\noindent One option is to use pointwise estimates. If
\begin{equation}
\|\bvar{\delta v_i}\|_{L^\infty(B_\ell)}
\le
P_{i,\ell}^{0}\mathfrak E_0+P_{i,\ell}^{k}\mathfrak H_k,
\qquad i\in\{\omega,r,z\},
\label{eq:localized-matrix-interior-pointwise-input}
\end{equation}
then \(|\bvar{d\xi}|^2=\sum_i\wgray{\Phi_i}|\bvar{\delta v_i}|^2\) and the identity \((a+b)^2\le2a^2+2b^2\) give \eqref{eq:localized-matrix-loc-concentration}, i.e.,
\begin{equation}
\int_{B_\ell}|\bvar{d\xi}|^2\,dr\,dz
\le
m_{\mathcal L}^{0,{\rm loc},\ell}\mathfrak E_0^2
+
m_{\mathcal L}^{k,{\rm loc},\ell}\mathfrak H_k^2,
\qquad \ell\in\mathcal I_{\mathcal L}^{\rm loc},
\end{equation}
with
\begin{align}
m_{\mathcal L}^{0,{\rm loc},\ell}
&:=
2|B_\ell|\sum_{i\in\{\omega,r,z\}}
\Bigl(\sup_{B_\ell}\wgray{\Phi_i}\Bigr)(P_{i,\ell}^{0})^2,
\label{eq:localized-matrix-interior-m0}\\
m_{\mathcal L}^{k,{\rm loc},\ell}
&:=
2|B_\ell|\sum_{i\in\{\omega,r,z\}}
\Bigl(\sup_{B_\ell}\wgray{\Phi_i}\Bigr)(P_{i,\ell}^{k})^2.
\label{eq:localized-matrix-interior-mk}
\end{align}
Thus shrinking \(B_\ell\) improves the constants, provided that the weight suprema and the pointwise constants remain uniformly controlled.\\

\paragraph{Local \(L^2\) estimate.}
\noindent Alternatively, suppose one has a local \(L^2\) estimate of the form
\begin{equation}
\int_{B_\ell}|\bvar{\delta v_i}|^2\,dr\,dz
\le
Q_{i,\ell}^{0}\mathfrak E_0^2
+
Q_{i,\ell}^{k}\mathfrak H_k^2,
\qquad i\in\{\omega,r,z\},
\label{eq:localized-matrix-local-L2-input}
\end{equation}
then \eqref{eq:localized-matrix-loc-concentration} holds with
\begin{equation}
m_{\mathcal L}^{0,{\rm loc},\ell}
:=
\sum_{i\in\{\omega,r,z\}}
\Bigl(\sup_{B_\ell}\wgray{\Phi_i}\Bigr)Q_{i,\ell}^{0},
\qquad
m_{\mathcal L}^{k,{\rm loc},\ell}
:=
\sum_{i\in\{\omega,r,z\}}
\Bigl(\sup_{B_\ell}\wgray{\Phi_i}\Bigr)Q_{i,\ell}^{k}.
\label{eq:localized-matrix-local-L2-constants}
\end{equation}

\hfill \\

\subsection{Additional details for the high-order linear estimate}
\label{app:high-order-linear}

\paragraph{Remainder structure.}
\noindent Section~\ref{sssec:complete-high-order-signed-matrix} defines the complete top-order signed matrix, including the transport contribution, every first transport commutator, and the local top-order component couplings. The high-order remainder begins exactly where that complete signed matrix ends.\\

\noindent To see the transport split directly, write \(\alpha=(q,k-q)\). In the Leibniz expansion of
\begin{equation}
D^\alpha\left(G_r\partial_r f+G_z\partial_z f\right),
\end{equation}
\noindent the terms in which exactly one derivative lands on \(G_r\) or \(G_z\) are the first transport commutators retained in the matrix. Equivalently, these are the \(|\beta|=1\) terms when \(\beta\) denotes the multi-index of derivatives falling on the transport coefficient. Each leaves exactly \(k\) derivatives on \(f\). Terms with \(|\beta|\ge2\) contain at most \(k-1\) derivatives of the perturbation and are therefore estimated as remainders.\\

\noindent The same principle applies to the local top-order coefficients. The undifferentiated coefficient multiplying a top derivative is retained in the signed matrix. Once a derivative lands on a spatially varying coefficient, the perturbation factor has order at most \(k-1\) and is treated outside the matrix. The nonlocal streamfunction terms and the remaining  high-order modulation terms are also kept outside the pointwise signed matrix and are controlled by the weighted elliptic, interpolation, and lower-order estimates.\\

\paragraph{Aggregate Euler estimate.}
\noindent All  high-order linear terms not represented by the complete signed matrix are required to satisfy
\begin{equation}
\left|
\sum_{i\in\{\omega,r,z\}}\sum_{|\alpha|=k}
\left\langle D^\alpha\mathcal L_i,D^\alpha\bvar{\delta v_i}\right\rangle_{\Phi_{i,k}}
-
\int_\mathbb{D} \mathbf V_k^{\top}\mathbb M_{\mathcal L}^{(k),\rm full}\mathbf V_k\,dr\,dz
\right|
\le
A_{\rm lin}^{k,\rm rem}\mathfrak H_k^2
+A_{\rm lin}^{0,\rm rem}\mathfrak E_0^2.
\label{eq:high-linear-full-remainder-bound-add}
\end{equation}

\hfill \\

\noindent Provided
\begin{equation}
A_{\rm lin}^{k,\rm rem}<\Lambda_{D^\alpha L}^{\rm full},
\label{eq:high-linear-full-euler-absorption-add}
\end{equation}
define the net high-order constants by
\begin{equation}
\DLamDLhigh
:=\Lambda_{D^\alpha L}^{\rm full}-A_{\rm lin}^{k,\rm rem},
\qquad
\DLamDLlow:=A_{\rm lin}^{0,\rm rem}.
\label{eq:high-linear-remainder-net-constants-reorg}
\end{equation}

\hfill 

\noindent Combining the remainder bound with the matrix certificate \eqref{eq:high-linear-full-matrix-certificate-add} gives the complete Euler high-order estimate
\begin{equation}
\sum_{i\in\{\omega,r,z\}}\sum_{|\alpha|=k}
\left\langle D^\alpha\mathcal L_i,D^\alpha\bvar{\delta v_i}\right\rangle_{\Phi_{i,k}}
\le
\DLamDLlow\,\mathfrak E_0^2-\DLamDLhigh\,\mathfrak H_k^2.
\label{eq:high-order-linear-final-full-add}
\end{equation}

\hfill \\

\noindent Thus the complete high-order matrix margin is reduced only by the top-order part of the remainder bound, while its low-order part becomes the coefficient \(\DLamDLlow\). This is the  high-order linear estimate used together with the low-order result when the final stability margin is assembled. 

\hfill

\hfill 

\paragraph{Possible refinements.}
\noindent The coupled low-order and high-order mechanism above is the first target.  If its certified margin is positive but smaller than desired, the natural next step is to improve the constants already present rather than change the proof strategy.  In particular, one can refine the steady profile, re-optimize the weights, shrink the localized sets \(B_\ell\), sharpen the concentration and interpolation bounds, retain additional favorable linear terms with their sign, and optimize \(k\) using the complete  high-order matrix together with its remainder.  These are quantitative improvements of the same damping mechanism.

\hfill \\
\noindent Only if the complete top-order signed matrix remained too weak on a substantial region after these refinements would a genuinely stronger top-order estimate be needed.  The low-order localization argument cannot simply be repeated for concentration of the top derivatives.  In that case one would refine the top-order energy or exploit additional forms of damping.

\clearpage

\section{Proofs of the Stability Theorems}\label{appx: nonlinear stability proof}

\hfill 

\noindent We now prove the single-radius stability Theorem~\ref{thm:main-stability}.\\

\begin{theorem}[Euler single-radius stability]\label{thm:main-stability_appx}
Assume that the following perturbation equations hold for admissible perturbations:
\begin{align}
\partial_t\bvar{\delta\omega}
&=\mathcal L_\omega(\bvar{\delta})+
\mathcal N_\omega(\bvar{\delta})+\mathcal R_\omega, \label{eq: Stability theorem perturbation 1_appx}\\
\partial_t\partial_r\bvar{\delta u}
&=\partial_r\mathcal L_u(\bvar{\delta})+
\mathcal N_{r}(\bvar{\delta})+\mathcal R_r,\label{eq: Stability theorem perturbation 2_appx}\\
\partial_t\partial_z\bvar{\delta u}
&=\partial_z\mathcal L_u(\bvar{\delta})+
\mathcal N_{z}(\bvar{\delta})+\mathcal R_z. \label{eq: Stability theorem perturbation 3_appx}
\end{align}
Assume that the following estimates have been certified with finite constants:

\begin{equation}
\begin{aligned}
&\langle \mathcal L_\omega,\bvar{\delta\omega}\rangle_{\wgray{\Phi_\omega}}
+\langle \partial_r\mathcal L_u,\partial_r\bvar{\delta u}\rangle_{\wgray{\Phi_r}}
+\langle \partial_z\mathcal L_u,\partial_z\bvar{\delta u}\rangle_{\wgray{\Phi_z}}
\le -\DLamLlow \, \mathfrak E_0^2+\DLamLhigh \, \mathfrak H_k^2.
\end{aligned}
\label{eq:ns-thm-low-linear_appx}
\end{equation}

\vspace{2mm}

\begin{equation}
\begin{aligned}
&\sum_{i\in\{\omega,r,z\}}\sum_{|\alpha|=k}
\langle D^\alpha\mathcal L_i,D^\alpha\bvar{\delta v_i}\rangle_{\wgray{\Phi_{i,k}}}
\le \DLamDLlow \, \mathfrak E_0^2-\DLamDLhigh \, \mathfrak H_k^2.
\end{aligned}
\label{eq:ns-thm-high-linear_appx}
\end{equation}

\vspace{2mm}

\begin{equation}
\begin{aligned}
&\left|\langle\mathcal N_\omega,\bvar{\delta\omega}\rangle_{\wgray{\Phi_\omega}}
+\langle\mathcal N_{r},\partial_r\bvar{\delta u}\rangle_{\wgray{\Phi_r}}
+\langle\mathcal N_{z},\partial_z\bvar{\delta u}\rangle_{\wgray{\Phi_z}}
\right|
\le \DLamNthree \, \mathfrak E_k^3.
\end{aligned}
\label{eq:ns-thm-low-nonlinear_appx}
\end{equation}

\vspace{2mm}

\begin{equation}
\begin{aligned}
&\mu_k\left|
\sum_{|\alpha|=k}\langle D^\alpha\mathcal N_\omega,D^\alpha\bvar{\delta\omega}\rangle_{\wgray{\Phi_{\omega,k}}}
+\sum_{|\alpha|=k}\langle D^\alpha\mathcal N_r,D^\alpha\partial_r\bvar{\delta u}\rangle_{\wgray{\Phi_{r,k}}}
+\sum_{|\alpha|=k}\langle D^\alpha\mathcal N_z,D^\alpha\partial_z\bvar{\delta u}\rangle_{\wgray{\Phi_{z,k}}}
\right|
\le \DLamDNthree \, \mathfrak E_k^3.
\end{aligned}
\label{eq:ns-thm-high-nonlinear_appx}
\end{equation}

\vspace{2mm}

\begin{equation}
\begin{aligned}
&\left|
\sum_{i\in\{\omega,r,z\}}\langle \mathcal R_i,\bvar{\delta v_i}\rangle_{\wgray{\Phi_i}}
+\mu_k\sum_{i\in\{\omega,r,z\}}\sum_{|\alpha|=k}
\langle D^\alpha\mathcal R_i,D^\alpha\bvar{\delta v_i}\rangle_{\wgray{\Phi_{i,k}}}
\right|
\le \DLamR \, \mathfrak E_k.
\end{aligned}
\label{eq:ns-thm-residual_appx}
\end{equation}

\hfill \\

\noindent Choose $\mu_k \in (0,1]$ and
\begin{equation}
\label{eq:ns-mu-choice_appx}
\DLamStab(\mu_k):=
\min\left\{
\DLamLlow-\mu_k\DLamDLlow,
\ \DLamDLhigh-\frac{\DLamLhigh}{\mu_k}
\right\}>0.
\end{equation}
After this choice of $\mu_k$ is fixed, write $\DLamStab$ for the positive number $\DLamStab(\mu_k)$. \\

\noindent Define the combined nonlinear constant by
\begin{equation}
\label{eq:ns-N2N3_appx}
\DLamThree:=\DLamNthree+\DLamDNthree.
\end{equation}

\noindent Assume that $\delta_*>0$ satisfies

\begin{equation}
\label{eq:ns-smallness_appx}
\DLamThree\delta_*+\frac{\DLamR}{\delta_*}<\DLamStab.
\end{equation}

\hfill 

\noindent Then every solution with $\mathfrak E_k(0)<\delta_*$ remains in the ball $\mathfrak E_k(t)<\delta_*$ for as long as the solution exists and the certified estimates apply.
\end{theorem}

\newpage

\begin{grayproof}  \Lean \\

\noindent We pair the three perturbation equations with $(\bvar{\delta\omega},\partial_r\bvar{\delta u},\partial_z\bvar{\delta u})$ in the low-order weights and add the $k$th-order pairings with prefactor $\mu_k$. \\

\noindent The two linear estimates in the theorem, with the high-order estimate multiplied by $\mu_k$, give
\begin{align}
\label{eq:ns-close-linear-add_appx}
&-\DLamLlow \, \mathfrak E_0^2+\DLamLhigh \, \mathfrak H_k^2
+\mu_k\bigl(-\DLamDLhigh \, \mathfrak H_k^2+\DLamDLlow \, \mathfrak E_0^2\bigr) \\ &  \qquad \qquad \qquad =-(\DLamLlow-\mu_k\DLamDLlow)\mathfrak E_0^2-(\mu_k\DLamDLhigh-\DLamLhigh)\mathfrak H_k^2. \nonumber
\end{align}

\hfill

\noindent By \eqref{eq:ns-mu-choice_appx},
\begin{equation}
\DLamLlow-\mu_k\DLamDLlow\ge \DLamStab,
\qquad
\mu_k\DLamDLhigh-\DLamLhigh\ge \mu_k\DLamStab.
\end{equation}

\hfill 

\noindent Therefore
\begin{equation}
\label{eq:ns-close-linear-gamma_appx}
-(\DLamLlow-\mu_k\DLamDLlow)\mathfrak E_0^2-(\mu_k\DLamDLhigh-\DLamLhigh)\mathfrak H_k^2
\le -\DLamStab\mathfrak E_0^2-\mu_k\DLamStab\mathfrak H_k^2
= -\DLamStab\mathfrak E_k^2,
\end{equation}
since $\mathfrak E_k^2=\mathfrak E_0^2+\mu_k\mathfrak H_k^2$.

\hfill  \\

\noindent The low-order nonlinear estimate \eqref{eq:ns-thm-low-nonlinear_appx}, the high-order nonlinear estimate \eqref{eq:ns-thm-high-nonlinear_appx}, and \nsgreen{\eqref{eq:ns-N2N3_appx} give}

\begin{equation}
\label{eq:ns-close-nonlinear_appx}
\left|
\sum_{i\in\{\omega,r,z\}}
\langle\mathcal N_i,\bvar{\delta v_i}\rangle_{\wgray{\Phi_i}}
+\mu_k\sum_{i\in\{\omega,r,z\}}\sum_{|\alpha|=k}
\langle D^\alpha\mathcal N_i,D^\alpha\bvar{\delta v_i}\rangle_{\wgray{\Phi_{i,k}}}
\right| \le \DLamThree\mathfrak E_k^3.
\end{equation}

\hfill \\ 

\noindent Combining \eqref{eq:ns-close-linear-gamma_appx}, \eqref{eq:ns-close-nonlinear_appx}, and the residual bound \eqref{eq:ns-thm-residual_appx} yields

\begin{equation}
\label{eq:ns-energy-squared-ineq_appx}
\frac12\frac{d}{dt}\mathfrak E_k^2
\le -\DLamStab\mathfrak E_k^2+\DLamThree\mathfrak E_k^3+
\DLamR \mathfrak E_k.
\end{equation}

\hfill 

\noindent It remains to show that the energy ball is forward invariant.  We suppose that \(\mathfrak E_k(0)<\delta_*\) and argue by contradiction. Define the first boundary-contact time by
\begin{equation}
t_*:=\inf\{t>0:\mathfrak E_k(t)=\delta_*\}
\label{eq:ns-closing-first-contact-time}
\end{equation}
 \\ 

\noindent Set \(F(t):=\mathfrak E_k(t)^2\). Then \(F(t)<\delta_*^2\) for \(t<t_*\) and \(F(t_*)=\delta_*^2\), so the left derivative satisfies
\begin{equation}
F'(t_*)=
\lim_{h\downarrow0}\frac{F(t_*)-F(t_*-h)}{h}\ge0.
\label{eq:ns-closing-first-contact-left-derivative}
\end{equation}

\hfill 

\noindent On the other hand, evaluating \eqref{eq:ns-energy-squared-ineq_appx} at \(t_*\) gives

\begin{align}
\frac12F'(t_*)
&\le -\DLamStab\delta_*^2
+\DLamThree\delta_*^3
+\DLamR \delta_*
=-\delta_*^2\left(
\DLamStab-\DLamThree\delta_*
-\frac{\DLamR}{\delta_*}\right)<0
\label{eq:euler-boundary-inward-derivative}
\end{align}
by \eqref{eq:ns-smallness_appx}, which contradicts \eqref{eq:ns-closing-first-contact-left-derivative}.  \\

\noindent Therefore the ball
\(\mathfrak E_k<\delta_*\) is forward invariant, which proves the stated
nonlinear stability. 
\end{grayproof}

\clearpage

\section{Proof of the Physical Reconstruction Proposition}
\label{appx:physical-reconstruction}

\noindent This appendix provides the proof of the reconstruction proposition stated in~\Cref{sec:stability-to-physical-blowup}.

\hfill 

\begin{proposition}[Finite-time physical blowup from a signed amplitude rate]
\label{thm:physical-reconstruction-signed-rate-appx}
Consider the Euler equations in the dynamic rescaling formulation above, with $\lambda=1/2$. Assume that the exact admissible rescaled solution exists for every $t\ge0$, remains in the corresponding stability ball, and satisfies the modulation bounds throughout. Assume in addition that
\begin{equation}
 c_u(t)\le -\gamma_u<0,
 \qquad t\ge0,
\label{eq:reconstruction-cu-sign-appx}
\end{equation}
for some $\gamma_u>0$. Then the physical time $\mathring t(t)$ converges to a finite limit $T$, with
\begin{equation}
0<T-\mathring t(t)
\le
\frac{e^{-\gamma_u t}}{\gamma_u\,\mathrm{s}_u(0)}.
\label{eq:reconstruction-time-tail-appx-statement}
\end{equation}
The traveling center $\mathring z_c(t)$ converges to a finite physical location. Moreover,
\begin{equation}
\left|\omega^z\bigl(0,\mathring z_c(t),\mathring t(t)\bigr)\right|
=2|\bar u_0|\,\mathrm{s}_u(t)
\longrightarrow\infty.
\label{eq:reconstruction-vorticity-appx-statement}
\end{equation}
Consequently,
\begin{equation}
\|\tilde\omega(\cdot,\mathring t(t))\|_{L^\infty(\mathbb R^3)}\longrightarrow\infty,
\end{equation}
and the reconstructed physical solution cannot be continued smoothly through $T$.
\end{proposition}

\begin{grayproof} \Lean \\

\noindent  The amplitude scaling law and \eqref{eq:reconstruction-cu-sign-appx} give
\begin{equation}
\partial_t\mathrm{s}_u=-c_u\mathrm{s}_u,
\qquad
\mathrm{s}_u(t)\ge \mathrm{s}_u(0) \, e^{\gamma_u t}.
\label{eq:reconstruction-su-growth-appx}
\end{equation}

\hfill 

\noindent Using \eqref{eq:reconstruction-physical-clock}, define
\begin{equation}
T=\mathring t(0)+\int_0^\infty\frac{d\tau}{\mathrm{s}_u(\tau)}.
\end{equation}

\hfill 

\noindent The lower bound in \eqref{eq:reconstruction-su-growth-appx} makes this integral finite and gives
\begin{equation}
0<T-\mathring t(t)
=\int_t^\infty\frac{d\tau}{\mathrm{s}_u(\tau)}
\le
\frac{e^{-\gamma_u t}}{\gamma_u\mathrm{s}_u(0)},
\end{equation}
which proves \eqref{eq:reconstruction-time-tail-appx-statement}. \\

\noindent Since $\mathrm{s}_u(t)>0$, the physical clock is strictly increasing and maps $0\le t<\infty$ onto $\mathring t(0)\le\mathring t<T$.\\

\noindent At the moving physical center, the rescaled coordinates are $(r,z)=(0,0)$. The normalization condition $u(0,0,t)=\bar u_0$ and the inverse amplitude rescaling give
\begin{equation}
 u_\bullet\bigl(0,\mathring z_c(t),\mathring t(t)\bigr)
=\mathrm{s}_u(t)\bar u_0.
\label{eq:reconstruction-u-center-appx}
\end{equation}

\hfill 

\noindent Since $u^\theta=\mathring r\,u_\bullet$, axisymmetry gives
\begin{equation}
\omega^z
=\frac{1}{\mathring r}\partial_{\mathring r}(\mathring r u^\theta)
=2u_\bullet+\mathring r\,\partial_{\mathring r}u_\bullet.
\label{eq:reconstruction-vorticity-axis-formula-appx}
\end{equation}

\hfill 

\noindent Axis regularity then yields
\begin{equation}
\omega^z(0,\mathring z,\mathring t)=2u_\bullet(0,\mathring z,\mathring t).
\end{equation}

\hfill 

\noindent Combining this identity with \eqref{eq:reconstruction-u-center-appx} and \eqref{eq:reconstruction-su-growth-appx} gives
\begin{equation}
\left|\omega^z\bigl(0,\mathring z_c(t),\mathring t(t)\bigr)\right|
=2|\bar u_0|\mathrm{s}_u(t)
\longrightarrow\infty.
\end{equation}

\hfill 

\noindent It remains only to locate the singularity in physical space.\\

\noindent The modulation bounds give $|C(t)|\le C_*$ for some finite $C_*$.  \\ 

\noindent Since $\lambda=1/2$,
\begin{equation}
\mathrm{s}_r(t)=\mathrm{s}_r(0) \, e^{-t/2}.
\label{eq:reconstruction-sr-appx}
\end{equation}

\hfill 

\noindent The center equation therefore gives
\begin{equation}
\int_0^\infty\left|\partial_t\mathring z_c(t)\right|dt
\le
C_*\mathrm{s}_r(0)\int_0^\infty e^{-t/2}dt
<\infty.
\end{equation}
Hence $\mathring z_c(t)$ converges to a finite limit $\mathring z_*$. 

\hfill 

\noindent For each finite rescaled time, the exact change of variables gives an equivalent smooth physical solution. A smooth continuation through $T$ would keep the vorticity finite near $(0,\mathring z_*,T)$, contradicting \eqref{eq:reconstruction-vorticity-appx-statement}. 
\end{grayproof}

\stopSupplementalTocCapture

\clearpage

\bibliographystyle{unsrtnat}
\bibliography{bib}

\end{document}